\documentclass[12pt, a4paper]{amsart}
\usepackage{amsmath,amsthm,amsfonts,amscd,amssymb,eucal,latexsym,mathrsfs}

\usepackage[all]{xy}
\usepackage{latexsym, amssymb, amscd, mathrsfs, graphics, fullpage}
\usepackage[backref=page,colorlinks,bookmarksopen=true]{hyperref}

\def\gG{\mathfrak{G}}
\def\gH{\mathfrak{H}}

\def\gM{\mathfrak{M}}
\def\gN{\mathfrak{N}}

\begin{document}

%%%%%%%%%%%%%%%%%%%%%TITLE
\title[Measured Quantum Groupoids in action] 
{Measured Quantum Groupoids in action}

\author{Michel Enock}
\address{Institut de Math\'ematiques de Jussieu, Unit\'{e} Mixte Paris 6 / Paris 7 /
CNRS de Recherche 7586 \\175, rue du Chevaleret, Plateau 7E, F-75013 Paris}
 \email{enock@math.jussieu.fr}

\begin{abstract}
Franck Lesieur had introduced in his thesis (now published in an expended and revised version in the {\it M\'emoires de la SMF} (2007)) a notion of measured quantum groupoid, in the setting of von Neumann algebras and a simplification of Lesieur's axioms is presented in an appendix of this article. We here develop the notions of actions, crossed-product, and obtain a biduality theorem, following what had been done by Stefaan Vaes for locally compact quantum groups.  Moreover, we prove that the inclusion of the initial algebra into its crossed-product is depth 2, which gives a converse of a result proved by Jean-Michel Vallin and the author. More precisely, to any action of a measured quantum groupoid, we associate  another measured quantum groupoid. In particular, starting from an action of a locally compact quantum group, we obtain a measured quantum groupoid canonically associated to this action;  when the action is outer, this measured quantum groupoid is the initial locally compact quantum group. \end{abstract}

\subjclass{46L55, 46L89}
\keywords{measured quantum groupoids, actions, crossed-product, biduality theorem, depth 2 inclusions}

\maketitle
\newpage
\tableofcontents
\newpage

 %%%%%%%%%%intro
\section{Introduction}
\label{intro}
\subsection{}
 In two articles (\cite{Val2}, \cite{Val3}), J.-M. Vallin has introduced two notions (pseudo-multiplicative
unitary, Hopf-bimodule), in order to generalize, up to the groupoid
case, the classical notions of multiplicative unitary \cite{BS} and of Hopf-von Neumann algebras \cite{ES2}
which were introduced to describe and explain duality of groups, and leaded to appropriate notions
of quantum groups (\cite{ES2}, \cite{W1}, \cite{W2}, \cite{BS}, \cite{MN}, \cite{W3}, \cite{KV1}, \cite{KV2}, \cite{MNW}). 
\\ In another article \cite{EV}, J.-M. Vallin and the author have constructed, from a depth 2 inclusion of
von Neumann algebras $M_0\subset M_1$, with an operator-valued weight $T_1$ verifying a regularity
condition, a pseudo-multiplicative unitary, which leaded to two structures of Hopf bimodules, dual
to each other. Moreover, we have then
constructed an action of one of these structures on the algebra $M_1$ such that $M_0$
is the fixed point subalgebra, the algebra $M_2$ given by the basic construction being then
isomorphic to the crossed-product. We construct on $M_2$ an action of the other structure, which
can be considered as the dual action.
\\  If the inclusion
$M_0\subset M_1$ is irreducible, we recovered quantum groups, as proved and studied in former papers
(\cite{EN}, \cite{E2}).
\\ Therefore, this construction leads to a notion of "quantum groupoid", and a construction of a
duality within "quantum groupoids". 
\subsection{}
In a finite-dimensional setting, this construction can be
mostly simplified, and is studied in \cite{NV1}, \cite{BSz1},
\cite{BSz2}, \cite{Sz},\cite{Val4}, \cite{Val5}, and examples are described. In \cite{NV2}, the link between these "finite quantum
groupoids" and depth 2 inclusions of
$II_1$ factors is given. 
\subsection{}
Franck Lesieur introduced \cite{L1} a notion of "measured quantum groupoids", in which a modular hypothesis on the basis is required. Mimicking in a wider setting the technics of Kustermans and Vaes (\cite{KV1}, \cite{KV2}), he obtained then a pseudo-multiplicative unitary, which, as in the locally compact quantum group case, "contains" all the information of the object (the von Neuman algebra, the coproduct, the antipod, the co-inverse). Unfortunately, the axioms chosen then by Lesieur don't fit perfectely with the duality (namely, the dual object does not fit the modular condition on the basis chosen in \cite{L1}), and, in order to get a perfect symmetry between his objects and their duals, Lesieur gave the name of "measured quantum groupoids" to a wider class \cite{L2}, whose axioms could be described as the analog of \cite{MNW}, in which a duality is defined and studied, the initial objects considered in \cite{L1} being denoted now "adapted measured quantum groupoids". In \cite{E3} had been shown that, with suitable conditions, the objects constructed in \cite{EV} from depth 2 inclusions, are "measured quantum groupoids" in this new setting. 
\subsection{}
Unfortunately, the axioms given in \cite{L2} are very complicated, and there was a serious need for simplification. This has been done in \cite{E6}, which is given in an appendix of this article. 
\subsection{}
Here are developped the notion of action (already introduced in \cite{EV}), crossed-product, etc, following what had been done for locally compact quantum groups in \cite{E1}, \cite{ES1}, \cite{V2}. Then we prove a Takesaki-like biduality theorem; moreover, we prove that the inclusion of the initial algebra into its crossed-product is depth 2, and, therefore, using (\cite{EV}, \cite{E3}), we can construct another measured quantum groupoid, canonically associated to the action. These results generalize both the case of locally compact quantum groups (\cite{V2}) and the case of measured groupoids (\cite{Y1}, \cite{Y2}, \cite{Y3}). 
 \subsection{}
This article is organized as follows : 
\newline
In chapter \ref{pr} are recalled all the definitions and constructions needed for that theory, namely Connes-Sauvageot's relative tensor product of Hilbert spaces, fiber product of von Neumann algebras, depth 2 inclusions, and Vaes' Radon-Nykodym theorem. \newline
The chapter \ref{MQG} is a r\'esum\'e of the construction of measured quantum groupoids; this theory is developped using the axiomatization given in \cite{E6}; the equivalence with Lesieur's axioms given in \cite{L2} is postponed in an appendix. 
\newline
Chapter \ref{strong} is a technical part, mostly inspired by \cite{KV2}, which will be needed in chapter \ref{action}. 
\newline
In chapter \ref{corep} is introduced the notion of (unitary) corepresentation of a measured quantum groupoid; this notion will be used throughout the paper. 
\newline
In chapter \ref{action} is defined the notion of an action of a measured quantum groupoid on a von Neumann algebra; we define, using chapter \ref{strong}, the notion of integrable action. 
\newline
In chapter \ref{tech}, we define some technical properties of actions (which will be proved to be always satisfied in chapter \ref{biduality}). 
\newline
Chapter \ref{deltainv} is very technical; we prove, when there exists a  $\delta$-invariant weight, that there is a standard implementation of an action. 
\newline
In chapter \ref{dual} is defined the crossed-product of a von Neumann algebra by a measured quantum groupoid, via an action, and the dual action of the dual measured quantum groupoid on this crossed-product. 
\newline
In chapter \ref{aux}, we construct on the crossed-product an auxilliary weight which will satisfy all the properties needed for applying chapter \ref{deltainv}. Therefore, we get that the dual action has a standard implementation. This auxilliary weight will be used in chapter \ref{dualw} to define the dual weight. 
\newline
In chapter \ref{biduality}, we obtain a Takesaki-like biduality theorem : the double crossed-product is, in a certain sense, equivalent to the initial algebra, and the bidual action to the initial action. 
\newline
In chapter \ref{char}, a measured quantum groupoid $\gG$ be given, we characterize the crossed products by an action of $\gG$ among von Neumann algebras on which acts the measured quantum groupoid $\widehat{\gG}^c$. 
\newline
In chapter \ref{dualw}, using the weight defined in chapter \ref{aux} and duality theorems proved in chapter \ref{biduality}, we define the dual weight on the crossed-product, and obtain some results on this weight, and then on the bidual weight, which allow us to study the basic construction made from the inclusion of the initial algebra into its crossed-product, and, using chapter \ref{char}, to prove that this inclusion is depth 2. 
\newline
In chapter \ref{mqga}, we apply the results of (\cite{EV}, \cite{E4}) to the depth 2 inclusion obtained in chapter \ref{dualw}, which allow us to construct a new measured quantum groupoid. In particular, starting with an action of a locally compact quantum group, we obtain a measured quantum groupoid; this gives a new example of a measured quantum groupoid which should look familiar to specialists of group actions and groupoids. 
\newline
In the appendix had been postponed the proof that our axioms, developped in \cite{E6}, are equivalent with Lesieur's axioms, given in \cite{L2}. Namely, in appendix \ref{coinverse}, with our axioms, we construct the co-inverse, the scaling group and the antipod, mimicking what Lesieur did for his "adapted measured quantum groupoids". Then, in appendix \ref{AGB}, we prove the equivalence of the two sets of axioms. 

%%%%%%%%%%Preliminaries 
\section{Preliminaries}
\label{pr}
In this chapter are mainly recalled definitions and notations about Connes' spatial
theory (\ref{spatial}, \ref{rel}) and the fiber product construction (\ref{fiber}, \ref{slice})
which are the main technical tools of the theory of measured quantum theory. We recall also some definitions and results about depth 2 inclusions (\ref{basic}), and Vaes' Radon-Nykodym theorem (\ref{Vaes}).

 %%%%%%%%%%%%%%%%%%%%%%%%%%%%%%%%%%%%%%spatial
\subsection{Spatial theory \cite{C1}, \cite{S2}, \cite{T}}
 \label{spatial}
 Let $N$ be a von Neumann algebra, 
and let $\psi$ be a faithful semi-finite normal weight on $N$; let $\gN_{\psi}$, 
$\gM_{\psi}$, $H_{\psi}$, $\pi_{\psi}$, $\Lambda_{\psi}$,$J_{\psi}$, 
$\Delta_{\psi}$,... be the canonical objects of the Tomita-Takesaki construction 
associated to the weight $\psi$. Let $\alpha$  be a non-degenerate normal representation of $N$ on a
Hilbert space
$\mathcal{H}$. We may as well consider $\mathcal{H}$ as a left $N$-module, and write it then
$_\alpha\mathcal{H}$. Following (\cite{C1}, definition 1), we define the set of 
$\psi$-bounded elements of $_\alpha\mathcal{H}$ as :
\[D(_\alpha\mathcal{H}, \psi)= \lbrace \xi \in \mathcal{H};\exists C < \infty ,\| \alpha (y) \xi\|
\leq C \| \Lambda_{\psi}(y)\|,\forall y\in \gN_{\psi}\rbrace\]
Then, for any $\xi$ in $D(_\alpha\mathcal{H}, \psi)$, there exists a bounded operator
$R^{\alpha,\psi}(\xi)$ from $H_\psi$ to $\mathcal{H}$,  defined, for all $y$ in $\gN_\psi$ by :
\[R^{\alpha,\psi}(\xi)\Lambda_\psi (y) = \alpha (y)\xi\]
If there is no ambiguity about the representation $\alpha$, we shall write
$R^{\psi}(\xi)$ instead of $R^{\alpha,\psi}(\xi)$. This operator belongs to $Hom_N (H_\psi , \mathcal{H})$;
therefore, for any
$\xi$, $\eta$ in
$D(_\alpha\mathcal{H}, \psi)$, the operator :
\[\theta^{\alpha,\psi} (\xi ,\eta ) =  R^{\alpha,\psi}(\xi)R^{\alpha,\psi}(\eta)^*\]
belongs to $\alpha (N)'$; moreover, $D(_\alpha\mathcal{H}, \psi)$ is dense (\cite{C1}, lemma 2), stable
under $\alpha (N)'$, and the linear span generated by the operators $\theta^{\alpha,\psi} (\xi ,\eta
)$ is a weakly dense ideal in $\alpha (N)'$. 
 \newline
 With the same hypothesis, the operator :
\[<\xi,\eta>_{\alpha,\psi} = R^{\alpha,\psi}(\eta)^* R^{\alpha,\psi}(\xi)\]
belongs to $\pi_{\psi}(N)'$. Using Tomita-Takesaki's theory, $\pi_{\psi}(N)'$ is equal to
$J_\psi
\pi_{\psi}(N)J_\psi$, and therefore anti-isomorphic to $N$ (or isomorphic to the opposite von
Neumann algebra $N^o$). We shall consider now $<\xi,\eta>_{\alpha,\psi}$ as an element of $N^o$, and
the linear span generated by these operators is a dense algebra in $N^o$. More precisely (\cite{C1}, lemma 4, and \cite{S1}, lemme 1.5), we get that $<\xi, \eta>_{\alpha, \psi}^o$ belongs to $\gM_\psi$, and that :
 \[\Lambda_{\psi}(<\xi, \eta>_{\alpha, \psi}^o)=J_\psi R^{\alpha, \psi}(\xi)^*\eta\]
 \newline
 If $y$ in $N$ is analytical with respect to $\psi$, and if $\xi\in D(_\alpha\mathcal{H}, \psi)$, then we get that $\alpha(y)\xi$ belongs to $D(_\alpha\mathcal{H}, \psi)$ and that :
 \[R^{\alpha,\psi}(\alpha(y)\xi)=R^{\alpha,\psi}(\xi)J_\psi\sigma_{-i/2}^{\psi}(y^*)J_\psi\]
  So, if $\eta$ is another $\psi$-bounded element of $_\alpha\mathcal H$, we get :
 \[<\alpha(y)\xi, \eta>_{\alpha, \psi}^o=\sigma_{i/2}^\psi(y)<\xi, \eta>_{\alpha, \psi}^o\]
There exists (\cite{C1}, prop.3) a family  $(e_i)_{i\in I}$ of 
$\psi$-bounded elements of $_\alpha\mathcal H$, such that
\[\sum_i\theta^{\alpha, \psi} (e_i ,e_i )=1\]
Such a family will be called an $(\alpha,\psi )$-basis of $\mathcal H$. 
 \newline
 It is possible (\cite{EN} 2.2) to construct 
an $(\alpha,\psi )$-basis of $\mathcal H$, $(e_i)_{i\in I}$, such that the
operators $R^{\alpha, \psi}(e_i)$ are partial isometries with final supports 
$\theta^{\alpha, \psi}(e_i ,e_i )$ 2 by 2 orthogonal, and such that, if $i\neq j$, then 
$<e_i ,e_j>_{\alpha, \psi}=0$. Such a family will be called an $(\alpha, \psi)$-orthogonal basis of $\mathcal H$.  We have then, for $\xi$, $\eta$ in $D(_\alpha\mathcal{H}, \psi)$ :
 \[R^{\alpha, \psi}(\xi)=\sum_i\theta^{\alpha, \psi}(e_i, e_i)R^{\alpha, \psi}(\xi)=\sum_iR^{\alpha, \psi}(e_i)<\xi, e_i>_{\alpha, \psi}\]
 \[<\xi, \eta>_{\alpha, \psi}=\sum_i<\eta, e_i>_{\alpha, \psi}^*<\xi, e_i>_{\alpha, \psi}\]
 \[\xi=\sum_iR^{\alpha, \psi}(e_i)J_\psi\Lambda_\psi(<\xi, e_i>^o_{\alpha, \psi})\]
 the sums being weakly convergent. 
 \newline
 Moreover, we get that, for all $n$ in $N$, $\theta^{\alpha, \psi}(e_i, e_i)\alpha(n)e_i=\alpha(n)e_i$, and $\theta^{\alpha, \psi}(e_i, e_i)$ is the orthogonal projection on the closure of the subspace $\{\alpha(n)e_i, n\in N\}$. 
 \newline
Let $\beta$ be a normal non-degenerate
anti-representation
 of
$N$ on
$\mathcal{H}$. We may then as well consider $\mathcal{H}$ as a right $N$-module, and write it $\mathcal{H}_\beta$, or
consider
$\beta$ as a normal non-degenerate representation of the opposite von Neumann algebra $N^o$, and
consider
$\mathcal{H}$ as a left $N^o$-module. 
\newline
We can then define on $N^o$ the opposite faithful
 semi-finite normal weight $\psi ^o$; we have $\gN_{\psi ^o}=\gN_\psi^*$, and 
 the Hilbert space $H_{\psi ^o}$ will be, as usual, identified with $H_\psi$, 
by the identification, for all $x$ in $\gN_\psi$, of 
$\Lambda_{\psi ^o}(x^*)$ with $J_\psi \Lambda_\psi (x)$.
\newline
From these remarks, we infer that the set of 
$\psi ^o$-bounded elements of
$\mathcal{H}_\beta$ is :
\[D(\mathcal{H}_\beta, \psi^o) = \lbrace\xi\in\mathcal{H} ;\exists C < \infty ,
\|\beta (y^*)\xi\|
\le C \| \Lambda_{\psi}(y)\|,\forall y\in \gN_{\psi}\rbrace\]
and, for any $\xi$ in $D(\mathcal{H}_\beta, \psi^o)$ and $y$ in $\gN_\psi$,
 the bounded operator $R^{\beta,\psi^{o}}(\xi )$ is given by the formula :
\[R^{\beta,\psi{^o}}(\xi)J_{\psi}\Lambda_\psi (y) = \beta (y^*)\xi\]
This operator belongs to $Hom_{N^{o}}(H_\psi ,\mathcal{H} )$.
Moreover, $D(\mathcal{H}_\beta, \psi^o)$ is dense, stable under $\beta( N)'=P$, 
and, for all $y$ in $P$, we have :
\[R^{\beta,\psi{^o}}(y\xi)= yR^{\beta,\psi{^o}}(\xi)\]
Then, for any $\xi$, $\eta$ in $D(\mathcal{H}_\beta, \psi^o)$, the operator 
\[\theta^{\beta, \psi^{o}}(\xi ,\eta )=R^{\beta,\psi{^o}}(\xi)R^{\beta,\psi{^o}}(\eta)^*\] 
belongs to $P$, and the linear span generated by these operators is a dense ideal in $P$;
moreover, the operator-valued product 
$<\xi ,\eta >_{\beta,\psi^o}= R^{\beta,\psi{^o}}(\eta)^*R^{\beta,\psi{^o}}(\xi)$
belongs to $\pi_\psi (N)$; we shall consider now, for simplification, that $<\xi ,\eta
>_{\beta,\psi^o}$ belongs to $N$, and the linear span generated by these operators is a dense algebra
in $N$, stable under multiplication by analytic elements with respect to $\psi$. More precisely, $<\xi ,\eta
>_{\beta,\psi^o}$ belongs to $\gM_\psi$ (\cite{C1}, lemma 4) and we have (\cite{S1}, lemme 1.5) :
\[\Lambda_\psi(<\xi, \eta>_{\beta,\psi^o})=R^{\beta,\psi{^o}}(\eta)^*\xi\]
A $(\beta,\psi^o)$-basis of $\mathcal H$ is a family 
$(e_i )_{i\in I}$ of $\psi^o$-bounded elements of $\mathcal H_\beta$, such that :
\[\sum_i\theta^{\beta, \psi^{o}}(e_i ,e_i )=1\]
We have then, for all $\xi$ in $D(\mathcal H_\beta, \psi^o)$ :
\[\xi=\sum_i R^{\beta,\psi^o}(e_i)\Lambda_\psi(<\xi, e_i>_{\beta,\psi^o})\]
It is possible to choose the $(e_i )_{i\in I}$ such that
the $R^{\beta, \psi{^o}}(e_i)$ are partial isometries, with final supports
$\theta^{\beta, \psi^{o}}(e_i ,e_i )$ 2 by 2 orthogonal, and $<e_i, e_j>_{\beta, \psi^o}=0$ if $i\neq j$; such
a family will be then called a $(\beta, \psi^o)$-orthogonal basis of $\mathcal H$. We have then 
\[R^{\beta, \psi{^o}}(e_i)=\theta^{\beta, \psi^{o}}(e_i ,e_i )R^{\beta, \psi{^o}}(e_i)=
R^{\beta, \psi{^o}}(e_i)<e_i, e_i>_{\beta, \psi^o}\]
 Moreover, we get that, for all $n$ in $N$, and for all $i$, we have $\theta^{\beta, \psi^{o}}(e_i ,e_i )\beta (n)e_i=\beta (n)e_i$, 
and therefore, $\theta^{\beta, \psi^{o}}(e_i ,e_i )$ is the orthogonal projection on the closure of the subspace $\{\beta (n)e_i, n\in N\}$.

 %%%%%%%%%%%%%%%%%%%%%%%%%%%%%%%%%basic
\subsection{Jones' basic construction and operator-valued weights; depth 2 inclusions}
\label{basic}
Let $M_0\subset M_1$ be an inclusion of von Neumann algebras  (for simplification, these algebras will be supposed to be $\sigma$-finite), equipped with a normal faithful semi-finite operator-valued weight $T_1$ from $M_1$ to $M_0$ (to be more precise, from $M_1^{+}$ to the extended positive elements of $M_0$ (cf. \cite{T} IX.4.12)). Let $\psi_0$ be a normal faithful semi-finite weight on $M_0$, and $\psi_1=\psi_0\circ T_1$; for $i=0,1$, let $H_i=H_{\psi_i}$, $J_i=J_{\psi_i}$, $\Delta_i=\Delta_{\psi_i}$ be the usual objects constructed by the Tomita-Takesaki theory associated to these weights. Following (\cite{J}, 3.1.5(i)), the von Neumann algebra $M_2=J_1M'_0J_1$ defined on the Hilbert space $H_1$ will be called the basic construction made from the inclusion $M_0\subset M_1$. We have $M_1\subset M_2$, and we shall say that the inclusion $M_0\subset M_1\subset M_2$ is standard.  
 \newline
Following (\cite{EN} 10.6), for $x$ in $\gN_{T_1}$, we shall define $\Lambda_{T_1}(x)$ by the following formula, for all $z$ in $\gN_{\psi_{0}}$ :
\[\Lambda_{T_1}(x)\Lambda_{\psi_{0}}(z)=\Lambda_{\psi_1}(xz)\]
This operator belongs to $Hom_{M_{0}^o}(H_{0}, H_1)$; if $x$, $y$ belong to $\gN_{T_1}$, then $\Lambda_{T_1}(x)\Lambda_{T_1}(y)^*$ belongs to $M_{2}$, and
$\Lambda_{T_1}(x)^*\Lambda_{T_1}(y)=T_1(x^*y)\in M_0$. 
\newline
 Using then Haagerup's construction (\cite{T}, IX.4.24), it is possible to construct a normal semi-finite faithful operator-valued weight $T_2$ from $M_2$ to $M_1$ (\cite{EN}, 10.7), which will be called the basic construction made from $T_1$. If $x$, $y$ belong to $\gN_{T_1}$, then $\Lambda_{T_1}(x)\Lambda_{T_1}(y)^*$
belongs to $\gM_{T_{2}}$, and $T_{2}(\Lambda_{T_1}(x)\Lambda_{T_1}(y)^*)=xy^*$.  As we have,  for all normal semi-finite faithful weight $\varphi$ on $M_1$ :
\[\frac {d\varphi}{d\varphi^o}=\Delta_{\varphi}\]
the operator-valued weight $T_2$ is charaterized by the equality (\cite{EN}, 10.3) :
\[\frac{d\psi_1\circ T_2}{d\psi_0^o}=\frac {d\psi_1}{d(\psi_0\circ T_1)^o}=\Delta_{\psi_1}\]
The operator-valued weight $T_2$ from $M_2$ to $M_1$ will be called the basic construction made from the operator-valued weight $T_1$ from $M_1$ to $M_0$. 
\newline
Repeating this construction, we obtain by recurrence successive basic constructions, which lead to Jones' tower $(M_i)_{i\in \mathbb{N}}$ of von Neumann algebras, which is the inclusion :
\[M_0\subset M_1\subset M_2\subset M_3\subset M_4\subset ...\]
which is equipped (for $i\geq 1$)  with normal faithful semi-finite operator-valued weights $T_i$ from $M_i$ to $M_{i-1}$. We define then, by recurrence, the weight $\psi_i=\psi_{i-1}\circ T_i$ on $M_i$, and we shall write $H_i$, $J_i$, $\Delta_i$ instead of $H_{\psi_i}$, etc. We shall define the mirroring $j_i$ on $\mathcal L(H_i)$ by $j_i(x)=J_ix^*J_i$, for all $x$ in $\mathcal L(H_i)$.
\newline
Following (\cite{EN} 10.6), for $x$ in $\gN_{T_i}$, we shall define $\Lambda_{T_i}(x)$ by the following formula, for all $z$ in $\gN_{\psi_{i-1}}$ :
\[\Lambda_{T_i}(x)\Lambda_{\psi_{i-1}}(z)=\Lambda_{\psi_i}(xz)\]
Then, $\Lambda_{T_i}(x)$ belongs to $Hom_{M_{i-1}^o}(H_{i-1}, H_i)$; if $x$, $y$ belong to $\gN_{T_i}$, then $\Lambda_{T_i}(x)^*\Lambda_{T_i}(y)=T_i(x^*y)$, and $\Lambda_{T_i}(x)\Lambda_{T_i}(y)^*$ belongs to $M_{i+1}$; more precisely, it belongs to $\gM_{T_{i+1}}$, and $T_{i+1}(\Lambda_{T_i}(x)\Lambda_{T_i}(y)^*)=xy^*$. 
\newline
Let $M_0\subset M_1$ be an inclusion of von Neumann algebras, equipped with a normal semi-finite faithful operator-valued weight $T_1$ from $M_1$ to $M_0$; following ([GHJ] 4.6.4), we shall say that the inclusion $M_0\subset M_1$ is depth 2 if the inclusion (called the derived tower) :
\[M'_0\cap M_1\subset M'_0\cap M_2\subset M'_0\cap M_3\]
 is standard, and, following (\cite{EN}, 11.12), we shall say that the operator-valued weight $T_1$ is regular if both restrictions $T_{2|M'_0\cap M_2}$ and $T_{3|M'_1\cap M_3}$ are semi-finite. 
\newline
By Tomita-Takesaki theory, the Hilbert space $H_1$ bears a natural structure of $M_1-M_1^o$-bimodule, and, therefore, by restriction, of $M_0-M_0^o$-bimodule. Let us write $r$ for the canonical representation of $M_0$ on $H_1$, and $s$ for the canonical antirepresentation given, for all $x$ in $M_0$, by $s(x)=J_1r(x)^*J_1$. Let us have now a closer look to the subspaces $D(H_{1s}, \psi_0^o)$ and $D(_rH_1, \psi_0)$. If $x$ belongs to $\gN_{T_1}\cap\gN_{\psi_1}$, we easily get that $J_1\Lambda_{\psi_1}(x)$ belongs to $D(_rH_1, \psi_0)$, with :
\[R^{r, \psi_0}(J_1\Lambda_{\psi_1}(x))=J_1\Lambda_{T_1}(x)J_0\]
and $\Lambda_{\psi_1}(x)$ belongs to $D(H_{1s}, \psi_0)$, with :
\[R^{s, \psi_0^o}(\Lambda_{\psi_1}(x))=\Lambda_{T_1}(x)\]
In (\cite{E4}, 2.3) was proved that the subspace $D(H_{1s}, \psi_0^o)\cap D(_rH_1, \psi_0)$ is dense in $H_1$; let us write down and precise this result :

%%%%%propbasic
\subsubsection{{\bf Proposition}}
\label{propbasic}
{\it  Let us keep on the notations of this paragraph; let $\mathcal T_{\psi_1, T_1}$ be the algebra made of elements $x$ in $\gN_{\psi_1}\cap\gN_{T_1}\cap\gN_{\psi_1}^*\cap\gN_{T_1}^*$, analytical with respect to $\psi_1$, and such that, for all $z$ in $\mathbb{C}$, $\sigma^{\psi_1}_z(x_n)$ belongs to $\gN_{\psi_1}\cap\gN_{T_1}\cap\gN_{\psi_1}^*\cap\gN_{T_1}^*$. Then :
\newline
(i) the algebra $\mathcal T_{\psi_1, T_1}$ is weakly dense in $M_1$; it will be called Tomita's algebra with respect to $\psi_1$ and $T_1$; 
\newline
(ii) for any $x$ in  $\mathcal T_{\psi_1, T_1}$, $\Lambda_{\psi_1}(x)$ belongs to $D(H_{1s}, \psi_0)\cap D(_rH_1, \psi_0)$;
\newline
(iii) for any $\xi$ in $D(H_{1s}, \psi_0^o))$, there exists a sequence $x_n$ in $\mathcal T_{\psi_1, T_1}$ such that $\Lambda_{T_1}(x_n)=R^{s, \psi_0^o}(\Lambda_{\psi_1}(x))$ is weakly converging to $R^{s, \psi_0^o}(\xi)$ and $\Lambda_{\psi_1}(x_n)$ is converging to $\xi$. }

\begin{proof}
The result (i) is taken from (\cite{EN}, 10.12); we get in (\cite{E4}, 2.3) an increasing sequence of projections $p_n$ in $M_1$, converging to $1$, and elements $x_n$ in $\mathcal T_{\psi_1, T_1}$ such that $\Lambda_{\psi_1}(x_n)=p_n\xi$. So, (i) and (ii) were obtained in (\cite{E4}, 2.3) from this construction. More precisely, we get that :
\begin{eqnarray*}
T_1(x_n^*x_n)&=&<R^{s, \psi^o}(\Lambda_{\psi_1}(x_n)), R^{s, \psi_0^o}(\Lambda_{\psi_1}(x_n))>_{s, \psi_0^o}\\
&=&<p_n\xi, p_n\xi>_{s, \psi_0^o}\\
&=&R^{s, \psi^o}(\xi)^*p_nR^{s, \psi^o}(\xi)
\end{eqnarray*}
which is increasing and weakly converging to $<\xi, \xi>_{s, \psi_0^o}$.  \end{proof}

We finish by writing a proof of this useful lemma, we were not able to find in litterature:   
%%%lemT
\subsubsection{{\bf Lemma}}
\label{lemT}
{\it Let $M_0\subset M_1$ be an inclusion of von neumann algebras, equipped with a normal faithful semi-finite operator-valued weight $T$ from $M_1$ to $M_0$. Let $\psi_0$ be a normal semi-finite faithful weight on $M_0$, and $\psi_1=\psi_0\circ T$; if $x$ is in $\gN_T$, and if $y$ is in $M'_0\cap M_1$, analytical with respect to $\psi_1$, then $xy$ belongs to $\gN_T$. }
\begin{proof}
Let $a$ be in $\gN_{\psi_0}$; then $xa$ belongs to $\gN_{\psi_1}$, and $xya=xay$ belongs to $\gN_{\psi_1}$; moreover, let us consider the element $T(y^*x^*xy)$ of the positive extended part of $M_0^+$; we have :
\begin{eqnarray*}
<T(y^*x^*xy), \omega_{\Lambda_{\psi_0}(a)}>&=&\psi_1(a^*y^*x^*xya)\\
&=&\|\Lambda_{\psi_1}(xay)\|^2\\
&=&\|J_{\psi_1}\sigma_{-i/2}^{\psi_1}(y^*)J_{\psi_1}\Lambda_{\psi_1}(xa)\|^2\\
&=&\|J_{\psi_1}\sigma_{-i/2}^{\psi_1}(y^*)J_{\psi_1}\Lambda_T(x)\Lambda_{\psi_0}(a)\|^2\\
&\leq&\|\sigma_{-i/2}^{\psi_1}(y^*)\|^2\|\Lambda_T(x)\Lambda_{\psi_0}(a)\|^2\\
&=&\|\sigma_{-i/2}^{\psi_1}(y^*)\|^2<T(x^*x), \omega_{\Lambda_{\psi_0}(a)}>
\end{eqnarray*}
from which we get that $T(y^*x^*xy)$ is bounded and 
\[T(y^*x^*xy)\leq \|\sigma_{-i/2}^{\psi_1}(y^*)\|^2T(x^*x)\] \end{proof}

%%%%%%%%%%%%%%%%rel

\subsection{Relative tensor product \cite{C1}, \cite{S2}, \cite{T}}
\label{rel}
Using the notations of \ref{spatial}, let now $\mathcal{K}$ be another Hilbert space on which there exists
a non-degenerate representation
$\gamma$ of
$N$. Following J.-L. Sauvageot (\cite{S2}, 2.1), we define
the relative tensor product $\mathcal{H}\underset{\psi}{_\beta\otimes_\gamma}\mathcal{K}$ as the
Hilbert space obtained from the algebraic tensor product $D(\mathcal{H}_\beta ,\psi^o )\odot
\mathcal{K} $ equipped with the scalar product defined, for $\xi_1$, $\xi_2$ in $D(\mathcal{H}_\beta
,\psi^o )$,
$\eta_1$, $\eta_2$ in $\mathcal{K}$, by 
\[(\xi_1\odot\eta_1 |\xi_2\odot\eta_2 )=(\gamma(<\xi_1 ,\xi_2 >_{\beta,\psi^o})\eta_1 |\eta_2 )\]
where we have identified $N$ with $\pi_\psi (N)$ to simplifly the notations.
\newline
The image of $\xi\odot\eta$ in $\mathcal{H}\underset{\psi}{_\beta\otimes_\gamma}\mathcal{K}$ will be
denoted by
$\xi\underset{\psi}{_\beta\otimes_\gamma}\eta$. We shall use intensively this construction; one
should bear in mind that, if we start from another faithful semi-finite normal weight $\psi '$, we
get another Hilbert space $\mathcal{H}\underset{\psi'}{_\beta\otimes_\gamma}\mathcal{K}$; there exists an isomorphism $U^{\psi, \psi'}_{\beta, \gamma}$ from $\mathcal{H}\underset{\psi}{_\beta\otimes_\gamma}\mathcal{K}$ to $\mathcal{H}\underset{\psi'}{_\beta\otimes_\gamma}\mathcal{K}$, which is unique up to some functorial property (\cite{S2}, 2.6) (but this isomorphism
does not send  $\xi\underset{\psi}{_\beta\otimes_\gamma}\eta$ on
$\xi\underset{\psi'}{_\beta\otimes_\gamma}\eta$ !). 
\newline
When no confusion is possible about the representation and the anti-representation, we shall write
$\mathcal{H}\otimes_{\psi}\mathcal{K}$ instead of
$\mathcal{H}\underset{\psi}{_\beta\otimes_\gamma}\mathcal{K}$, and $\xi\otimes_\psi\eta$ instead
of
$\xi\underset{\psi}{_\beta\otimes_\gamma}\eta$.
\newline
In (\cite{S1} 2.1), the relative tensor product
$\mathcal{H}\underset{\psi}{_\beta\otimes_\gamma}\mathcal{K}$ is defined also, if
$\xi_1$, $\xi_2$ are in $\mathcal{H}$, $\eta_1$, $\eta_2$ are in $D(_\gamma\mathcal{K},\psi)$, by the
following formula :
\[(\xi_1\odot\eta_1 |\xi_2\odot\eta_2 )= (\beta(<\eta_1, \eta_2>_{\gamma,\psi})\xi_1 |\xi_2)\]
which leads to the the definition of a relative flip $\sigma_\psi$ which will be an isomorphism from
$\mathcal{H}\underset{\psi}{_\beta\otimes_\gamma}\mathcal{K}$ onto
$\mathcal{K}\underset{\psi^o}{_\gamma\otimes _\beta}\mathcal{H}$, defined, for any 
$\xi$ in $D(\mathcal{H}_\beta ,\psi^o )$, $\eta$ in $D(_\gamma \mathcal{K},\psi)$, by :
\[\sigma_\psi (\xi\otimes_\psi\eta)=\eta\otimes_{\psi^o}\xi\]
This allows us to define a relative flip $\varsigma_\psi$ from
$\mathcal{L}(\mathcal{H}\underset{\psi}{_\beta\otimes_\gamma}\mathcal{K})$ to $\mathcal{L}
(\mathcal{K}\underset{\psi^o}{_\gamma\otimes_\beta}\mathcal{H})$ which sends
$X$ in 
$\mathcal{L}(\mathcal{H}\underset{\psi}{_\beta\otimes_\gamma}\mathcal{K})$ onto
$\varsigma_\psi(X)=\sigma_\psi X\sigma_\psi^*$. Starting from another faithful semi-finite normal
weight $\psi'$, we get a von Neumann algebra
$\mathcal{L}(\mathcal{H}\underset{\psi'}{_\beta\otimes_\gamma}\mathcal{K})$ which is isomorphic to
$\mathcal{L}(\mathcal{H}\underset{\psi}{_\beta\otimes_\gamma}\mathcal{K})$, and a von Neumann
algebra $\mathcal{L} (\mathcal{K}\underset{\psi^{'o}}{_\gamma\otimes_\beta}\mathcal{H})$ which is
isomorphic to
$\mathcal{L} (\mathcal{K}\underset{\psi^o}{_\gamma\otimes_\beta}\mathcal{H})$; as we get that :
 \[\sigma_{\psi'}\circ U^{\psi, \psi'}_{\beta, \gamma}=U^{\psi^o, \psi'^o}_{\gamma, \beta}\]
 we see that these isomorphisms exchange $\varsigma_\psi$ and
$\varsigma_{\psi '}$. Therefore, the homomorphism $\varsigma_{\psi}$ can be denoted $\varsigma_N$
without any reference to a specific weight.
\newline
We recall, following
(\cite{S2}, 2.2b) that, for all $\xi$ in $\mathcal{H}$, $\eta$ in $D(_\gamma\mathcal{K},\psi)$, $y$ in
$N$, analytic with respect to $\psi$, we have :
\[\beta (y)\xi \otimes_\psi\eta=\xi\otimes_\psi\gamma(\sigma^\psi_{-i/2}(y))\eta\]
 With the notations of \ref{spatial}, let $(e_i)_{i\in I}$ a $(\beta, \psi^o)$-orthogonal basis of $\mathcal H$; 
 let us remark that, for all $\eta$ in $\mathcal K$, we have :
 \[e_i\underset{\psi}{_\beta\otimes_\gamma}\eta=e_i\underset{\psi}{_\beta\otimes_\gamma}\gamma(<e_i, e_i>_{\beta, \psi^o})\eta\]
On the other hand, $\theta^{\beta, \psi^o}(e_i, e_i)$ is an orthogonal projection, and so is $\theta^{\beta, \psi^o}(e_i, e_i)\underset{N}{_\beta\otimes_\gamma}1$; this last operator is the projection on the subspace $e_i\underset{\psi}{_\beta\otimes_\gamma}\gamma(<e_i, e_i>_{\beta, \psi^o})\mathcal K$ (\cite{E3}, 2.3) and, therefore, we get that $\mathcal H\underset{\psi}{_\beta\otimes_\gamma}\mathcal K$ is the orthogonal sum of the subspaces 
\[e_i\underset{\psi}{_\beta\otimes_\gamma}\gamma(<e_i, e_i>_{\beta, \psi^o})\mathcal K\]
 For any $\Xi$ in $\mathcal{H}\underset{\psi}{_\beta\otimes_\gamma}
\mathcal{K}$, there exist $\xi_i$ in $\mathcal K$, such that $\gamma(<e_i, e_i>_{\beta, \psi^o})\xi_i=\xi_i$ and $\Xi=\sum_i e_i\underset{\psi}{_\beta\otimes_\gamma}\xi_i$, from which we get that $\sum_i\|\xi_i\|^2=\|\Xi\|^2$. 
\newline
Let us suppose now that $\mathcal{K}$ is a $N-P$ bimodule; that means that there exists a von
Neumann algebra $P$, and a non-degenerate normal anti-representation $\epsilon$ of $P$ on
$\mathcal{K}$, such that
$\epsilon (P)\subset\gamma (N)'$. We shall write then $_\gamma\mathcal{K}_\epsilon$. If $y$ is in $P$, we
have seen that it is possible to define then the operator
$1_{\mathcal{H}}\underset{\psi}{_\beta\otimes_\gamma}\epsilon (y)$ on
$\mathcal{H}\underset{\psi}{_\beta\otimes_\gamma}\mathcal{K}$, and we define this way a
non-degenerate normal antirepresentation of $P$ on
$\mathcal{H}\underset{\psi}{_\beta\otimes_\gamma}\mathcal{K}$, we shall call again $\epsilon$ for
simplification. If $\mathcal H$ is a $Q-N$ bimodule, then
$\mathcal{H}\underset{\psi}{_\beta\otimes_\gamma}\mathcal{K}$ becomes a $Q-P$ bimodule (Connes'
fusion of bimodules).
\newline
Taking a faithful semi-finite normal weight
$\nu$  on $P$, and a left $P$-module $_{\zeta}\mathcal{L}$ (i.e. a Hilbert space $\mathcal{L}$ and a normal
non-degenerate representation $\zeta$ of $P$ on $\mathcal{L}$), it is possible then to define
$(\mathcal{H}\underset{\psi}{_\beta\otimes_\gamma}\mathcal{K})\underset{\nu}{_\epsilon\otimes_\zeta}\mathcal{L}$.
Of course, it is possible also to consider the Hilbert space
$\mathcal{H}\underset{\psi}{_\beta\otimes_\gamma}(\mathcal{K}\underset{\nu}{_\epsilon\otimes_\zeta}\mathcal{L})$.
It can be shown that these two Hilbert spaces are isomorphics as $\beta (N)'-\zeta
(P)^{'o}$-bimodules. (In (\cite{Val1} 2.1.3), the proof, given for $N=P$ abelian can be used, without
modification, in that wider hypothesis). We shall write then
$\mathcal{H}\underset{\psi}{_\beta\otimes_\gamma}\mathcal{K}\underset{\nu}{_\epsilon\otimes_\zeta}\mathcal{L}$
without parenthesis, to emphazise this coassociativity property of the relative tensor
product.
\newline
Dealing now with that Hilbert space $\mathcal{H}\underset{\psi}{_\beta\otimes_\gamma}\mathcal{K}\underset{\nu}{_\epsilon\otimes_\zeta}\mathcal{L}$, there exist different flips, and it is necessary to be careful with notations. For instance, $1\underset{\psi}{_\beta\otimes\sigma_\nu}$ is the flip from this Hilbert space onto $\mathcal{H}\underset{\psi}{_\beta\otimes_\gamma}(\mathcal L\underset{\nu^o}{_\zeta\otimes_\epsilon}\mathcal K)$, where $\gamma$ is here acting on the second leg of $\mathcal L\underset{\nu^o}{_\zeta\otimes_\epsilon}\mathcal K$ (and should therefore be written $1\underset{\nu^o}{_\zeta\otimes_\epsilon}\gamma$, but this will not be done for obvious reasons). Here, the parentheses remain, because there is no associativity rule, and to remind that $\gamma$ is not acting on $\mathcal L$. The adjoint of $1\underset{\psi}{_\beta\otimes\sigma_\nu}$ is $1\underset{\psi}{_\beta\otimes\sigma_{\nu^o}}$. 
\newline
The same way, we can consider $\sigma_\psi\underset{\nu}{_\epsilon\otimes_\zeta}1$ from $\mathcal{H}\underset{\psi}{_\beta\otimes_\gamma}\mathcal{K}\underset{\nu}{_\epsilon\otimes_\zeta}\mathcal{L}$ onto $(\mathcal K\underset{\psi^o}{_\gamma\otimes_\beta}\mathcal H)\underset{\nu}{_\epsilon\otimes_\zeta}\mathcal L$. 
 \newline
 Another kind of flip sends $\mathcal{H}\underset{\psi}{_\beta\otimes_\gamma}(\mathcal L\underset{\nu^o}{_\zeta\otimes_\epsilon}\mathcal K)$ onto $\mathcal L\underset{\nu^o}{_\zeta\otimes_\epsilon}(\mathcal{H}\underset{\psi}{_\beta\otimes_\gamma}\mathcal{K})$. We shall denote this application $\sigma^{1,2}_{\gamma, \epsilon}$ (and its adjoint $\sigma^{1,2}_{\epsilon, \gamma}$), in order to emphasize that we are exchanging the first and the second leg, and the representations $\gamma$ and $\epsilon$ on the third leg. 
 \newline
If $\pi$ denotes the canonical left representation of $N$ on the Hilbert space $L^2(N)$, we verify that the application which sends, for all $\xi$ in $\mathcal H$, $\psi$ normal semi-finite faithful weight on $N$, and $x$ in $\gN_\psi$, the vector $\xi\underset{\psi}{_\beta\otimes_\pi}J_\psi\Lambda_\psi(x)$ on $\beta(x^*)\xi$, gives an isomorphism of $\mathcal H\underset{\psi}{_\beta\otimes_\pi}L^2(N)$, which sends the antirepresentation of $N$ given by $n\mapsto 1_\mathcal H\underset{N}{_\beta\otimes_\pi}L^2(N)$ on $\beta$. The same way, there exists a canonical identification, as left $N$-modules, of $L^2(N)\otimes_\psi\mathcal K$ with $\mathcal K$. 
For any $\xi$ in $D(\mathcal{H}_\beta,
\psi^o)$, we define the bounded linear application $\lambda_\xi^{\beta, \gamma}$ from $\mathcal
K$ to
$\mathcal{H}\underset{\psi}{_\beta\otimes_\gamma}\mathcal{K}$ by, for all $\eta$ in $\mathcal K$,
$\lambda_\xi^{\beta, \gamma} (\eta)=\xi\underset{\psi}{_\beta\otimes_\gamma}\eta$. We shall write
$\lambda_\xi$ if no confusion is possible. We get (\cite{EN}, 3.10) :
\[\lambda_\xi^{\beta, \gamma}=R^{\beta, \psi^o}(\xi)\otimes_\psi 1_\mathcal K\]
We have :
\[(\lambda_\xi^{\beta, \gamma})^*\lambda_\xi^{\beta, \gamma}=\gamma(<\xi, \xi>_{\beta, \psi^o})\]
We may define, for any $\eta$ in $D(_\gamma\mathcal{K}, \psi)$, an application
$\rho_\eta^{\beta, \gamma}$ from $\mathcal H$ to
$\mathcal{H}\underset{\psi}{_\beta\otimes_\gamma}\mathcal{K}$ by
$\rho_\eta^{\beta, \gamma} (\xi)=\xi\underset{\psi}{_\beta\otimes_\gamma}\eta$. We shall write
$\rho_\eta$ if no confusion is possible. We get that :
\[(\rho_\eta^{\beta, \gamma})^*\rho_\eta^{\beta, \gamma}=\beta(<\eta, \eta>_{\gamma, \psi})\]
\newline
Let $x$ be an element
of
$\mathcal{L}(\mathcal{H})$, commuting with the right action of $N$ on $\mathcal{H}_\beta$ (i.e. $x$ belongs to $\beta(N)'$). It
is possible to define an operator $x\underset{\psi}{_\beta\otimes_\gamma} 1_{\mathcal{K}}$ on
$\mathcal{H}\underset{\psi}{_\beta\otimes_\gamma}
\mathcal{K}$. We can easily evaluate $\|x\underset{\psi}{_\beta\otimes_\gamma} 1_{\mathcal{K}}\|$ : for any finite $J\subset I$, for any $\eta_i$ in $\mathcal K$, we have :
 \begin{eqnarray*}
 ((x^*x\underset{\psi}{_\beta\otimes_\gamma} 1_{\mathcal{K}})(\Sigma_{i\in J}e_i\underset{\psi}{_\beta\otimes_\gamma}\eta_i)|(\Sigma_{i\in J}e_i\underset{\psi}{_\beta\otimes_\gamma}\eta_i))
 &=&
 \Sigma_{i\in J}(\gamma(<xe_i, xe_i>_{\beta, \psi^o})\eta_i|\eta_i)\\
 &\leq& \|x\|^2\Sigma_{i\in J}(\gamma(<e_i, e_i>_{\beta, \psi^o})\eta_i|\eta_i)\\
 &=&
 \|x\|^2\|\Sigma_{i\in J}e_i\underset{\psi}{_\beta\otimes_\gamma}\eta_i\|
 \end{eqnarray*}
 from which we get $\|x\underset{\psi}{_\beta\otimes_\gamma} 1_{\mathcal{K}}\|\leq\|x\|$. 
 \newline
 By the same way, if $y$ commutes with the left action of $N$ on
$_\gamma\mathcal{K}$ (i.e. $y\in\gamma(N)'$), it is possible to define
$1_{\mathcal{H}}\underset{\psi}{_\beta\otimes_\gamma}y$ on
$\mathcal{H}\underset{\psi}{_\beta\otimes_\gamma} \mathcal{K}$, and by composition, it is possible
to define then
$x\underset{\psi}{_\beta\otimes_\gamma} y$. If we start from another faithful semi-finite normal
weight $\psi '$, the canonical isomorphism $U^{\psi, \psi'}_{\beta, \gamma}$ from $\mathcal{H}\underset{\psi}{_\beta\otimes_\gamma}
\mathcal{K}$ to $\mathcal{H}\underset{\psi'}{_\beta\otimes_\gamma} \mathcal{K}$ sends
$x\underset{\psi}{_\beta\otimes_\gamma} y$ on $x\underset{\psi'}{_\beta\otimes_\gamma} y$ (\cite{S2},
2.3 and 2.6); therefore, this operator can be denoted $x\underset{N}{_\beta\otimes_\gamma} y$
without any reference to a specific weight, and we get $\|x\underset{N}{_\beta\otimes_\gamma} y\|\leq\|x\|\|y\|$. 
\newline
 If $\mathcal K$ is a Hilbert space on which there exists a non-degenerate representation $\gamma$ of $N$, then $\mathcal K$ is a $N-\gamma(N)'^o$ bimodule, and the conjugate Hilbert space $\overline{\mathcal K}$ is a $\gamma(N)'-N^o$ bimodule, and, (\cite{S2}), for any normal faithful semi-finite weight $\phi$ on $\gamma(N)'$, the fusion $_\gamma\mathcal K\underset{\phi^o}{\otimes}\overline{\mathcal K}_\gamma$ is isomorphic to the standard space $L^2(N)$, equipped with its standard left and right representation. 
 \newline
 Using that remark, the associativity rule,  and the identification (as right $N$-modules) of 
 $\mathcal H\otimes_\psi L^2(N)$ with $\mathcal H$, one gets for any $x\in \beta(N)'$  :
  \[\|x\underset{N}{_\beta\otimes_\gamma}1_{\mathcal K}\|\leq\|x\underset{N}{_\beta\otimes_\gamma}1_{\mathcal K}\underset{\gamma(N)'^o}{\otimes}1_{\overline{\mathcal K}}\|=\|x\underset{N}{_\beta\otimes}1_{L^2(N)}\|=\|x\|\]
from which we get $\|x\underset{N}{_\beta\otimes_\gamma}1_{\mathcal K}\|=\|x\|$ \vspace{2mm}. 
\newline
If $\mathcal H$ and $\mathcal K$ are finite-dimensional Hilbert spaces, the relative tensor product 
$\mathcal{H}\underset{\psi}{_\beta\otimes_\gamma}
\mathcal{K}$ can be identified with a subspace of the tensor Hilbert space 
$\mathcal{H}\otimes
\mathcal{K}$ (\cite{EV} 2.4), the projection on which belonging to $\beta (N)\otimes\gamma (N)$.

%%%%%%%%%%%%%fiber
\subsection{Fiber product \cite{V1}, \cite{EV}} 
\label{fiber}
Let us follow the notations of \ref{rel}; let now
$M_1$ be a von Neumann algebra on $\mathcal{H}$, such that $\beta (N)\subset
M_1$, and $M_2$ be a von Neumann algebra on $\mathcal{K}$, such that $\gamma (N)\subset
M_2$. The von Neumann algebra generated by all elements $x\underset{N}{_\beta\otimes_\gamma} y$,
where
$x$ belongs to $M'_1$, and $y$ belongs to $M'_2$, will be denoted
$M'_1\underset{N}{_\beta\otimes_\gamma} M'_2$ (or $M'_1\otimes_N M'_2$ if no confusion if
possible), and will be called the relative tensor product of
$M'_1$ and $M'_2$ over $N$. The commutant of this algebra will be denoted 
$M_1\underset{N}{_\beta *_\gamma} M_2$ (or $M_1*_N M_2$ if no confusion is possible) and called the
fiber product of $M_1$ and
$M_2$, over
$N$. It is straightforward to verify that, if $P_1$ and $P_2$ are two other von Neumann
algebras satisfying the same relations with $N$, we have :
\[M_1*_N M_2\cap P_1*_N P_2=(M_1\cap P_1)*_N (M_2\cap P_2)\]
Moreover, we get that $\varsigma_N (M_1\underset{N}{_\beta *_\gamma}
M_2)=M_2\underset{N^o}{_\gamma *_\beta}M_1$.
\newline
In particular, we have :
\[(M_1\cap \beta (N)')\underset{N}{_\beta\otimes_\gamma} (M_2\cap \gamma (N)')\subset
M_1\underset{N}{_\beta *_\gamma} M_2\] and :
\[M_1\underset{N}{_\beta *_\gamma} \gamma(N)=(M_1\cap\beta (N)')\underset{N}{_\beta\otimes_\gamma} 1\]
More generally, if
$\beta$ is a non-degenerate normal involutive antihomomorphism from
$N$ into a von Neumann algebra
$M_1$, and
$\gamma$ a non-degenerate normal involutive homomorphism from $N$ into a von Neumann
algebra
$M_2$, it is possible
to define, without any reference to a specific Hilbert space, a von Neumann algebra
$M_1\underset{N}{_\beta *_ \gamma}M_2$. 
\newline
Moreover, if now $\beta '$ is a non-degenerate normal involutive antihomomorphism from $N$ into
another von Neumann algebra
$P_1$,
$\gamma '$ a non-degenerate normal involutive homomorphism from $N$ into another
von Neumann algebra $P_2$, $\Phi$ a normal involutive homomorphism from $M_1$ into $P_1$ such that
$\Phi\circ\beta =\beta '$, and $\Psi$ a normal involutive homomorphism from $M_2$ into $P_2$ such that
$\Psi\circ\gamma=\gamma'$, it is possible then to define a normal involutive homomorphism (the proof
given in (\cite{S1} 1.2.4) in the case when $N$ is abelian can be extended without modification in the
general case) :
\[\Phi\underset{N}{_\beta *_\gamma}\Psi 
: M_1\underset{N}{_\beta
*_\gamma}M_2\mapsto P_1\underset{N}{_{\beta '}*_{\gamma '}}P_2\]
In the case when $_\gamma\mathcal{K}_\epsilon$ is a $N-P^o$ bimodule as explained in \ref{rel} and
$_\zeta\mathcal{L}$ a $P$-module, if
$\gamma (N)\subset M_2$ and $\epsilon (P)\subset M_2$, and if $\zeta (P)\subset M_3$, where $M_3$ is
a von Neumann algebra on $\mathcal{L}$, it is possible to consider then $(M_1\underset{N}{_\beta
*_\gamma}M_2)\underset{P}{_\epsilon *_\zeta}M_3$ and $M_1\underset{N}{_\beta
*_\gamma}(M_2\underset{P}{_\epsilon *_\zeta}M_3)$. The coassociativity property for relative tensor
products leads then to the isomorphism of these von Neumann algebra we shall write now 
$M_1\underset{N}{_\beta
*_\gamma}M_2\underset{P}{_\epsilon *_\zeta}M_3$ without parenthesis.
\newline
If $M_1$ and $M_2$ are finite-dimensional, the fiber product $M_1\underset{N}{_\beta *_ \gamma}M_2$
can be identified to a reduced algebra of $M_1\otimes M_2$ (reduced by a projector which belongs to
$\beta (N)\otimes \gamma (N)$) (\cite{EV} 2.4).

%%%%%%%%%%%slice
\subsection{Slice maps \cite{E4}}
\label{slice}
Let $A$ be in $M_1\underset{N}{_\beta *_\gamma}M_2$, $\psi$ a normal faithful semi-finite weight on $N$, $\mathcal{H}$ an Hilbert space on which $M_1$ is acting, $\mathcal{K}$ an Hilbert space on which $M_2$ is acting, and let $\xi_1$, $\xi_2$ be in
$D(\mathcal{H}_\beta,\psi^o)$; let us define :
\[(\omega_{\xi_1, \xi_2}\underset{\psi}{_\beta*_\gamma}id)(A)=(\lambda^{\beta, \gamma}_{\xi_2})^*A\lambda^{\beta, \gamma}_{\xi_1}\]
We define this way $(\omega_{\xi_1, \xi_2}\underset{\psi}{_\beta*_\gamma}id)(A)$ as a bounded operator on $\mathcal{K}$,
which belongs to $M_2$, such that :
\[((\omega_{\xi_1, \xi_2}\underset{\psi}{_\beta*_\gamma}id)(A)\eta_1|\eta_2)=
(A(\xi_1\underset{\psi}{_\beta\otimes_\gamma}\eta_1)|
\xi_2\underset{\psi}{_\beta\otimes_\gamma}\eta_2)\]
One should note that $(\omega_{\xi_1, \xi_2}\underset{\psi}{_\beta*_\gamma}id)(1)=\gamma (<\xi_1, \xi_2 >_{\beta, \psi^o})$. 
\newline
Let us define the same way, for any $\eta_1$, $\eta_2$ in
$D(_\gamma\mathcal{K}, \psi)$ :
\[(id\underset{\psi}{_\beta*_\gamma}\omega_{\eta_1, \eta_2})(A)=(\rho^{\beta, \gamma}_{\eta_2})^*A\rho^{\beta, \gamma}_{\eta_1}\]
which belongs to $M_1$. 
\newline
We therefore have a Fubini formula for these slice maps : for any $\xi_1$, $\xi_2$ in
$D(\mathcal{H}_\beta,\psi^o)$, $\eta_1$, $\eta_2$ in $D(_\gamma\mathcal{K}, \psi)$, we have :
\[<(\omega_{\xi_1, \xi_2}\underset{\psi}{_\beta*_\gamma}id)(A), \omega_{\eta_1, \eta_2}>=<(id\underset{\psi}{_\beta*_\gamma}\omega_{\eta_1,
\eta_2})(A),\omega_{\xi_1, \xi_2}>\]
Let $\phi_1$ be a normal semi-finite weight on
$M_1^+$, and $A$ be a positive element of the fiber product
$M_1\underset{N}{_\beta*_\gamma}M_2$, then we may define an element of the extended positive part
of $M_2$, denoted
$(\phi_1\underset{\psi}{_\beta*_\gamma}id)(A)$, such that, for all $\eta$ in $D(_\gamma L^2(M_2), \psi)$, we have :
\[\|(\phi_1\underset{\psi}{_\beta*_\gamma}id)(A)^{1/2}\eta\|^2=\phi_1(id\underset{\psi}{_\beta*_\gamma}\omega_\eta)(A)\]
Moreover, then, if $\phi_2$ is a normal semi-finite weight on $M_2^+$, we have :
\[\phi_2(\phi_1\underset{\psi}{_\beta*_\gamma}id)(A)=\phi_1(id\underset{\psi}{_\beta*_\gamma}\phi_2)(A)\]
and if $\omega_i$ be in $M_{1*}$ such that $\phi_1=sup_i\omega_i$, we
have $(\phi_1\underset{\psi}{_\beta*_\gamma}id)(A)=sup_i(\omega_i\underset{\psi}{_\beta*_\gamma}id)(A)$.
\newline
Let now $P_1$ be a von Neuman algebra such that :
\[\beta(N)\subset P_1\subset M_1\]
and let $\Phi_i$ ($i=1,2$)
be a normal faithful semi-finite operator valued weight from $M_i$ to $P_i$; for any positive
operator $A$ in the fiber product
$M_1\underset{N}{_\beta*_\gamma}M_2$, there exists an element $(\Phi_1\underset{N}{_\beta*_\gamma}id)(A)$
of the extended positive part
of $P_1\underset{N}{_\beta*_\gamma}M_2$, such that (\cite{E4}, 3.5), for all
$\eta$ in $D(_\gamma L^2(M_2), \psi)$, and $\xi$ in $D(L^2(P_1)_\beta, \psi^o)$, we have :
\[\|(\Phi_1\underset{N}{_\beta*_\gamma}id)(A)^{1/2}(\xi\underset{\psi}{_\beta\otimes_\gamma}\eta)\|^2=
\|\Phi_1(id\underset{\psi}{_\beta*_\gamma}\omega_\eta)(A)^{1/2}\xi\|^2\]
If $\phi$ is a normal semi-finite weight on $P$, we have :
\[(\phi\circ\Phi_1\underset{\psi}{_\beta*_\gamma}id)(A)=(\phi\underset{\psi}{_\beta*_\gamma}id)
(\Phi_1\underset{N}{_\beta*_\gamma}id)(A)\]
We define the same way an element $(id\underset{N}{_\beta*_\gamma}\Phi_2)(A)$ of the extended positive part
of
$M_1\underset{N}{_\gamma*_\beta}P_2$, and we have :
\[(id\underset{N}{_\beta*_\gamma}\Phi_2)((\Phi_1\underset{N}{_\beta*_\gamma}id)(A))=
(\Phi_1\underset{N}{_\beta*_\gamma}id)((id\underset{N}{_\beta*_\gamma}\Phi_2)(A))\]
Let $\pi$ denotes the canonical left representation of $N$ on the Hilbert space $L^2(N)$; let 
 $x$ be an element of $M_1{}_\beta\underset{N}{*}{}_\pi \pi(N)$, which can be identified 
(\ref{fiber}) to $M_1\cap\beta(N)'$, we get that, for $e$ in $\gN_\psi$, we have : \[(id_\beta\underset{\psi}{*}{}_\pi\omega_{J_\psi \Lambda_{\psi}(e)})(x)=\beta(ee^*)x\]
 Therefore, by increasing limits, we get that $(id_\beta\underset{\psi}{*}{}_\pi\psi)$ is the injection of $M_1\cap\beta(N)'$ into $M_1$.  More precisely, if $x$ belongs to $M_1\cap\beta(N)'$, we have :
 \[(id_\beta\underset{\psi}{*}{}_\pi\psi)(x{}_\beta\underset{N}{\otimes}{}_\pi 1)=x\]
 \newline
 Therefore, if $\Phi_2$ is a normal faithful semi-finite operator-valued weight from $M_2$ onto $\gamma(N)$, we get that, for all $A$ positive in $M_1\underset{N}{_\beta*_\gamma}M_2$, we have :
 \[(id_\beta\underset{\psi}{*}{}_\gamma\psi\circ\Phi_2)(A){}_\beta\underset{N}{\otimes}{}_\gamma 1=
 (id_\beta\underset{\psi}{*}{}_\gamma\Phi_2)(A)\]

With the notations of \ref{spatial}, let $(e_i)_{i\in I}$ be a $(\beta, \psi^o)$-orthogonal basis of $\mathcal H$; using the fact (\ref{rel}) that, for all $\eta$ in $\mathcal K$, we have :\[e_i\underset{\psi}{_\beta\otimes_\gamma}\eta=e_i\underset{\psi}{_\beta\otimes_\gamma}\gamma(<e_i, e_i>_{\beta, \psi^o})\eta\]
we get that, for all $X$ in $M_1\underset{N}{_\beta*_\gamma}M_2$, $\xi$ in $D(\mathcal H_\beta, \psi^o)$, we have 
\[(\omega_{\xi, e_i}\underset{\psi}{_\beta *_\gamma}id)(X)=\gamma(<e_i, e_i>_{\beta, \psi^o})(\omega_{\xi, e_i}\underset{\psi}{_\beta *_\gamma}id)(X)\]

%%%%%Vaes
\subsection{Vaes' Radon-Nikodym theorem}
\label{Vaes}
In \cite{V1} is proved a very nice Radon-Nikodym theorem for two normal faithful semi-finite weights on a von Neumann algebra $M$. If $\Phi$ and $\Psi$ are such weights, then are equivalent :
\newline
- the two modular automorphism groups $\sigma^\Phi$ and $\sigma^\Psi$ commute;
\newline
- the Connes' derivative $[D\Psi : D\Phi]_t$ is of the form :
\[[D\Psi : D\Phi]_t=\lambda^{it^2/2}\delta^{it}\]
where $\lambda$ is a non-singular positive operator affiliated to $Z(M)$, and $\delta$ is a non-singular positive operator affiliated to $M$. 
\newline
It is then easy to verify that $\sigma^\Phi_t(\delta^{is})=\lambda^{ist}\delta^{is}$, and that 
\[[D\Phi\circ\sigma^\Psi_t:D\Phi]_s=\lambda^{ist}\]
\[[D\Psi\circ\sigma^\Phi_t:D\Psi]_s=\lambda^{-ist}\]
Moreover, we have also, for any $x\in M^+$ :
\[\Psi(x)=lim_n\Phi((\delta^{1/2}e_n)x(\delta^{1/2}e_n))\]
where the $e_n$ are self-adjoint elements of $M$ given by the formula :
\[e_n=\frac{2n^2}{\Gamma(1/2)\Gamma(1/4)}\int_{\mathbb{R}^2}e^{-n^2x^2-n^4y^4}\lambda^{ix}\delta^{iy}dxdy\]
The operators $e_n$ are analytic with respect to $\sigma^\Phi$ and such that, for any $z\in \mathbb{C}$,  the sequence $\sigma_z^\Phi(e_n)$ is bounded and strongly converges to $1$. 
\newline
In that situation, we shall write $\Psi=\Phi_\delta$ and call $\delta$ the modulus of $\Psi$ with respect to $\Phi$; $\lambda$ will be called the scaling operator of $\Psi$ with respect to $\Phi$. 
\newline
Moreover, if $a\in M$ is such that $a\delta^{1/2}$ is bounded and its closure $\overline{a\delta^{1/2}}$ belongs to $\gN_\Phi$, then $a$ belongs to $\gN_\Psi$. We may then identify $\Lambda_\Psi(a)$ with $\Lambda_\Phi(\overline{a\delta^{1/2}})$, which leads to the identifications of $J_\Psi$ with $\lambda^{i/4}J_\Phi$, and of $\Delta_\Psi$ with $\overline{J_\Phi\delta^{-1}J_\Phi\delta\Delta_\Phi}$.

%%%%%%%%%%%MQG

\section{Measured quantum groupoids}
\label{MQG}
In this chapter, we first recall the definition of Hopf-bimodules
(\ref{Hbimod}). We then give
(\ref{defmult}) the definition of a pseudo-multiplicative unitary, give the fundamental example given by groupoids (\ref{gd}), and construct the algebras and the Hopf-bimodules "generated by the left (resp. right) leg" of a pseudo-multiplicative unitary (\ref{AW}). We recall the definition of left-(resp. right-) invariant operator-valued weights on a Hopf-bimodule; we then give the definition of a measured  quantum groupoid (\ref{defMQG}). We give the essential theorems which found the theory of measured  quantum groupoids and its duality (\ref{thL2}, \ref{thL3}, \ref{thL4}, \ref{thL5}). 

%%%%%%%%%Hbimod
\subsection{Definition}
\label{Hbimod}
A quintuplet $(N, M, \alpha, \beta, \Gamma)$ will be called a Hopf-bimodule, following (\cite{Val2}, \cite{EV} 6.5), if
$N$,
$M$ are von Neumann algebras, $\alpha$ a faithful non-degenerate representation of $N$ into $M$, $\beta$ a
faithful non-degenerate anti-representation of
$N$ into $M$, with commuting ranges, and $\Gamma$ an injective involutive homomorphism from $M$
into
$M\underset{N}{_\beta *_\alpha}M$ such that, for all $X$ in $N$ :
\newline
(i) $\Gamma (\beta(X))=1\underset{N}{_\beta\otimes_\alpha}\beta(X)$
\newline
(ii) $\Gamma (\alpha(X))=\alpha(X)\underset{N}{_\beta\otimes_\alpha}1$ 
\newline
(iii) $\Gamma$ satisfies the co-associativity relation :
\[(\Gamma \underset{N}{_\beta *_\alpha}id)\Gamma =(id \underset{N}{_\beta *_\alpha}\Gamma)\Gamma\]
This last formula makes sense, thanks to the two preceeding ones and
\ref{fiber}. The von Neumann algebra $N$ will be called the basis of $(N, M, \alpha, \beta, \Gamma)$\vspace{5mm}.\newline
If $(N, M, \alpha, \beta, \Gamma)$ is a Hopf-bimodule, it is clear that
$(N^o, M, \beta, \alpha,
\varsigma_N\circ\Gamma)$ is another Hopf-bimodule, we shall call the symmetrized of the first
one. (Recall that $\varsigma_N\circ\Gamma$ is a homomorphism from $M$ to
$M\underset{N^o}{_r*_s}M$).
\newline
If $N$ is abelian, $\alpha=\beta$, $\Gamma=\varsigma_N\circ\Gamma$, then the quadruplet $(N, M, \alpha, \alpha,
\Gamma)$ is equal to its symmetrized Hopf-bimodule, and we shall say that it is a symmetric
Hopf-bimodule\vspace{5mm}.\newline
Let $\mathcal G$ be a measured groupo\"{\i}d, with $\mathcal G^{(0)}$ as its set of units, and let us denote
by $r$ and $s$ the range and source applications from $\mathcal G$ to $\mathcal G^{(0)}$, given by
$xx^{-1}=r(x)$ and $x^{-1}x=s(x)$. As usual, we shall denote by $\mathcal G^{(2)}$ (or $\mathcal
G^{(2)}_{s,r}$) the set of composable elements, i.e. 
\[\mathcal G^{(2)}=\{(x,y)\in \mathcal G^2; s(x)=r(y)\}\]
Let $(\lambda^u)_{u\in \mathcal G^{(0)}}$ be a Haar system on $\mathcal G$ and $\nu$ a measure $\nu$ on $\mathcal G^{(0)}$. Let us write $\mu$ the measure on $\mathcal G$ given by integrating $\lambda^u$ by $\nu$ :
\[\mu=\int_{{\mathcal G}^{(0)}}\lambda^ud\nu\]
By definition,  $\nu$ is said quasi-invariant if $\mu$ is equivalent to its image under the inverse $x\mapsto x^{-1}$ of $\mathcal G$ (see \cite{Ra}, \cite{R1},
\cite{R2}, \cite{C2} II.5, \cite{P} and \cite{AR} for more details, precise definitions and examples of groupo\"{\i}ds). 
\newline
In \cite{Y1}, \cite{Y2}, \cite{Y3} and \cite{Val2} were associated to a measured groupo\"{\i}d $\mathcal G$, equipped with a Haar system $(\lambda^u)_{u\in \mathcal G ^{(0)}}$ and a quasi-invariant measure $\nu$ on $\mathcal G ^{(0)}$ two
Hopf-bimodules : 
\newline
The first one is $(L^\infty (\mathcal G^{(0)}, \nu), L^\infty (\mathcal G, \mu), r_{\mathcal G}, s_{\mathcal G}, \Gamma_{\mathcal
G})$, where we define $r_{\mathcal G}$ and $s_{\mathcal G}$ by writing , for $g$ in $L^\infty (\mathcal G^{(0)})$ :
\[r_{\mathcal G}(g)=g\circ r\]
\[s_{\mathcal G}(g)=g\circ s\]
 and where
$\Gamma_{\mathcal G}(f)$, for $f$ in $L^\infty (\mathcal G)$, is the function defined on $\mathcal G^{(2)}$ by $(s,t)\mapsto f(st)$;
$\Gamma_{\mathcal G}$ is then an involutive homomorphism from $L^\infty (\mathcal G)$ into $L^\infty
(\mathcal G^2_{s,r})$ (which can be identified to
$L^\infty (\mathcal G){_s*_r}L^\infty (\mathcal G)$).
\newline
The second one is symmetric; it is $(L^\infty (\mathcal G^{(0)}, \nu), \mathcal L(\mathcal G), r_{\mathcal G}, r_{\mathcal G},
\widehat{\Gamma_{\mathcal G}})$, where
$\mathcal L(\mathcal G)$ is the von Neumann algebra generated by the convolution algebra associated to the
groupo\"{\i}d
$\mathcal G$, and $\widehat{\Gamma_{\mathcal G}}$ has been defined in \cite{Y3} and
\cite{Val2}\vspace{5mm}.\newline
If $(N,M,r,s,\Gamma)$ be a Hopf-bimodule with a finite-dimensional algebra $M$, then, the
identification of $M\underset{N}{_\beta*_\alpha}M$ with a reduced algebra $(M\otimes M)_e$ (\ref{fiber})
leads to an injective homomorphism $\widetilde{\Gamma}$ from $M$ to $M\otimes M$ such that
$\widetilde{\Gamma}(1)=e\not= 1$ and $(\widetilde{\Gamma}\otimes
id)\widetilde{\Gamma}=(id\otimes\widetilde{\Gamma})\widetilde{\Gamma}$ (\cite{EV} 6.5). Then $(M,
\widetilde{\Gamma})$ is a weak Hopf
$\mathbb{C}^*$-algebra in the sense of (\cite{BSz1}, \cite{BSz2}, \cite{Sz}). 

%%%%%%%%%%defmult
\subsection{Definition}
\label{defmult}
Let $N$ be a von Neumann algebra; let
$\gH$ be a Hilbert space on which $N$ has a non-degenerate normal representation $\alpha$ and two
non-degenerate normal anti-representations $\hat{\beta}$ and $\beta$. These 3 applications
are supposed to be injective, and to commute two by two.  Let $\nu$ be a normal semi-finite faithful weight on
$N$; we can therefore construct the Hilbert spaces
$\gH\underset{\nu}{_\beta\otimes_\alpha}\gH$ and
$\gH\underset{\nu^o}{_\alpha\otimes_{\hat{\beta}}}\gH$. A unitary $W$ from
$\gH\underset{\nu}{_\beta\otimes_\alpha}\gH$ onto
$\gH\underset{\nu^o}{_\alpha\otimes_{\hat{\beta}}}\gH$
will be called a pseudo-multiplicative unitary over the basis $N$, with respect to the
representation $\alpha$, and the anti-representations $\hat{\beta}$ and $\beta$ (we shall write it is an $(\alpha, \hat{\beta}, \beta)$-pseudo-multiplicative unitary), if :
\newline
(i) $W$ intertwines $\alpha$, $\hat{\beta}$, $\beta$  in the following way : we have, for all $X\in N$ :
\[W(\alpha
(X)\underset{N}{_\beta\otimes_\alpha}1)=
(1\underset{N^o}{_\alpha\otimes_{\hat{\beta}}}\alpha(X))W\]
\[W(1\underset{N}{_\beta\otimes_\alpha}\beta
(X))=(1\underset{N^o}{_\alpha\otimes_{\hat{\beta}}}\beta (X))W\]
\[W(\hat{\beta}(X) \underset{N}{_\beta\otimes_\alpha}1)=
(\hat{\beta}(X)\underset{N^o}{_\alpha\otimes_{\hat{\beta}}}1)W\]
\[W(1\underset{N}{_\beta\otimes_\alpha}\hat{\beta}(X))=
(\beta(X)\underset{N^o}{_\alpha\otimes_{\hat{\beta}}}1)W\]
(ii) The operator $W$ satisfies :
\[(1_\gH\underset{N^o}{_\alpha\otimes_{\hat{\beta}}}W)
(W\underset{N}{_\beta\otimes_\alpha}1_{\gH})
=(W\underset{N^o}{_\alpha\otimes_{\hat{\beta}}}1_{\gH})
\sigma^{2,3}_{\alpha, \beta}(W\underset{N}{_{\hat{\beta}}\otimes_\alpha}1)
(1_{\gH}\underset{N}{_\beta\otimes_\alpha}\sigma_{\nu^o})
(1_{\gH}\underset{N}{_\beta\otimes_\alpha}W)\]
Here, $\sigma^{2,3}_{\alpha, \beta}$
goes from $(H\underset{\nu^o}{_\alpha\otimes_{\hat{\beta}}}H)\underset{\nu}{_\beta\otimes_\alpha}H$ to $(H\underset{\nu}{_\beta\otimes_\alpha}H)\underset{\nu^o}{_\alpha\otimes_{\hat{\beta}}}H$, 
and $1_{\gH}\underset{N}{_\beta\otimes_\alpha}\sigma_{\nu^o}$ goes from $H\underset{\nu}{_\beta\otimes_\alpha}(H\underset{\nu^o}{_\alpha\otimes_{\hat{\beta}}}H)$ to $H\underset{\nu}{_\beta\otimes_\alpha}H\underset{\nu}{_{\hat{\beta}}\otimes_\alpha}H$. 
\newline
All the properties supposed in (i) allow us to write such a formula, which will be called the
"pentagonal relation". 
\newline
One should note that this definition is different from the definition introduced in \cite{EV} (and repeated afterwards). It is in fact the same formula, the new writing :
\[\sigma^{2,3}_{\alpha, \beta}
(W\underset{N}{_{\hat{\beta}}\otimes_\alpha}1)
(1_{\gH}\underset{N}{_\beta\otimes_\alpha}\sigma_{\nu^o})\]
is here replacing the rather akward writing :
\[(\sigma_{\nu^o}\underset{N^o}{_\alpha\otimes_{\hat{\beta}}}1_{\gH})
(1_{\gH}\underset{N^o}{_\alpha\otimes_{\hat{\beta}}}W)
\sigma_{2\nu}
(1_{\gH}\underset{N}{_\beta\otimes_\alpha}\sigma_{\nu^o})\]
but denotes the same operator, and we suggest the reader to convince himself of this easy fact. 
\newline
If we start from another normal semi-finite faithful weight $\nu'$ on $N$, we may define, using \ref{rel}, another unitary $W^{\nu'}=U^{\nu^o,
\nu^{'o}}_{\alpha, {\hat{\beta}}}WU^{\nu', \nu}_{\beta, \alpha}$ from $\gH\underset{\nu'}{_\beta\otimes_\alpha}\gH$ onto
$\gH\underset{\nu^{'o}}{_\alpha\otimes_{\hat{\beta}}}\gH$. The formulae which link these isomorphims between relative product Hilbert spaces and the
relative flips allow us to check that this operator $W^{\nu'}$ is also pseudo-multiplicative; which can be resumed in saying that a
pseudo-multiplicative unitary does not depend on the choice of the weight on $N$. 
\newline
If $W$ is an $(\alpha, \hat{\beta}, \beta)$-pseudo-multiplicative unitary, then the unitary $\sigma_\nu W^*\sigma_\nu$ from $\gH\underset{\nu}{_{\hat{\beta}}\otimes_\alpha}\gH$ to $\gH\underset{\nu^o}{_\alpha\otimes_\beta}\gH$ is an $(\alpha, \beta, \hat{\beta})$-pseudo-multiplicative unitary, called the dual of $W$.

%%%%%%%%%AW
\subsection{Algebras and Hopf-bimodules associated to a pseudo-multiplicative unitary}
\label{AW}
For $\xi_2$ in $D(_\alpha\gH, \nu)$, $\eta_2$ in $D(\gH_{\hat{\beta}}, \nu^o)$, the operator $(\rho_{\eta_2}^{\alpha,
\hat{\beta}})^*W\rho_{\xi_2}^{\beta, \alpha}$ will be written $(id*\omega_{\xi_2, \eta_2})(W)$; we have, therefore, for all
$\xi_1$, $\eta_1$ in $\gH$ :
\[((id*\omega_{\xi_2, \eta_2})(W)\xi_1|\eta_1)=(W(\xi_1\underset{\nu}{_\beta\otimes_\alpha}\xi_2)|
\eta_1\underset{\nu^o}{_\alpha\otimes_{\hat{\beta}}}\eta_2)\]
and, using the intertwining property of $W$ with $\hat{\beta}$, we easily get that $(id*\omega_{\xi_2, \eta_2})(W)$ belongs
to $\hat{\beta} (N)'$. 
\newline
If $x$ belongs to $N$, we have :
\[(id*\omega_{\xi_2, \eta_2})(W)\alpha (x)=(id*\omega_{\xi_2, \alpha(x^*)\eta_2})(W)\]
\[\beta(x)(id*\omega_{\xi_2, \eta_2})(W)=(id*\omega_{\hat{\beta}(x)\xi_2, \eta_2})(W)\]
We shall write $A_w(W)$ the weak closure of the linear span of these operators, which are right $\alpha(N)$-modules and left $\beta(N)$-modules. Applying (\cite{E3} 3.6), we get that $A_w(W)$ and $A_w(W)^*$ are non-degenerate algebras (one should note that the notations of (\cite{E3}) had been changed in order to fit with Lesieur's notations). 
We shall write $\mathcal A(W)$ the von Neumann algebra generated by $A_w(W)$ .
We then have $\mathcal A(W)\subset\hat{\beta}(N)'$.
\newline
For $\xi_1$ in $D(\gH_\beta,\nu^o)$, $\eta_1$ in $D(_\alpha\gH, \nu)$, we shall write $(\omega_{\xi_1,
\eta_1}*id)(W)$ for the operator $(\lambda_{\eta_1}^{\alpha,
\hat{\beta}})^*W\lambda_{\xi_1}^{\beta, \alpha}$; we have,
therefore, for all
$\xi_2$,
$\eta_2$ in
$\gH$ :
\[((\omega_{\xi_1,\eta_1}*id)(W)\xi_2|\eta_2)=(W(\xi_1\underset{\nu}{_\beta\otimes_\alpha}\xi_2)|
\eta_1\underset{\nu^o}{_\alpha\otimes_{\hat{\beta}}}\eta_2)\]
and, using the intertwining property of $W$ with $\beta$, we easily get that $(\omega_{\xi_1,
\eta_1}*id)(W)$ belongs to $\beta(N)'$. 
\newline
We shall write $\widehat{A_w(W)}$ the weak closure of the linear span of these operators. It is clear that this weakly closed subspace is a non degenerate algebra; following (\cite{EV} 6.1 and 6.5), we shall write $\widehat{\mathcal A(W})$ the von Neumann algebra generated by  $\widehat{A_w(W)}$. We then have $\widehat{\mathcal A(W)}\subset\beta(N)'$. 
\newline
In (\cite{EV} 6.3 and 6.5), using the pentagonal equation, we got
that
$(N,\mathcal A(W),\alpha,\beta,\Gamma)$, and
$(N,\widehat{\mathcal A(W)}, \alpha, \hat{\beta}, \widehat{\Gamma})$ are Hopf-bimodules, where $\Gamma$ and
$\widehat{\Gamma}$ are defined, for any $x$ in $\mathcal A(W)$ and $y$ in $\widehat{\mathcal
A(W)}$, by :
\[\Gamma(x)=W^*(1\underset{N^o}{_\alpha\otimes_{\hat{\beta}}}x)W\]
\[\widehat{\Gamma}(y)=\sigma_{\nu^o}W(y\underset{N}{_\beta\otimes_\alpha}1)W^*\sigma_\nu\]
(Here also, we have changed the notations of \cite{EV}, in order to fit with Lesieur's notations). 
In (\cite{EV} 6.1(iv)), we had obtained that $x$ in $\mathcal L(\gH)$ belongs to $\mathcal A(W)'$ if and only if $x$ belongs to $\alpha(N)'\cap
\beta(N)'$ and verifies $(x\underset{N^o}{_\alpha\otimes_{\hat{\beta}}}1)W=W(x\underset{N}{_\beta\otimes_\alpha}1)$. 
We obtain the same way
that $y$ in $\mathcal L(\gH)$ belongs to $\widehat{\mathcal A(W)}'$ if and only if $y$ belongs to $\alpha(N)'\cap
\hat{\beta}(N)'$ and verifies $(1\underset{N^o}{_\alpha\otimes_{\hat{\beta}}}y)W=W(1\underset{N}{_\beta\otimes_\alpha}y)$.  \newline
Moreover, we get that $\alpha(N)\subset\mathcal A\cap\widehat{\mathcal A}$, $\beta(N)\subset\mathcal A$,
$\hat{\beta}(N)\subset\widehat{\mathcal A}$, and, for all $x$ in $N$ :
\[\Gamma (\alpha (x))=\alpha (x)\underset{N}{_\beta\otimes_\alpha}1\]
\[\Gamma (\beta (x))=1\underset{N}{_\beta\otimes_\alpha}\beta (x)\]
\[\widehat{\Gamma}(\alpha(x))=\alpha (x)\underset{N}{_{\hat{\beta}}\otimes_\alpha}1\]
\[\widehat{\Gamma}(\hat{\beta}(x))=1\underset{N}{_{\hat{\beta}}\otimes_\alpha}\hat{\beta}(x)\]

%%%%%gd
\subsection{Fundamental example}
\label{gd}
Let $\mathcal G$ be a measured groupoid; let's use all notations introduced in \ref{Hbimod}. Let us
note :
\[\mathcal G^2_{r,r}=\{(x,y)\in \mathcal G^2, r(x)=r(y)\}\]
Then, it has been shown \cite{Val2} that the formula $W_{\mathcal G}f(x,y)=f(x,x^{-1}y)$, where $x$, $y$ are
in
$\mathcal G$, such that $r(y)=r(x)$, and $f$ belongs to $L^2(\mathcal G^{(2)})$ (with respect to an
appropriate measure, constructed from $\lambda^u$ and $\nu$), is a unitary from $L^2(\mathcal G^{(2)})$ to $L^2(\mathcal G^2_{r,r})$ (with respect also to another
appropriate measure, constructed from $\lambda^u$ and $\nu$). 
\newline
Let us define $r_{\mathcal G}$ and $s_{\mathcal G}$ from
$L^\infty (\mathcal G^{(0)}, \nu)$ to $L^\infty (\mathcal G, \mu)$ (and then considered as representations on $\mathcal L(L^2(\mathcal
G, \mu))$, for any
$f$ in
$L^\infty (\mathcal G^{(0)}, \nu)$, by
$r_{\mathcal G}(f)=f\circ r$ and $s_{\mathcal G}(f)=f\circ s$.
\newline
We shall identify (\cite{Y3}, 3.2.2) the Hilbert space $L^2(\mathcal G^{(2)})$ with the relative Hilbert tensor product $L^2(\mathcal G, \mu)\underset{L^{\infty}(\mathcal G^{(0)}, \nu)}{_{s_{\mathcal G}}\otimes_{r_{\mathcal G}}}L^2(\mathcal G, \mu)$, and the Hilbert space $L^2(\mathcal G^2_{r,r})$ with the relative Hilbert tensor product $L^2(\mathcal G, \mu)\underset{L^{\infty}(\mathcal G^{(0)}, \nu)}{_{r_{\mathcal G}}\otimes_{r_{\mathcal G}}}L^2(\mathcal G, \mu)$. Moreover, the unitary $W_{\mathcal G}$ can be then interpreted \cite{Val3} as a pseudo-multiplicative unitary over the basis
$L^\infty (\mathcal G^{(0)}, \nu)$, with respect to the representation $r_{\mathcal G}$, and anti-representations
$s_{\mathcal G}$ and
$r_{\mathcal G}$ (as here the basis is abelian, the notions of representation and anti-representations are the same, and the commutation property is fulfilled). So, we get that $W_{\mathcal G}$ is a $(r_\mathcal G, s_\mathcal G, r_\mathcal G)$ pseudo-multiplicative unitary. 
\newline
Let us take the notations of \ref{AW}; the von Neumann algebra $\mathcal A(W_{\mathcal G})$ is equal to the von Neumann algebra $L^{\infty}(\mathcal G, \nu)$ (\cite{Val3}, 3.2.6 and 3.2.7); using (\cite{Val3}
3.1.1), we get that the Hopf-bimodule homomorphism
$\Gamma$ defined on
$L^{\infty}(\mathcal G, \mu)$ by $W_{\mathcal G}$ is equal to the usual Hopf-bimodule homomorphism $\Gamma_{\mathcal G}$ studied in \cite{Val2}, and
recalled in
\ref{Hbimod}.
Moreover, the von Neumann algebra $\widehat{\mathcal A(W_{\mathcal G})}$ is equal to the von Neumann algebra $\mathcal
L(\mathcal G)$ (\cite{Val3}, 3.2.6 and 3.2.7); using (\cite{Val3} 3.1.1), we get that the Hopf-bimodule homomorphism $\widehat{\Gamma}$ defined on $\mathcal L(\mathcal
G)$ by
$W_{\mathcal G}$ is the usual Hopf-bimodule homomorphism $\widehat{\Gamma_{\mathcal G}}$ studied in \cite{Y3} and \cite{Val2}. 
\newline
Let us suppose now that the groupoid $\mathcal G$ is locally compact in the sense of \cite{R1}; it has been proved in (\cite{E3} 4.8) that $W_\mathcal G$ satisfies a strong condition of regularity, called "norm regularity".

%%%%%%LW
\subsection{Definitions (\cite{L1}, \cite{L2})}
\label{LW}
Let $(N, M, \alpha, \beta, \Gamma)$ be a Hopf-bimodule, as defined in \ref{Hbimod}; a normal, semi-finite, faithful operator valued weight $T$ from $M$ to $\alpha (N)$ is said to be left-invariant if, for all $x\in \gM_T^+$, we have :
\[(id\underset{N}{_\beta*_\alpha}T)\Gamma (x)=T(x)\underset{N}{_\beta\otimes_\alpha}1\]
or, equivalently (\ref{slice}), if we write $\Phi=\nu\circ\alpha^{-1}\circ T$ :
\[(id\underset{N}{_\beta*_\alpha}\Phi)\Gamma (x)=T(x)\]
A normal, semi-finite, faithful operator-valued weight $T'$ from $M$ to $\beta (N)$ will be said to be right-invariant if it is left-invariant with respect to the symmetrized Hopf-bimodule, i.e., if, for all $x\in\gM_{T'}^+$, we have :
\[(T'\underset{N}{_\beta*_\alpha}id)\Gamma (x)=1\underset{N}{_\beta\otimes_\alpha}T'(x)\]
or, equivalently, if we write $\Psi=\nu\circ\beta^{-1}\circ T'$ : 
\[(\Psi\underset{N}{_\beta*_\alpha}id)\Gamma (x)=T'(x)\]

%%%%thL1
\subsection{Theorem(\cite{L1}, \cite{L2})}
\label{thL1}
{\it Let $(N, M, \alpha, \beta, \Gamma)$ be a Hopf-bimodule, as defined in \ref{Hbimod}, and let $T$ be a left-invariant normal, semi-finite, faithful operator valued weight from $M$ to $\alpha (N)$; let us choose a normal, semi-finite, faithful weight $\nu$ on $N$, and let us write $\Phi=\nu\circ\alpha^{-1}\circ T$, which is a normal, semi-finite, faithful weight on $M$; let us write $H_\Phi$, $J_\Phi$, $\Delta_\Phi$ for the canonical objects of the Tomita-Takesaki theory associated to the weight $\Phi$, and let us define, for $x$ in $N$, $\hat{\beta}(x)=J_\Phi\alpha(x^*)J_\Phi$.
\newline
(i) There exists an unique isometry $U$ from $H_\Phi\underset{\nu^o}{_\alpha\otimes_{\hat{\beta}}}H_\Phi$ to $H_\Phi\underset{\nu}{_\beta\otimes_\alpha}H_\Phi$, such that, for any $(\beta, \nu^o)$-orthogonal 
basis $(\xi_i)_{i\in I}$ of  $(H_\Phi)_\beta$, for any $a$ in $\gN_T\cap\gN_\Phi$ and for any $v$ in $D((H_\Phi)_\beta, \nu^o)$, we have 
\[U(v\underset{\nu^o}{_\alpha\otimes_{\hat{\beta}}}\Lambda_\Phi (a))=\sum_{i\in I} \xi_i\underset{\nu}{_\beta\otimes_\alpha}\Lambda_{\Phi}((\omega_{v, \xi_i}\underset{\nu}{_\beta*_\alpha}id)(\Gamma(a)))\]
(ii) Let us suppose there exists a right-invariant normal, semi-finite, faithful operator valued weight $T'$ from $M$ to $\beta (N)$; then this isometry is a unitary, and $W=U^*$ is an $(\alpha, \hat{\beta}, \beta)$-pseudo-multiplicative unitary from $H_\Phi\underset{\nu}{_\beta\otimes_\alpha}H_\Phi$ to $H_\Phi\underset{\nu^o}{_\alpha\otimes_{\hat{\beta}}}H_\Phi$ which verifies, for any $x$, $y_1$, $y_2$ in $\gN_T\cap\gN_\Phi$ :
\[(i*\omega_{J_\Phi\Lambda_\Phi (y_1^*y_2), \Lambda_\Phi (x)})(W)=
(id\underset{N}{_\beta*_\alpha}\omega_{J_\Phi\Lambda_\Phi(y_2), J_\Phi\Lambda_\Phi(y_1)})\Gamma (x^*)\]
Clearly, the pseudo-multplicative unitary $W$ does not depend upon the choice of the right-invariant operator-valued weight $T'$, and, for any $y$ in $M$, we have : }
\[\Gamma(y)=W^*(1\underset{N^o}{_\alpha\otimes_{\hat{\beta}}}y)W\]

\begin{proof} This is \cite{L2} 3.51 and 3.52. \end{proof}
%%%%defMQG
\subsection{Definitions}
\label{defMQG}

Let us take the notations of \ref{thL1}; let us write $\Psi=\nu\circ\beta^{-1}\circ T'$. We shall say that $\nu$ is relatively invariant with respect to $T$ and $T'$ if the two modular automorphism groups associated to the two weights $\Phi$ and $\Psi$ commute; we then write down\vspace{5mm}: \newline 
{\bf Definition}
\newline
A measured quantum groupoid is an octuplet $\gG=(N, M, \alpha, \beta, \Gamma, T, T', \nu)$ such that :
\newline
(i) $(N, M, \alpha, \beta, \Gamma)$ is a Hopf-bimodule, as defined in \ref{Hbimod}, 
\newline
(ii) $T$ is a left-invariant normal, semi-finite, faithful operator valued weight $T$ from $M$ to $\alpha (N)$, as defined in \ref{LW}, 
\newline
(iii) $T'$ is a right-invariant normal, semi-finite, faithful operator-valued weight $T'$ from $M$ to $\beta (N)$, as defined in \ref{LW}, 
\newline
(iv) $\nu$ is normal semi-finite faithful weight on $N$, which is relatively invariant with respect to $T$ and $T'$\vspace{5mm}. \newline 
{\bf Remark}
These axioms are not Lesieur's axioms, given in (\cite{L2}, 4.1). The equivalence of these axioms with Lesieur's axioms had been written down in \cite{E6}, and is recalled in the appendix (\ref{th}). 

%%%%%thL2
\subsection{Theorem (\cite{L2}, \cite{E6})}
\label{thL2}
{\it Let $\gG=(N, M, \alpha, \beta, \Gamma, T, T', \nu)$ be a measured quantum groupoid in the sense of \ref{defMQG}. Let us write $\Phi=\nu\circ\alpha^{-1}\circ T$, which is a normal, semi-finite faithful weight on $M$. Then 
\newline
(i) there exists a $*$-antiautomorphism $R$ on $M$, such that $R^2=id$, $R(\alpha(n))=\beta(n)$ for all $n\in N$, and 
\[\Gamma\circ R=\varsigma_{N^o}(R\underset{N}{_\beta*_\alpha}R)\Gamma\]
$R$ will be called the coinverse; 
\newline
(ii) there exists a one-parameter group $\tau_t$ of automorphisms of $M$, such that $R\circ\tau_t=\tau_t\circ R$ for all $t\in\mathbb{R}$, and, for all $t\in\mathbb{R}$ and $n\in N$, $\tau_t(\alpha(n))=\alpha(\sigma^\nu_t(n))$, $\tau_t(\beta(n))=\beta(\sigma^\nu_t(n))$ and :
\[\Gamma\circ\tau_t=(\tau_t\underset{N}{_\beta*_\alpha}\tau_t)\Gamma\]
\[\Gamma\circ\sigma_t^\Phi=(\tau_t\underset{N}{_\beta*_\alpha}\sigma_t^\Phi)\Gamma\]
$\tau_t$ will be called the scaling group;
\newline
(iii) the weight $\nu$ is relatively invariant with respect to $T$ and $RTR$; moreover, $R$ and $\tau_t$ are still the co-inverse and the scaling group of this new measured quantum groupoid, we shall denote 
\[\underline{\gG}=(N, M, \alpha, \beta, \Gamma, T, RTR, \nu) \]
\newline
(iv) for any $\xi$, $\eta$ in $D(_\alpha H_\Phi, \nu)\cap D((H_\Phi)_{\hat{\beta}}, \nu^o)$, $(id*\omega_{\xi, \eta})(W)$ belongs to $D(\tau_{i/2})$, and, if we define $S=R\tau_{-i/2}$, we have :
\[S((id*\omega_{\xi, \eta})(W))=(id*\omega_{\eta, \xi})(W)^*\]
More generally, for any $x$ in $D(S)=D(\tau_{-i/2})$, we get that $S(x)^*$ belongs to $D(S)$ and $S(S(x)^*)^*=x$; 
$S$ will be called the antipod of $\gG$ (or $\underline{\gG}$), and, therefore, the co-inverse and the scaling group, given by polar decomposition of the antipod, rely only upon the pseudo-multiplicative $W$. 
\newline
(v) there exists a one-parameter group $\gamma_t$ of automorphisms of $N$ such that, for all $t\in\mathbb{R}$ and $n\in N$, we have :
\[\sigma_t^{T}(\beta(n))=\beta(\gamma_t(n))\]
Moreover, we get that $\nu\circ\gamma_t=\nu$. 
\newline
(vi) there exists a positive non-singular operator $\lambda$ affiliated to $Z(M)$, and a positive non singular operator $\delta$ affiliated to $M$, such that :
\[(D\Phi\circ R: D\Phi)_t=\lambda^{it^2/2}\delta^{it}\]
and, therefore, we have :
\[(D\Phi\circ\sigma_s^{\Phi\circ R}:D\Phi)_t=\lambda^{ist}\]
The operator $\lambda$ will be called the scaling operator, and there exists a positive non-singular operator $q$ affiliated to $N$ such that $\lambda=\alpha(q)=\beta(q)$. We have $R(\lambda)=\lambda$. 
\newline
The operator $\delta$ will be called the modulus; we have $R(\delta)=\delta^{-1}$, and $\tau_t(\delta)=\delta$, for all $t\in\mathbb{R}$, and we can define a one-parameter group of unitaries $\delta^{it}\underset{N}{_\beta\otimes_\alpha}\delta^{it}$ which acts naturally on elementary tensor products, which verifies, for all $t\in\mathbb{R}$ :
\[\Gamma(\delta^{it})=\delta^{it}\underset{N}{_\beta\otimes_\alpha}\delta^{it}\]
(vii) we have $(D\Phi\circ\tau_t : D\Phi)_s=\lambda^{-ist}$, which leads to define a one-parameter group of unitaries by, for any $x\in\gN_\Phi$ :
\[P^{it}\Lambda_\Phi(x)=\lambda^{t/2}\Lambda_\Phi(\tau_t(x))\]
Moreover, for any $y$ in $M$, we get :
\[\tau_t(y)=P^{it}yP^{-it}\]
and it is possible to define one parameter groups of unitaries $P^{it}\underset{N}{_\beta\otimes_\alpha}P^{it}$ and $P^{it}\underset{N^o}{_\alpha\otimes_{\hat{\beta}}}P^{it}$ such that :
\[W(P^{it}\underset{N}{_\beta\otimes_\alpha}P^{it})=(P^{it}\underset{N^o}{_\alpha\otimes_{\hat{\beta}}}P^{it})W\]
Moreover, for all $v\in D(P^{-1/2})$, $w\in D(P^{1/2})$, $p$, $q$ in $D(_\alpha H_\Phi, \nu)\cap D((H_\Phi)_{\hat{\beta}}, \nu^o)$, we haveÊ:
\[(W^*(v\underset{\nu^o}{_\alpha\otimes_{\hat{\beta}}}q)|w\underset{\nu}{_\beta\otimes_\alpha}p)=
(W(P^{-1/2}v\underset{\nu}{_\beta\otimes_\alpha}J_\Phi p)|P^{1/2}w\underset{\nu^o}{_\alpha\otimes_{\hat{\beta}}}J_\Phi q)\]
We shall say that the pseudo-multiplicative unitary $W$ is "manageable", with "managing operator" $P$, which implies it is weakly regular in the sense of \cite{E3}, 4.1, which implies (with the notations of \ref{AW}) that $A_w(W)=\mathcal A(W)=M$ and $\widehat{A_w(W)}=\widehat{{\mathcal A}(W)}=\widehat{M}$
\newline
As $\tau_t\circ\sigma_t^\Phi=\sigma_t^\Phi\circ\tau_t$, we get that $J_\Phi PJ_\Phi=P$. }

\begin{proof} Result (i) is \ref{R}, result (ii) is \ref{tau}, result (iii) is \ref{thcentral}(iv) and \ref{def}; result (iv) is \ref{thS}, result (v) is \ref{gamma} and \ref{thcentral}(iii); result (vi) is given by \ref{Vaes} and (iii), \ref{thcentral}(iv), (v), and, using \ref{th}, \cite{L2}, 5.6 and 5.20; result (vii), using \ref{th}, is \cite{L2} 7.3. \end{proof}
%%%%VRN
\subsection{Notations}
\label{VRN}
We shall use the notations of \ref{Vaes} about Vaes' Radon-Nykodym theorem for the weights $\Phi$ and $\Phi\circ R$. Let $\lambda$ be the scaling operator of $\Phi\circ R$ with respect to $\Phi$, and $\delta$ be the modulus of $\Phi\circ R$ with respect to $\Phi$, defined in \ref{thL2}(vi). Then, the self-adjoint elements $e_n$ of $M$ defined in \ref{Vaes} verify $R(e_n)=e_n$. 
\newline
If $x\in M$ is such that $x\delta^{1/2}$ is bounded and its closure $\overline{x\delta^{1/2}}$ belongs to $\gN_\Phi$, then $x$ belongs to $\gN_{\Phi\circ R}$ (\ref{Vaes}), and, then $R(x^*)$ belongs to $\gN_\Phi$. In particular, if $x\in\gN_\Phi$, then $xe_n$ belongs to $\gN_\Phi\cap\gN_{\Phi\circ R}$ (and $R(x^*)e_n$ belongs also to $\gN_\Phi\cap\gN_{\Phi\circ R}$), and $\lim_n\Lambda_\Phi(xe_n)=\Lambda_\Phi(x)$.

%%%%thL3
\subsection{Theorem (\cite{L2})}
\label{thL3}
{\it Let $\gG=(N, M, \alpha, \beta, \Gamma, T, T', \nu)$ be a measured quantum groupoid in the sense of \ref{defMQG}. Let $\mathcal T$ be the subset of $\gN_T\cap\gN_{RTR}\cap\gN_\Phi\cap\gN_{\Phi\circ R}$ made of elements $x$ which are analytic with respect to both $\Phi$ and $\Phi\circ R$, and such that, for all $z$, $z'$ in $\mathbb{C}$, $\sigma_z^\Phi\circ\sigma_{z'}^{\Phi\circ R}(x)$ belongs to $\gN_T\cap\gN_{RTR}\cap\gN_\Phi\cap\gN_{\Phi\circ R}$. Let $E$ be the linear subspace generated by the elements of the form $e_nx$, where $e_n$ had been defined in \ref{VRN}, and $x\in \mathcal T$, and let us denote $\mathcal E$ the linear subspace $J_\Phi \Lambda_\Phi(E)$. Let us denote $M_*^{\alpha, \beta}$ the subspace of $M_*$ spanned by the positive elements $\omega$ in $M_*^+$ such that there exists $k\in\mathbb{R}^+$ such that $\omega\circ\alpha\leq k\nu$ and $\omega\circ\beta\leq k\nu$. Then :
\newline
(i)  $\mathcal T$ and $E$ are dense in $M$; $\Lambda_\Phi(\mathcal T)$ and $\mathcal E$ are dense in $H_\Phi$. As 
\[\mathcal E\subset D(_\alpha H_\Phi, \nu)\cap D((H_\Phi)_\beta, \nu^o)\]
we get that $D(_\alpha H_\Phi, \nu)\cap D((H_\Phi)_\beta, \nu^o)$ is dense in $H_\Phi$; therefore $M_*^{\alpha, \beta}$ is dense in $M_*$. Moreover, if $\eta$, $\eta'$ are in $\mathcal E$, then $<\eta, \eta'>_{\alpha, \nu}$ and $<\eta, \eta'>_{\beta, \nu^o}$ are analytic with respect to both $\sigma_t^\nu$ and $\gamma_t$. 
\newline
(ii) If $\omega_1$, $\omega_2$ are in $M_*^{\alpha, \beta}$, let us define $\omega_1\omega_2=(\omega_1\underset{N}{_\beta*_\alpha}\omega_2)\Gamma$. If $\omega\in M_*^{\alpha, \beta}$ and $x\in M$, let us define $\omega ^*(x)=\overline{\omega(R(x^*))}$;   with such product and involution, $M_*^{\alpha, \beta}$ is an involutive algebra, and we can define a faithful representation of this involutive algebra by :
\[\hat{\pi}(\omega)=(\omega*id)(W)\]
 Let $\widehat{M}$ be the weak closure of the algebra generated (which is a von Neumann algebra, by \ref{thL2} (vii) ) and is included in $\beta(N)'$. More generally, for any $v$ in $D(_\alpha H_\phi, \nu)$, $w$ in $D((H_\Phi)_\beta, \nu^o)$, let us write $\hat{\pi}(\omega_{v, w})=(\omega_{v, w}*id)(W)$; then $\hat{\pi}(\omega_{v, w})$ belongs to $\widehat{M}$. If, moreover, $v$ belongs also to $D((H_\Phi)_\beta, \nu^o)$, we have, for all $x\in\gN_\Phi$ :
\[\Lambda_\Phi(\omega_{w,v}\underset{N}{_\beta*_\alpha}id)\Gamma(x))=\hat{\pi}(\omega_{v,w})^*\Lambda_\Phi(x)\]
(iii) If $y$ is in $\widehat{M}$, let us define :
\[\widehat{\Gamma}(y)=\sigma_{\nu^o}W(y\underset{N}{_\beta\otimes_\alpha}1)W^*\sigma_\nu\] 
with these notations, $(N, \widehat{M}, \alpha, \hat{\beta}, \widehat{\Gamma})$ is a Hopf-bimodule. 
\newline
(iv) there exists a unique $*$-anti-homomorphism $\widehat{R}$ of $\widehat{M}$, such that $\widehat{R}\hat{\pi}(\omega)=\hat{\pi}(\omega\circ R)$, for all $\omega\in M_*^{\alpha, \beta}$. Moreover, we have $\widehat{R}(y)=J_\Phi y^*J_\Phi$ for all $y\in \widehat{M}$, $\widehat{R}(\alpha(n))=\hat{\beta}(n)$, for all $n\in N$, and :
\[\widehat{\Gamma}\circ \widehat{R}=\varsigma_{N^o}(\widehat{R}\underset{N}{_{\hat{\beta}}*_\alpha}\widehat{R})\widehat{\Gamma}\]
and, therefore, $\widehat{R}$ is a co-inverse of the Hopf-bimodule $(N, \widehat{M}, \alpha, \hat{\beta}, \widehat{\Gamma})$. 
\newline
(v) let $I_\Phi\subset M_*$ be defined as :
\[I_\Phi=\{\omega\in M_*; \exists k'\in\mathbb{R}^+, |\omega(x^*)|\leq k'\|\Lambda_\Phi(x)\|, \forall x\in\gN_\Phi\}\]
and, for any $\omega\in I_\Phi$, let $a_\Phi(\omega)$ be the unique vector in $H_\Phi$ such that, for all $x\in\gN_\Phi$ :
\[(a_\Phi(\omega)|\Lambda_\Phi(x))=\omega(x^*)\]
Then $I_\Phi\cap M_*^{\alpha, \beta}$ is a left ideal of $M_*^{\alpha, \beta}$, and 
there exists a unique normal semifinite faithful weight $\widehat{\Phi}$ on $\widehat{M}$, such that 
$\widehat{\pi}(I_\Phi\cap M_*^{\alpha, \beta})$ is a core for $\Lambda_{\widehat{\Phi}}$. Moreover, we may identify $H_{\widehat{\Phi}}$ with $H_\Phi$ in such a way that we have, for any $\omega\in M_*^{\alpha, \beta}\cap I_\Phi$ :
\[\Lambda_{\widehat{\Phi}}(\hat{\pi}(\omega))=a_\Phi(\omega)\]
We have also, with the notations of \ref{VRN}, for all $x\in\gN_\Phi$ and $n\in\mathbb{N}$ :
\[J_{\widehat{\Phi}}\Lambda_\Phi(xe_n)=\Lambda_\Phi(R(x^*)\overline{e_n\delta^{1/2}})\]
and, for any $y\in M$, we get $R(y)=J_{\widehat{\Phi}}y^*J_{\widehat{\Phi}}$. 
\newline
Moreover, there exists a unique 
left-invariant normal faithful semi-finite operator-valued weight $\widehat{T}$ from $\widehat{M}$ to $\alpha (N)$ such that 
\[\widehat{\Phi}=\nu\circ\alpha^{-1}\circ \widehat{T}\]
(vi)  the octuplet $(N, \widehat{M}, \alpha, \hat{\beta}, \widehat{\Gamma}, \widehat{T}, \widehat{R}\circ\widehat{T}\circ\widehat{R}, \nu)$ is a measured quantum groupoid, denoted by $\widehat{\gG}$ and called the dual of $\gG$ (or of $\underline{\gG}$). 
\newline
(vii) The pseudo-multiplicative unitary $\widehat{W}$ constructed from $\widehat{\gG}$ is :
\[\widehat{W}=\sigma_\nu W^*\sigma_\nu\]
The scaling group $\hat{\tau}_t$ of $\mathcal{\gG}$ is given, for all $y\in\widehat{M}$, by :
\[\hat{\tau}_t(y)=P^{it}yP^{-it}\]
The scaling operator $\hat{\lambda}$ of $\mathcal{\gG}$ satifies $\hat{\lambda}=\lambda^{-1}$.  The managing operator $\hat{P}$ of $\widehat{W}$ is equal to $P$. 
\newline
The modular operator $\Delta_{\widehat{\Phi}}$ is equal to the closure of $PJ_\Phi\delta^{-1}J_\Phi$, and, therefore, for all $x\in M$, we have :
\[\tau_t(x)=\Delta_{\widehat{\Phi}}^{it}x\Delta_{\widehat{\Phi}}^{-it}\]
The automorphism group $\hat{\gamma}_t$ is equal to $\gamma_{-t}$. }
\begin{proof} Result (i) is \cite{L2} 6.7, result (ii) is \cite{L2} 8.4,8.5, 3.18 and 8.1; result (iii) is \cite{L2} 8.2, (iv) is \cite{L2} 8.6, (v) is \cite{L2} 8.7, 8.14, 8.18 and 8.28, and results (vi) and (vii) are \cite{L2} 8.24. \end{proof}
%%%%thL4
\subsection{Theorem (\cite{L2})}
\label{thL4}
{\it Let $\gG=(N, M, \alpha, \beta, \Gamma, T, T', \nu)$ be a measured quantum groupoid in the sense of \ref{defMQG}, $\widehat{\gG}$ its dual in the sense of \ref{thL3};  then :
\newline
(i) the bidual $\widehat{\widehat{\gG}}$ is equal to $\underline{\gG}$; 
\newline
(ii) the modulus $\hat{\delta}$ is equal to the closure of $P^{-1}J_\Phi\delta J_\Phi\delta^{-1}\Delta_\Phi^{-1}$; therefore, we have $\overline{\delta\Delta_{\widehat{\Phi}}}=\overline{\hat{\delta}\Delta_\Phi}^{-1}$. 
\newline
(iii) it is possible to define one-parameter groups of unitaries $\Delta_{\widehat{\Phi}}^{it}\underset{\nu}{_\beta\otimes_\alpha}\Delta_\Phi^{it}$ and $\Delta_{\widehat{\Phi}}^{it}\underset{\nu^o}{_\alpha\otimes_{\hat{\beta}}}\Delta_\Phi^{it}$ and we have :
\[W(\Delta_{\widehat{\Phi}}^{it}\underset{\nu}{_\beta\otimes_\alpha}\Delta_\Phi^{it})=
(\Delta_{\widehat{\Phi}}^{it}\underset{\nu^o}{_\alpha\otimes_{\hat{\beta}}}\Delta_\Phi^{it})W\]
\[W(J_{\widehat{\Phi}}\underset{\nu^o}{_\alpha\otimes_{\hat{\beta}}}J_\Phi)=(J_{\widehat{\Phi}}\underset{\nu}{_\beta\otimes_\alpha}J_\Phi)W^*\]
(iv) for all $s$, $t$ in $\mathbb{R}$, we have :
\[\Delta_{\widehat{\Phi}}^{it}\Delta_\Phi^{is}=\lambda^{ist}\Delta_\Phi^{is}\Delta_{\widehat{\Phi}}^{it}\]
and we have $J_{\widehat{\Phi}}J_\Phi=\lambda^{i/4}J_\Phi J_{\widehat{\Phi}}$. 
\newline
(v) we have the following "Heisenberg-type" relations :
\[M\cap\widehat{M}=\alpha(N)\]
\[M\cap\widehat{M}'=\beta(N)\]
\[M'\cap\widehat{M}=\hat{\beta}(N)\]
\[M'\cap\widehat{M}'=\hat{\alpha}(N)\] 
where we put, for all $n\in N$, $\hat{\alpha}(n)=J_\Phi\beta(n)^*J_\Phi=J_{\widehat{\Phi}}\hat{\beta}(n^*)J_{\widehat{\Phi}}$. }
\begin{proof}
Result (i) is \cite{L2} 8.25, (ii) is \cite{L2} 8.24 (iii); results (iii) and (iv) are \cite{L2} 8.28, and (v) is \cite{L2} 8.30. \end{proof}

%%%%%thL5
\subsection{Theorem(\cite{L2})}
\label{thL5}
{\it Let $\gG=(N, M, \alpha, \beta, \Gamma, T, T', \nu)$ be a measured quantum groupoid in the sense of \ref{defMQG}, $\widehat{\gG}$ its dual in the sense of \ref{thL3}. Then :
\newline
(i) The octuplet $(N^o, M, \beta, \alpha, \varsigma_N\circ \Gamma, RTR, T, \nu^o)$ is a measured quantum groupoid we shall write $\gG^o$ and call the opposite measured quantum groupoid. We have $(\gG^o)^o=\underline{\gG}$. 
\newline
(ii) For any $x$ in $M$, let us write $j(x)=J_\Phi x^*J_\Phi$ the canonical $*$-anti-homomorphism from $M$ onto $M'$ and let us define then $\Gamma^c=(j\underset{N}{_\beta*_\alpha}j)\Gamma\circ j$, $R^c=j\circ R\circ j$, $T^c=j\circ T\circ j$. The octuplet $(N^o, M', \hat{\beta}, \hat{\alpha}, \Gamma^c, T^c, R^cT^cR^c, \nu^o)$
 is a measured quantum groupoid we shall write $\gG^c$ and call the commutant of $\gG$. We have $(\gG^c)^c=\underline{\gG}$. 
\newline
(iii) we have ${\gG^o}^c={\gG^c}^o$ and this last measured quantum groupoid is isomorphic to $\underline{\gG}$, via the isomorphism implemented by $J_\Phi J_{\widehat{\Phi}}$.
\newline
(iv) We have $\widehat{\gG^o}=(\widehat{\gG})^c$ and $\widehat{\gG^c}=(\widehat{\gG})^o$
\newline
(v) The pseudo-multiplicative unitaries $W^o$ and $W^c$ of $\gG^o$ and $\gG^c$ are given by :  
\[W^o=(J_{\widehat{\Phi}}\underset{N^o}{_\alpha\otimes_{\hat{\beta}}}J_{\widehat{\Phi}})W(J_{\widehat{\Phi}}\underset{N^o}{_\alpha\otimes_\beta}J_{\widehat{\Phi}})\]
\[W^c=(J_\Phi\underset{N^o}{_\alpha\otimes_{\hat{\beta}}}J_\Phi)W(J_\Phi\underset{N^o}{_{\hat{\alpha}}\otimes_{\hat{\beta}}}J_\Phi)\]
where $W^o$ is a unitary from $H_\Phi\underset{\nu^o}{_\alpha\otimes_\beta}H_\Phi$ onto $H_\Phi\underset{\nu}{_\beta\otimes_{\hat{\alpha}}}H_\Phi$ and $W^c$ is a unitary from $H_\Phi\underset{\nu^o}{_{\hat{\alpha}}\otimes_{\hat{\beta}}}H_\Phi$ onto $H_\Phi\underset{\nu}{_{\hat{\beta}}\otimes_\alpha}H_\Phi$. 
\newline
Moreover, for any $x\in M$, $y\in\widehat{M}$, we have :
\[\Gamma(x)=(\sigma_\nu W^o\sigma_\nu)^*(x\underset{N^o}{_{\hat{\alpha}}\otimes_\beta} 1)(\sigma_\nu W^o\sigma_\nu)\]
\[\widehat{\Gamma}(y)=W^c(y\underset{N^o}{_{\hat{\alpha}}\otimes_{\hat{\beta}}}1)(W^c)^*\]
(vi) The pseudo-multplicative unitary $W^{oc}$ of $\gG^{oc}$ is given by :
\[W^{oc}=(J_\Phi J_{\widehat{\Phi}}\underset{N^o}{_\alpha\otimes_{\hat{\beta}}}J_\Phi J_{\widehat{\Phi}})W(J_{\widehat{\Phi}}J_\Phi\underset{N}{_{\hat{\beta}}\otimes_{\hat{\alpha}}}J_{\widehat{\Phi}}J_\Phi)\]
and is a unitary from $H_\Phi\underset{\nu}{_{\hat{\beta}}\otimes_{\hat{\alpha}}}H_\Phi$ onto $H_\Phi\underset{\nu^o}{_{\hat{\alpha}}\otimes_\beta}H_\Phi$. 
\newline
Let us denote $\mathcal I_{\gG}$ the isomorphism between $M$ and $M'$ given, for $x\in M$, by :
\[\mathcal I_{\gG}(x)=J_\Phi J_{\widehat{\Phi}}xJ_{\widehat{\Phi}}J_\Phi\]
Then $\mathcal I_{\gG}$ is an isomorphism of measured quantum groupoids between $\gG$ and $\gG^{oc}$. 
\newline
Moreover, for any $y\in\widehat{M}'$, we have :}
\[\widehat{\Gamma}^c(y)=W^{oc}(y\underset{N}{_{\hat{\beta}}\otimes_{\hat{\alpha}}}1)(W^{oc})^*\]
\begin{proof}
Result (i) is \cite{L2} 17.8, (ii) is \cite{L2} 17.10, (iii) and (iv) are \cite{L2} 17.14; result (v) is \cite{L2} 17.11, and (vi) is a straightforward corollary of the preceeding statements. \end{proof}
%%%fd2
\subsection{Example}
\label{gd2}
Let $\mathcal G$ be a measured groupoid; let's use all notations introduced in \ref{Hbimod} and \ref{gd}. If $f\in L^\infty(\mathcal G, \mu)^+$, let's consider the function on $\mathcal G^{(0)}$, $u\mapsto \int_{\mathcal G}fd\lambda^u$, which belongs to $L^\infty (\mathcal G^{(0)}, \nu)$; the image of this function by the homomorphism $r_\mathcal G$ is the function on $\mathcal G$, $\gamma\mapsto \int_{\mathcal G}fd\lambda^{r(\gamma)}$; the application which sends $f$ on this function can be considered as an operator-valued weight from $L^\infty(\mathcal G, \mu)$ to $r_{\mathcal G}(L^\infty (\mathcal G^{(0)}, \nu))$ which is normal, semi-finite, faithful; by definition of the Haar system $(\lambda^u)_{u\in\mathcal G^{(0)}}$, it is left-invariant, in the sense of \ref{LW}; we shall denote by $T_\mathcal G$ this operator-valued weight from $L^\infty(\mathcal G, \mu)$ to $r_{\mathcal G}(L^\infty (\mathcal G^{(0)}, \nu))$.
 If we write $\lambda_u$ for the image of $\lambda^u$ under the inverse $x\mapsto x^{-1}$ of the groupoid $\mathcal G$, starting from the application which sends $f$ to the function on $\mathcal G^{(0)}$ defined by $u\mapsto \int_{\mathcal G}fd\lambda_u$, we define a normal semifinite faithful operator-valued weight from $L^\infty(\mathcal G, \mu)$ to $s_{\mathcal G}(L^\infty (\mathcal G^{(0)}, \nu))$, which is right-invariant in the sense of \ref{LW}, we shall denote by $T_{\mathcal G}^{(-1)}$.
 \newline
  We then get that :
\[(L^\infty(\mathcal G^{(0)}, \nu), L^\infty(\mathcal G, \mu), r_\mathcal G, s_\mathcal G, \Gamma_\mathcal G, T_\mathcal G, T_{\mathcal G}^{(-1)}, \nu)\] is a measured quantum groupoid, we shall denote again $\mathcal G$. 
\newline
It can be proved (\cite{E5}) that any measured quantum groupoid, whose underlying von Neumann algebra is abelian, is of that type. 
\newline
The constructions obtained by applying \ref{thL2} to this measured quantum groupoid are the following : the pseudo-multiplicative unitary is the operator $W_\mathcal G$ defined in \ref{gd}; the co-inverse is given by dualizing the inverse $x\mapsto x^{-1}$ of $\mathcal G$ to $L^\infty (\mathcal G, \mu)$; the scaling group is trivial, and the co-inverse is equal to the antipod; the automorphism group $\gamma_t$ is trivial, the scaling operator $\lambda$ is $1$, and the modulus $\delta$ is the usual Radon-Nikodym derivative between $\mu$ and its image under the inverse $x\mapsto x^{-1}$ of $\mathcal G$. 
\newline
The underlying Hopf bimodule of the dual $\widehat{\mathcal G}$ constructed in \ref{thL3} is the symmetric Hopf-bimodule $(L^\infty (\mathcal G^{(0)}, \nu), \mathcal L(\mathcal G), r_\mathcal G, r_\mathcal G, \widehat{\Gamma_\mathcal G})$ constructed in \cite{Y3} (\ref{Hbimod}).  
\newline
The constructions obtained by applying \ref{thL5} are the following : $(\mathcal G)^o$ is given by the opposite groupoid (in which we return the product, and exchange source and range applications). 
As $L^\infty(\mathcal G, \mu)$ is abelian, we have clearly $(\mathcal G)^c=\mathcal G$.

%%%quantumgroups
\subsection{Example}
\label{lcqg}
Let ${\bf G}=(M, \Gamma, \Phi, \Psi)$ be a locally compact quantum group, in the sense of Kustermanns and Vaes (\cite{KV1}, or, preferably, the von Neumann version given in \cite{KV2}). Here $M$ is a von Neumann algebra, $\Gamma$ a coassociative coproduct $M\mapsto M\otimes M$ (where $\otimes$ is the von Neumann tensor product), and $\Phi$ (resp. $\Psi$) is a left-invariant (resp. a right-invariant) normal semi-finite faithful weight on $M^+$. Then, thanks to (\cite{KV1}, prop. 6.1), we see that :
\[(\mathbb{C}, M, id, id, \Gamma, \Phi, \Psi, \omega)\]
(where $id$ is the canonical homomorphism from $\mathbb{C}$ into $M$, and $\omega$ is the canonical state on $\mathbb{C}$), is a measured quantum groupoid, we shall still denote by ${\bf G}$.  
\newline
Conversely, it is clear that any measured quantum groupoid whose basis $N$ is equal to $\mathbb{C}$ is of that type. 

%%%%d2
\subsection{Theorem(\cite{EV}, \cite{E4} 8.2 and 8.3)}
\label{d2}
{\it Let $M_0\subset M_1$ be a depth 2 inclusion of $\sigma$-finite von Neumann algebras, equipped with a regular (in the sense of \ref{basic}) normal semi-finite faithful operator-valued weight $T_1$. Then :
\newline
 (i) there exists an application $\tilde{\Gamma}$ from $M'_0\cap M_2$ to 
\[(M'_0\cap M_2)_{j_1}\underset{M'_0\cap M_1}{*}{}_{id}(M'_0\cap M_2)\]
such that $(M'_0\cap M_1, M'_0\cap M_2, id, j_1, \tilde{\Gamma})$ is a Hopf-bimodule, (where $id$ means here the injection of $M'_0\cap M_1$ into $M'_0\cap M_2$, and $j_1$ means here the restriction of $j_1$ to $M'_0\cap M_1$, considered then as an anti-representation of $M'_0\cap M_1$ into $M'_0\cap M_2$). Moreover, the anti-automorphism $j_1$ of $M'_0\cap M_2$ is a co-inverse for this Hopf-bimodule structure. 
\newline 
 (ii) Let us write $\tilde{T_2}$ for the restriction of $T_2$ to $M'_0\cap M_2$ (which is semi-finite, by definition of the regularity of $T_1$(\ref{basic})); we then get that $\tilde{T_2}$ is a left-invariant operator-valued weight from $M'_0\cap M_2$ to $M'_0\cap M_1$, and, therefore, that $j_1\circ\tilde{T_2}\circ j_1$ is a right-invariant operator-valued weight from $M'_0\cap M_2$ to $M'_1\cap M_2$. 
\newline 
(iii) Let us suppose moreover that there exists on $M'_0\cap M_1$ a normal faithful semi-finite weight $\chi$ invariant under the modular automorphism group of $T_1$and let $\chi_2$ be the weight $\chi\circ \tilde{T_2}$; the modular automorphism groups $\sigma_t^{\chi_2}$ and $\sigma_s^{\chi_2\circ j}$ commute. 
\newline
Therefore, $(M'_0\cap M_1, M'_0\cap M_2, id, j_1, \tilde{\Gamma}, \tilde{T_2}, j_1\circ\tilde{T_2}\circ j_1, \chi)$ is a measured quantum groupoid $\gG_1$ in the sense of \ref{defMQG}. 
\newline
Moreover, the inclusion $M_1\subset M_2$ satisfies the same hypothesis, and leads to another measured quantum groupoid $\gG_2$, which is isomorphic to $\widehat{\gG_1}^o$. }

Let us describe how the coproduct $\tilde{\Gamma}$ is constructed : let $\psi_0$ be a normal semi-finite faithful weight on $M_0$, $\psi_1=\psi_0\circ T_1$ the normal semi-finite faithful weight constructed then on $M_1$, $\psi_2=\psi_1\circ T_2$ the normal semi-finite faithful weight constructed then on $M_2$; let $\chi$ be a normal semi-finite faithful weight on $M'_0\cap M_1$, $\chi_2=\chi_\circ\tilde{T_2}$ the normal semi-finite faithful weight constructed then on $M'_0\cap M_2$. 
\newline
Then, the application defined, for $x\in\gN_{\tilde{T_2}}\cap\gN_{\chi_2}$ and $y\in\gN_{\psi_1}$ by :
\[U(\Lambda_{\chi_2}(x)\otimes_\chi\Lambda_{\psi_1}(y))=\Lambda_{\psi_2}(xy)\]
is a unitary from $H_{\chi_2}\otimes_\chi H_{\psi_1}$ onto $H_{\psi_2}$ (\cite{EV}, 4.2 and 4.6). 
\newline
Moreover, using th fact that the inclusion $M_0\subset M_2\subset M_4$ is standard (\cite{EN}, 4.5), we get (\cite{EV}, 4.6) :
\[U^*[M'_0\cap M_4\cap (M'_0\cap M_2)']U=\pi_{\chi_2}(M'_0\cap M_2)'*_{M'_0\cap M_1}M'_0\cap M_2\]
and, on the other hand, composing the mirrorings (with the notations of (\ref{basic})), we get that :
\[j_2\circ j_1(M'_0\cap M_2)=M'_2\cap M_4\subset M'_0\cap M_4\cap (M'_0\cap M_2)'\]
and, therefore, $x\mapsto U^*j_2\circ j_1(x)U$ is an injective $*$-homomorphism from $M'_0\cap M_2$ into $\pi_{\chi_2}(M'_0\cap M_2)'*_{M'_0\cap M_1}M'_0\cap M_2$. 
\newline
On the other hand, using $U$, it is possible to define a pseudo-multiplicative unitary $\mathcal W$ (\cite{EV}, 5.3), which leads to a coproduct $\Gamma$ on $\pi_{\chi_2}(M'_0\cap M_2)'$ by the formula ($y\in \pi_{\chi_2}(M'_0\cap M_2)'$) :
\[\Gamma(y)=\mathcal W (y\otimes_\chi 1)\mathcal W^*\]
and we have (\cite{E4}, 5.1), for all $x\in M'_0\cap M_2$ :
\[\mathcal W^*(\pi_{\chi_2}(x)\otimes_\chi1)\mathcal W=\varsigma_{M'_0\cap M_1}(id*_{M'_0\cap M_1}\pi_{\chi_2})(U^*j_2\circ j_1(x)U)\]
which shows that, if we write $W=(J_{\chi_2}\hat{J}\otimes_\chi 1)\sigma_{\chi^o}\mathcal W\sigma_{\chi^o}(\hat{J}J_{\chi_2}\otimes_\chi1)$ (where $\hat{J}$ is the standard implementation of $j_{1|M'_0\cap M_2}$ on $H_{\chi_2}$), or, equivalently $\mathcal W=\widehat{W^c}$, we obtain a coproduct $\tilde{\Gamma}$ on $M'_0\cap M_2$ which satisfies 
\[\tilde{\Gamma}(x)=W^*(1\otimes_{\chi^o}x)W=(J_{\chi_2}\hat{J}\otimes_\chi 1)U^*j_2\circ j_1(x)U(\hat{J}J_{\chi_2}\otimes_\chi1)\]
(in which, for simplification, we have identified $M'_0\cap M_2$ with $\pi_{\chi_2}(M'_0\cap M_2)$.)
%%%%%%%strong
\section{Left invariance revisited}
\label{strong}
In the definition of a left-invariant operator-valued weight, recalled in \ref{LW}, we have, for any $x\in\gM_T^+$, $(id\underset{\nu}{_\beta*_\alpha}T)\Gamma (x)=T(x)\underset{N}{_\beta\otimes_\alpha}1$ (or, equivalently, $T(x)=(id\underset{\nu}{_\beta*_\alpha}\Phi)\Gamma(x)$, where $\Phi=\nu\circ\alpha^{-1}\circ T$). In this chapter, we extend this formula to any positive $x$ in $M$, and even, any positive $x$ in a fiber product $\mathcal L(\gH)\underset{N}{_b*_\alpha}M$, where $b$ is a faithful antireprentation of $N$ on $\gH$ (\ref{thinv}(i) and (ii)). On our way, we prove a Plancherel-like theorem (\ref{Plancherel}), which extends the construction of the dual left-invariant operator-valued weight recalled in \ref{thL3}(iv). 
The result \ref{thinv}(ii) will be used in \ref{defovw}, and I am indebted to S. Vaes who pointed in \cite{V2} how this kind of result was necessary in order to get correct proofs.

%%%%defI
\subsection{Definitions and notations}
\label{defI}
Let $N$ be a von Neumann algebra, $\psi$ a normal semi-finite faithful weight on $N$; following (\cite{ES2}, 2.1.6), we define, for all $\omega\in M_*$ :
\[\|\omega\|_\psi=\sup\{|\omega(x^*)|, x\in\gN_\psi, \psi(x^*x)\leq 1\}\]
and $I_\psi=\{\omega\in M_*, \|\omega\|_\psi <\infty\}$. Then $\omega\in M_*$ belongs to $I_\psi$ if and only if there exists a vector $a_\psi(\omega)$ in $H_\psi$ such that :
\[\omega(x^*)=(a_\psi(\omega)|\Lambda_\psi(x))\]
and we have then $\|a_\psi(\omega)\|=\|\omega\|_\psi$. Moreover, the application $a$ here defined is clearly linear from $I_\psi$ into $H_\psi$, and, $\psi$ being semi-finite, injective. 
\newline
Let us identify $N$ and $\pi_\psi(N)$; then, any element $\omega\in M_*$ is of the form $\omega_{\xi, \eta}$, where $\xi$, $\eta$ belong to $H_\psi$. For any $\eta$ in $H_\psi$, let us define the following operator $\pi'(\eta)$, with $\mathcal D(\pi'(\eta))=\Lambda_\psi(\gN_\psi)$, and, for any $x\in\gN_\psi$, by $\pi'(\eta)\Lambda_\psi(x)=x\eta$. We then easily get that $\omega_{\xi, \eta}$ belongs to $I_\psi$ if and only if $\xi$ belongs to $\mathcal D(\pi'(\eta)^*)$, and we have then :
\[a_\psi(\omega_{\xi, \eta})=\pi'(\eta)^*\xi\]
Let us write $K_\eta=\overline{\mathcal D(\pi'(\eta)^*)}$; then $(\pi'(\eta)^*)_{|K_\eta}$ is a densely defined operator from $K_\eta$ into $H_\psi$, and we may consider $((\pi'(\eta)^*)_{|K_\eta})^*(\pi'(\eta)^*)_{|K_\eta}$ which is a positive self-adjoint operator on $K_\eta$; we shall write, for simplicity, $\mathcal P_\eta$ for $((\pi'(\eta)^*)_{|K_\eta})^*(\pi'(\eta)^*)_{|K_\eta}+\infty (1-p_{K_\eta})$, which belongs to the positive extension of $N'$. We have then, for any $\xi$, $\eta$ in $H_\psi$ :
\[\|\omega_{\xi, \eta}\|_\psi^2=<\mathcal P_\eta, \omega_\xi>\]

%%%%lema
\subsection{Lemma}
\label{lema}
{\it  Let $\gG=(N, M, \alpha, \beta, \Gamma, T, T', \nu)$ be a measured quantum groupoid; $\widehat{\gG}$ its dual measured quantum groupoid; let $u$, $v$ in $H_\Phi$ such that $\omega_{u,v}$ belongs to $I_\Phi$, and let $e_n$ the self-adjoint elements introduced in \ref{VRN}. Then, we have, for all $x\in\gN_\Phi$ :
\[(J_{\widehat{\Phi}}a_\Phi(\omega_{u,v})|\Lambda_\Phi(xe_n))=(R(x^*)\overline{e_n\delta^{1/2}}v|u)\]
where $\overline{e_n\delta^{1/2}}$ is the closure of the bounded operator $e_n\delta^{1/2}$. 
\newline
If $v$ belongs to $D(\delta^{1/2})$ and $\omega_{u,v}$ belongs to $I_\Phi$ (for instance, if $v\in\mathcal E$, and any $u\in H_\Phi$), we have : 
\[(J_{\widehat{\Phi}}a_\Phi(\omega_{u,v})|\Lambda_\Phi(x))=(R(x^*)\delta^{1/2}v|u)=<x^*, \omega_{\delta^{1/2}v, u}\circ R>\]
and, therefore, $\omega_{\delta^{1/2}v, u}\circ R$ belongs to $I_\Phi$, and $a_\Phi(\omega_{\delta^{1/2}v, u}\circ R)=J_{\widehat{\Phi}}a_\Phi(\omega_{u,v})$. }
\begin{proof}
We have, using \ref{thL3}(v), \ref{VRN} and \ref{defI} :
\begin{eqnarray*}
(J_{\widehat{\Phi}}a_\Phi(\omega_{u,v})|\Lambda_\Phi(xe_n))
&=&
(J_{\widehat{\Phi}}\Lambda_\Phi(xe_n)|a_\Phi(\omega_{u,v}))\\
&=&(\Lambda_\Phi(R(x^*)\overline{e_n\delta^{1/2}}))|a_\Phi(\omega_{u,v}))\\
&=&(R(x^*)\overline{e_n\delta^{1/2}}v|u)
\end{eqnarray*}
which is the first result. The second result is clear by continuity. \end{proof}

%%%%lemE
\subsection{Lemma}
\label{lemE}
{\it Let $\gG=(N, M, \alpha, \beta, \Gamma, T, T', \nu)$ be a measured quantum groupoid; let $\mathcal E$ be the subspace introduced in \ref{thL3}; then :
\newline
(i) $\mathcal E\subset\{\xi\in D(_\alpha H_\Phi, \nu)\cap \mathcal D(\delta^{1/2}), \delta^{1/2}\xi\in D((H_\Phi)_\beta, \nu^o)\}$
\newline
(ii) $\mathcal E\subset\{\xi\in D((H_\Phi)_\beta, \nu^o)\cap \mathcal D(\delta^{1/2}), \delta^{1/2}\xi\in D(_\alpha H_\Phi, \nu)\}$}

\begin{proof}
Let us recall that $\mathcal E$ is the linear subspace generated by all the elements of the form $J_\Phi\Lambda_\Phi(e_nx)$, where the self-adjoint elements $e_n$ had been introduced in \ref{Vaes}, and $x$ belongs to the subset $\mathcal T$ of elements in $\gN_T\cap\gN_{RTR}\cap\gN_\Phi\cap\gN_{\Phi\circ R}$ which are analytic with respect both to $\Phi$ and $\Phi\circ R$, and such that, for any $z$, $z'$ and $\mathbb{C}$, $\sigma_z^\Phi\circ\sigma_{z'}^{\Phi\circ R}(x)$ belongs to $\gN_T\cap\gN_{RTR}\cap\gN_\Phi\cap\gN_{\Phi\circ R}$. Using then the identifications made in \ref{Vaes}, we have :
\[J_\Phi\Lambda_\Phi(e_nx)=\delta^{-1/2}J_{\Phi\circ R}\Lambda_{\Phi\circ R}(e_nx)\]
and, therefore, $J_\Phi\Lambda_\Phi(e_nx)$ belongs to $D(\delta^{1/2})$ and $\delta^{1/2}J_\Phi\Lambda_\Phi(e_nx)=J_{\Phi\circ R}\Lambda_{\Phi\circ R}(e_nx)$ which belongs to $D((H_\Phi)_\beta, \nu^o)$. So we have (i). 
\newline
On the other hand, for all $t\in\mathbb{R}$, as $\delta^{-t}x\delta^{t}\subset \sigma_{-it}^\Phi\sigma_{it}^{\Phi\circ R}(x)$ is bounded, we get that $e_nx\delta^{t}\subset \overline{e_n\delta^{t}} \overline{\delta^{-t}x\delta^{t}}$ is bounded also, and belongs to $\gN_\Phi$, and that :
\[J_\Phi\Lambda_\Phi(\overline{e_nx\delta^{t}})=\delta^{t}J_\Phi\Lambda_\Phi(\lambda^{-it/2}e_nx)\]
So,  $\delta^{t}J_\Phi\Lambda_\Phi(e_nx)=J_\Phi\Lambda_\Phi(\lambda^{-it/2}\overline{e_nx\delta^{t}})$ which belongs to $D(_\alpha H_\Phi, \nu)$. So, we get (ii). \end{proof}

%%%lemE1
\subsection{Definitions and lemma}
\label{lemE1}
{\it Let $\gG=(N, M, \alpha, \beta, \Gamma, T, T', \nu)$ be a measured quantum groupoid; let $\mathcal T$ and $\mathcal E$ be the subspaces introduced in \ref{thL3}; let us define $\mathcal T_\tau$ as the subspace of $\mathcal T$ made of elements which are analytic with respect to $\tau_t$, and such that, for any $z\in\mathbb{C}$, $\tau_z(x)$ belongs to $\mathcal T$, and $\mathcal E_\tau$ as the linear subspace of $\mathcal E$ generated by the elements of the form $J_\Phi\Lambda_\Phi(e_nx)$, where the self-adjoint elements $e_n$ had been introduced in \ref{Vaes}, and $x$ belongs to the subset $\mathcal T_\tau$. Then :
\newline
(i) $\mathcal E_\tau$ is dense in $H_\Phi$; moreover, for any $\xi$ in $D(_\alpha H_\Phi, \nu)$, there exists $\xi_n$ in $\mathcal E_\tau$ such that $R^{\alpha, \nu}(\xi_n)$ weakly converges to $R^{\alpha, \nu}(\xi)$. 
\newline
(ii) $\mathcal E_\tau\subset\{\xi\in \mathcal D(\hat{\delta}^{1/2}), \hat{\delta}^{1/2}\xi\in D((H_\Phi)_{\hat{\beta}}, \nu^o)\}$.
\newline
(iii) for all $s$, $t$ in $\mathbb{R}$, $\xi\in\mathcal E_\tau$, $\xi$ belongs to $\mathcal D((\delta\Delta_{\hat{\Phi}})^t)$, and $(\delta\Delta_{\hat{\Phi}})^t\xi$ belongs to $\mathcal E_\tau$; moreover, $\xi$ belongs also to $\mathcal D(P^t\delta^s)$, and $P^t\delta^s\xi$ belongs to $\mathcal E_\tau$. 
\newline
(iv) let us define $\hat{\mathcal E}_{\hat{\tau}}$ the linear subspace associated to the dual measured quantum groupoid $\widehat{\gG}$. Then :
\[\hat{\mathcal E}_{\hat{\tau}}\subset \{\xi\in D(_\alpha H_\Phi, \nu)\cap \mathcal D(\delta^{1/2}), \delta^{1/2}\xi\in D((H_\phi)_\beta, \nu^o)\}\]}

\begin{proof}
There exists an invertible positive operator $q$ affiliated to $Z(N)$ such that $\lambda=\alpha(q)=\beta(q)$ (\ref{thL2}(vi)). Let us write $q=\int_0^\infty \mu de_\mu$, and let $f_r=\int_{1/r}^rde_\mu$. 
\newline
Let $x\in\mathcal T$; as the automorphisms groups $\tau_t$, $\sigma_t^\Phi$ and $\sigma_t^{\Phi\circ R}$ are two by two commuting, we can check that, for any $x\in\mathcal T$, the operators :
\[x_{n,r}=\sqrt{\frac{n}{\pi}}\int_{-\infty}^\infty e^{-nt^2}f_r\tau_t(x)dt\]
belong to $\mathcal T_\tau$, the strong limit $lim_{n, r}x_{n, r}=x$, and $lim_{n, r}\Lambda_\Phi(x_{n,r})=\Lambda_\Phi(x)$. From which we get that $\mathcal E_\tau$ is dense in $\mathcal E$, and, therefore, in $H_\Phi$; moreover, as we can prove that all the operators $T(x_{n,r}^*e_p^2x_{n,r})$ are uniformly bounded, we get also that $R^{\alpha, \nu}(J_\Phi\Lambda_\Phi(e_px_{n,r}))$ is strongly converging to $R^{\alpha, \nu}(J_\Phi\Lambda_\Phi(e_px))$, from which we get (i). 
\newline
Let now $x$ be in $\mathcal T_\tau$. Then, using \ref{thL4}(ii), we get that :
\begin{eqnarray*}
\hat{\delta}^{1/2}J_\Phi\Lambda_\Phi(e_nx)
&=&P^{-1/2}J_\Phi\delta^{1/2}J_\Phi\delta^{-1/2}\Delta_\Phi^{-1/2}J_\Phi\Lambda_\Phi(e_nx)\\
&=&P^{-1/2}\lambda^{i/4}J_\Phi\Lambda_\Phi(\sigma_{-i/2}(e_n)\sigma_{-i/2}^{\Phi\circ R}(x))\\
&=&\lambda^{i/4}J_\Phi\Lambda_\Phi(\sigma_{-i/2}(e_n)\tau_{-i/2}(\sigma_{-i/2}^{\Phi\circ R}(x)))
\end{eqnarray*}
We get also that $J_\Phi\hat{\delta}^{1/2}J_\Phi\Lambda_\Phi(e_nx)=\lambda^{-i/4}J_\Phi\Lambda_\Phi(\sigma_{-i/2}\tau_{i/2}\sigma_{i/2}^{\Phi\circ R}(x^*)e_n)$, which belongs to $D(_\alpha H_\Phi, \nu)$, and, therefore, $\hat{\delta}^{1/2}J_\Phi\Lambda_\Phi(e_nx)$ belongs to $D((H_\Phi)_{\hat{\beta}}, \nu^o)$, which is (ii). 
\newline
If $x$ belongs to $\mathcal T_\tau$, and $t\in\mathbb{R}$, we get that :
\[P^tJ_\Phi\Lambda_\Phi(e_nx)=J_\Phi(\lambda^{-t/2}e_n\tau_{-it}(x))\in\mathcal E_\tau\]
and, using \ref{thL2}(vii) :
\[(\delta \Delta_{\widehat{\Phi}})^tJ_\Phi\Lambda_\Phi(e_nx)=J_\Phi\Lambda_\Phi(\lambda^{-it}e_n\sigma_{-it}^\Phi\sigma_{it}^{\Phi\circ R}\tau_{-it}(x))\]
and, therefore, $\mathcal E_\tau\subset D((\delta\Delta_{\widehat{\Phi}})^t)$, and $(\delta\Delta_{\widehat{\Phi}})^t\mathcal E_\tau \subset \mathcal E_\tau$
\newline
Moreover, we have also :
\[\delta^sJ_\Phi\Lambda_\Phi(e_nx)=J_\Phi\Lambda_\Phi(\lambda^{-is/2}e_nx\delta^s)\in\mathcal E_\tau\]
from which we get (iii). 
\newline
If we consider now the dual measured quantum groupoid $\widehat{\gG}$, and the linear subset $\widehat{\mathcal E_{\hat{\tau}}}$ associated with, we get, using (iii), that $\delta^{1/2}\widehat{\mathcal E_{\hat{\tau}}}\subset D((H_\Phi)_\beta, \nu^o)$, which is (iv). \end{proof}

%%%%propW
\subsection{Proposition}
\label{propW}
{\it Let $\gG=(N, M, \alpha, \beta, \Gamma, T, T', \nu)$ be a measured quantum groupoid in the sense of \ref{defMQG}, $\gG^o$ the opposite measured quantum groupoid in the sense of \ref{thL5}(i), $\eta$ in $D(_\alpha H_\Phi, \nu)\cap \mathcal D(\delta^{1/2})$ such that $\delta^{1/2}\eta$ belongs to $D((H_\Phi)_\beta, \nu^o)$, $(e_i)_{i\in I}$ an $(\alpha, \nu)$-orthonormal basis of $H_\Phi$; then, for all $x\in\gN_\Phi$, $\xi\in D(_\alpha H_\Phi, \nu)$, $(id\underset{N}{_\beta*_\alpha}\omega_{\eta, \xi})\Gamma(x)$ belongs to $\gN_\Phi$, and the pseudo-multiplicative $\widehat{W^o}$ satisfies :
\[\widehat{W^o}(\Lambda_\Phi(x)\underset{\nu^o}{_{\hat{\alpha}}\otimes_\beta}\delta^{1/2}\eta)=
\sum_i\Lambda_\Phi[(id\underset{N}{_\beta*_\alpha}\omega_{\eta, e_i})\Gamma(x)]\underset{\nu}{_\beta\otimes_\alpha}e_i\]
\begin{eqnarray*}
\Lambda_\Phi[((id\underset{N}{_\beta*_\alpha}\omega_{\eta, \xi})\Gamma(x)]
&=&
(id*\omega_{\delta^{1/2}\eta, \xi})(\widehat{W^o})\Lambda_\Phi(x)\\
&=&
J_{\widehat{\Phi}}\hat{\pi}(\omega_{J_{\widehat{\Phi}}\xi, J_{\widehat{\Phi}}\delta^{1/2}\eta})^*J_{\widehat{\Phi}}\Lambda_\Phi(x)
\end{eqnarray*}
If $u$ belongs to $D((H_\Phi)_\beta, \nu^o)$ and $v$ belongs to $D(_\alpha H_\Phi, \nu)\cap\mathcal D(\delta^{1/2})$ such that $\delta^{1/2}v$ belongs to $D((H_\Phi)_\beta, \nu^o)$, we get that $(id\underset{N}{_\beta*_\alpha}\omega_{\delta^{1/2}u, v}\circ R)\Gamma(x)$ belongs to $\gN_\Phi$, and : }
\[\Lambda_\Phi((id\underset{N}{_\beta*_\alpha}\omega_{u, \delta^{1/2}v}\circ R)\Gamma(x))=J_{\widehat{\Phi}}\hat{\pi}(\omega_{u,v})^*J_{\widehat{\Phi}}\Lambda_\Phi(x)\]
\begin{proof}
Let's apply the definition of $W^*$ (\ref{thL1}) to the measured quantum groupoid $\gG^o$, with the identification made in \ref{VRN} of $\Lambda_\Phi(xe_n)$ with $\Lambda_{\Phi\circ R}(x\overline{e_n\delta^{-1/2}})$ (where $e_n$ has been defined in \ref{Vaes}) :
\[\widehat{W^o}(\Lambda_\Phi(xe_n)\underset{\nu^o}{_{\hat{\alpha}}\otimes_\beta}\delta^{1/2}\eta)=
\sum_i\Lambda_\Phi[(id\underset{N}{_\beta*_\alpha}\omega_{\eta, e_i})\Gamma(xe_n)]\underset{\nu}{_\beta\otimes_\alpha}e_i\]
from which we get :
\[\Lambda_\Phi[((id\underset{N}{_\beta*_\alpha}\omega_{\eta, \xi})\Gamma(xe_n)]=
(id*\omega_{\delta^{1/2}\eta, \xi})(\widehat{W^o})\Lambda_\Phi(xe_n)\]
As $\Lambda_\Phi$ is $\sigma$-strong $*$-norm closed, we get that $(id\underset{N}{_\beta*_\alpha}\omega_{\eta, \xi})\Gamma(x)$ belongs to $\gN_\Phi$, and that :
\[\Lambda_\Phi[((id\underset{N}{_\beta*_\alpha}\omega_{\eta, \xi})\Gamma(x)]=
(id*\omega_{\delta^{1/2}\eta, \xi})(\widehat{W^o})\Lambda_\Phi(x)\]
Now, we get easily the first formula. The second equality is just an application of \ref{thL3}(ii) and \ref{thL5}(v). Then, writing $u=J_{\widehat{\Phi}}\delta^{1/2}\eta$, and $v=J_{\widehat{\Phi}}\xi$, and using \ref{thL3}(v) and \ref{thL2}(vi), we get the last formula. \end{proof}

%%%%%propleftw
\subsection{Proposition}
\label{propleftw}
{\it  Let $\gG=(N, M, \alpha, \beta, \Gamma, T, T', \nu)$ be a measured quantum groupoid; $\widehat{\gG}$ its dual measured quantum groupoid; let $\xi$ in $D((H_\Phi)_\beta, \nu^o)$, $\eta$ in $D(_\alpha H_\Phi, \nu)$ such that $\omega_{\xi, \eta}$ belongs to $I_\Phi$; then, the operator $\hat{\pi}(\omega_{\xi, \eta})$ belongs to $\gN_{\widehat{\Phi}}$, and we have then :}
\[\Lambda_{\widehat{\Phi}}(\hat{\pi}(\omega_{\xi, \eta}))=a(\omega_{\xi, \eta})\]

\begin{proof}
Let $\xi\in D(_\alpha H_\Phi, \nu)\cap D((H_\Phi)_\beta, \nu^o)$, $\eta\in D(_\alpha H_\Phi, \nu)$ such that $\omega_{\xi, \eta}$ belongs to $I_\Phi$; let $u\in D(_\alpha H_\Phi, \nu)$, $v\in \mathcal D(\delta^{1/2})$  such that $\delta^{1/2}v$ is in $D((H_\Phi)_\beta, \nu)$ and $\omega_{u,v}$ is in $I_\Phi$; we have then, for any $x\in\gN_\Phi$, using successively \ref{thL3}(ii), \ref{lema} and \ref{propW} :
\begin{eqnarray*}
(\hat{\pi}(\omega_{\xi, \eta})J_{\widehat{\Phi}}a_\Phi(\omega_{u,v})|\Lambda_\Phi(x))
&=&(J_{\widehat{\Phi}}a_\Phi(\omega_{u,v})|\hat{\pi}(\omega_{\xi, \eta})^*\Lambda_\Phi(x))\\
&=&(J_{\widehat{\Phi}}a_\Phi(\omega_{u,v})|\Lambda_\Phi((\omega_{\eta, \xi}\underset{N}{_\beta*_\alpha}id)\Gamma(x))\\
&=&(a_\Phi(\omega_{\delta^{1/2}v,u}\circ R)|\Lambda_\Phi((\omega_{\eta, \xi}\underset{N}{_\beta*_\alpha}id)\Gamma(x))\\
&=&<[(\omega_{\eta, \xi}\underset{N}{_\beta*_\alpha}id)\Gamma(x)]^*, \omega_{\delta^{1/2}v,u}\circ R>\\
&=&<(\omega_{\xi, \eta}\underset{N}{_\beta*_\alpha}id)\Gamma(x^*), \omega_{\delta^{1/2}v,u}\circ R>\\
&=&<[(id\underset{N}{_\beta*_\alpha}\omega_{u, \delta^{1/2}v}\circ R)\Gamma(x)]^*, \omega_{\xi, \eta}>\\
&=&(a_\Phi(\omega_{\xi, \eta})|\Lambda_\Phi((id\underset{N}{_\beta*_\alpha}\omega_{u, \delta^{1/2}v}\circ R)\Gamma(x))\\
&=&(a_\Phi(\omega_{\xi, \eta})|J_{\widehat{\Phi}}\hat{\pi}(\omega_{u,v})^*J_{\widehat{\Phi}}\Lambda_\Phi(x))\\
&=&(J_{\widehat{\Phi}}\hat{\pi}(\omega_{u,v})J_{\widehat{\Phi}}a_\Phi(\omega_{\xi, \eta})|\Lambda_\Phi(x))
\end{eqnarray*}
From which we get, by density :
\[\hat{\pi}(\omega_{\xi, \eta})J_{\widehat{\Phi}}a_\Phi(\omega_{u,v})=J_{\widehat{\Phi}}\hat{\pi}(\omega_{u,v})J_{\widehat{\Phi}}a_\Phi(\omega_{\xi, \eta})\]
Let us take $v=J_\Phi\Lambda_\Phi(e_nx)$, with $e_n$ defined in \ref{VRN} and $x\in \mathcal T$ (\ref{thL3}). Then, this formula gives :
\[\hat{\pi}(\omega_{\xi, \eta})J_{\widehat{\Phi}}J_\Phi (e_nx)^*J_\Phi u=J_{\widehat{\Phi}}\hat{\pi}(\omega_{u,J_\Phi\Lambda_\Phi(e_nx)})J_{\widehat{\Phi}}a_\Phi(\omega_{\xi, \eta})\]
which, by continuity, gives, for any $u\in D(_\alpha H_\Phi, \nu)$, $v\in J_\Phi\Lambda_\Phi(\gN_\Phi)\cap D((H_\Phi)_\beta, \nu^o)$ :
\[\hat{\pi}(\omega_{\xi, \eta})J_{\widehat{\Phi}}\pi'(v)^* u=J_{\widehat{\Phi}}\hat{\pi}(\omega_{u,v})J_{\widehat{\Phi}}a_\Phi(\omega_{\xi, \eta})\]
Let us suppose now that $u$, $v$ belong to $D(_\alpha H_\Phi, \nu)\cap D((H_\Phi)_\beta, \nu^o)$ and are such that $\omega_{u,v}$ belongs to $I_\Phi$; with the notations of \ref{defI}, let us write :
\[\mathcal P_v=\int_0^\infty\lambda de_\lambda +(1-p)\infty\]
As $\omega_{u,v}$ belongs to $I_\Phi$, we get that $pu=u$. Let us define $v_\mu=e_\mu v$ which belongs to $J_\Phi\Lambda_\Phi(\gN_\Phi)\cap D((H_\Phi)_\beta, \nu^o)$. We have :
\[lim_\mu \pi'(v_\mu)^*u=\lim_\mu \pi'(v)^*e_\mu u=\pi'(v)^*pu=\pi'(v)^*u=a_\Phi(\omega_{u,v})\]
and we get also, for the weak limits :
\[lim_\mu R^{\beta, \nu^o}(v_\mu)=lim_\mu e_\mu R^{\beta, \nu^o}(v)=pR^{\beta, \nu^o}(v)=R^{\beta, \nu^o}(pv)\]
and, as $p\in M'$, we get for the weak limits :
\begin{multline*}
lim_\mu \hat{\pi}(\omega_{u, v_\mu})=lim_\mu (\omega_{u, v_\mu}*id)(W)=(\omega_{u, pv}*id)(W)=\\=(\omega_{pu, v}*id)(W)=(\omega_{u,v}*id)(W)=\hat{\pi}(\omega_{u,v})
\end{multline*}
and, taking the limits of the equality :
\[\hat{\pi}(\omega_{\xi, \eta})J_{\widehat{\Phi}}\pi'(v_\mu)^* u=J_{\widehat{\Phi}}\hat{\pi}(\omega_{u,v_\mu})J_{\widehat{\Phi}}a_\Phi(\omega_{\xi, \eta})\]
we get :
\[\hat{\pi}(\omega_{\xi, \eta})J_{\widehat{\Phi}}a_\Phi(\omega_{u,v})=J_{\widehat{\Phi}}\hat{\pi}(\omega_{u,v})J_{\widehat{\Phi}}a_\Phi(\omega_{\xi, \eta})\]
for all $u$, $v$ in $D(_\alpha H_\Phi, \nu)\cap D((H_\Phi)_\beta, \nu^o)$ and are such that $\omega_{u,v}$ belongs to $I_\Phi$; therefore, by linearity, we get, for all $\omega\in M_*^{\alpha, \beta}\cap I_\Phi$, that  :
\[\hat{\pi}(\omega_{\xi, \eta})J_{\widehat{\Phi}}a(\omega)=J_{\widehat{\Phi}}\hat{\pi}(\omega)J_{\widehat{\Phi}}a(\omega_{\xi, \eta})=J_{\widehat{\Phi}}\hat{\pi}(\omega)J_{\widehat{\Phi}}\pi'(\eta)^*\xi\]
By the closedness of $\pi'(\eta)^*_{|K_\eta}$ (\ref{defI}), this formula remains true for $\xi$ in $D((H_\Phi)_\beta, \nu^o)$ such that $\omega_{\xi, \eta}$ belongs to $I_\Phi$.  From which we get that $a(\omega_{\xi, \eta})$ is left-bounded with respect to $\widehat{\Phi}$, and the result, by standard arguments of Tomita-Takesaki's theory. \end{proof}

%%%%thleftw
\subsection{Proposition}
\label{thleftw}
{\it Let $\gG=(N, M, \alpha, \beta, \Gamma, T, T', \nu)$ be a measured quantum groupoid; $\widehat{\gG}$ its dual measured quantum groupoid; let's take $\xi$ in $D((H_\Phi)_\beta, \nu^o)$ and $\eta$ in $D(_\alpha H_\Phi, \nu)$, and let us suppose that $\hat{\pi}(\omega_{\xi, \eta})$ belongs to $\gN_{\widehat{\Phi}}$. Then $\omega_{\xi, \eta}$ belongs to $I_\Phi$ and : }
\[\Lambda_{\widehat{\Phi}}(\hat{\pi}(\omega_{\xi, \eta}))=a_\Phi(\omega_{\xi, \eta})\]
\begin{proof}
Let $u$, $v$ be in $D(_\alpha H_\Phi, \nu)\cap D((H_\Phi)_{\hat{\beta}}, \nu^o)$; then the linear form $\hat{\omega}_{u,v}$ on $\widehat{M}$ defined, for $y$ in $\widehat{M}$ by $\hat{\omega}_{u,v}(y)=(yu|v)$ belongs to $\widehat{M}_*^{\alpha, \hat{\beta}}$ (with the notations of \ref{thL3}(i) applied to $\widehat{\gG}$). Let us suppose that $\hat{\omega}_{u,v}$ belongs to $I_{\widehat{\Phi}}$. Using  \ref{thL3}(v) and (ii), applied to $\widehat{\gG}$, we get :
\begin{eqnarray*}
(a_{\widehat{\Phi}}(\hat{\omega}_{u,v})|\Lambda_{\widehat{\Phi}}(\hat{\pi}(\omega_{\xi, \eta}))
&=& <\hat{\pi}(\omega_{\xi, \eta})^*, \hat{\omega}_{u,v}>\\
&=& ((\omega_{\xi, \eta}*id)(W)^*u|v)\\
&=&(\eta\underset{\nu^o}{_\alpha\otimes_{\hat{\beta}}}u|W(\xi\underset{\nu}{_\beta\otimes_\alpha}v))\\
&=&(\widehat{W}(u\underset{\nu}{_{\hat{\beta}}\otimes_\alpha}\eta)|v\underset{\nu^o}{_\alpha\otimes_\beta}\xi)\\
&=&(\hat{\hat{\pi}}(\hat{\omega}_{u,v})\eta|\xi)
\end{eqnarray*}
and, by linearity, we obtain, for all $\hat{\omega}$ in $\widehat{M}_*^{\alpha, \hat{\beta}}\cap I_{\widehat{\Phi}}$ :
\[(\Lambda_{\widehat{\Phi}}(\hat{\pi}(\omega_{\xi, \eta}))|a_{\widehat{\Phi}}(\hat{\omega}))=
<\hat{\hat{\pi}}(\hat{\omega})^*, \omega_{\xi, \eta}>\] 
But, by  \ref{thL3}(v) applied to $\widehat{\gG}$, and \ref{thL4}(i), we have $a_{\widehat{\Phi}}(\hat{\omega})=\Lambda_\Phi(\hat{\hat{\pi}}(\hat{\omega}))$, and, by \ref{thL3}(v) applied to $\widehat{\gG}$ and again \ref{thL4}(i), we know that the algebra $\hat{\hat{\pi}}(\widehat{M}_*^{\alpha, \hat{\beta}}\cap I_{\widehat{\Phi}})$ is a core for $\Lambda_\Phi$; so, by the closedness of $\Lambda_\Phi$, we get, for all $x\in\gN_\Phi$ :
\[(\Lambda_{\widehat{\Phi}}(\hat{\pi}(\omega_{\xi, \eta}))|\Lambda_\Phi(x))=
<x^*, \omega_{\xi, \eta}>\]
from which we get that $|<x^*, \omega_{\xi, \eta}>|\leq\|\Lambda_{\widehat{\Phi}}(\hat{\pi}(\omega_{\xi, \eta}))\|\|\Lambda_\phi(x)\|$, and that $\omega_{\xi, \eta}$ belongs to $I_\Phi$. Then, using \ref{thL3}(v), the proof is finished. \end{proof}

%%%%Plancherel
\subsection{Theorem}
\label{Plancherel}
{\it  Let $\gG=(N, M, \alpha, \beta, \Gamma, T, T', \nu)$ be a measured quantum groupoid; $\widehat{\gG}$ its dual measured quantum groupoid; let's take $\xi$ in $D((H_\Phi)_\beta, \nu^o)$ and $\eta$ in $D(_\alpha H_\Phi, \nu)$; then, are equivalent :
\newline
(i) $\hat{\pi}(\omega_{\xi, \eta})$ belongs to $\gN_{\widehat{\Phi}}$;
\newline
(ii) $\omega_{\xi, \eta}$ belongs to $I_\Phi$. 
\newline
We have then $\Lambda_{\widehat{\Phi}}(\hat{\pi}(\omega_{\xi, \eta}))=a_\Phi(\omega_{\xi, \eta})$. }

\begin{proof}
This is just \ref{propleftw} and \ref{thleftw}. \end{proof}

%%%%%S
\subsection{Proposition}
\label{S}
{\it Let $\gG=(N, M, \alpha, \beta, \Gamma, T, T', \nu)$ be a measured quantum groupoid; $\widehat{\gG}$ its dual measured quantum groupoid; for any $X\in\hat{\alpha}(N)'$, let us define :
\[\mathcal S(X)=(id\underset{N}{_{\hat{\beta}}*_\alpha}\widehat{\Phi})(W^c(X\underset{\nu^o}{_{\hat{\alpha}}\otimes_{\hat{\beta}}}1)W^{c*})\]
Then :
\newline
(i) $\mathcal S(X)$ belongs to the positive extended part of $\hat{\alpha}(N)'$, and,  for any $\xi\in D((H_\Phi)_\beta, \nu^o)$, $\eta\in D(_\alpha H_\Phi, \nu)$, we have :
\[<\mathcal S(\theta^{\hat{\alpha}}(J_{\widehat{\Phi}}J_\Phi\eta, J_{\widehat{\Phi}}J_\Phi\eta)),\omega_{J_{\widehat{\Phi}}J_\Phi\xi}>=\widehat{\Phi}(\hat{\pi}(\omega_{\xi, \eta})^*\hat{\pi}(\omega_{\xi, \eta}))=<\mathcal P_\eta, \omega_\xi>\]
(i) for all $\eta\in D(_\alpha H_\Phi, \nu)$, we have :
\[\mathcal S(\theta^{\hat{\alpha}}(J_{\widehat{\Phi}}J_\Phi\eta, J_{\widehat{\Phi}}J_\Phi\eta))=J_{\widehat{\Phi}}J_\Phi\mathcal P_\eta J_\Phi J_{\widehat{\Phi}}\]
(ii) $\mathcal S(X)$ belongs to the positive extended part of $M$. }

\begin{proof}
Using \ref{thL5}(v), we get that $W^{c*}=(J_\Phi J_{\widehat{\Phi}}\underset{\nu^o}{_\alpha\otimes_{\hat{\beta}}} 1)W (J_{\widehat{\Phi}}J_\Phi\underset{\nu}{_{\hat{\beta}}\otimes_\alpha} 1)$; so, with the intertwining properties of $W$, we easily get that $S(X)$ is an element of the extended positive part of 
$\hat{\alpha}(N)'$. 
\newline
Let us notice that $J_{\widehat{\Phi}}J_\Phi\xi$ belongs to $D((H_\Phi)_{\hat{\beta}}, \nu^o)$, and 
$J_{\widehat{\Phi}}J_\Phi\eta$ to $D(_{\hat{\alpha}}H_\Phi, \nu)$. We then get that $<\mathcal S(\theta^{\hat{\alpha}}(J_{\widehat{\Phi}}J_\Phi\eta, J_{\widehat{\Phi}}J_\Phi\eta)),\omega_{J_{\widehat{\Phi}}J_\Phi\xi}>$ is equal to :
\begin{multline*}
\widehat{\Phi}((\omega_{J_{\widehat{\Phi}}J_\Phi\xi, J_{\widehat{\Phi}}J_\Phi\eta}*id)(W^{c*})^*(\omega_{J_{\widehat{\Phi}}J_\Phi\xi, J_{\widehat{\Phi}}J_\Phi\eta}*id)(W^{c*}))=\\
=\widehat{\Phi}((\omega_{\xi, \eta}*id)(W)^*(\omega_{\xi, \eta}*id)(W))
=\widehat{\Phi}(\hat{\pi}(\omega_{\xi, \eta})^*\hat{\pi}(\omega_{\xi, \eta}))
\end{multline*}
If it is finite, it is equal to $\|\Lambda_{\widehat{\Phi}}(\hat{\pi}(\omega_{\xi, \eta})\|^2=\|a_\Phi(\omega_{\xi, \eta})\|^2=\|\omega_{\xi, \eta}\|_\Phi^2$, thanks to \ref{Plancherel}. If it is infinite, using \ref{Plancherel} again, $\|\omega_{\xi, \eta}\|_\Phi$ also is infinite. So, in both cases, it is equal to $\|\omega_{\xi, \eta}\|_\Phi^2$, which is equal, using \ref{defI}, to $<\mathcal P_\eta, \omega_\xi>$, which finishes the proof of (i). Then (ii) (resp. (iii)) is a straightforward corollary of (i) (resp. (ii)). \end{proof}

%%%%%stronginv
\subsection{Proposition}
\label{stronginv}
{\it  Let $\gG=(N, M, \alpha, \beta, \Gamma, T, T', \nu)$ be a measured quantum groupo-id; $\widehat{\gG}$ its dual measured quantum groupoid. Let $X$ be a positive element in $\widehat{M}$; then, we have: }
\[(id\underset{\nu}{_{\hat{\beta}}*_\alpha}\widehat{\Phi})\widehat{\Gamma}(X)=\widehat{T}(X)\]

\begin{proof}
Using \ref{thL5}(v), we get that :
\[\widehat{\Gamma}(X)=W^c(X\underset{N^o}{_{\hat{\alpha}}\otimes_{\hat{\beta}}}1)W^{c*}\]
and, therefore, with the notations of \ref{S}, $(id\underset{\nu}{_{\hat{\beta}}*_\alpha}\widehat{\Phi})\widehat{\Gamma}(X)=\mathcal S(X)$. Using \ref{S}(iii), it belongs to the positive extended part of $M$, and, by construction, it belongs to the positive extended part of $\widehat{M}$. So, it belongs to the positive extended part of $\alpha (N)$. If $u$ is a unitary in $N$, we easily get that $\mathcal S(\alpha(u^*)X\alpha(u))=\alpha(u^*)\mathcal S(X)\alpha(u)$, and, therefore, the application $X\rightarrow \mathcal S(X)$ is an operator-valued weight from $\widehat{M}$ on $\alpha(N)$. It is clearly normal and faithful, and, as $\mathcal S(X)=\widehat{T}(X)$ for any $X\in\gM_{\widehat{T}}^+$, we get it is also semi-finite, and that, for any positive $X$ in $\widehat{M}$, we have $\mathcal S(X)\leq\widehat{T}(X) $
\newline
Let us define $\tilde{\Phi}(X)=\nu\circ\alpha^{-1}\circ \mathcal S(X)$; we have $\tilde{\Phi}\leq\widehat{\Phi}$, and these two weights are equal on $\gM_{\widehat{T}}$; moreover, as, for all $t\in\mathbb{R}$, we have :\[\widehat{\Gamma}\sigma_t^{\widehat{\Phi}}=(\hat{\tau_t}\underset{N}{_{\hat{\beta}}*_\alpha}\sigma_t^{\widehat{\Phi}})\widehat{\Gamma}\]
we get that $\mathcal S\sigma_t^{\widehat{\Phi}}=\hat{\tau_t}\mathcal S$, and that $\tilde{\Phi}$ is invariant under $\sigma_t^{\widehat{\Phi}}$, by \ref{thL2}(ii) applied to $\widehat{\gG}$; therefore, we have $\tilde{\Phi}=\widehat{\Phi}$, from which we deduce that $\mathcal S=\widehat{T}$, which gives the result. \end{proof}

%%%%%%vstronginv
\subsection{Proposition}
\label{vstronginv}
{\it  Let $\gG=(N, M, \alpha, \beta, \Gamma, T, T', \nu)$ be a measured quantum groupo-id; $\widehat{\gG}$ its dual measured quantum groupoid, $\gH$ an Hilbert space with $b$ a faithful normal anti-representation of $N$ on $\gH$. Let $X$ be a positive element in $\mathcal L(\gH)\underset{N}{_b*_\alpha}\widehat{M}$; then, we have:}
\[(id\underset{N}{_b*_\alpha}id\underset{\nu}{_{\hat{\beta}}*_\alpha}\widehat{\Phi})(id\underset{N}{_b*_\alpha}\widehat{\Gamma})(X)=(id\underset{N}{_b*_\alpha}\widehat{\Phi})(X)\underset{N}{_b\otimes_\alpha}1\]

\begin{proof}
Using \ref{thL5}(v), we get that, for $X\in \mathcal L(\gH)\underset{N}{_b*_\alpha}\widehat{M}$ :
\[(id\underset{N}{_b*_\alpha}\widehat{\Gamma})(X)=(1\underset{N}{_b\otimes_\alpha}W^c)(X\underset{N^o}{_{\hat{\alpha}}\otimes_{\hat{\beta}}}1)(1\underset{N}{_b\otimes_\alpha}W^c)^*\]
For any $Y$ in $\mathcal L(\gH\underset{\nu}{_b\otimes_\alpha}H_\Phi)^+$, commuting with $1\underset{N}{_b\otimes_\alpha}\hat{\alpha}(N)$, let us write 
\[\tilde{S}(Y)=(id\underset{N}{_b*_\alpha}id\underset{\nu}{_{\hat{\beta}}*_\alpha}\widehat{\Phi})[(1\underset{N}{_b\otimes_\alpha}W^c)(Y\underset{N^o}{_{\hat{\alpha}}\otimes_{\hat{\beta}}}1)(1\underset{N}{_b\otimes_\alpha}W^c)^*]\]
which is an element of the positive extended part of $\mathcal L(\gH\underset{\nu}{_b\otimes_\alpha}H_\Phi)$. 
\newline
Let's take  $\xi$ in $D((\gH\underset{\nu}{_b\otimes_\alpha}H_\Phi)_{1\underset{N}{_b\otimes_\alpha}\hat{\beta}}, \nu^o)$ and $\eta$ in $D(_{1\underset{N}{_b\otimes_\alpha}\hat{\alpha}}(\gH\underset{\nu}{_b\otimes_\alpha}H_\Phi), \nu)$. Then, we get that $<S(\theta^{1\underset{N}{_b\otimes_\alpha}\hat{\alpha}}(\eta, \eta)), \omega_\xi>$ is equal to :
\[\widehat{\Phi}[(\omega_\xi\underset{\nu}{_{\hat{\beta}}*_\alpha}id)((1\underset{N}{_b\otimes_\alpha}W^c)(\theta^{1\underset{N}{_b\otimes_\alpha}\hat{\alpha}}(\eta, \eta)\underset{N^o}{_{\hat{\alpha}}\otimes_{\hat{\beta}}}1)(1\underset{N}{_b\otimes_\alpha}W^c)^*)]\]
which can be written also :
\[\widehat{\Phi}[((\omega_{\xi, \eta}*id)((1\underset{N}{_b\otimes_\alpha}W^{c*}))^*(\omega_{\xi, \eta}*id)((1\underset{N}{_b\otimes_\alpha}W^{c*})]\]
Using now \ref{thL5}(v), we get that $W^{c*}=(J_\Phi J_{\widehat{\Phi}}\underset{\nu^o}{_\alpha\otimes_{\hat{\beta}}} 1)W (J_{\widehat{\Phi}}J_\Phi\underset{\nu}{_{\hat{\beta}}\otimes_\alpha} 1)$ and, therefore, $<S(\theta^{1\underset{N}{_b\otimes_\alpha}\hat{\alpha}}(\eta, \eta)), \omega_\xi>$ is equal to 
$\widehat{\Phi}(A^*A)$, where :
\[A=(\omega_{(1\underset{N}{_b\otimes_\alpha}J_{\widehat{\Phi}}J_\Phi)\xi, (1\underset{N}{_b\otimes_\alpha}J_\Phi J_{\widehat{\Phi}})\eta}*id)(1\underset{N^o}{_b\otimes_\alpha}W)\]
Let $(e_i)_{i\in I}$ be a $(b, \nu^o)$-orthogonal basis of $\gH$. Using (\cite{E3}, 2.3(ii)), there exist $\xi_i$ in $D((H_\Phi)_{\hat{\beta}}, \nu^o)$ and $\eta_i$ in $D(_{\hat{\alpha}}H_\Phi, \nu)$ such that :
\[\xi=\sum_ie_i\underset{\nu}{_b\otimes_\alpha}\xi_i\]
\[\eta=\sum_ie_i\underset{\nu}{_b\otimes_\alpha}\eta_i\]
and we get that $A=\sum_i\alpha(<e_i, e_i>_{b, \nu^o})(\omega_{J_{\widehat{\Phi}}J_\Phi\xi_i, J_{\widehat{\Phi}}J_\Phi\eta_i}*id)(W)$
because $J_{\widehat{\Phi}}J_\Phi  \xi_i$ belongs to $D((H_\Phi)_\beta, \nu^o)$, and $J_{\widehat{\Phi}}J_\Phi\eta_i$ belongs to $D(_\alpha H_\Phi, \nu)$;  so, $<S(\theta^{1\underset{N}{_b\otimes_\alpha}\hat{\alpha}}(\eta, \eta)), \omega_\xi>$ is equal to :
\[
\sum_{i,j}\widehat{\Phi}(\hat{\pi}(\omega_{J_{\widehat{\Phi}}J_\Phi\xi_i, J_{\widehat{\Phi}}J_\Phi\eta_i})^*\alpha(<e_i, e_i>_{b, \nu^o}<e_j, e_j>_{b, \nu^o})\hat{\pi}(\omega_{J_{\widehat{\Phi}}J_\Phi\xi_j, J_{\widehat{\Phi}}J_\Phi\eta_j}))\]
If it is bounded, then so are, for all $i\in I$ :
\[\widehat{\Phi}(\hat{\pi}(\omega_{J_{\widehat{\Phi}}J_\Phi\xi_i, J_{\widehat{\Phi}}J_\Phi\eta_i})^*\alpha(<e_i, e_i>_{b, \nu^o})\hat{\pi}(\omega_{J_{\widehat{\Phi}}J_\Phi\xi_i, J_{\widehat{\Phi}}J_\Phi\eta_i}))\]
which implies (\ref{thleftw}) that, for all $i\in I$, $\omega_{J_{\widehat{\Phi}}J_\Phi\xi_i, J_{\widehat{\Phi}}J_\Phi\eta_i}$ belongs to $I_\Phi$, and, therefore, using \ref{defI}, that $J_{\widehat{\Phi}}J_\Phi\xi_i$ belongs to $D(\pi'(J_{\widehat{\Phi}}J_\Phi\eta_i)^*)$; moreover, we then get :
\begin{eqnarray*}
<S(\theta^{1\underset{N}{_b\otimes_\alpha}\hat{\alpha}}(\eta, \eta)), \omega_\xi>
&=&
\|\sum_i\alpha(<e_i, e_i>_{b, \nu^o})\pi'(J_{\widehat{\Phi}}J_\Phi\eta_i)^*J_{\widehat{\Phi}}J_\Phi\xi_i\|^2\\
&=&\|\sum_i e_i\underset{\nu}{_b\otimes_\alpha}\pi'(J_{\widehat{\Phi}}J_\Phi\eta_i)^*J_{\widehat{\Phi}}J_\Phi\xi_i\|^2
\end{eqnarray*}
By definition of the extended positive part, there exists a closed subspace $K\subset \gH\underset{\nu}{_b\otimes_\alpha}H_\Phi$ on which $S(\theta^{1\underset{N}{_b\otimes_\alpha}\hat{\alpha}}(\eta, \eta))$ is a positive self-adjoint operator, and it is $\infty$ on the orthogonal of $K$.
\newline
So, $K=\oplus_i K_i$, where the $K_i$ are the closed orthogonal subspaces $e_i\underset{\nu}{_b\otimes_\alpha}J_\Phi J_{\widehat{\Phi}}K_{J_{\widehat{\Phi}}J_\Phi\eta_i}$ (with the notations of \ref{defI}). 
\newline
Let $u$ be a unitary in $M'$; we get that each $K_i$ is invariant under $1\underset{N}{_b\otimes_\alpha}u$, and so, we get that $p_K$ belongs to $\mathcal L(\gH)\underset{N}{_b*_\alpha}M$; as $(1\underset{N}{_b\otimes_\alpha}u)\xi$ belongs to $D((\gH\underset{\nu}{_b\otimes_\alpha}H_\Phi)_{1\underset{N}{_b\otimes_\alpha}\hat{\beta}}, \nu^o)$, we get that :
\begin{eqnarray*}
<S(\theta^{1\underset{N}{_b\otimes_\alpha}\hat{\alpha}}(\eta, \eta)), \omega_{u\xi}>
&=&\|\sum_i e_i\underset{\nu}{_b\otimes_\alpha}\pi'(J_{\widehat{\Phi}}J_\Phi\eta_i)^*J_{\widehat{\Phi}}J_\Phi u\xi_i\|^2\\
&=&\|\sum_i e_i\underset{\nu}{_b\otimes_\alpha}R(J_\Phi u^*J_\Phi)\pi'(J_{\widehat{\Phi}}J_\Phi\eta_i)^*J_{\widehat{\Phi}}J_\Phi \xi_i\|^2\\
&=&<S(\theta^{1\underset{N}{_b\otimes_\alpha}\hat{\alpha}}(\eta, \eta)), \omega_\xi>
\end{eqnarray*}
So, we get that $(1\underset{N}{_b\otimes_\alpha}u^*)S(\theta^{1\underset{N}{_b\otimes_\alpha}\hat{\alpha}}(\eta, \eta))(1\underset{N}{_b\otimes_\alpha}u)$ is equal to $\theta^{1\underset{N}{_b\otimes_\alpha}\hat{\alpha}}(\eta, \eta)$ from which we get that $S(\theta^{1\underset{N}{_b\otimes_\alpha}\hat{\alpha}}(\eta, \eta))$ belongs to the positive extended part of $\mathcal L(\gH)\underset{N}{_b*_\alpha}M$. 
\newline
Let $(\eta_j)_{j\in J}$ be an $((1\underset{N}{_b\otimes_\alpha}\hat{\alpha}), \nu)$ orthogonal basis of $\gH\underset{\nu}{_b\otimes_\alpha}H_\Phi$; for any $Y\in \mathcal L(\gH\underset{\nu}{_b\otimes_\alpha}H_\Phi)^+$, commuting with $1\underset{N}{_b\otimes_\alpha}\hat{\alpha}(N)$. We have :
\[Y=\sum_j\theta^{1\underset{N}{_b\otimes_\alpha}\hat{\alpha}}(Y^{1/2}\eta_j, Y^{1/2}\eta_j)\]
and, therefore :
\[S(Y)=\sum_j S(\theta^{1\underset{N}{_b\otimes_\alpha}\hat{\alpha}}(Y^{1/2}\eta_j, Y^{1/2}\eta_j))\]
from which we get that $S(Y)$ belongs to the positive extended part of $\mathcal L(\gH)\underset{N}{_b*_\alpha}M$. 
\newline
In particular, if $X$ is a positive element in $\mathcal L(\gH)\underset{N}{_b*_\alpha}\widehat{M}$; then 
\[(id\underset{N}{_b*_\alpha}id\underset{\nu}{_{\hat{\beta}}*_\alpha}\widehat{\Phi})(id\underset{N}{_b*_\alpha}\widehat{\Gamma})(X)=S(X)\]
belongs to the positive extended part of $\mathcal L(\gH)\underset{N}{_b*_\alpha}M$.  As, by definition,  it belongs to the positive extended part of $\mathcal L(\gH)\underset{N}{_b*_\alpha}\widehat{M}$, we get that it belongs to the positive extended part of $\mathcal L(\gH)\underset{N}{_b*_\alpha}\alpha(N)$, from which we get that there exists $X'$ in the positive extended part of $b(N)'$ such that :
\[(id\underset{N}{_b*_\alpha}id\underset{\nu}{_{\hat{\beta}}*_\alpha}\widehat{\Phi})(id\underset{N}{_b*_\alpha}\widehat{\Gamma})(X)=X'\underset{N}{_b\otimes_\alpha}1\]
Let $\xi$ be in $D(\gH_b, \nu^o)$, $\eta$ in $D(_\alpha H_\Phi, \nu)\cap D((H_\Phi)_{\hat{\beta}}, \nu^o)$; we get :
\begin{eqnarray*}
<b(<\eta, \eta>_{\alpha, \nu})X',\omega_\xi>
&=&(\omega_\xi\underset{N}{_b*_\alpha}\omega_\eta\underset{\nu}{_{\hat{\beta}}*_\alpha}\widehat{\Phi})(id\underset{N}{_b*_\alpha}\widehat{\Gamma})(X)\\
&=&(\omega_\eta\underset{\nu}{_{\hat{\beta}}*_\alpha}\widehat{\Phi})\widehat{\Gamma}[(\omega_\xi\underset{N}{_b*_\alpha}id)(X)]
\end{eqnarray*}
which, using \ref{stronginv}, is equal to :
\begin{eqnarray*}
<\widehat{T}(\omega_\xi\underset{N}{_b*_\alpha}id)(X), \omega_\eta>
&=&(\omega_\xi\underset{N}{_b*_\alpha}\omega_\eta)(id\underset{N}{_b*_\alpha}\widehat{T})(X)\\
&=&<b(<\eta, \eta>_{\alpha, \nu})(id\underset{\nu}{_b*_\alpha}\widehat{\Phi})(X), \omega_\xi>
\end{eqnarray*}
from which we infer that $X'=(id\underset{\nu}{_b*_\alpha}\widehat{\Phi})(X)$, which gives the result. \end{proof}

%%%%%%%%%%thinv
\subsection{Theorem}
\label{thinv}
{\it Let $\gG=(N, M, \alpha, \beta, \Gamma, T, T', \nu)$ be a measured quantum groupoid; let $\gH$ be an Hilbert space with $b$ a faithful normal anti-representation of $N$ on $\gH$. Let $X$ be a positive element in $M$, $\tilde{X}$ a positive element in $\mathcal L(\gH)\underset{N}{_b*_\alpha}M$; then :
\newline
(i) we have $(id\underset{N}{_\beta*_\alpha}T)\Gamma(X)=T(X)\underset{N}{_\beta\otimes_\alpha}1$;
\newline
(ii) we have $(id\underset{N}{_b*_\alpha}id\underset{N}{_\beta*_\alpha}T)(id\underset{N}{_b*_\alpha}\Gamma)(\tilde{X})=(id\underset{N}{_b*_\alpha}T)(\tilde{X})\underset{N}{_\beta\otimes_\alpha}1$}

\begin{proof}
Let us apply \ref{stronginv} and \ref{vstronginv} to the dual measured quantum groupoid $\widehat{\gG}$ . \end{proof}

 %%%%%corepresentation
 \section{Corepresentations of measured quantum groupoids}
 \label{corep}
 In this chapter, we define corepresentations of measured quantum groupoids (\ref{defcorep}), and we prove a fundamental property of these with respect to the antipod $S$ (\ref{thS1}). 
 
 %%%%defcorep
 \subsection{Definition}
 \label{defcorep}
 Let $\gG=(N, M, \alpha, \beta, \Gamma, T, T', \nu)$ be a measured quantum groupoid; we shall use all the notations of chapter \ref{MQG}; let $\gH$ be a $N-N$-bimodule, i.e. an Hilbert space equipped with a normal faithful representation $a$ of $N$ on $\gH$ and a normal faithful anti-representation $b$ on $\gH$, such that $b(N)\subset a(N)'$. Let now $V$ be a unitary from $\gH\underset{\nu^o}{_a\otimes_\beta}H_\Phi$ onto 
 $\gH\underset{\nu}{_b\otimes_\alpha}H_\Phi$, satisfying all the following properties, for all $x\in N$ :
 \[V(b(x)\underset{N^o}{_a\otimes_\beta}1)=(1\underset{N}{_b\otimes_\alpha}\beta(x))V\]
 \[V(1\underset{N^o}{_a\otimes_\beta}\alpha(x))=(a(x)\underset{N}{_b\otimes_\alpha}1)V\]
 \[V(1\underset{N^o}{_a\otimes_\beta}\hat{\beta}(x))=(1\underset{N}{_b\otimes_\alpha}\hat{\beta}(x))V\]
Thanks to these intertwining properties, the following operator has a meaning :
\[(1\underset{N}{_b\otimes_\alpha}W^*)(1\underset{N}{_b\otimes_\alpha}\sigma_\nu)(V\underset{N^o}{_{\hat{\beta}}\otimes_\alpha}1)(1\underset{N^o}{_a\otimes_\beta}\sigma_{\nu^o})(1\underset{N^o}{_a\otimes_\beta}W)\]
and is a unitary from $\gH\underset{\nu^o}{_a\otimes_\beta}(H_\Phi\underset{\nu}{_\beta\otimes_\alpha}H_\Phi)$ onto $\gH\underset{\nu^o}{_b\otimes_\alpha}H_\Phi\underset{\nu}{_\beta\otimes_\alpha}H_\Phi$. 
Let us recall that the parenthesis in the first Hilbert space means that the antirepresentation $\beta$ is here acting on the second leg of $H_\Phi\underset{\nu}{_\beta\otimes_\alpha}H_\Phi$. Moreover, here, $\sigma_\nu$ is the flip from $H_\Phi\underset{\nu}{_{\hat{\beta}}\otimes_\alpha}H_\Phi$ onto 
$H_\Phi\underset{\nu^o}{_\alpha\otimes_{\hat{\beta}}}H_\Phi$, and $\sigma_{\nu^o}$ is the flip from $H_\Phi\underset{\nu}{_\beta\otimes_\alpha}H_\phi$ onto $H_\Phi\underset{\nu^o}{_\alpha\otimes_\beta}H_\Phi$. 
\newline
We may consider also the following unitary from $\gH\underset{\nu^o}{_a\otimes_\beta}(H_\Phi\underset{\nu}{_\beta\otimes_\alpha}H_\Phi)$ onto $\gH\underset{\nu^o}{_b\otimes_\alpha}H_\Phi\underset{\nu}{_\beta\otimes_\alpha}H_\Phi$:
\[(V\underset{N}{_b\otimes_\alpha}1)\sigma^{2,3}_{b, a}(V\underset{N^o}{_\alpha\otimes_\beta}1)(1\underset{N^o}{_a\otimes_\beta}\sigma_\nu)\]
Now, in that formula, $\sigma_\nu$ is the flip from $H_\Phi\underset{\nu}{_{\beta}\otimes_\alpha}H_\Phi$ onto 
$H_\Phi\underset{\nu^o}{_\alpha\otimes_{\beta}}H_\Phi$, and $\sigma^{2,3}_{b, a}$ is the flip from $(\gH\underset{\nu}{_b\otimes_\alpha}H_\Phi)\underset{\nu^o}{_a\otimes_\beta}H_\Phi$ onto $(\gH\underset{\nu^o}{_a\otimes_\beta}H_\Phi)\underset{\nu}{_b\otimes_\alpha}H_\Phi$. 
\newline
We shall say that $V$ is a corepresentation of $\gG$ on the bimodule $_a\gH_b$ if we have :
\begin{multline*}
(1\underset{N}{_b\otimes_\alpha}W^*)(1\underset{N}{_b\otimes_\alpha}\sigma_\nu)(V\underset{N^o}{_{\hat{\beta}}\otimes_\alpha}1)(1\underset{N^o}{_a\otimes_\beta}\sigma_{\nu^o})(1\underset{N^o}{_a\otimes_\beta}W)\\
=(V\underset{N}{_b\otimes_\alpha}1)\sigma^{2,3}_{b, a}(V\underset{N^o}{_\alpha\otimes_\beta}1)(1\underset{N^o}{_a\otimes_\beta}\sigma_\nu)
\end{multline*}
 
%%%%%thV
\subsection{Theorem}
\label{thV}
{\it Let $\gG=(N, M, \alpha, \beta, \Gamma, T, T', \nu)$ be a measured quantum groupoid, and $_a\gH_b$ a $N-N$-bimodule; let $V$ be a unitary from $\gH\underset{\nu^o}{_a\otimes_\beta}H_\Phi$ onto  $\gH\underset{\nu}{_b\otimes_\alpha}H_\Phi$, satisfying all the following properties, for all $x\in N$ :
 \[V(b(x)\underset{N^o}{_a\otimes_\beta}1)=(1\underset{N}{_b\otimes_\alpha}\beta(x))V\]
 \[V(1\underset{N^o}{_a\otimes_\beta}\alpha(x))=(a(x)\underset{N}{_b\otimes_\alpha}1)V\]
 \[V(1\underset{N^o}{_a\otimes_\beta}\hat{\beta}(x))=(1\underset{N}{_b\otimes_\alpha}\hat{\beta}(x))V\]
 Then, are equivalent :
 \newline
 (i) $V$ is a corepresentation of $\gG$ on the bimodule $_a\gH_b$;
 \newline
 (ii) for any $\xi\in D(_a\gH, \nu)$ and $\eta\in D(\gH_b, \nu^o)$, the operator $(\omega_{\xi, \eta}*id)(V)$ belongs to $M$ and is such that :}
\[\Gamma((\omega_{\xi, \eta}*id)(V))=(\omega_{\xi, \eta}*id*id)[(V\underset{N}{_b\otimes_\alpha}1)\sigma^{2,3}_{b, a}(V\underset{N^o}{_\alpha\otimes_\beta}1)(1\underset{N^o}{_a\otimes_\beta}\sigma_\nu)]\]

\begin{proof}
Let $V$ be a corepresentation on $_a\gH_b$; let us take now a unitary $u$ in $M'$; as we have then :
\[(u^*\underset{\nu^o}{_\alpha\otimes_{\hat{\beta}}}1)W(u\underset{\nu}{_\beta\otimes_\alpha}1)=W\]
it is easy to get, because $V$ is a unitary, that :
\[(1\underset{\nu}{_b\otimes_\alpha}u^*)V(1\underset{\nu^o}{_a\otimes_\beta}u)=V\]
from which we get, for any $\xi\in D(_a\gH, \nu)$ and $\eta\in D(\gH_b, \nu^o)$, that the operator $(\omega_{\xi, \eta}*id)(V)$ belongs to $M$; from the definition of a corepresentation, we then get :
\[\Gamma((\omega_{\xi, \eta}*id)(V))=(\omega_{\xi, \eta}*id*id)[(V\underset{N}{_b\otimes_\alpha}1)\sigma^{2,3}_{b, a}(V\underset{N^o}{_\alpha\otimes_\beta}1)(1\underset{N^o}{_a\otimes_\beta}\sigma_\nu)]\]
The converse is trivial. \end{proof}

%%%%corV*
\subsection{Corollary}
\label{corV*}
{\it Let $\gG=(N, M, \alpha, \beta, \Gamma, T, T', \nu)$ be a measured quantum groupoid, and $_a\gH_b$ be a $N-N$-bimodule; let $V$ be unitary from $\gH\underset{\nu^o}{_a\otimes_\beta}H_\Phi$ onto  $\gH\underset{\nu}{_b\otimes_\alpha}H_\Phi$; then, are equivalent :
\newline
(i) $V$ is a corepresentation of $\gG$ on the $N-N$-bimodule $_a\gH_b$;
\newline
(ii) $V^*$ is a corepresentation of $\gG^o$ on the $N^o-N^o$-bimodule $_b\gH_a$.}

\begin{proof}
Let us suppose that $V$ is a corepresentation of $\gG$ on the $N-N$-bimodule $_a\gH_b$; then, for any $\xi\in D(_a\gH, \nu)$ and $\eta\in D(\gH_b, \nu^o)$, the operator $(\omega_{\xi, \eta}*id)(V)$ belongs to $M$ and we have :
\[\Gamma((\omega_{\xi, \eta}*id)(V))=(\omega_{\xi, \eta}*id*id)[(V\underset{N}{_b\otimes_\alpha}1)\sigma^{2,3}_{b, a}(V\underset{N^o}{_\alpha\otimes_\beta}1)(1\underset{N^o}{_a\otimes_\beta}\sigma_\nu)]\]
Therefore, taking the adjoint, we get :
\[\varsigma_\nu\Gamma((\omega_{\eta, \xi}*id)(V^*))=(\omega_{\eta, \xi}*id*id)[(V^*\underset{N^o}{_a\otimes_\beta}1)\sigma^{2,3}_{a, b}(V^*\underset{N}{_b\otimes_\alpha}1)(1\underset{N^o}{_a\otimes_\beta}\sigma_\nu)]\]
which gives, thanks to \ref{thV}, that $V^*$ is a corepresentation of $\gG^o$ on the $N^o-N^o$-bimodule $_b\gH_a$. The proof of the converse is the same, applied to $\gG^o$. \end{proof}

%%%%isom
\subsection{Proposition}
\label{isom}
{\it Let $\gG=(N, M, \alpha, \beta, \Gamma, T, T', \nu)$ be a measured quantum groupoid, and $_a\gH_b$ a $N-N$-bimodule; let $V$ be an isometry from $\gH\underset{\nu^o}{_a\otimes_\beta}H_\Phi$ to  $\gH\underset{\nu}{_b\otimes_\alpha}H_\Phi$, satisfying all the following properties, for all $x\in N$ :
 \[V(b(x)\underset{N^o}{_a\otimes_\beta}1)=(1\underset{N}{_b\otimes_\alpha}\beta(x))V\]
 \[V(1\underset{N^o}{_a\otimes_\beta}\alpha(x))=(a(x)\underset{N}{_b\otimes_\alpha}1)V\]
 \[V(1\underset{N^o}{_a\otimes_\beta}\hat{\beta}(x))=(1\underset{N}{_b\otimes_\alpha}\hat{\beta}(x))V\]
and let us suppose that for any $\xi\in D(_a\gH, \nu)$ and $\eta\in D(\gH_b, \nu^o)$, that the operator $(\omega_{\xi, \eta}*id)(V)$ belongs to $M$ and is such that :
\[\Gamma((\omega_{\xi, \eta}*id)(V))=(\omega_{\xi, \eta}*id*id)[(V\underset{N}{_b\otimes_\alpha}1)\sigma^{2,3}_{b, a}(V\underset{N^o}{_\alpha\otimes_\beta}1)(1\underset{N^o}{_a\otimes_\beta}\sigma_\nu)]\]
then $V$ is a unitary and a corepresentation of $\gG$ on $_a\gH_b$. }

\begin{proof}
We easily have :
\[(1\underset{N}{_b\otimes_\alpha}\sigma_{\nu^o}W^*\sigma_\nu) (V\underset{N}{_{\hat{\beta}}\otimes_\alpha}1)(1\underset{N^o}{_a\otimes_\beta}\sigma_{\nu^o}W\sigma_\nu)=
(1\underset{N}{_b\otimes_\alpha}\sigma_{\nu^o})(V\underset{N}{_b\otimes_\alpha}1)\sigma^{2,3}_{b, a}(V\underset{N^o}{_\alpha\otimes_\beta}1)\]
Then, taking the adjoints, we get :
\[(1\underset{N^o}{_a\otimes_\beta}\sigma_{\nu^o}W^*\sigma_\nu)(V^*\underset{N}{_{\hat{\beta}}\otimes_\alpha}1)(1\underset{N}{_b\otimes_\alpha}\sigma_{\nu^o}W\sigma_\nu) =
(V^*\underset{N^o}{_a\otimes_\beta}1)\sigma^{2,3}_{a, b}(V^*\underset{N}{_\beta\otimes_\alpha}1)(1\underset{N}{_b\otimes_\alpha}\sigma_{\nu})\]
and, by multiplication, we get, because $Q=VV^*$ is a projection in $\mathcal L(\gH)\underset{N}{_b*_\alpha}M$ :
\begin{multline*}
(1\underset{N}{_b\otimes_\alpha}\sigma_{\nu^o}W^*\sigma_\nu) (Q\underset{N}{_{\hat{\beta}}\otimes_\alpha}1)(1\underset{N^o}{_a\otimes_\beta}\sigma_{\nu^o}W\sigma_\nu)\\
=(1\underset{N}{_b\otimes_\alpha}\sigma_{\nu^o})(V\underset{N}{_b\otimes_\alpha}1)\sigma^{2,3}_{b, a}(Q\underset{N}{_\beta\otimes_\alpha}1)\sigma^{2,3}_{a, b}(V^*\underset{N^o}{_a\otimes_\beta}1)(1\underset{N^o}{_a\otimes_\beta}\sigma_{\nu})\end{multline*}
which can be written :
\[(id\underset{N}{_b*_\alpha}\Gamma)(Q)= 
(V\underset{N}{_b\otimes_\alpha}1)\sigma^{2,3}_{b, a}(Q\underset{N}{_\beta\otimes_\alpha}1)\sigma^{2,3}_{a, b}(V^*\underset{N^o}{_a\otimes_\beta}1)\]
from which we get that $(id\underset{N}{_b*_\alpha}\Gamma)(Q)=(Q\underset{\nu}{_\beta\otimes_\alpha}1)(id\underset{N}{_b*_\alpha}\Gamma)(Q)$. 
\newline
But, using \ref{thL5}(v), we know that $(id\underset{N}{_b*_\alpha}\Gamma)(Q)$ is equal to 
\[(1\underset{N}{_b\otimes_\alpha}\sigma_\nu W^o\sigma_\nu)^*(Q\underset{N^o}{_{\hat{\alpha}}\otimes_\beta}1)(1\underset{N}{_b\otimes_\alpha}\sigma_\nu W^o\sigma_\nu)\]
from which we infer that :
\[(1\underset{N}{_b\otimes_\alpha}\sigma_\nu W^o\sigma_\nu)^*(Q\underset{N^o}{_{\hat{\alpha}}\otimes_\beta}1)=(Q\underset{N^o}{_\beta\otimes_\alpha}1)(1\underset{N}{_b\otimes_\alpha}\sigma_\nu W^o\sigma_\nu)^*(Q\underset{N^o}{_{\hat{\alpha}}\otimes_\beta}1)\]
Let's take now $\xi\in D((H_\Phi)_\beta, \nu^o)$ and $\eta\in D(_\alpha H_\Phi, \nu)$, and let us apply $(id*id*\omega_{\xi, \eta})$ to this last equality. We find :
\[(1\underset{N}{_b\otimes_\alpha}(\omega_{\eta, \xi}*id)(W^o)^*)Q=Q(1\underset{N}{_b\otimes_\alpha}(\omega_{\eta, \xi}*id)(W^o)^*)Q\]
and, by density, we get $(1\underset{N}{_b\otimes_\alpha}y)Q=Q(1\underset{N}{_b\otimes_\alpha}y)Q$ for all $y\in \widehat{M}'$, from which we get that $Q$ belongs also to $\mathcal L(\gH)\underset{N}{_b*_\alpha}\widehat{M}$.
\newline
Therefore $Q$ belongs to $\mathcal L(\gH)\underset{N}{_b*_\alpha}\alpha(N)=b(N)'\underset{N}{_b\otimes_\alpha}1$. So, let $q\in b(N)'$ such that $Q=q\underset{N}{_b\otimes_\alpha}1$, and $(id\underset{N}{_b*_\alpha}\Gamma)(Q)=q\underset{N}{_b\otimes_\alpha}1\underset{N}{_\beta\otimes_\alpha}1$
\newline
So, we get :
\[(q\underset{N}{_b\otimes_\alpha}1\underset{N}{_\beta\otimes_\alpha}1)=(V\underset{N}{_b\otimes_\alpha}1)(q\underset{N}{_b\otimes_\alpha}1\underset{N}{_\beta\otimes_\alpha}1)(V^*\underset{N^o}{_a\otimes_\beta}1)\]
and, taking the adjoints, we get : 
\[(V^*\underset{N}{_\beta\otimes_\alpha}1)(q\underset{N}{_b\otimes_\alpha}1\underset{N}{_\beta\otimes_\alpha}1)(V\underset{N^o}{_\alpha\otimes_\beta}1)=(q\underset{N}{_b\otimes_\alpha}1\underset{N}{_\beta\otimes_\alpha}1)\]
but, as, by definition of $q$, we have $(q\underset{N}{_b\otimes_\alpha}1\underset{N}{_\beta\otimes_\alpha}1)(V\underset{N^o}{_\alpha\otimes_\beta}1)=V\underset{N^o}{_\alpha\otimes_\beta}1$, and, as $V$ is an isometry, we finally get that $q=1$; so, $V$ is a unitary, which finishes the proof. \end{proof}

%%%%rep
\subsection{Proposition}
\label{rep}
{\it Let $\gG=(N, M, \alpha, \beta, \Gamma, T, T', \nu)$ be a measured quantum groupoid, and $_a\gH_b$ a $N-N$-bimodule; let $V$ be a unitary from $\gH\underset{\nu^o}{_a\otimes_\beta}H_\Phi$ onto 
 $\gH\underset{\nu}{_b\otimes_\alpha}H_\Phi$, satisfying all the following properties, for all $x\in N$ :
 \[V(b(x)\underset{N^o}{_a\otimes_\beta}1)=(1\underset{N}{_b\otimes_\alpha}\beta(x))V\]
 \[V(1\underset{N^o}{_a\otimes_\beta}\alpha(x))=(a(x)\underset{N}{_b\otimes_\alpha}1)V\]
 \[V(1\underset{N^o}{_a\otimes_\beta}\hat{\beta}(x))=(1\underset{N}{_b\otimes_\alpha}\hat{\beta}(x))V\]
and let us suppose that, for any $\xi\in D(_a\gH, \nu)$ and $\eta\in D(\gH_b, \nu^o)$, the operator $(\omega_{\xi, \eta}*id)(V)$ belongs to $M$; therefore, it is possible to define $(id*\theta)(V)$, for any $\theta$ in $M_*^{\alpha, \beta}$. Moreover, are equivalent :
\newline
(i) $V$ is a corepresentation of $\gG$ on $_a\gH_b$; 
\newline
(ii) the application $\theta\rightarrow (id*\theta)(V)$ from $M_*^{\alpha, \beta}$ into $\mathcal L(\gH)$ is multiplicative. }

\begin{proof} Clear. \end{proof}

%%%%Example
\subsection{Example}
\label{Example}
Let $\gG=(N, M, \alpha, \beta, \Gamma, T, T', \nu)$ be a measured quantum groupoid; then the pseudo-multiplicative unitary $W$ verifies :
\[(1_{\gH}\underset{N^o}{_\alpha\otimes_{\hat{\beta}}}W)
(W\underset{N}{_\beta\otimes_\alpha}1_{\gH})(1_{\gH}\underset{N}{_\beta\otimes_\alpha}W^*)
=(W\underset{N^o}{_\alpha\otimes_{\hat{\beta}}}1_{\gH})
\sigma^{2,3}_{\alpha, \beta}
(W\underset{N}{_{\hat{\beta}}\otimes_\alpha}1)
(1_{\gH}\underset{N}{_\beta\otimes_\alpha}\sigma_{\nu^o})\]
from which one gets, using the definition of $\widehat{\Gamma}$, that we have, for any $\xi\in D((H_\Phi)_\beta, \nu^o)$ and $\eta\in D(_\alpha H_\Phi, \nu)$ :
\[\varsigma_N\widehat{\Gamma}((\omega_{\xi, \eta}*id)(W))=(\omega_{\xi, \eta}*id*id)[(W\underset{N^o}{_\alpha\otimes_{\hat{\beta}}}1_{\gH})
\sigma^{2,3}_{\alpha, \beta}
(W\underset{N}{_{\hat{\beta}}\otimes_\alpha}1)
(1_{\gH}\underset{N}{_\beta\otimes_\alpha}\sigma_{\nu^o})]\]
which means, thanks to \ref{isom}, that $W$ is a corepresentation of $(\widehat{\gG})^o$ on the $N^o-N^o$ bimodule $_\beta (H_\Phi)_\alpha$. So, $(\sigma_{\nu^o}W\sigma_{\nu^o})^*$ is a corepresentation of $\gG^o$ on the $N^o-N^o$ bimodule $_{\hat{\beta}}(H_\Phi)_\alpha$, and, using \ref{corV*}, $\sigma_{\nu^o}W\sigma_{\nu^o}$ is a corepresentation of $\gG$ on the $N-N$ bimodule $_\alpha(H_\Phi)_{\hat{\beta}}$. 
\newline
We get also that $(\sigma_\nu W^o\sigma_\nu)^*$ is a corepresentation of $\gG$ on the $N-N$ bimodule $_{\hat{\alpha}}(H_\Phi)_\beta$ and that $W^{oc}$ is a corepresentation of $\widehat{\gG}^c$ on the $N^o-N^o$ bimodule $_{\hat{\beta}}(H_\Phi)_{\hat{\alpha}}$.

%%%%propS
\subsection{Proposition}
\label{propS}
{\it Let $\gG=(N, M, \alpha, \beta, \Gamma, T, T', \nu)$ be a measured quantum groupoid, and $_a\gH_b$ a $N-N$-bimodule; let $V$ be a corepresentation of $\gG$ on $_a\gH_b$. Then, we have :}
\[(V\underset{N^o}{_\alpha\otimes_\beta}1)(1\underset{N^o}{_a\otimes_\beta}\sigma_{\nu^o}W^*\sigma_\nu)(V^*\underset{N}{_{\hat{\beta}}\otimes_\alpha}1)=
\sigma^{2,3}_{a, b}(V^*\underset{N}{_\beta\otimes_\alpha}1)(1\underset{N}{_b\otimes_\alpha}\sigma_\nu)(1\underset{N}{_b\otimes_\alpha}\sigma_{\nu^o}W^*\sigma_\nu)\]

\begin{proof}
From \ref{defcorep}, we easily have :
\[(1\underset{N}{_b\otimes_\alpha}\sigma_{\nu^o}W^*\sigma_\nu) (V\underset{N}{_{\hat{\beta}}\otimes_\alpha}1)(1\underset{N^o}{_a\otimes_\beta}\sigma_{\nu^o}W\sigma_\nu)=
(1\underset{N}{_b\otimes_\alpha}\sigma_{\nu^o})(V\underset{N}{_b\otimes_\alpha}1)\sigma^{2,3}_{b, a}(V\underset{N^o}{_\alpha\otimes_\beta}1)\]
Then, taking the adjoints, we get :
\[(1\underset{N^o}{_a\otimes_\beta}\sigma_{\nu^o}W^*\sigma_\nu)(V^*\underset{N}{_{\hat{\beta}}\otimes_\alpha}1)(1\underset{N}{_b\otimes_\alpha}\sigma_{\nu^o}W\sigma_\nu) =
(V^*\underset{N^o}{_a\otimes_\beta}1)\sigma^{2,3}_{a, b}(V^*\underset{N}{_\beta\otimes_\alpha}1)(1\underset{N}{_b\otimes_\alpha}\sigma_{\nu})\]
from which we get the result, because $V$ is an unitary.  \end{proof}

%%%%%prop1S
\subsection{Proposition}
\label{prop1S}
{\it Let $\gG=(N, M, \alpha, \beta, \Gamma, T, T', \nu)$ be a measured quantum groupoid, and $_a\gH_b$ a $N-N$-bimodule; let $V$ be a corepresentation of $\gG$ on $_a\gH_b$. Let $\xi_1$, $\xi_2$ be in $D(_a\gH, \nu)\cap D(\gH_b, \nu^o)$, $\eta_1$, $\eta_2$ be in $D(_\alpha H_\Phi, \nu)\cap D((H_\Phi)_{\hat{\beta}}, \nu^o)$. Then, the element :
\[(\omega_{V^*(\xi_1\underset{\nu}{_b\otimes_\alpha}\eta_1), V^*(\xi_2\underset{\nu}{_b\otimes_\alpha}\eta_2)}*id)(1\underset{N^o}{_a\otimes_\beta}\sigma_{\nu^o}W\sigma_\nu)\]
 belongs to $D(S)$, and :
 \begin{multline*}
S((\omega_{V^*(\xi_1\underset{\nu}{_b\otimes_\alpha}\eta_1), V^*(\xi_2\underset{\nu}{_b\otimes_\alpha}\eta_2)}*id)(1\underset{N^o}{_a\otimes_\beta}\sigma_{\nu^o}W\sigma_\nu))=\\
(\omega_{V^*(\xi_1\underset{\nu}{_b\otimes_\alpha}\eta_1), V^*(\xi_2\underset{\nu}{_b\otimes_\alpha}\eta_2)}*id)(1\underset{N^o}{_a\otimes_\beta}\sigma_{\nu^o}W^*\sigma_\nu)
\end{multline*}
}
\begin{proof}
By linearity, it suffices to prove that the element of $\widehat{M}_*$ defined, for $y$ in $\widehat{M}$ by :
\[y\mapsto \omega_{V^*(\xi_1\underset{\nu}{_b\otimes_\alpha}\eta_1)}(1\underset{N^o}{_a\otimes_\beta}y)\]
belongs to $\widehat{M}_*^{\alpha, \hat{\beta}}$. For $n\in N$, we have :
\[(1\underset{N^o}{_a\otimes_\beta}\alpha(n))V^*(\xi_1\underset{\nu}{_b\otimes_\alpha}\eta_1)=
V^*(a(n)\xi_1\underset{\nu}{_b\otimes_\alpha}\eta_1)\]
\[(1\underset{N^o}{_a\otimes_\beta}\hat{\beta}(n))V^*(\xi_1\underset{\nu}{_b\otimes_\alpha}\eta_1)=
V^*(\xi_1\underset{\nu}{_b\otimes_\alpha}\hat{\beta}(n)\eta_1)\]
and, therefore, if $n\in\gN_\nu$, 
\[\|(1\underset{N^o}{_a\otimes_\beta}\alpha(n))V^*(\xi_1\underset{\nu}{_b\otimes_\alpha}\eta_1)\|
\leq \|R^\alpha(\eta_1)\|\|R^a(\xi_1)\|\|\Lambda_\nu(n)\|\]
\[\|(1\underset{N^o}{_a\otimes_\beta}\hat{\beta}(n))V^*(\xi_1\underset{\nu}{_b\otimes_\alpha}\eta_1)\|\leq\|R^b(\xi_1)\|\|R^\alpha(\eta_1)\|\|\Lambda_\nu(n)\|\]
which, therefore, gives the result. \end{proof}

%%%%%%%prop2S
\subsection{Proposition}
\label{prop2S}
{\it Let $\gG=(N, M, \alpha, \beta, \Gamma, T, T', \nu)$ be a measured quantum groupoid, and $_a\gH_b$ a $N-N$-bimodule; let $V$ be a corepresentation of $\gG$ on $_a\gH_b$. Let $\xi_1$, $\xi_2$ be in $D(_a\gH, \nu)\cap D(\gH_b, \nu^o)$, $\eta_1$, $\eta_2$ be in $D(_\alpha H_\Phi, \nu)\cap D((H_\Phi)_{\hat{\beta}}, \nu^o)$. Then, the element :
\[(id*\omega_{\eta_1, \eta_2})(W)(\omega_{\xi_1, \xi_2}*id)(V)\]
belongs to $D(S)$, and :}
\[S[(id*\omega_{\eta_1, \eta_2})(W)(\omega_{\xi_1, \xi_2}*id)(V)]=(\omega_{\xi_1, \xi_2}*id)(V^*)(id*\omega_{\eta_1, \eta_2})(W^*)\]

\begin{proof}
Using \ref{propS}, we get that :
\[(id*\omega_{\eta_1, \eta_2})(W)(\omega_{\xi_1, \xi_2}*id)(V)=(\omega_{V^*(\xi_1\underset{\nu}{_b\otimes_\alpha}\eta_1), V^*(\xi_2\underset{\nu}{_b\otimes_\alpha}\eta_2)}*id)(1\underset{N^o}{_a\otimes_\beta}\sigma_{\nu^o}W\sigma_\nu)\]
\[(\omega_{\xi_1, \xi_2}*id)(V^*)(id*\omega_{\eta_1, \eta_2})(W^*)=(\omega_{V^*(\xi_1\underset{\nu}{_b\otimes_\alpha}\eta_1), V^*(\xi_2\underset{\nu}{_b\otimes_\alpha}\eta_2)}*id)(1\underset{N^o}{_a\otimes_\beta}\sigma_{\nu^o}W^*\sigma_\nu)\]
We then get the result by \ref{prop1S}. 
\end{proof}

%%%%%thS1
\subsection{Theorem}
\label{thS1}
{\it Let $\gG=(N, M, \alpha, \beta, \Gamma, T, T', \nu)$ be a measured quantum groupoid, and $_a\gH_b$ a $N-N$-bimodule; let $V$ be a corepresentation of $\gG$ on $_a\gH_b$. Let $\xi_1$, $\xi_2$ in $D(_a\gH, \nu)\cap D(\gH_b, \nu^o)$; then the element $(\omega_{\xi_1, \xi_2}*id)(V^*)$ belongs to $D(S)$, and :}
\[S((\omega_{\xi_1, \xi_2}*id)(V))=(\omega_{\xi_1, \xi_2}*id)(V^*)\]

\begin{proof}
Using \ref{prop2S}, we get that, for any $x$ of the form $(id*\omega_{\eta_1, \eta_2})(W^*)$, where $\eta_1$, $\eta_2$ are in $D(_\alpha H_\Phi, \nu)\cap D((H_\Phi)_{\hat{\beta}}, \nu^o)$, the operator
 $x(\omega_{\xi_1, \xi_2}*id)(V)$ belongs to $D(S)$, and $S(x(\omega_{\xi_1, \xi_2}*id)(V^*))=(\omega_{\xi_1, \xi_2}*id)(V^*)S(x)$, \newline
If, moreover, $\eta_1$ belongs to $D(P^{-1/2})$ and is such that  $P^{-1/2}\eta_1$ belongs to $D((H_\Phi)_{\hat{\beta}}, \nu^o)$, and if $\eta_2$ belongs to $D(P^{1/2})$ and is such that  $P^{1/2}\eta_2$ belongs to $D(_\alpha(H_\Phi, \nu)$, then using (\cite{E5}, 3.8), we know that :
\[(id*\omega_{J_\Phi\eta_1, J_\Phi\eta_2})(W^*)=(id*\omega_{P^{-1/2}\eta_1, P^{1/2}\eta_2})(W^*)^*\]
and we know (\cite{E5}, 3.9) that the linear space generated by such elements is weakly dense in $M$. Let $\mathcal A$ be the algebra generated by these elements, which is therefore involutive, and dense in $M$. 
\newline
Therefore, we obtain that $S(x(\omega_{\xi_1, \xi_2}*id)(V^*))=(\omega_{\xi_1, \xi_2}*id)(V^*)S(x)$ for any $x\in\mathcal A$. Then, by Kaplanski's theorem, it is therefore possible to find a sequence $a_n$ in $D(S)$ increasing to $1$ such that :
\[S(a_n(\omega_{\xi_1, \xi_2}*id)(V^*))=(\omega_{\xi_1, \xi_2}*id)(V^*)S(a_n)\]
Let's take now 
\[e_n=1/\sqrt{\pi}\int_{-\infty}^\infty e^{-t^2}\tau_t(a_n)dt\]
As $S$ is closed, we have again 
\[S(e_n(\omega_{\xi_1, \xi_2}*id)(V^*))=(\omega_{\xi_1, \xi_2}*id)(V^*)S(e_n)\]
Moreover, $e_n$ is increasing to $1$, and, as $\tau_{-i/2}(e_n)$ is bounded, $\tau_{-i/2}(e_n)$ converges also to $1$, and so does $S(e_n)$. So, we get the result, using again the closedness of $S$.  \end{proof}

%%%%%thS2
\subsection{Theorem}
\label{thS2}
{\it Let  $\gG=(N, M, \alpha, \beta, \Gamma, T, T', \nu)$ be a measured quantum groupoid, and $_a\gH_b$ a $N-N$-bimodule, and let us suppose that $D(_a\gH, \nu)\cap D(\gH_b, \nu^o)$ is dense in $\gH$. Let $V$ be a corepresentation of $\gG$ on $_a\gH_b$. Then, for any $\zeta_1$, $\zeta_2$ in $\gH$, for any $\xi\in D((H_\Phi)_\beta, \nu^o)\cap\mathcal D(P^{1/2})$ such that $P^{1/2}\xi$ belongs to $D(_\alpha H_\Phi, \nu)$, and for any $\eta\in D(_\alpha H_\Phi, \nu)\cap D(P^{-1/2})$ such that $P^{-1/2}\eta$ belongs to $D((H_\Phi)_\beta, \nu^o)$, we have : }
\[(V(\zeta_1\underset{\nu^o}{_a\otimes_\beta}\xi)|\zeta_2\underset{\nu}{_b\otimes_\alpha}\eta)
=(V^*(\zeta_1\underset{\nu}{_b\otimes_\alpha}J_{\widehat{\Phi}}P^{-1/2}\eta)|
\zeta_2\underset{\nu^o}{_a\otimes_\beta}J_{\widehat{\Phi}}P^{1/2}\xi)
\]
\begin{proof}
It is clear that $J_{\widehat{\Phi}}P^{-1/2}\eta$ belongs to $D(_\alpha H_\Phi, \nu)$ and that $J_{\widehat{\Phi}}P^{1/2}\xi$ belongs to $D((H_\Phi)_\beta, \nu^o)$. Moreover, we have, if we take $\zeta_1$, $\zeta_2$ in $D(_a\gH, \nu)\cap D(\gH_b, \nu^o)$ :
\begin{eqnarray*}
(V(\zeta_1\underset{\nu^o}{_a\otimes_\beta}\xi)|\zeta_2\underset{\nu}{_b\otimes_\alpha}\eta)
&=&(\zeta_1\underset{\nu^o}{_a\otimes_\beta}\xi|V^*(\zeta_2\underset{\nu}{_b\otimes_\alpha}\eta))\\
&=&(\xi|(\omega_{\zeta_2, \zeta_1}*id)(V^*)\eta)
\end{eqnarray*}
which, using \ref{thS1}, and the fact that $J_{\widehat{\Phi}}P^{1/2}=P^{-1/2}J_{\widehat{\Phi}}$, is equal to :
\begin{eqnarray*}
(\xi|S[(\omega_{\zeta_2, \zeta_1}*id)(V)]\eta)&=&
(\xi|J_{\widehat{\Phi}}P^{-1/2}(\omega_{\zeta_2, \zeta_1}*id)(V)^*P^{1/2}J_{\widehat{\Phi}}\eta)\\
&=&((\omega_{\zeta_2, \zeta_1}*id)(V)^*J_{\widehat{\Phi}}P^{-1/2}\eta|J_{\widehat{\Phi}}P^{1/2}\xi)\\
&=&(J_{\widehat{\Phi}}P^{-1/2}\eta|(\omega_{\zeta_2, \zeta_1}*id)(V)J_{\widehat{\Phi}}P^{1/2}\xi)\\
&=&(\zeta_1\underset{\nu}{_b\otimes_\alpha}J_{\widehat{\Phi}}P^{-1/2}\eta|V(\zeta_2\underset{\nu^o}{_a\otimes_\beta}J_{\widehat{\Phi}}P^{1/2}\xi))\\
&=&(V^*(\zeta_1\underset{\nu}{_b\otimes_\alpha}J_{\widehat{\Phi}}P^{-1/2}\eta)|
\zeta_2\underset{\nu^o}{_a\otimes_\beta}J_{\widehat{\Phi}}P^{1/2}\xi)
\end{eqnarray*}
Therefore, the result is proved for $\zeta_1$, $\zeta_2$ in $D(_a\gH, \nu)\cap D(\gH_b, \nu^o)$, and it remains true, by continuity, for any $\zeta_1$, $\zeta_2$ in $\gH$. \end{proof}

%%%%%action
\section{Actions of measured quantum groupoids}
\label{action}
In this chapter, we define actions of measured quantum groupoids on von Neumann algebras (\ref{defaction}) and actions implemented by a corepresentation (\ref{PropV}). We define also invariant elements by an action (\ref{definvariant}), the canonical operator-valued weight on the invariants (\ref{defovw1}), and integrable actions (\ref{defintegrable})

%%%%%defaction
\subsection{Definition}
\label{defaction}
Let $\gG=(N, M, \alpha, \beta, \Gamma, T, T', \nu)$ be a measured quantum groupoid, and let $A$ be a von Neumann algebra. 
\newline
An action of $\gG$ on $A$ is a couple $(b, \mathfrak a)$, where :
\newline
(i) $b$ is an injective $*$-antihomomorphism from $N$ into $A$; 
\newline
(ii) $\mathfrak a$ is an injective $*$-homomorphism from $A$ into $A\underset{N}{_b*_\alpha}M$; 
\newline
(iii) $b$ and $\mathfrak a$ are such that, for all $n$ in $N$:
\[\mathfrak a (b(n))=1\underset{N}{_b\otimes_\alpha}\beta(n)\]
(which allow us to define $\mathfrak a\underset{N}{_b*_\alpha}id$ from $A\underset{N}{_b*_\alpha}M$ into $A\underset{N}{_b*_\alpha}M\underset{N}{_\beta*_\alpha}M$)
and such that :
\[(\mathfrak a\underset{N}{_b*_\alpha}id)\mathfrak a=(id\underset{N}{_b*_\alpha}\Gamma)\mathfrak a\]
If there is no ambiguity, we shall say that $\mathfrak a$ is the action. 
\newline
So, a measured quantum groupoid $\gG$ can act only on a von Neumann algebra $A$ which is a right module over the basis $N$.
\newline
Moreover, if $M$ is abelian, then, as, for all $n\in N$, $\mathfrak a (b(n))=1\underset{N}{_b\otimes_\alpha}\beta(n)$ commutes with $\mathfrak a (x)$, for all $x\in A$, we see that $b(N)$ is in the center of $A$. As in that case (\ref{gd2}) the measured quantum groupoid comes then from a measured groupoid $\mathcal G$, we have $N=L^\infty(\mathcal G^{(0)}, \nu)$, and $A$ can be decomposed as $A=\int_{\mathcal G^{(0)}}A^xd\nu (x)$. 

%%%%triaction
\subsection{Example}
\label{triaction}
Let $\gG=(N, M, \alpha, \beta, \Gamma, T, T', \nu)$ be a measured quantum groupoid; then $(id, \beta)$ is an action of $\gG$ on $N^o$, we shall call the trivial action of $\gG$ on $N^o$. 

%%%%graction
\subsection{Example}
\label{graction}
Let $\mathcal G$ be a measured groupoid, and let $\mathfrak a$ be an action of $\mathcal G$ on a von Neumann algebra $A=\int_{\mathcal G^{(0)}}^\oplus A^xd\nu (x)$, as defined in (\cite{Y3}, def. 3.1(i)); let us define then 
$b$ as the homomorphism which sends $L^\infty(\mathcal G^{(0)})$ into the diagonalizable operators in $A$. Then, $(b, \mathfrak a)$ is an action (in the sense of \ref{defaction}) of the measured quantum groupoid $\mathcal G$  defined in \ref{gd2} on $A$ (\cite{OR}, Th. 3.2), i.e. for all $g\in\mathcal G$, there exists a family of $*$-isomorphisms $\mathfrak a_g$ from $A^{s(g)}$ onto $A^{r(g)}$, such that, if $(g_1, g_2)\in\mathcal G^{(2)}$, we have $\mathfrak a_{g_1g_2}=\mathfrak a_{g_1}\mathfrak a_{g_2}$, and such that, for any normal positive functional $\omega=\int_{\mathcal G^{(0)}}^\oplus \omega^x d\nu(x)$, and any $y=\int_{\mathcal G^{(0)}}^\oplus y^x d\nu (x)$, the function $g\mapsto \omega^{r(g)}(\mathfrak a_g(y^{s(g)}))$ is $\mu$-measurable. We then get :
\[\mathfrak a (\int_{\mathcal G^{(0)}}^\oplus y^x d\nu (x))=\int_{\mathcal G}^\oplus \mathfrak a_g(y^{s(g)})d\mu(g)\]
\newline
Let now $\mathfrak b$ be a coaction of $\mathcal G$ on a von Neumann algebra $B$ which is a $L^\infty(\mathcal G^{(0)}, \nu)$-module, as defined in (\cite{Y3}, def. 3.1(ii)); let us define then $b$ as the homomorphism which sends $L^\infty(\mathcal G^{(0)}, \nu)$ into $B$; then, $(b, \mathfrak b)$ is an action of $\widehat{\mathcal G}$ (as defined in \ref{gd2}) on $B$. (\cite{Y3}, 3.1 (ii) and \cite{OR}, 2.3). 

%%%%lcqgaction
\subsection{Example}
\label{lcqgaction}
Let ${\bf G}=(M, \Gamma, \Phi, \Psi)$ be a locally compact quantum group; an action of ${\bf G}$ on a von Neumann algebra $A$ is an injective $*$-homomorphism $\mathfrak a$ from $M$ into $A\otimes M$, such that 
\[(\mathfrak a\otimes id)\mathfrak a=(id\otimes\Gamma)\mathfrak a\]
(where $id$ means the identity of $A$ or of $M$). (\cite{V2}, 1.1)
\newline
Writing now $id$ for the canonical $*$-homomorphism from $\mathbb{C}$ into $A$, we get that $(id, \mathfrak a)$ is an action of the measured quantum groupoid ${\bf G}$ (\ref{lcqg}) on $A$.

%%%d2action
\subsection{Example}
\label{d2action}
Let $M_0\subset M_1$ a depth 2 inclusion of von Neumann algebras, with a regular operator-valued 
weight $T_1$ from $M_1$ to $M_0$, as defined in \ref{basic}. Let $\gG_1$ be the measured quantum groupoid constructed (\ref{d2}) from this inclusion, and $\gG_2$ the measured quantum groupoid constructed from the inclusion $M_1\subset M_2$, which is isomorphic to $\widehat{\gG_1}^o$; we shall use the notaions introduced in \ref{basic} and \ref{d2}. Then, there exists (\cite{EV}, 7.3) a canonical action of $\gG_2$ on $M_1$, which can be described as follows : the anti-representation of the basis $M'_1\cap M_2$ (which, using $j_1$, is anti-isomorphic to $M'_0\cap M_1$), is given by the natural inclusion of $M'_0\cap M_1$ into $M_1$, and the homomorphism from $M_1$ is given by the natural inclusion of $M_1$ into $M_3$ (which is, thanks to (\cite{EV}, 4.6), isomorphic to $M_1*\mathcal L(H_{\chi_2})$).

%%%%PropV
\subsection{Proposition}
\label{PropV}
{\it Let $\gG=(N, M, \alpha, \beta, \Gamma, T, T', \nu)$ be a measured quantum groupoid, $_a\gH_b$ be a $N-N$ bimodule, and $V$ be a corepresentation of $\gG$ on $\gH$; for any $x\in a(N)'$, let us write $\mathfrak a(x)=V(x\underset{N^o}{_a\otimes_\beta}1)V^*$ and let $A$ be a von Neumann algebra on $\gH$, such that $b(N)\subset A\subset a(N)'$, and $\mathfrak a (A)\subset A\underset{N}{_b*_\alpha}M$; then : 
\newline
(i) $(b, \mathfrak a)$ is an action of $\gG$ on $a(N)'$; 
\newline
(ii) $(b, \mathfrak a_{|A})$ is an action  of $\gG$ on $A$. 
\newline
In such a situation, we shall say that $\mathfrak a$ and $\mathfrak a_{|A}$ are implemented by $V$. }

\begin{proof}
The intertwining properties of $V$ imply that, for all $n\in N$ :
\[\mathfrak a (b(n))=1\underset{N}{_b\otimes_\alpha}\beta (n)\]
\[(a(n)\underset{N}{_b\otimes_\alpha}1)\mathfrak a (x)=\mathfrak a(x)(a(n)\underset{N}{_b\otimes_\alpha}1)\]
and this last property gives that $\mathfrak a(a(N)')\subset a(N)'\underset{N}{_b*_\alpha}M$. Therefore, we can consider $(\mathfrak a \underset{N}{_b*_\alpha}id)\mathfrak a (x)$, which will be equal to :
\[(V\underset{N}{_b\otimes_\alpha}1)\sigma^{2,3}_{b, a}(V\underset{N^o}{_\alpha\otimes_\beta}1)(x\underset{N^o}{_a\otimes_\beta}1\underset{N^o}{_\alpha\otimes_\beta}1)
(V^*\underset{N^o}{_\alpha\otimes_\beta}1)\sigma^{2,3}_{a, b}(V^*\underset{N}{_b\otimes_\alpha}1)\]
which, thanks to \ref{defcorep}, is equal to :
\[(1\underset{N}{_b\otimes_\alpha}W^*)(1\underset{N}{_b\otimes_\alpha}\sigma_\nu)(\mathfrak a (x)\underset{N}{_{\hat{\beta}}\otimes_\alpha}1)(1\underset{N}{_b\otimes_\alpha}\sigma_{\nu^o})(1\underset{N}{_b\otimes_\alpha}W)\]
which is equal to $(id\underset{N}{_b*_\alpha}\Gamma)\mathfrak a (x)$.  \end{proof}

%%%%%thimple
\subsection{Theorem}
\label{thimple}
{\it Let $\gG=(N, M, \alpha, \beta, \Gamma, T, T', \nu)$ be a measured quantum groupoid, $A$ a von Neumann algebra, $(b, \mathfrak a)$ an action of $\gG$; then, there exist an Hilbert space $\gH$, a faithful normal representation $\pi$ of $A$ on $\gH$, a normal faithful representation $a$ of $N$ on $\gH$, such that $\pi (A)\subset a(N)'$, a corepresentation $V$ of $\gG$ on the $N-N$ bimodule $_a\gH_{\pi\circ b}$, such that, for all $x$ in $A$ :}
\[(\pi\underset{N}{_b*_\alpha}id)\mathfrak a (x)=V(\pi (x)\underset{N^o}{_a\otimes_\beta}1)V^*\]

\begin{proof}
Let's take $\gH=L^2(A)\underset{N}{_b\otimes_\alpha}H_\Phi$, let us write  $a(n)=1_{L^2(A)}\underset{N}{_b\otimes_\alpha}\hat{\alpha}(n)$, for all $n\in N$, and take $V=1_{L^2(A)}\underset{N^o}{_b\otimes_\alpha}(\sigma_\nu W^o\sigma_\nu)^*$ which is (\ref{Example}) a corepresentation of $\gG$ on $_a\gH_{(1_{L^2(A)}\underset{N^o}{_b\otimes_\alpha}\beta)}$ and choose $\pi=\mathfrak a$ to obtain the result. \end{proof}

%%%%%propVV*
\subsection{Proposition}
\label{propVV*}
{\it Let $\gG=(N, M, \alpha, \beta, \Gamma, T, T', \nu)$ be a measured quantum groupoid, $_a\gH_b$ be a $N-N$ bimodule, and $V$ be a corepresentation of $\gG$ on $\gH$; let $A$ be a von Neumann algebra on $\gH$, such that $b(N)\subset A\subset a(N)'$; then, are equivalent :
\newline
(i) the corepresentation $V$ of $\gG$ implements an action $\mathfrak a$ of $\gG$ on $A$; 
\newline
(ii) the corepresentation $V^*$ of $\gG^o$ implements an action $\mathfrak a '$ of $\gG^o$ on $A'$. }

\begin{proof}
Let us suppose (i); then, we have $V(A\underset{N^o}{_a\otimes_b}1)V^*\subset A\underset{N}{_b*_\alpha}M$; taking the commutants, we get :
\[A'\underset{N}{_b\otimes_\alpha}M'\subset V(A'\underset{N^o}{_a*_\beta}\mathcal L(H_\Phi))V^*\]
and, therefore, $A'\underset{N}{_b\otimes_\alpha}1\subset V(A'\underset{N^o}{_a*_\beta}\mathcal L(H_\Phi))V^*$, from which have :
\[V^*(A'\underset{N}{_b\otimes_\alpha}1)V\subset A'\underset{N^o}{_a*_\beta}\mathcal L(H_\Phi)\]
but, from \ref{corV*} and \ref{PropV} applied to $V^*$, we get that :
\[V^*(b(N)'\underset{N}{_b\otimes_\alpha}1)V\subset b(N)'\underset{N^o}{_a*_\beta} M\]
and, therefore, we have :
\[V^*(A'\underset{N}{_b\otimes_\alpha}1)V\subset A'\underset{N^o}{_a*_\beta}M\]
from which we get (ii), using \ref{PropV} again. If we apply this proof to $\gG^o$, we obtain that (ii) implies (i). \end{proof}

%%%%standard
\subsection{Definition}
\label{standard}
Let $\gG=(N, M, \alpha, \beta, \Gamma, T, T', \nu)$ be a measured quantum groupoid, $A$ a von Neumann algebra, $(b, \mathfrak a)$ an action of $\gG$ on $A$; let $\psi$ a normal semi-finite faithful weight on $A$, and let us define, for all $n\in N$, $a(n)=J_\psi b(n^*)J_\psi$, which gives to $_a(H_\psi)_b$ a canonical structure of $N-N$ bimodule; let $V$ be a coreprentation of $\gG$ on $H_\psi$; let us suppose that $V$ implements $\mathfrak a$, i.e. that, for any $x\in A$, we have :
\[\mathfrak a (x)=V(x\underset{N^o}{_a\otimes_\beta}1)V^*\]
We shall say that $V$ is a standard implementation of $\mathfrak a$ if, moreover, we have :
\[V^*=(J_\psi\underset{\nu^o}{_\alpha\otimes_\beta}J_{\widehat{\Phi}})V(J_\psi\underset{\nu}{_b\otimes_\alpha} J_{\widehat{\Phi}})\]

%%%Gamma
\subsection{Example}
\label{Gamma}
Let $\gG=(N, M, \alpha, \beta, \Gamma, T, T', \nu)$ be a measured quantum groupoid; then $(\beta, \Gamma)$ is an action of $\gG$ on $M$. 
\newline
Moreover, $(\sigma_\nu W^o\sigma_\nu)^*$ is a corepresentation of $\gG$ on $_{\hat{\alpha}}(H_\Phi)_\beta$ which is a standard implementation of $\Gamma$.

%%%%definvariant
\subsection{Definition}
\label{definvariant}
Let $\gG=(N, M, \alpha, \beta, \Gamma, T, T', \nu)$ be a measured quantum groupo-id, and let $A$ be a von Neumann algebra, $(b, \mathfrak a)$ an action of $\gG$ on $A$; we shall call $A^\mathfrak a$ the invariant subalgebra of $A$ defined by :
\[A^\mathfrak a=\{x\in A\cap b(N)'; \mathfrak a(x)=x\underset{N}{_b\otimes_\alpha}1\}\]
In the example \ref{triaction}, the invariant subalgebra is $Z(N)$; in the example \ref{d2action}, the invariant subalgebra is $M_0$ (\cite{EV}, 7.5); in the example \ref{Gamma}, the invariant subalgebra is $\alpha (N)$(\ref{thL4}(v)). 
%%%%defovw
\subsection{Proposition}
\label{defovw}
{\it Let $\gG=(N, M, \alpha, \beta, \Gamma, T, T', \nu)$ be a measured quantum groupo-id, and let $A$ be a von Neumann algebra, $(b, \mathfrak a)$ an action of $\gG$ on $A$. Let $\Phi=\nu\circ\alpha^{-1}\circ T$; for any $x\in A^+$, the extended positive element of $A$
\[T_\mathfrak a(x)=(id\underset{\nu}{_b*_\alpha}\Phi)\mathfrak a(x)\]
is an extended positive element of $A^\mathfrak a$. }

\begin{proof}
Let $\omega\in (A\underset{N}{_b*_\alpha}M)_*^+$; then, we have :
\begin{eqnarray*}
\omega\circ\mathfrak a(T_\mathfrak a(x))
&=&
\omega(id\underset{N}{_b*_\alpha}id\underset{\nu}{_\beta*_\alpha}\Phi)
(\mathfrak a\underset{N}{_b*_\alpha}id)\mathfrak a(x)\\
&=&
\omega(id\underset{N}{_b*_\alpha}id\underset{N}{_\beta*_\alpha}\Phi)
(id\underset{N}{_b*_\alpha}\Gamma)\mathfrak a(x)
\end{eqnarray*}
which, thanks to \ref{thinv}(ii), is equal to :
\[\omega((id\underset{\nu}{_b*_\alpha}\Phi)\mathfrak a(x)\underset{N}{_\beta\otimes_\alpha}1)\]
from which we get that $\mathfrak a(T_\mathfrak a(x))=T_\mathfrak a(x)\underset{N}{_b\otimes_\alpha}1$, which is the result.  \end{proof}

%%%%%%defovw1
\subsection{Proposition}
\label{defovw1}
{\it Let $(N, M, \alpha, \beta, \Gamma, T, T', \nu)$ be a measured quantum groupoid, and let $A$ be a von Neumann algebra, $(b, \mathfrak a)$ an action of $(N, M, \alpha, \beta, \Gamma, T, T', \nu)$ on $A$. Let $\Phi=\nu\circ\alpha^{-1}\circ T$, and, for any $x\in A^+$, let $T_\mathfrak a(x)=(id\underset{N}{_b*_\alpha}\Phi)\mathfrak a(x)$ be the extended positive element of $A^\mathfrak a$ defined in \ref{defovw}. Then $T_\mathfrak a$ is a normal faithful operator-valued weight from $A$ on $A^\mathfrak a$. }

\begin{proof}
Let $y$ in $A^\mathfrak a$; let $\omega\in A_*^+$, such that there exists $k>0$ for which $\omega\circ b\leq k\nu$; as $y\in b(N)'$, the element $\omega_y$ of $A^+_*$ defined, for all $x\in A$ by $\omega_y(x)=\omega(y^*xy)$ satisfies the same property, and we have, for $x\in A^+$:
\begin{eqnarray*}
<T_\mathfrak a(y^*xy), \omega>&=&
\omega((id\underset{N}{_b*_\alpha}\Phi)\mathfrak a(y^*xy))\\
&=&\omega((id\underset{N}{_b*_\alpha}\Phi)(y^*\underset{N}{_b\otimes_\alpha}1)\mathfrak a(x)(y\underset{N}{_b\otimes_\alpha}1))\\
&=&\omega_y(T_\mathfrak a(x))\\
&=&<y^*T_\mathfrak a(x)y, \omega>
\end{eqnarray*}
which proves that $T_\mathfrak a$ is indeed an operator-valued weight. Normality and faithfulness are trivial. \end{proof}

%%%%defintegrable
\subsection{Definition}
\label{defintegrable}
Let $\gG=(N, M, \alpha, \beta, \Gamma, T, T', \nu)$ be a measured quantum groupo-id, and let $A$ be a von Neumann algebra, $(b, \mathfrak a)$ an action of $\gG$ on $A$; we shall say that the action is integrable if the operator-valued weight $T_\mathfrak a$ introduced in \ref{defovw} and \ref{defovw1} is semi-finite.

\section{Some technical properties of actions}
\label{tech}

In this chapter, we define technical properties (property A in \ref{defA}, saturation property in \ref{defsat}) and prove technical results (\ref{propsat}, \ref{propstandard}) about these properties which will be necessary in the sequel. We shall prove in chapter \ref{biduality} that these properties are always fulfilled.

%%%%defpropA
\subsection{Definition}
\label{defA}
Let $\gG=(N, M, \alpha, \beta, \Gamma, T, T', \nu)$ be a measured quantum groupoid, $(b, \mathfrak a)$ an action of $\gG$ on a von Neumann algebra $A$; following (\cite{E1}, III2), we shall say that the action $\mathfrak a$ satisfies the property (A) if :
\[(\mathfrak a (A)\cup 1\underset{N}{_b\otimes_\alpha}\alpha (N)')''=A\underset{N}{_b*_\alpha}\mathcal L (H_\Phi)\]

%%%%GammaA
\subsection{Proposition}
\label{GammaA}
{\it Let $\gG=(N, M, \alpha, \beta, \Gamma, T, T', \nu)$ be a measured quantum groupoid; then, we have :
\[(\Gamma (M)\cup 1\underset{N}{_\beta\otimes_\alpha}\alpha(N)')''=M\underset{N}{_\beta*_\alpha}\mathcal L(H_\Phi)\]
In other words, the action $(\beta, \Gamma)$ of $\gG$ on $M$ (\ref{Gamma}) satisfies the property (A) (\ref{defA}). }

\begin{proof}
Let $y$ be in $\beta(N)'$ such that $W(y\underset{N}{_\beta\otimes_\alpha}1)W^*$ belongs to $\mathcal L(H_\Phi)\underset{N^o}{_\alpha*_{\hat{\beta}}}M'$; then, for any $u$ unitary in $M$, we have :
\[(1\underset{N^o}{_\alpha\otimes_{\hat{\beta}}}u)W(y\underset{N}{_\beta\otimes_\alpha}1)W^*(1\underset{N^o}{_\alpha\otimes_{\hat{\beta}}}u^*)=W(y\underset{N}{_\beta\otimes_\alpha}1)W^*\]
which can be written :
\[\Gamma(u)(y\underset{N}{_\beta\otimes_\alpha}1)\Gamma(u^*)=y\underset{N}{_\beta\otimes_\alpha}1\]
and, therefore, we have $\Gamma(u)(y\underset{N}{_\beta\otimes_\alpha}1)=(y\underset{N}{_\beta\otimes_\alpha}1)\Gamma(u)$; therefore, $y$ commutes with all elements of the form $(id\underset{N}{_\alpha*_\beta}\omega)\Gamma(u)$, for all $u$ unitary in $M$ and all $\omega\in M_*$ for which there exists $k>0$ such that $\omega\circ\alpha\leq k\nu$, and, using \ref{dens}, we get that $y$ belongs to $M'$. So, we have :
\[(\beta(N)'\underset{N}{_\beta\otimes_\alpha}1)\cap W^*(\mathcal L(H_\Phi)\underset{N^o}{_\alpha*_{\hat{\beta}}}M')W=M'\underset{N}{_\beta\otimes_\alpha}1\]
and, as $(\beta(N)'\underset{N}{_\beta\otimes_\alpha}1)=\mathcal L(H_\Phi)\underset{N}{_\beta*_\alpha}\alpha (N)=(1\underset{N}{_\beta\otimes_\alpha}\alpha (N)')'$, we get the result, taking the commutants. \end{proof}

%%%%%defsat
\subsection{Definition}
\label{defsat}
Let $\gG=(N, M, \alpha, \beta, \Gamma, T, T', \nu)$ be a measured quantum groupoid, $(b, \mathfrak a)$ an action of $\gG$ on a von Neumann algebra $A$; following (\cite{ES1}, II8), we shall denote:
\[Sat\mathfrak a=\{X\in A\underset{N}{_b*_\alpha}M, (\mathfrak a \underset{N}{_b*_\alpha}id)(X)=(id\underset{N}{_b*_\alpha}\Gamma)(X)\}\]
By definition, we have $\mathfrak a (A)\subset Sat \mathfrak a$, and we shall say that $\mathfrak a$ is saturated if $\mathfrak a (A)= Sat \mathfrak a$

%%%%lemV
\subsection{Lemma}
\label{lemV}
{\it Let $\gG=(N, M, \alpha, \beta, \Gamma, T, T', \nu)$ be a measured quantum groupoid, $_a\gH_b$ a $N-N$ bimodule, $V$ a corepresentation of $\gG$ on $_a\gH_b$; let $Sat V$ be 
the set of elements $X$ in $a(N)'\underset{N}{_b*_\alpha}\hat{\beta}(N)'$ such that :
\begin{multline*}
(V\underset{N}{_b\otimes_\alpha}1)\sigma^{2,3}_{b,a}(X\underset{N^o}{_a\otimes_\beta}1)\sigma^{2,3}_{a,b}(V^*\underset{N}{_\beta\otimes_\alpha}1)
\\
=(1\underset{N}{_b\otimes_\alpha}W^*)(1\underset{N}{_b\otimes_\alpha}\sigma_\nu)(X\underset{N}{_{\hat{\beta}}\otimes_\alpha}1)(1\underset{N}{_b\otimes_\alpha}\sigma_\nu)(1\underset{N}{_b\otimes_\alpha}W)
\end{multline*}
Then, we have :
\newline
(i) $Sat V=V(\mathcal L(\gH)\underset{N^o}{_a*_\beta}\widehat{M}')V^*$; 
\newline
(ii) $V^*(1\underset{N}{_b\otimes_\alpha}\widehat{M}')V\subset \mathcal L(\gH)\underset{N^o}{_a*_\beta}\widehat{M}'$}

\begin{proof}
As $V$ is a corepresentation, we have, by definition : 
\begin{multline*}
(1\underset{N}{_b\otimes_\alpha}W^*)(1\underset{N}{_b\otimes_\alpha}\sigma_\nu)(V\underset{N^o}{_{\hat{\beta}}\otimes_\alpha}1)(1\underset{N^o}{_a\otimes_\beta}\sigma_{\nu^o})(1\underset{N^o}{_a\otimes_\beta}W)\\
=(V\underset{N}{_b\otimes_\alpha}1)\sigma^{2,3}_{b,a}(V\underset{N^o}{_\alpha\otimes_\beta}1)(1\underset{N^o}{_a\otimes_\beta}\sigma_\nu)
\end{multline*}
and, therefore : 
\begin{multline*}
(1\underset{N}{_b\otimes_\alpha}\sigma_{\nu^o})(1\underset{N}{_b\otimes_\alpha}W)(V\underset{N}{_b\otimes_\alpha}1)\sigma^{2,3}_{b,a}\\
=(V\underset{N^o}{_{\hat{\beta}}\otimes_\alpha}1)(1\underset{N^o}{_a\otimes_\beta}\sigma_{\nu^o})(1\underset{N^o}{_a\otimes_\beta}W)(1\underset{N^o}{_\alpha\otimes_\beta}\sigma_{\nu^o})(V^*\underset{N^o}{_\alpha\otimes_\beta}1)
\end{multline*}
So, $X$ belongs to $Sat V$ if and only if we have :
\[(1\underset{N^o}{_\alpha\otimes_\beta}\sigma_{\nu^o}W\sigma_{\nu^o})(V^*XV\underset{N^o}{_\alpha\otimes_\beta}1)=(V^*XV\underset{N}{_{\hat{\beta}}\otimes_\alpha}1)(1\underset{N^o}{_\alpha\otimes_\beta}\sigma_{\nu^o}W\sigma_{\nu^o})\]
from which we get that $X$ belongs to $Sat V$ if and only if $V^*XV$ belongs to $\mathcal L(\gH)\underset{N^o}{_a*_\beta}\widehat{M}'$, which finishes the proof of (i).
\newline
As it is clear that $1\underset{N}{_b\otimes_\alpha}\widehat{M}'\subset Sat V$, we get (ii) from (i). \end{proof}

%%%%%Vsat
\subsection{Lemma}
\label{Vsat}
{\it Let $\gG=(N, M, \alpha, \beta, \Gamma, T, T', \nu)$ be a measured quantum groupoid, $_a\gH_b$ a $N-N$ bimodule, $A$ a von Neumann algebra such that $b(N)\subset A\subset a(N)'$, $V$ a corepresentation of $\gG$ on $_a\gH_b$ which implements an action $(b, \mathfrak a)$ on $A$; then we have :
\[Sat\mathfrak a = A\underset{N}{_b*_\alpha}M\cap V(\mathcal L(\gH)\underset{N^o}{_a*_\beta}\widehat{M}')V^*\]
and also :}
\[Sat\mathfrak a = A\underset{N}{_b*_\alpha}\mathcal L(H_\Phi)\cap V(a(N)'\underset{N^o}{_a\otimes_\beta}1)V^*\]
\begin{proof}
Let $X\in A\underset{N}{_b*_\alpha}M$; as we have :
\[(\mathfrak a\underset{N}{_b*_\alpha}id)(X)=(V\underset{N}{_b\otimes_\alpha}1)\sigma^{2,3}_{b,a}(X\underset{N^o}{_a\otimes_\beta}1)\sigma^{2,3}_{a,b}(V^*\underset{N}{_\beta\otimes_\alpha}1)\]
and :
\[(id\underset{N}{_b*_\alpha}\Gamma)(X)=(1\underset{N}{_b\otimes_\alpha}W^*)(1\underset{N}{_b\otimes_\alpha}\sigma_\nu)(X\underset{N}{_{\hat{\beta}}\otimes_\alpha}1)(1\underset{N}{_b\otimes_\alpha}\sigma_\nu)(1\underset{N}{_b\otimes_\alpha}W)\]
we get that $Sat\mathfrak a=A\underset{N}{_b*_\alpha}M\cap Sat V$, which, thanks to \ref{lemV}(i), gives the first formula. 
\newline
Starting from this formula, we get, because $V$ is a corepresentation of $\gG$ :
\begin{eqnarray*}
Sat\mathfrak a &=&V(V^*( A\underset{N}{_b*_\alpha}M)V\cap \mathcal L(\gH)\underset{N^o}{_a*_\beta}\widehat{M}')V^*\\
&\subset&V(\mathcal L(\gH)\underset{N}{_a*_\beta}M\cap \mathcal L(\gH)\underset{N^o}{_a*_\beta}\widehat{M}')V^*\\
&=&V(\mathcal L(\gH)\underset{N}{_a*_\beta}\beta(N))V^*\\
&=&V(a(N)'\underset{N^o}{_a\otimes_\beta}1)V^*
\end{eqnarray*}
and, therefore, $Sat\mathfrak a\subset A\underset{N}{_b*_\alpha}M\cap V(a(N)'\underset{N^o}{_a\otimes_\beta}1)V^*$. 
\newline
As $a(N)'\underset{N^o}{_a\otimes_\beta}1\subset \mathcal L(\gH)\underset{N^o}{_a*_\beta}\widehat{M}'$, using the first formula, we get :
\[Sat\mathfrak a=A\underset{N}{_b*_\alpha}M\cap V(a(N)'\underset{N^o}{_a\otimes_\beta}1)V^*\]
And, as $V(a(N)'\underset{N^o}{_a\otimes_\beta}1)V^*\subset \mathcal L(\gH)\underset{N}{_b*_\alpha}M$, we get the second formula. \end{proof}

%%%%%propsat
\subsection{Proposition}
\label{propsat}
{\it Let $\gG=(N, M, \alpha, \beta, \Gamma, T, T', \nu)$ be a measured quantum groupoid, $_a\gH_b$ a $N-N$ bimodule, $A$ a von Neumann algebra such that $b(N)\subset A\subset a(N)'$, $V$ a corepresentation of $\gG$ on $_a\gH_b$ which implements an action $(b, \mathfrak a)$ of $\gG$ on $A$ and an action $(a, \mathfrak a')$ of $\gG^o$ on $A'$. Then, are equivalent :
\newline
(i) $\mathfrak a$ is saturated;
\newline
(ii) $\mathfrak a'$ satisfies property (A). }

\begin{proof}
Using \ref{Vsat}, we get that $\mathfrak a$ is saturated if and only if :
\[V(A\underset{N^o}{_a\otimes_\beta}1)V^*=A\underset{N}{_b*_\alpha}\mathcal L(H_\Phi)\cap V(a(N)'\underset{N^o}{_a\otimes_\beta}1)V^*\]
which, taking the commutants, is equivalent to :
\[V(A'\underset{N^o}{_a*_\beta}\mathcal L(H_\Phi))V^*=(A'\underset{N}{_b\otimes_\alpha}1\cup V(a(N)\underset{N^o}{_a*_\beta}\mathcal L(H_\Phi))V^*)''\]
or, to :
\[A'\underset{N^o}{_a*_\beta}\mathcal L(H_\Phi)=(V^*(A'\underset{N}{_b\otimes_\alpha}1)V\cup a(N)\underset{N^o}{_a*_\beta}\mathcal L(H_\Phi))''\]
which is :
\[A'\underset{N^o}{_a*_\beta}\mathcal L(H_\Phi)=(\mathfrak a'(A')\cup (1\underset{N^o}{_a\otimes_\beta}\beta(N)'))''\]
which means that $\mathfrak a'$ satisfies property (A). The converse is proved the same way. \end{proof}

%%%%propstandard
\subsection{Proposition}
\label{propstandard}
{\it Let $\gG=(N, M, \alpha, \beta, \Gamma, T, T', \nu)$ be a measured quantum groupoid, $A$ a von Neumann algebra, $(b, \mathfrak a)$ an action of $\gG$ on $A$; let $\psi$ be a normal semi-finite faithful weight on $A$, and let us suppose that there exists a corepresentation $V$ on $H_\psi$ which is a standard implementation for $\mathfrak a$ in the sense of \ref{standard}. Then, are equivalent :
\newline
(i) $\mathfrak a$ is saturated;
\newline
(ii) $\mathfrak a$ satisfies property (A). }

\begin{proof}
By \ref{standard}, for all $n\in N$, we write $a(n)=J_\psi b(n^*)J_\psi$, and $V$ is a corepresentation of 
$\gG$ on $_a(H_\psi)_b$, which implements $\mathfrak a$, and verifies :
\[V^*=(J_\psi\underset{\nu^o}{_\alpha\otimes_\beta}J_{\widehat{\Phi}})V(J_\psi\underset{\nu}{_b\otimes_\alpha} J_{\widehat{\Phi}})\]
from which we get that the action $\mathfrak a'$ implemented by $V^*$ satisfies :
\[(J_\psi\underset{\nu^o}{_\alpha\otimes_\beta}J_{\widehat{\Phi}})\mathfrak a'(A')(J_\psi\underset{\nu}{_b\otimes_\alpha} J_{\widehat{\Phi}})=\mathfrak a (A)\]
As $(J_\psi\underset{\nu^o}{_\alpha\otimes_\beta}J_{\widehat{\Phi}})(A\underset{N^o}{_a\otimes_\beta}1)(J_\psi\underset{\nu}{_b\otimes_\alpha} J_{\widehat{\Phi}})=A'\underset{N}{_b\otimes_\alpha}1$, we get that :
\[(J_\psi\underset{\nu^o}{_\alpha\otimes_\beta}J_{\widehat{\Phi}})(A'\underset{N^o}{_a*_\beta}\mathcal L(H_\Phi))(J_\psi\underset{\nu}{_b\otimes_\alpha} J_{\widehat{\Phi}})=A\underset{N}{_b*_\alpha}\mathcal L(H_\Phi)\]
Thanks to \ref{propsat}, $\mathfrak a$ is saturated if and only if $\mathfrak a'$ satisfies property (A), which means that :
\[A'\underset{N^o}{_a*_\beta}\mathcal L(H_\Phi)=(\mathfrak a'(A')\cup (1\underset{N^o}{_a\otimes_\beta}\beta(N)'))''\]
from which we get that :
\[A\underset{N}{_b*_\alpha}\mathcal L(H_\Phi)=(\mathfrak a (A)\cup (1\underset{N}{_b\otimes_\alpha}\alpha(N)')''\]
which means that $\mathfrak a$ satisfies property (A). \end{proof}

%%%%corGamma
\subsection{Corollary}
\label{corGamma}
{\it Let $\gG=(N, M, \alpha, \beta, \Gamma, T, T', \nu)$ be a measured quantum groupoid; then, we have :}
\[\Gamma (M)=\{X\in M\underset{N}{_\beta*_\alpha}M; (\Gamma\underset{N}{_\beta*_\alpha}id)(X)=(id\underset{N}{_\beta*_\alpha}\Gamma)(X)\}\]

\begin{proof}
We have seen in \ref{Gamma} that $\Gamma$ can be considered as an action of $\gG$ on $M$, which has a standard implementation; moreover, by \ref{GammaA}, it satisfies property (A), and, therefore, by \ref{propstandard}, this action is saturated, which finishes the proof. \end{proof}

%%%%%deltainv
\section{The standard implementation of an action : 
the case of a $\delta$-invariant weight}
\label{deltainv}
In this chapter, we define $\delta$-invariant weights on a von Neumann algebra on which is acting a measured quantum groupoid (\ref{defdeltainv}); in that situation, we construct (\ref{thV2}) a standard implementation of the action. As a corollary, we obtain in that case that the property A and the saturation property introduced in chapter \ref{tech} are equivalent (\ref{corstandardV}). 

%%%%%defdeltainv
\subsection{Definition}
\label{defdeltainv}
Let $\gG=(N, M, \alpha, \beta, \Gamma, T, T', \nu)$ be a measured quantum groupoid; let $A$ be a von Neumann algebra, $(b, \mathfrak a)$ an action of $\gG$ on $A$ and let $\psi$ be a normal faithful semi-finite weight on $A$. We shall identify $A$ and $\pi_\psi (A)$, for simplification. Let $a$ be the representation of $N$ on $H_\Psi$ given by :
\[a(n)=J_\psi b(n^*)J_\psi\]
We shall say that $\psi$ is $\delta$-invariant if, for all $\eta\in D(_\alpha H_\Phi, \nu)\cap\mathcal D(\delta^{1/2})$, such that $\delta^{1/2}\eta$ belongs to $D((H_\Phi)_\beta, \nu^o)$, and for all $x\in\gN_\psi$, we have:
\[\psi((id\underset{N}{_b*_\alpha}\omega_\eta)\mathfrak a(x^*x))=\|\Lambda_\psi(x)\underset{\nu^o}{_a\otimes_\beta}\delta^{1/2}\eta\|^2\]

One should remark that $(id\underset{N}{_b*_\alpha}\omega_\eta)\mathfrak a(x^*x)$ belongs then to $\gM_\psi^+$. 
\newline
Moreover, if the subset $D((H_\psi)_b, \nu^o)\cap D(_aH_\psi, \nu)$ is dense in $H_\psi$, we shall say that $\psi$ bears the density property, 

%%%%Phideltainv
\subsection{Example}
\label{Phideltainv}
Let us consider (\ref{Gamma}) the action $(\beta, \Gamma)$ of $\gG$ on $M$; then, by \ref{propW}, we get that $\Phi$ is $\delta$-invariant. Namely, taking $(e_i)_{i\in I}$ an orthonormal $(\alpha, \nu)$ basis of $H_\Phi$, we get, for $x\in\gN_\Phi$ and convenient $\eta$, that :
\begin{eqnarray*}
\Phi[(id\underset{N}{_\beta*_\alpha}\omega_\eta)\Gamma(x^*x)]
&=&\sum_i\Phi((id\underset{N}{_\beta*_\alpha}\omega_{\eta, e_i})\Gamma(x)^*(id\underset{N}{_\beta*_\alpha}\omega_{\eta, e_i})\Gamma(x))\\
&=&\sum_i\|(id*\omega_{\delta^{1/2}\eta, e_i})(\widehat{W^o})\Lambda_\Phi(x)\|^2\\
&=&\|\sum_i(1\underset{N}{_\beta\otimes_\alpha}\theta^{\alpha, \nu}(e_i, e_i))\widehat{W^o}(\Lambda_\Phi(x)\underset{\nu^o}{_{\hat{\alpha}}\otimes_\beta}\delta^{1/2}\eta)\|^2\\
&=&\|\Lambda_\Phi(x)\underset{\nu^o}{_{\hat{\alpha}}\otimes_\beta}\delta^{1/2}\eta\|^2
\end{eqnarray*}
Moreover, as $D((H_\Phi)_\beta,\nu^o)\cap D(_{\hat{\alpha}}H_\Phi, \nu)$ is dense in $H_\Phi$, $\Phi$ bears the density property.

%%%%%lemV1
\subsection{Lemma}
\label{lemV1}
{\it Let $\gG=(N, M, \alpha, \beta, \Gamma, T, T', \nu)$ be a measured quantum groupoid; let $A$ be a von Neumann algebra, $(b, \mathfrak a)$ an action of $\gG$ on $A$ and let $\psi$ be a $\delta$-invariant normal faithful semi-finite weight on $A$. Let $(e_i)_{i\in I}$ be an $(\alpha, \nu)$-orthogonal basis of $H_\Phi$. Let $x\in \gN_\psi$, and $\eta$ as in \ref{defdeltainv}. Then:
\newline
(i) for all $\xi\in D(_\alpha H_\Phi, \nu)$, $(id\underset{N}{_b*_\alpha}\omega_{\eta, \xi})\mathfrak a(x)$ belongs to $\gN_\psi$; 
\newline
(ii) the sum $\sum_i\Lambda_\psi((id\underset{N}{_b*_\alpha}\omega_{\eta, e_i})\mathfrak a(x))\underset{\nu}{_b\otimes_\alpha}e_i$ is strongly converging, the limit does not depend upon the choice of the $(\alpha, \nu)$-orthogonal basis of $H_\Phi$, and :}
\[\|\sum_i\Lambda_\psi((id\underset{N}{_b*_\alpha}\omega_{\eta, e_i})\mathfrak a(x))\underset{\nu}{_b\otimes_\alpha}e_i\|^2=\psi((id\underset{N}{_b*_\alpha}\omega_\eta)\mathfrak a(x^*x))\]

\begin{proof}
We have :
\[(id\underset{N}{_b*_\alpha}\omega_{\eta, \xi})\mathfrak a(x)^*(id\underset{N}{_b*_\alpha}\omega_{\eta, \xi})\mathfrak a(x)
=(id\underset{N}{_b*_\alpha}\omega_\eta)(\mathfrak a(x^*)(1\underset{N}{_b\otimes_\alpha}\theta^{\alpha, \nu}(\xi, \xi))\mathfrak a(x))\]
which is less than $\|R^{\alpha, \nu}(\xi)\|^2(id\underset{N}{_b*_\alpha}\omega_\eta)\mathfrak a(x^*x)$, 
which belongs to $\gM_\psi^+$, by hypothesis; so we get (i). 
\newline
As the vectors $\Lambda_\psi((id\underset{N}{_b*_\alpha}\omega_{\eta, e_i})\mathfrak a(x))\underset{\nu}{_b\otimes_\alpha}e_i$ are two by two orthogonal, for any finite $J\subset I$, we have :
\begin{multline*}
\|\sum_{i\in J}\Lambda_\psi((id\underset{N}{_b*_\alpha}\omega_{\eta, e_i})\mathfrak a(x))\underset{\nu}{_b\otimes_\alpha}e_i\|^2\\
=\sum_{i\in J}\|\Lambda_\psi((id\underset{N}{_b*_\alpha}\omega_{\eta, e_i})\mathfrak a(x))\underset{\nu}{_b\otimes_\alpha}e_i\|^2\\
=\psi(id\underset{N}{_b*_\alpha}\omega_\eta)(\mathfrak a(x^*)(1\underset{N}{_b\otimes_\alpha}\sum_{i\in J}\theta^{\alpha, \nu}(e_i, e_i))\mathfrak a(x))
\end{multline*}
which is less than (and uprising to) $\psi(id\underset{N}{_b*_\alpha}\omega_\eta)\mathfrak a(x^*x)$. 
If $\xi\in D(_\alpha H_\Phi, \nu)$, we have :
\[(\rho^{b, \alpha}_\xi)^*(\sum_i\Lambda_\psi((id\underset{N}{_b*_\alpha}\omega_{\eta, e_i})\mathfrak a(x))\underset{\nu}{_b\otimes_\alpha}e_i)=\Lambda_\psi((id\underset{N}{_b*_\alpha}\omega_{\eta, \xi})\mathfrak a(x))\] which does not depend upon the choice of the basis, and finishes the proof of (ii). 
\end{proof}

%%%%propV2
\subsection{Proposition}
\label{propV2}
{\it Let $\gG=(N, M, \alpha, \beta, \Gamma, T, T', \nu)$ be a measured quantum groupoid; let $A$ be a von Neumann algebra, $(b, \mathfrak a)$ an action of $\gG$ on $A$ and let $\psi$ be a $\delta$-invariant normal faithful semi-finite weight on $A$. Let us define the representation $a$ of $N$ on $H_\psi$ by, for $n\in N$ :
\[a(n)=J_\Psi b(n^*)J_\psi\]
 Then :
\newline
(i) there exists an isometry $V_\psi$ from $H_\psi\underset{\nu^o}{_a\otimes_\beta}H_\Phi$ to $H_\psi\underset{\nu}{_b\otimes_\alpha}H_\Phi$ such that:
\[V_\psi(\Lambda_\psi (x)\underset{\nu^o}{_a\otimes_\beta}\delta^{1/2}\eta)=\sum_i\Lambda_\psi((id\underset{N}{_b*_\alpha}\omega_{\eta, e_i})\mathfrak a(x))\underset{\nu}{_b\otimes_\alpha}e_i\]
for all $x\in \gN_\psi$, $\eta$ as in \ref{defdeltainv} and for all $(\alpha, \nu)$-orthogonal basis $(e_i)_{i\in I}$ of $H_\Phi$.
\newline
(ii) for all $e\in\gN_\psi$, we have :
\[(J_\psi eJ_\psi\underset{\nu}{_b\otimes_\alpha}1)V_\psi(\Lambda_\psi (x)\underset{\nu^o}{_a\otimes_\beta}\delta^{1/2}\eta)=\mathfrak a(x)(J_\psi\Lambda_\psi (e)\underset{\nu}{_b\otimes_\alpha}\eta)\]
(iii) for all $\xi\in D(_\alpha H_\Phi, \nu)$, we have :
\[\Lambda_\psi((id\underset{N}{_b*_\alpha}\omega_{\eta, \xi})\mathfrak a(x))=(id*\omega_{\delta^{1/2}\eta, \xi})(V_\psi)\Lambda_\psi (x)\]
(iv) for all $y\in A$, $z\in M'$, we have :
\[\mathfrak a(y)V_\psi=V_\psi(y\underset{\nu^o}{_a\otimes_\beta}1)\]
\[(1\underset{N}{_b\otimes_\alpha}z)V_\psi=V_\psi(1\underset{N^o}{_a\otimes_\beta}z)\]
(v) for all $n\in N$; we have :}
\[(a(n)\underset{N}{_b\otimes_\alpha}1)V_\psi=V_\psi(1\underset{N^o}{_a\otimes_\beta}\alpha(n))\]
\[(1\underset{N}{_b\otimes_\alpha}\beta(n))V_\psi=V_\psi(b(n)\underset{N^o}{_a\otimes_\beta}1)\]
\[(1\underset{N}{_b\otimes_\alpha}\hat{\beta}(n))V_\psi=V_\psi(1\underset{N^o}{_a\otimes_\beta}\hat{\beta}(n))\]

\begin{proof}
Defining $V_\psi$ by the formula given in (i), we easily obtain, using \ref{lemV1}(ii), that for all $x'\in\gN_\psi$, $\eta'$ as in \ref{defdeltainv}, we get :
\[(V_\psi(\Lambda_\psi (x)\underset{\nu^o}{_a\otimes_\beta}\delta^{1/2}\eta)|V_\psi(\Lambda_\psi (x')\underset{\nu^o}{_a\otimes_\beta}\delta^{1/2}\eta'))=
\psi((id\underset{N}{_b*_\alpha}\omega_{\eta, \eta'})\mathfrak a(x'^*x))\]
which, by polarization of the definition of $\delta$-invariance, is equal to:
\[(\Lambda_\psi (x)\underset{\nu^o}{_a\otimes_\beta}\delta^{1/2}\eta|\Lambda_\psi (x')\underset{\nu^o}{_a\otimes_\beta}\delta^{1/2}\eta')\]
which implies that this formula defines an isometry which can be extended by continuity to $H_\psi\underset{\nu^o}{_a\otimes_\beta}H_\Phi$, and does not depend upon the choice of the basis, which is (i). We then get :
\begin{eqnarray*}
(J_\psi eJ_\psi\underset{\nu}{_b\otimes_\alpha}1)V_\psi(\Lambda_\psi (x)\underset{\nu^o}{_a\otimes_\beta}\delta^{1/2}\eta)
&=&
\sum_iJ_\psi eJ_\psi\Lambda_\psi((id\underset{N}{_b*_\alpha}\omega_{\eta, e_i})\mathfrak a(x))\underset{\nu}{_b\otimes_\alpha}e_i\\
&=&\sum_i(id\underset{N}{_b*_\alpha}\omega_{\eta, e_i})\mathfrak a(x))J_\psi\Lambda_\psi(e)\underset{\nu}{_b\otimes_\alpha}e_i\\
&=&\mathfrak a(x)(J_\psi\Lambda_\psi (e)\underset{\nu}{_b\otimes_\alpha}\eta)
\end{eqnarray*}
which is (ii). Using (ii), we then get :
\begin{eqnarray*}
J_\psi eJ_\psi\Lambda_\psi((id\underset{N}{_b*_\alpha}\omega_{\eta, \xi})\mathfrak a(x))
&=&(id\underset{N}{_b*_\alpha}\omega_{\eta, \xi})\mathfrak a(x)J_\psi\Lambda_\psi (e)\\
&=&J_\psi eJ_\psi(id*\omega_{\delta^{1/2}\eta, \xi})(V_\psi)\Lambda_\psi (x)
\end{eqnarray*}
from which we get (iii), by taking the limit when $e$ goes to $1$. Using (ii) again, we have also : 
\begin{eqnarray*}
(J_\psi eJ_\psi\underset{N}{_b\otimes_\alpha}1)\mathfrak a(y)V_\psi(\Lambda_\psi (x)\underset{\nu^o}{_a\otimes_\beta}\delta^{1/2}\eta)
&=&
\mathfrak a(y)(J_\psi eJ_\psi\underset{\nu}{_b\otimes_\alpha}1)V_\psi(\Lambda_\psi (x)\underset{\nu^o}{_a\otimes_\beta}\delta^{1/2}\eta)\\
&=&\mathfrak a(yx)(J_\psi\Lambda_\psi (e)\underset{\nu}{_b\otimes_\alpha}\eta)
\end{eqnarray*}
which, thanks to (ii) again, is equal to :
\[(J_\psi eJ_\psi\underset{N}{_b\otimes_\alpha}1)V_\psi(\Lambda_\psi (yx)\underset{\nu^o}{_a\otimes_\beta}\delta^{1/2}\eta)
=(J_\psi eJ_\psi\underset{N}{_b\otimes_\alpha}1)V_\psi(y\underset{N^o}{_a\otimes_\beta}1)(\Lambda_\psi (x)\underset{\nu^o}{_a\otimes_\beta}\delta^{1/2}\eta)\]
from which we get the first formula of (iv), by taking the limit when $e$ goes to $1$. 
We have also :
\begin{eqnarray*}
(J_\psi eJ_\psi\underset{N}{_b\otimes_\alpha}z)V_\psi(\Lambda_\psi (x)\underset{\nu^o}{_a\otimes_\beta}\delta^{1/2}\eta)
&=&
(1\underset{N}{_b\otimes_\alpha}z)\mathfrak a(x)(J_\psi\Lambda_\psi (e)\underset{\nu}{_b\otimes_\alpha}\eta)\\
&=&
\mathfrak a(x)(J_\psi\Lambda_\psi (e)\underset{\nu}{_b\otimes_\alpha}z\eta)
\end{eqnarray*}
and, thanks to (ii) again, we get :
\begin{eqnarray*}
(J_\psi eJ_\psi\underset{N}{_b\otimes_\alpha}1)V_\psi(\Lambda_\psi (x)\underset{\nu^o}{_a\otimes_\beta}\delta^{1/2}z\eta)
&=&
(J_\psi eJ_\psi\underset{N}{_b\otimes_\alpha}1)V_\psi(\Lambda_\psi (x)\underset{\nu^o}{_a\otimes_\beta}z\delta^{1/2}\eta)\\
&=&(J_\psi eJ_\psi\underset{N}{_b\otimes_\alpha}1)V_\psi(1\underset{N^o}{_a\otimes_\beta}z)(\Lambda_\psi (x)\underset{\nu^o}{_a\otimes_\beta}\delta^{1/2}\eta)
\end{eqnarray*}
from which we get the second formula of (iv), by taking the limit when $e$ goes to $1$. 
\newline
Let's take now $n\in N$, analytic with respect both to $\sigma^\nu$ and $\gamma$. We have :
\begin{eqnarray*}
(J_\psi eJ_\psi a(n)\underset{N}{_b\otimes_\alpha}1)V_\psi(\Lambda_\psi (x)\underset{\nu^o}{_a\otimes_\beta}\delta^{1/2}\eta)
&=&(J_\psi eb(n^*)J_\psi\underset{N}{_b\otimes_\alpha}1)V_\psi(\Lambda_\psi (x)\underset{\nu^o}{_a\otimes_\beta}\delta^{1/2}\eta)\\
&=&\mathfrak a(x)(J_\psi\Lambda_\psi(eb(n^*))\underset{\nu}{_b\otimes_\alpha}\eta)
\end{eqnarray*}
It is equal to :
\[\mathfrak a(x)(\sigma_{-i/2}^\psi(b(n))J_\psi\Lambda_\psi(e)\underset{\nu}{_b\otimes_\alpha}\eta)
=
 \mathfrak a(x)(b(\gamma_{-i/2}(n))J_\psi\Lambda_\psi(e)\underset{\nu}{_b\otimes_\alpha}\eta)\]
which is equal to $\mathfrak a(x)(J_\psi\Lambda_\psi(e)\underset{\nu}{_b\otimes_\alpha}\alpha(\sigma_{-i/2}^\nu\gamma_{-i/2}(n))\eta)$, and, thanks to (ii) again, to :
 \begin{multline*}
(J_\psi eJ_\psi\underset{N}{_b\otimes_\alpha}1)V_\psi(\Lambda_\psi (x)\underset{\nu^o}{_a\otimes_\beta}\delta^{1/2}\alpha(\sigma_{-i/2}^\nu\gamma_{-i/2}(n))\eta)\\ 
=(J_\psi eJ_\psi\underset{N}{_b\otimes_\alpha}1)V_\psi(\Lambda_\psi (x)\underset{\nu^o}{_a\otimes_\beta}\alpha(n)\delta^{1/2}\eta)\\ 
=(J_\psi eJ_\psi\underset{N}{_b\otimes_\alpha}1)V_\psi(1\underset{N^o}{_a\otimes_\beta}\alpha(n))(\Lambda_\psi (x)\underset{\nu^o}{_a\otimes_\beta}\delta^{1/2}\eta)
\end{multline*}
from which we get the first formula of (v), by taking the limit when $e$ goes to $1$. 
\newline
The other formulae of (v) are special cases of (iv), taking $y=b(n)$, and $z=\hat{\beta}(n)$. \end{proof}

%%%%%propV3
\subsection{Proposition}
\label{propV3}
{\it Let $\gG=(N, M, \alpha, \beta, \Gamma, T, T', \nu)$ be a measured quantum groupoid; let $A$ be a von Neumann algebra, $(b, \mathfrak a)$ an action of $\gG$ on $A$ and let $\psi$ be a $\delta$-invariant normal faithful semi-finite weight on $A$. Let us define the representation $a$ of $N$ on $H_\psi$ by, for $n\in N$ :
\[a(n)=J_\psi b(n^*)J_\psi\]
and let $V_\psi$ the isometry from $H_\psi\underset{\nu^o}{_a\otimes_\beta}H_\Phi$ to $H_\psi\underset{\nu}{_b\otimes_\alpha}H_\Phi$ constructed in \ref{propV2}, which satisfies the intertwining properties proved in \ref{propV2}(v). Then we have :}
\[(1\underset{N}{_b\otimes_\alpha}\sigma_{\nu^o}W^*\sigma_\nu) (V_\psi\underset{N}{_{\hat{\beta}}\otimes_\alpha}1)(1\underset{N^o}{_a\otimes_\beta}\sigma_{\nu^o}W\sigma_\nu)=
(1\underset{N}{_b\otimes_\alpha}\sigma_{\nu^o})(V_\psi\underset{N}{_b\otimes_\alpha}1)\sigma^{2,3}_{b,a}(V_\psi\underset{N^o}{_\alpha\otimes_\beta}1)\]

\begin{proof}
Let $x_1$, $x_2$ in $\gN_\psi$, $\xi_1$, $\xi_2$ in $D(_\alpha H_\Phi, \nu)$, and $\eta_1$, $\eta_2$ as in \ref{defdeltainv} $(i=1,2)$; then, the scalar product of the vector :
\[(V_\psi\underset{N}{_b\otimes_\alpha}1)(1\underset{N^o}{_a\otimes_\beta}\sigma_{\nu^o})(V_\psi\underset{N^o}{_\alpha\otimes_\beta}1)(1\underset{N^o}{_\alpha\otimes_\beta}\sigma_\nu)(\Lambda_\psi(x_1)\underset{\nu^o}{_a\otimes_\beta}(\delta^{1/2}\eta_1\underset{\nu}{_\beta\otimes_\alpha}\delta^{1/2}\eta_2))\]
with $\Lambda_{\psi}(x_2)\underset{\nu}{_b\otimes_\alpha}\xi_1\underset{\nu}{_\beta\otimes_\alpha}\xi_2$ is equal to :
\[(V_\psi((id*\omega_{\delta^{1/2}\eta_2, \xi_2})(V_\psi)\Lambda_\psi(x_1)\underset{\nu^o}{_a\otimes_\beta}\delta^{1/2}\eta_1)|\Lambda_\psi(x_2)\underset{\nu}{_b\otimes_\alpha}\xi_1)\]
which, thanks to \ref{propV2}(iii), is equal to :
\[(V_\psi(\Lambda_\psi[(id\underset{N}{_b*_\alpha}\omega_{\eta_2, \xi_2})\mathfrak a(x_1)]\underset{\nu^o}{_a\otimes_\beta}\delta^{1/2}\eta_1)|\Lambda_\psi(x_2)\underset{\nu}{_b\otimes_\alpha}\xi_1)\]
and, using \ref{propV2}(iii) again, is equal to :
\[(\Lambda_\psi[(id\underset{N}{_b*_\alpha}\omega_{\eta_1, \xi_1})\mathfrak a((id\underset{N}{_b*_\alpha}\omega_{\eta_2, \xi_2})\mathfrak a(x_1))]|\Lambda_\psi(x_2))\]
which, using \ref{defaction}, is equal to :
\[(\Lambda_\psi((id\underset{N}{_b*_\alpha}\omega_{\eta_1, \xi_1}\underset{N}{_\beta*_\alpha}\omega_{\eta_2, \xi_2})(id\underset{N}{_\beta*_\alpha}\Gamma)\mathfrak a(x_1))|\Lambda_\psi(x_2))\]
Let us define $\Omega\in M_*$ such that, for any $x\in M$, $\Omega (x)=(\omega_{\eta_1, \xi_1}\underset{N}{_\beta*_\alpha}\omega_{\eta_2, \xi_2})\Gamma(x)$. 
\newline 
Let $(e_i)_{i\in I}$ be an $(\alpha, \nu)$-orthogonal basis of $H_\Phi$. There exists a family of vectors $\xi_i$ in $H$, such that $\xi_i=\hat{\beta}(<e_i, e_i>_{\alpha, \nu})\xi_i$, and :
\[W(\xi_1\underset{\nu}{_\beta\otimes_\alpha}\xi_2)=\sum_i e_i\underset{\nu^o}{_\alpha\otimes_{\hat{\beta}}}\xi_i\]
Moreover, as, for $n\in\gN_\nu$ :
\begin{eqnarray*}
\sum_i\|\alpha(n)\xi_i\|^2
&=&\|(1\underset{N^o}{_\alpha\otimes_{\hat{\beta}}}\alpha (n))W(\xi_1\underset{\nu}{_\beta\otimes_\alpha}\xi_2)\|^2\\
&=&\|W(\alpha(n)\xi_1\underset{\nu}{_\beta\otimes_\alpha}\xi_2)\|^2\\
&\leq&\|R^{\alpha, \nu}(\xi_2)\|^2\|R^{\alpha, \nu}(\xi_1)\|^2\|\Lambda_\nu(n)\|^2
\end{eqnarray*}
we get that the vectors $\xi_i$ belong to $D(_\alpha H_\Phi, \nu)$. 
The same way, there exist vectors $\eta_i=\hat{\beta}(<e_i, e_i>_{\alpha, \nu})\eta_i\in D(_\alpha H_\Phi, \nu)$, such that :
\[W(\eta_1\underset{\nu}{_\beta\otimes_\alpha}\eta_2)=\sum_i e_i\underset{\nu^o}{_\alpha\otimes_{\hat{\beta}}}\eta_i\]
Moreover, as we have $\Gamma(\delta^{it})=\delta^{it}\underset{N}{_\beta\otimes_\alpha}\delta^{it}$, we get that the vectors $\eta_i$ belongs to $\mathcal D(\delta^{1/2})$, and that :
\[W(\delta^{1/2}\eta_1\underset{\nu}{_\beta\otimes_\alpha}\delta^{1/2}\eta_2)=\sum_i e_i\underset{\nu^o}{_\alpha\otimes_{\hat{\beta}}}\delta^{1/2}\eta_i\]
We have also :
\begin{eqnarray*}
\sum_i\|\beta(n)\delta^{1/2}\eta_i\|^2
&=&\|(1\underset{N^o}{_\alpha\otimes_{\hat{\beta}}}\beta (n))W(\delta^{1/2}\eta_1\underset{\nu}{_\beta\otimes_\alpha}\delta^{1/2}\eta_2)\|^2\\
&\leq&\|R^{\beta, \nu^o}(\delta^{1/2}\eta_1)\|^2\|R^{\beta, \nu^o}(\delta^{1/2}\eta_2)\|^2|\Lambda_\nu(n)\|^2
\end{eqnarray*}
from which we get that the vectors $\delta^{1/2}\eta_i$ belong to $D((H_\Phi)_\beta, \nu^o)$. 
From these results, we get that, for any $x\in M$, we have $\Omega (x)=\sum_i\omega_{\eta_i, \xi_i}(x)$, and, therefore, we have :
\[(\Lambda_\psi((id\underset{N}{_b*_\alpha}\omega_{\eta_1, \xi_1}\underset{N}{_\beta*_\alpha}\omega_{\eta_2, \xi_2})(id\underset{N}{_\beta*_\alpha}\Gamma)\mathfrak a(x_1))|\Lambda_\psi(x_2))=\sum_i(\Lambda_\psi((id\underset{N}{_b*_\alpha}\omega_{\eta_i, \xi_i})\mathfrak a (x_1))|\Lambda_\psi (x_2))\]
The properties of the $\xi_i$s and the $\eta_i$s allow us to use again \ref{propV3}(iii), and this is equal to :
\[\sum_i((id*\omega_{\delta^{-1/2}\eta_i, \xi_i})(V_\psi)\Lambda_\psi (x_1)|\Lambda_\psi (x_2))\]
On the other hand, if $\zeta_1$ belongs to $D(_a H_\psi, \nu)$ and $\zeta_2$ to $D((H_\psi)_b, \nu^o)$, we get, by \ref{propV3}(iv), that $(\omega_{\zeta_1, \zeta_2}*id)(V)$ belongs to $M$, and that :
\begin{eqnarray*}
\sum_i((id*\omega_{\delta^{1/2}\eta_i, \xi_i})(V_\psi)\zeta_1|\zeta_2)
&=&
\sum_i \omega_{\delta^{1/2}\eta_i, \xi_i}(\omega_{\zeta_1, \zeta_2}*id)(V_\psi)\\
&=&(\omega_{\delta^{1/2}\eta_1, \xi_1}\underset{N}{_\beta*_\alpha}\omega_{\delta^{1/2}\eta_2, \xi_2})\Gamma[(\omega_{\zeta_1, \zeta_2}*id)(V_\psi)]
\end{eqnarray*}
is equal to :
\[(W^*((1\underset{N^o}{_\alpha\otimes_{\hat{\beta}}}(\omega_{\zeta_1, \zeta_2}*id)(V_\psi))W(\delta^{1/2}\eta_1\underset{\nu}{_\beta\otimes_\alpha}\delta^{1/2}\eta_2)|\xi_1\underset{\nu}{_\beta\otimes_\alpha}\xi_2)\]
which is equal to :
\[((1\underset{N}{_b\otimes_\alpha}\sigma_{\nu^o}W^*\sigma_\nu) (V_\psi\underset{N}{_{\hat{\beta}}\otimes_\alpha}1)(1\underset{N^o}{_a\otimes_\beta}\sigma_{\nu^o}W\sigma_\nu)(\zeta_1\underset{\nu^o}{_a\otimes_\beta}\delta^{1/2}\eta_2\underset{\nu^o}{_\alpha\otimes_\beta}\delta^{1/2}\eta_1)|\zeta_2\underset{\nu^o}{_a\otimes_\beta}\xi_2\underset{\nu^o}{_\alpha\otimes_\beta}\xi_1)\]
From which, by continuity, we get that the scalar product of the vector :
\[(V_\psi\underset{N}{_b\otimes_\alpha}1)\sigma^{2,3}_{b,a}(V_\psi\underset{N^o}{_\alpha\otimes_\beta}1)(1\underset{N^o}{_\alpha\otimes_\beta}\sigma_\nu)(\Lambda_\psi(x_1)\underset{\nu^o}{_a\otimes_\beta}(\delta^{1/2}\eta_1\underset{\nu}{_\beta\otimes_\alpha}\delta^{1/2}\eta_2))\]
with $\Lambda_{\psi}(x_2)\underset{\nu}{_b\otimes_\alpha}\xi_1\underset{\nu}{_\beta\otimes_\alpha}\xi_2$ is equal to the scalar product of the vector :
\[(1\underset{N}{_b\otimes_\alpha}\sigma_{\nu^o}W^*\sigma_\nu) (V_\psi\underset{N}{_{\hat{\beta}}\otimes_\alpha}1)(1\underset{N^o}{_a\otimes_\beta}\sigma_{\nu^o}W\sigma_\nu)(\Lambda_\psi(x_1)\underset{\nu^o}{_a\otimes_\beta}\delta^{1/2}\eta_2\underset{\nu^o}{_\alpha\otimes_\beta}\delta^{1/2}\eta_1)\]
with $\Lambda_{\psi}(x_2)\underset{\nu}{_b\otimes_\alpha}\xi_1\underset{\nu}{_\beta\otimes_\alpha}\xi_2$, from which, by continuity, we get the result. \end{proof}

%%%%thV1
\subsection{Theorem}
\label{thV1}
{\it Let $\gG=(N, M, \alpha, \beta, \Gamma, T, T', \nu)$ be a measured quantum groupoid; let $A$ be a von Neumann algebra, $(b, \mathfrak a)$ an action of $\gG$ on $A$ and let $\psi$ be a $\delta$-invariant normal faithful semi-finite weight on $A$. Let us define the representation $a$ of $N$ on $H_\psi$ by, for $n\in N$ :
\[a(n)=J_\psi b(n^*)J_\psi\]
and let $V_\psi$ the isometry from $H_\psi\underset{\nu^o}{_a\otimes_\beta}H_\Phi$ to $H_\psi\underset{\nu}{_b\otimes_\alpha}H_\Phi$ constructed in \ref{propV2}; then :
\newline
(i) $V_\psi$ is a unitary, and, for all $x$ in $A$, we have $\mathfrak a(x)=V_\psi(x\underset{N^o}{_a\otimes_\beta}1)V_\psi^*$; 
\newline
(ii) $V_\psi$ is a corepresentation of $\gG$ on the bimodule $_a(H_\psi)_b$ which implements $\mathfrak a$. }

\begin{proof}
Thanks to \ref{propV2}(i), (v), (iv) and \ref{propV3}, we can apply \ref{isom}, which says that $V$ is a unitary and a corepresentation, which gives (ii), and (i) by \ref{propV2}(iv) again.  \end{proof}

%%%lemV2
\subsection{Lemma}
\label{lemV2}
{\it Let $\gG=(N, M, \alpha, \beta, \Gamma, T, T', \nu)$ be a measured quantum groupoid; let $A$ be a von Neumann algebra, $(b, \mathfrak a)$ an action of $\gG$ on $A$ and let $\psi$ be a $\delta$-invariant normal faithful semi-finite weight on $A$. Let us define the representation $a$ of $N$ on $H_\psi$ such that, for $n\in N$:
\[a(n)=J_\psi b(n^*)J_\psi\]
and let us suppose that $\psi$ bears the density property, as defined in \ref{defdeltainv}.
Let $V_\psi$ be the corepresentation of $\gG$ on the bimodule $_a(H_\psi)_b$ constructed in \ref{thV1}. 
Then, for all $x$, $y$ in $\gN_\psi\cap \gN_\psi^*$, analytic with respect to $\psi$, and $\xi$ in $D((H_\Phi)_\beta, \nu^o)\cap\mathcal D((P\delta)^{-1})$ such that $(P\delta)^{-1}\xi$ belongs to $D((H_\Phi)_\beta, \nu^o)$ and $\eta$ in $D(_\alpha H_\Phi, \nu)\cap\mathcal D((P\delta)^{-1})$ such that $(P\delta)^{-1}\eta$ belongs to $D(_\alpha H_\Phi, \nu)$,we have :}
\[(V_\psi(\Delta_\psi\Lambda_\psi(x)\underset{\nu^o}{_a\otimes_\beta}(P\delta)^{-1}\xi)|
\Lambda_\psi(y)\underset{\nu}{_b\otimes_\alpha}\eta)
=(V_\psi(\Lambda_\psi(x)\underset{\nu^o}{_a\otimes_\beta}\xi)|\Delta_\psi
\Lambda_\psi (y)\underset{\nu}{_b\otimes_\alpha}(P\delta)^{-1}\eta)\]

\begin{proof}
Thanks to the density property, we can use \ref{thS2}; therefore, with $\xi$, $\eta$ in $\mathcal E_\tau$, using \ref{lemE1}(iii) :
\begin{multline*}
(V_\psi(\Delta_\psi\Lambda_\psi(x)\underset{\nu^o}{_a\otimes_\beta}P^{-1}\delta^{-1}\xi)|
\Lambda_\psi(y)\underset{\nu}{_b\otimes_\alpha}\eta)\\
=
(V_\psi^*(\Delta_\psi\Lambda_\psi(x)\underset{\nu}{_b\otimes_\alpha}J_{\widehat{\Phi}}P^{-1/2}\eta)|
\Lambda_\psi(y)\underset{\nu^o}{_a\otimes_\beta}J_{\widehat{\Phi}}P^{-1/2}\delta^{-1}\xi)\\
=
(\Delta_\psi\Lambda_\psi(x)\underset{\nu}{_b\otimes_\alpha}J_{\widehat{\Phi}}P^{-1/2}\eta|
V_\psi(\Lambda_\psi(y)\underset{\nu^o}{_a\otimes_\beta}J_{\widehat{\Phi}}P^{-1/2}\delta^{-1}\xi))
\end{multline*}
which, thanks to \ref{propV2}(iii), is equal to :
\begin{multline*}
(\Delta_\psi\Lambda_\psi(x)|(id*\omega_{J_{\widehat{\Phi}}P^{-1/2}\delta^{-1}\xi, J_{\widehat{\Phi}}P^{-1/2}\eta})(V_\psi)\Lambda_\psi(y))\\
=(\Delta_\psi\Lambda_\psi(x)|\Lambda_\psi[(id*\omega_{J_{\widehat{\Phi}}P^{-1/2}\delta^{-1/2}\xi, J_{\widehat{\Phi}}P^{-1/2}\eta})\mathfrak a(y)])\\
=(\Lambda_\psi[(id*\omega_{J_{\widehat{\Phi}}P^{-1/2}\eta, J_{\widehat{\Phi}}P^{-1/2}\delta^{-1/2}\xi})\mathfrak a(y^*)]|\Lambda_\psi(x^*))\\
\end{multline*}
which, using again \ref{propV2}(iii), is equal to :
\begin{multline*}
((id*\omega_{J_{\widehat{\Phi}}P^{-1/2}\delta^{-1/2}\eta, J_{\widehat{\Phi}}P^{-1/2}\delta^{-1/2}\xi})(V_\psi)\Lambda_\psi(y^*)|\Lambda_\psi(x^*))\\ 
=(V_\psi(\Lambda_\psi(y^*)\underset{\nu^o}{_a\otimes_\beta}J_{\widehat{\Phi}}P^{-1/2}\delta^{-1/2}\eta)|\Lambda_\psi(x^*)\underset{\nu}{_b\otimes_\alpha}J_{\widehat{\Phi}}P^{-1/2}\delta^{-1/2}\xi)\\
=(\Lambda_\psi(y^*)\underset{\nu^o}{_a\otimes_\beta}J_{\widehat{\Phi}}P^{-1/2}\delta^{-1/2}\eta|V_\psi^*(\Lambda_\psi(x^*)\underset{\nu}{_b\otimes_\alpha}J_{\widehat{\Phi}}P^{-1/2}\delta^{-1/2}\xi))
\end{multline*}
which, using again \ref{thS2}, is equal to :
\[(\Lambda_\psi(y^*)\underset{\nu}{_b\otimes_\alpha}\delta^{-1/2}\xi|
V_\psi(\Lambda_\psi(x^*)\underset{\nu^o}{_a\otimes_\beta}P^{-1}\delta^{-1/2}\eta))
=(\Lambda_\psi(y^*)|(id*\omega_{P^{-1}\delta^{-1/2}\eta, \delta^{-1/2}\xi})(V_\psi)\Lambda_\psi(x^*))\]
which, using again twice \ref{propV2}(iii), is equal to :
\begin{eqnarray*}
(\Lambda_\psi(y^*)|\Lambda_\psi[(id\underset{N}{_b*_\alpha}\omega_{(P\delta)^{-1}\eta, \delta^{-1/2}\xi})\mathfrak a (x^*)])
&=&(S_\psi\Lambda_\psi(y)|S_\psi\Lambda_\psi[(id\underset{N}{_b*_\alpha}\omega_{\delta^{-1/2}\xi, (P\delta)^{-1}\eta})\mathfrak a (x)])\\
&=&(\Delta_\psi\Lambda_\psi[(id\underset{N}{_b*_\alpha}\omega_{\delta^{-1/2}\xi, (P\delta)^{-1}\eta})\mathfrak a (x)]|\Lambda_\psi(y))\\
&=&(\Delta_\psi(id*\omega_{\xi, (P\delta)^{-1}\eta})(V_\psi)\Lambda_\psi(x)|\Lambda_\psi(y))\\
&=&(V_\psi(\Lambda_\psi(x)\underset{\nu^o}{_a\otimes_\beta}\xi)|\Delta_\psi\Lambda_\psi(y)\underset{\nu}{_b\otimes_\alpha}(P\delta)^{-1}\eta)
\end{eqnarray*}
from which we get the result, thanks to \ref{lemE1}(i). \end{proof}

%%%%%thV2
\subsection{Theorem}
\label{thV2}
{\it Let $\gG=(N, M, \alpha, \beta, \Gamma, T, T', \nu)$ be a measured quantum groupoid; let $A$ be a von Neumann algebra, $(b, \mathfrak a)$ an action of $\gG$ on $A$ and let $\psi$ be a $\delta$-invariant normal faithful semi-finite weight on $A$. Let us define the representation $a$ of $N$ on $H_\psi$ such that, for $n\in N$:
\[a(n)=J_\psi b(n^*)J_\psi\]
and let us suppose that $\psi$ bears the density property, as defined in \ref{defdeltainv}.
Let $V_\psi$ the corepresentation of $\gG$ on the bimodule $_a(H_\psi)_b$ constructed in \ref{thV1}. 
Then :
\newline
(i) for all $t\in\mathbb{R}$ and $n\in N$, we have $\sigma_t^\psi(b(n))=b(\gamma_t(n))$. 
\newline
(ii) it is possible to define two one parameter groups of unitaries $\Delta_\psi^{it}\underset{N^o}{_a\otimes_\beta}\delta^{-it}P^{-it}$ and $\Delta_\psi^{it}\underset{N}{_b\otimes_\alpha}\delta^{-it}P^{-it}$ with natural values on elementary tensor products, and we have :
\[V_\psi(\Delta_\psi^{it}\underset{N^o}{_a\otimes_\beta}\delta^{-it}P^{-it})=(\Delta_\psi^{it}\underset{N}{_b\otimes_\alpha}\delta^{-it}P^{-it})V_\psi\]
(iii) for every $x$ in $A$, $t\in\mathbb{R}$, we have :
\[\mathfrak a(\sigma_t^\psi(x))=(\sigma_t^\psi\underset{N}{_b*_\alpha}\tau_{-t}\sigma_{-t}^{\Phi\circ R}\sigma_{t}^\Phi)\mathfrak a (x)\]
(iv) we have :
\[V_\psi=(J_\psi\underset{\nu}{_b\otimes_\alpha}J_{\widehat{\Phi}})V_\psi^*(J_\psi\underset{\nu}{_b\otimes_\alpha}J_{\widehat{\Phi}})\]
and, therefore, $V_\psi$ is a standard implementation of $\mathfrak a$. }
\begin{proof}
Let $x$, $y$ in $\gN_\psi\cap\gN_\psi^*$, analytic with respect to $\psi$; let $\xi$, $\eta$ as in \ref{lemV2}, $n$ in $N$ analytic with respect both to $\sigma_t^\nu$ and $\gamma_t$; we have :
\begin{multline*}
(V_\psi(\Delta_\psi\Lambda_\psi(b(n)x)\underset{\nu^o}{_a\otimes_\beta}(P\delta)^{-1}\xi)|\Lambda_\psi(y)\underset{\nu}{_b\otimes_\alpha}\eta)\\
= (V_\psi(\Lambda_\psi(b(n)x)\underset{\nu^o}{_a\otimes_\beta}\xi)|\Delta_\psi\Lambda_\psi(y)\underset{\nu}{_b\otimes_\alpha}(P\delta)^{-1}\eta)\\
= (V_\psi(b(n)\Lambda_\psi(x)\underset{\nu^o}{_a\otimes_\beta}\xi)|\Delta_\psi\Lambda_\psi(y)\underset{\nu}{_b\otimes_\alpha}(P\delta)^{-1}\eta)
\end{multline*}
which, using \ref{propV2}(v), is equal to :
\begin{multline*}
((1\underset{N}{_b\otimes_\alpha}\beta(n))V_\psi(\Lambda_\psi(x)\underset{\nu^o}{_a\otimes_\beta}\xi)|\Delta_\psi\Lambda_\psi(y)\underset{\nu}{_b\otimes_\alpha}(P\delta)^{-1}\eta)\\
=(V_\psi(\Lambda_\psi(x)\underset{\nu^o}{_a\otimes_\beta}\xi)|\Delta_\psi\Lambda_\psi(y)\underset{\nu}{_b\otimes_\alpha}\beta(n^*)(P\delta)^{-1}\eta)
\end{multline*}
but, we know that $\beta(n^*)(P\delta)^{-1}\eta=(P\delta)^{-1}\beta(\gamma_i(n^*))\eta$, and $\beta(\gamma_i(n^*))\eta$ bears the property necessary to apply \ref{lemV2};  therefore, it is equal to :
\begin{multline*}
(V_\psi(\Lambda_\psi(x)\underset{\nu^o}{_a\otimes_\beta}\xi)|\Delta_\psi\Lambda_\psi(y)\underset{\nu}{_b\otimes_\alpha}(P\delta)^{-1}\beta(\gamma_i(n^*))\eta)\\
= (V_\psi(\Delta_\psi\Lambda_\psi(x)\underset{\nu^o}{_a\otimes_\beta}(P\delta)^{-1}\xi)|\Lambda_\psi(y)\underset{\nu}{_b\otimes_\alpha}\beta(\gamma_i(n^*))\eta)\\
=((1\underset{N}{_b\otimes_\alpha}\beta(\gamma_{-i}(n))V_\psi(\Delta_\psi\Lambda_\psi(x)\underset{\nu^o}{_a\otimes_\beta}(P\delta)^{-1}\xi)|\Lambda_\psi(y)\underset{\nu}{_b\otimes_\alpha}\eta)
\end{multline*}
which, using again \ref{propV2}(v) again, is equal to :
\[(V_\psi(b(\gamma_{-i}(n))\Delta_\psi\Lambda_\psi(x)\underset{\nu^o}{_a\otimes_\beta}(P\delta)^{-1}\xi)|\Lambda_\psi(y)\underset{\nu}{_b\otimes_\alpha}\eta)\]
from which we get that :
\[\Delta_\psi\Lambda_\psi(b(n)x)=b(\gamma_{-i}(n))\Delta_\psi\Lambda_\psi(x)\]
and we see that $b(n)$ belongs to $\mathcal D(\sigma_{-i}^\psi)$ , and that $\sigma_{-i}^\psi(b(n))=b(\gamma_{-i}(n))$. By classical approximations, we can prove this result for any $n$ in $\mathcal D(\gamma_{-i})$, from which, using (\cite{H}, 4.3 and 4.4), one get (i). 
\newline
Let $\xi\in D(_aH_\psi)$; we have, for all $n\in\gN_\nu$, $t\in\mathbb{R}$ :
\begin{eqnarray*}
a(n)\Delta^{it}_\psi\xi&=&J_\psi b(n^*)J_\psi\Delta^{it}_\psi\xi\\
&=&J_\psi b(n^*)\Delta^{it}_\psi J_\psi\xi\\
&=&J_\psi \Delta^{it}_\psi\sigma_{-t}^\psi(b(n^*))J_\psi\xi\\
&=&J_\psi \Delta^{it}_\psi b(\gamma_{-t}(n^*))J_\psi\xi\\
&=&\Delta^{it}_\psi a(\gamma_{-t}(n))\xi\\
&=&\Delta^{it}_\psi R^{\alpha, \nu}(\xi)\Lambda_\nu(\gamma_{-t}(n))
\end{eqnarray*}
If we define $H^{it}$ by $H^{it}\Lambda_\nu(n)=\Lambda_\nu(\gamma_t(n))$, we get that $\Delta^{it}_\psi\xi$ belongs to $D(_aH_\psi, \nu)$, and that $R^{\alpha, \nu}(\Delta^{it}_\psi\xi)=\Delta^{it}_\psi R^{\alpha, \nu}(\xi)H^{-it}$; then, we get, if $\xi' \in D(_aH_\psi, \nu)$, that :
\[<\Delta^{it}_\psi\xi,\Delta^{it}_\psi\xi'>_{\alpha, \nu^o}^o=\gamma_{t}(<\xi, \xi'>_{\alpha, \nu^o}^o)\]
But, for all $n\in N$, we have:
\begin{eqnarray*}
\beta(\gamma_t(n))&=&\sigma_t^\Phi(\beta(n))\\
&=&\delta^{-it}\sigma_t^{\Phi\circ R}(\beta(n))\delta^{it}\\
&=&\delta^{-it}\beta(\sigma_{-t}^\nu(n))\delta^{it}\\
&=&\delta^{-it}\tau_{-t}(\beta(n))\delta^{it}\\
&=&(\delta P)^{-it}\beta(n)(\delta P)^{it}
\end{eqnarray*}
Therefore,  it is possible to define a one-parameter group of unitaries $\Delta_\psi^{it}\underset{N^o}{_a\otimes_\beta}\delta^{-it}P^{-it}$ on $H_\psi\underset{N^o}{_a\otimes_\beta}H_\Phi$, with natural values on elementary tensor products. 
\newline
The construction of $\Delta_\psi^{it}\underset{N}{_b\otimes_\alpha}\delta^{-it}P^{-it}$ is done in a similar way. 
\newline
Let $\Delta_\psi\underset{\nu^o}{_a\otimes_\beta}(P\delta)^{-1}$ and 
$\Delta_\psi
\underset{\nu}{_b\otimes_\alpha}(P\delta)^{-1}$ be the analytic generators of these two one parameter automorphims groups. We had got in \ref{lemV2} that :
\[V_\psi(\Delta_\psi\underset{\nu^o}{_a\otimes_\beta}(P\delta)^{-1})\subset (\Delta_\psi
\underset{\nu}{_b\otimes_\alpha}(P\delta)^{-1})V_\psi\]
from which we finish the proof of (ii). Then, we get  (iii) from (ii) and \ref{thV1}(i). 
\newline
We can verify that it is possible to define an antilinear unitary $J_\psi\underset{\nu}{_b\otimes_\alpha}J_{\widehat{\Phi}}$ from $H_\psi\underset{\nu}{_b\otimes_\alpha}H_{\Phi}$ onto $H_\psi\underset{\nu^o}{_a\otimes_\beta}H_{\Phi}$, with natural values on elementary tensors, whose inverse is $J_\psi\underset{\nu^o}{_a\otimes_\beta}J_{\widehat{\Phi}}$; moreover, we have, using \ref{thS2} :
\begin{multline*}
((J_\psi\underset{\nu^o}{_a\otimes_\beta}J_{\widehat{\Phi}})V_\psi^*(J_\psi\underset{\nu^o}{_a\otimes_\beta}J_{\widehat{\Phi}})(\Lambda_\psi(x)\underset{\nu^o}{_a\otimes_\beta}\xi)
|\Lambda_\psi(y)\underset{\nu}{_b\otimes_\alpha}\eta)\\
=\overline{(V_\psi^*(J_\psi\Lambda_\psi(x)\underset{\nu}{_b\otimes_\alpha}J_{\widehat{\Phi}}\xi)|
J_\psi\Lambda_\psi(y)\underset{\nu^o}{_a\otimes_\beta}J_{\widehat{\Phi}}\eta)}\\
=\overline{(V_\psi(J_\psi\Lambda_\psi (x)\underset{\nu^o}{_a\otimes_\beta}P^{-1/2}\eta)|
J_\psi\Lambda_\psi(y)\underset{\nu}{_b\otimes_\alpha}P^{1/2}\xi)}
\end{multline*}
which, using (ii), is equal to the scalar product of $J_\psi\Lambda_\psi(y)\underset{\nu}{_b\otimes_\alpha}P^{1/2}\xi$ with :
\[(\Delta_\psi^{1/2}\underset{\nu}{_b\otimes_\alpha}\delta^{-1/2}P^{-1/2})V_\psi(\Delta_\psi^{-1/2}\underset{\nu^o}{_a\otimes_\beta}\delta^{1/2}P^{1/2})(J_\psi\Lambda_\psi (x)\underset{\nu^o}{_a\otimes_\beta}P^{-1/2}\eta)\]
which is :
\[(F_\psi\Lambda_\psi(y)\underset{\nu}{_b\otimes_\alpha}\delta^{-1/2}\xi| 
V_\psi(\Lambda_\psi(x^*)\underset{\nu^o}{_a\otimes_\beta}\delta^{1/2}\eta))
=(F_\psi\Lambda_\psi(y)|(id*\omega_{\delta^{1/2}\eta, \delta^{-1/2}\xi})(V_\psi)\Lambda_\psi(x))\]
and, using twice \ref{propV2}(iii), is equal to :
\begin{multline*}
(F_\psi\Lambda_\psi(y)|\Lambda_\psi((id\underset{N}{_b*_\alpha}\omega_{\eta, \delta^{-1/2}\xi})\mathfrak a (x^*))\\
=(\Lambda_\psi ((id\underset{N}{_b*_\alpha}\omega_{\delta^{-1/2}\xi, \eta})\mathfrak a(x)|\Lambda_\psi(y))\\
=((id*\omega_{\xi, \eta})(V_\psi)\Lambda_\psi (x)|\Lambda_\psi (y))
\end{multline*}
from which we get (iv). \end{proof}

%%%%corV2
\subsection{Corollary}
\label{corV2}
{\it Let $\gG=(N, M, \alpha, \beta, \Gamma, T, T', \nu)$ be a measured quantum groupoid; let $A$ be a von Neumann algebra, $(b, \mathfrak a)$ an action of $\gG$ on $A$ and let $\psi$ be a $\delta$-invariant normal faithful semi-finite weight on $A$. Let us define the representation $a$ of $N$ on $H_\psi$ such that, for $n\in N$:
\[a(n)=J_\psi b(n^*)J_\psi\]
and let us suppose that $\psi$ bears the density property, as defined in \ref{defdeltainv}.
Let $V_\psi$ be the corepresentation of $\gG$ on the bimodule $_a(H_\psi)_b$ constructed in \ref{thV1}. Then :
\newline
(i) it is possible to define the one parameter group of unitaries $\Delta_\psi^{it}\underset{N^o}{_a\otimes_\beta}\delta^{-it}\Delta_{\widehat{\Phi}}^{-it}$ and another one $\Delta_\psi^{it}\underset{N}{_b\otimes_\alpha}\delta^{-it}\Delta_{\widehat{\Phi}}^{-it}$, with natural values on elementary tensor products.
\newline
(ii) we have :
\[V_\psi(\Delta_\psi^{it}\underset{N^o}{_a\otimes_\beta}\delta^{-it}\Delta_{\widehat{\Phi}}^{-it})=(\Delta_\psi^{it}\underset{N}{_b\otimes_\alpha}\delta^{-it}\Delta_{\widehat{\Phi}}^{-it})V_\psi\]
(iii) for any $x$ in $A$, $t$ in $\mathbb{R}$, we have :}
\[\mathfrak a(\sigma_t^\psi(x))=(\Delta_\psi^{it}\underset{N}{_b\otimes_\alpha}\delta^{-it}\Delta_{\widehat{\Phi}}^{-it})\mathfrak a(x)(\Delta_\psi^{-it}\underset{N}{_b\otimes_\alpha}\delta^{it}\Delta_{\widehat{\Phi}}^{it})\]
\begin{proof}
Using \ref{thL3}(vii), $\Delta_{\widehat{\Phi}}$ is the closure of $PJ_\Phi\delta^{-1}J_\Phi$, and, therefore, $\overline{\delta\Delta_{\widehat{\Phi}}}$ is equal to $\overline{P\delta J_\Phi\delta^{-1}J_\Phi}$, and, using \ref{thV2}(i), we can define $\Delta_\psi^{it}\underset{N^o}{_a\otimes_\beta}\delta^{-it}\Delta_{\widehat{\Phi}}^{-it}$ and $\Delta_\psi^{it}\underset{N}{_b\otimes_\alpha}\delta^{-it}\Delta_{\widehat{\Phi}}^{-it}$ by writing :
\[\Delta_\psi^{it}\underset{N^o}{_a\otimes_\beta}\delta^{-it}\Delta_{\widehat{\Phi}}^{-it}=
(\Delta_\psi^{it}\underset{N^o}{_a\otimes_\beta}\delta^{-it}P^{-it})(1\underset{N^o}{_a\otimes_\beta}J_\Phi\delta^{-it}J_\Phi)\]
\[\Delta_\psi^{it}\underset{N}{_b\otimes_\alpha}\delta^{-it}\Delta_{\widehat{\Phi}}^{-it}=
(\Delta_\psi^{it}\underset{N}{_b\otimes_\alpha}\delta^{-it}P^{-it})(1\underset{N}{_b\otimes_\alpha}J_\Phi\delta^{-it}J_\Phi)\]
We then obtain (ii), thanks to \ref{thV2}(i) and \ref{thV1}(ii). Then (iii) is an easy corollary of (ii) and \ref{thV1}(i). \end{proof}

%%%%Corollary
\subsection{Corollary}
\label{corstandardV}
{\it Let $\gG=(N, M, \alpha, \beta, \Gamma, T, T', \nu)$ be a measured quantum groupoid; let $A$ be a von Neumann algebra, $(b, \mathfrak a)$ an action of $\gG$ on $A$; let $\psi$ a $\delta$-invariant normal semi-finite faithful weight on $A$, and let us define the representation $a$ of $N$ on $H_\psi$ such that, for $n\in N$:
\[a(n)=J_\psi b(n^*)J_\psi\]
and let us suppose that $\psi$ bears the density property, as defined in \ref{defdeltainv}. Then :
\newline
(i) $\mathfrak a$ has a standard implementation; 
\newline
(ii) $\mathfrak a$ satifies property (A) if and only if $\mathfrak a$ is saturated. }

\begin{proof}
We obtain (i) by \ref{thV2}. We obtain then (ii) by using (i) and \ref{propstandard}. 
\end{proof}

%%%%cocycle
\subsection{Corollary}
\label{cocycle}
{\it Let $\gG=(N, M, \alpha, \beta, \Gamma, T, T', \nu)$ be a measured quantum groupoid; let $A$ be a von Neumann algebra, $(b, \mathfrak a)$ an action of $\gG$ on $A$ and let $\psi_1$, $\psi_2$ be two $\delta$-invariant normal faithful semi-finite weights on $A$, and let us suppose that both $\psi_1$ and $\psi_2$ bear the density property. Then, for all $t$ in $\mathbb{R}$, $(D\psi_1:D\psi_2)_t$ belongs to $A^\mathfrak a$.}

\begin{proof}
Let us consider the von Neumann algebra $B=A\otimes M_2(\mathbb{C})$; for $n$ in $N$, let us define $\tilde{b}(n)=b(n)\otimes (e_{11}+e_{22})$, which is an injective $*$-antihomomorphism from $N$ into $B$, and, for $X=(x_{i,j})_{i,j}$ in $B$; let $\sigma$ denotes the flip from $H_\Phi\otimes\mathbb{C}^2$ onto $\mathbb{C}^2\otimes H_\Phi$, and let us define :
\[\mathfrak b(X)=(1\underset{N}{_b\otimes_\alpha}\sigma)(\mathfrak a(x_{i,j}))_{i,j}(1\underset{N}{_b\otimes_\alpha}\sigma^*)\]
which belongs to $B\underset{N}{_{\tilde{b}}*_\alpha}M$.
\newline
It is easy to check that $(\tilde{b}, \mathfrak b)$ is an action of $\gG$ on $B$. Moreover, let $\psi$ be the normal semi-finite faithful weight defined on $B$ by, for $X=(x_{i,j})_{i,j}$ in $B^+$, $\psi(X)=\psi_1(x_{11})+\psi_2(x_{22})$; we then can verify that $\psi$ is a $\delta$-invariant weight with respect to $\mathfrak b$, and that $\psi$ bears the density property. 
\newline
Therefore, we have, by definition of $(D\psi_1:D\psi_2)_t$ :
\begin{eqnarray*}
(1\underset{N}{_b\otimes_\alpha}\sigma)[\mathfrak a(D\psi_1:D\psi_2)_t)\otimes e_{21}](1\underset{N}{_b\otimes_\alpha}\sigma^*)
&=&
\mathfrak b[(D\psi_1:D\psi_2)_t\otimes e_{21}]\\
&=&
\mathfrak b\sigma_t^\psi(1\otimes e_{21})
\end{eqnarray*}
which, using \ref{thV2}(ii), is equal to :
\begin{eqnarray*}
(\sigma_t^\psi\underset{N}{_{\tilde{b}}*_\alpha}\tau_{-t}\sigma_{-t}^{\Phi\circ R}\sigma_{t}^\Phi)\mathfrak b(1\otimes e_{21})
&=&
(\sigma_t^\psi\underset{N}{_{\tilde{b}}*_\alpha}\tau_{-t}\sigma_{-t}^{\Phi\circ R}\sigma_{t}^\Phi)[(1\otimes e_{21})\underset{N}{_{\tilde{b}}\otimes_\alpha}1]\\
&=&[(D\psi_1:D\psi_2)_t\otimes e_{21}]\underset{N}{_{\tilde{b}}\otimes_\alpha}1)
\end{eqnarray*}
from which we get the result.  \end{proof}

%%%%%dual
\section{Crossed-product and dual actions}
\label{dual}
In that chapter, we define the crossed product of a von Neumann algebra by a measured quantum groupoid (\ref{defcross}), and the dual action on this new von Neumann algebra (\ref{dualaction}). We prove that the dual action is integrable (\ref{Propdual2}), and satisfies (\ref{Propdual}) the property (A) introduced in chapter \ref{action}.

%%%%defcross
\subsection{Definition}
\label{defcross}
Let $\gG=(N, M, \alpha, \beta, \Gamma, T, T', \nu)$ be a measured quantum groupoid; let $A$ be a von Neumann algebra acting on a Hilbert space $\gH$, $(b, \mathfrak a)$ an action of $\gG$ on $A$; we define the crossed-product of $A$ by $\gG$ (via the action $\mathfrak a$) as the von Neumann algebra generated by $\mathfrak a(A)$ and $1\underset{N}{_b\otimes_\alpha}\widehat{M}'$ on the Hilbert space $\gH\underset{N}{_b\otimes_\alpha}H_\Phi$. This von Neumann algebra will be denoted by $A\rtimes_{\mathfrak a}\gG$. Clearly, this von Neumann algebra is included in $A\underset{N}{_b*_\alpha}\mathcal L(H_\Phi)$. 

%%%gcross
\subsection{Example}
\label{gcross}
Let $\mathcal G$ be a measured groupoid, and $\mathfrak a$ an action of $\mathcal G$ on a von Neumann algebra $A$, in the sense of \cite{Y1}, \cite{Y2}, \cite{Y3}. We have seen in \ref{graction} that $\mathfrak a$ can be as well considered as an action of the abelian measured quantum groupoid constructed from $\mathcal G$ (\ref{gd2}) on $A$. In \cite{Y1} and \cite{Y3} is given a construction of the crossed product of $A$ by $\mathcal G$; using (\cite{Y1}, 2.14), we get that the crossed product defined in \ref{defcross} is the same. 

%%%d2cross
\subsection{Example}
\label{d2cross}
Let $M_0\subset M_1$ be a depth 2 inclusion of von Neumann algebras, with a regular operator-valued 
weight $T_1$ from $M_1$ to $M_0$, as defined in \ref{basic}. Let $\gG_1$ be the measured quantum groupoid, whose underlying von Neumann algebra is $M'_0\cap M_2$ and whose basis is $M'_0\cap M_1$, constructed (\ref{d2}) from this inclusion, and $\gG_2$ the measured quantum groupoid constructed from the inclusion $M_1\subset M_2$, which is isomorphic to $\widehat{\gG_1}^o$ (therefore, we get that $\widehat{\gG_2}^c=\widehat{\gG_1}$). Then (\ref{d2action}), there exists (\cite{EV}, 7.3) a canonical action $\mathfrak a_2$ of $\gG_2$ on $M_1$, which has the following properties : 
\newline
(i) the basis of $\gG_2$ is $M'_1\cap M_2$ (which, using $j_1$, is anti-isomorphic to $M'_0\cap M_1$), the underlying von Neumann algebra of $\gG_2$ is $M'_1\cap M_3$ and, for $x$ in $M_1$,  $\mathfrak a_2(x)$ belongs to $M_1\underset{M'_1\cap M_2}{_{j_1}*_{j_2}}M'_1\cap M_3$; in fact $\mathfrak a_2$ is given by the natural inclusion of $M_1$ into $M_3$;  
\newline
(ii) we have $M_0=M_1^{\mathfrak a_2}$;
\newline
(iii) there is an isomorphism $I_{\mathfrak a_2}$ from $M_1\rtimes_{\mathfrak a_2}\gG_2$ onto $M_2$(\cite{EV}, 7.6), such that, for any $x\in M_1$, we have $I_{\mathfrak a_2}(\mathfrak a_2(x))=x$, and, for any $y\in M'_0\cap M_2$, we have $I_{\mathfrak a_2}(1\underset{M'_1\cap M_2}{_{j_1}\otimes_{j_2}} y)=y$.

%%%%dualaction
\subsection{Theorem}
\label{dualaction}
{\it Let $\gG=(N, M, \alpha, \beta, \Gamma, T, T', \nu)$ be a measured quantum groupoid; let $A$ be a von Neumann algebra acting on a Hilbert space $\gH$, $(b, \mathfrak a)$ an action of $\gG$ on $A$, $A\rtimes_{\mathfrak a}\gG$ the crossed product  of $A$ by $\gG$ via the action $\mathfrak a$. Then :
\newline
(i) the operator $1_\gH\underset{N}{_b\otimes_\alpha}W^{co}$ is a corepresentation of $(\widehat{\gG})^c$ on the $N^o-N^o$ bimodule $_{1\underset{N}{_b\otimes_\alpha}\hat{\beta}}(\gH\underset{N}{_b\otimes_\alpha}H_\Phi)_{1\underset{N}{_b\otimes_\alpha}\hat{\alpha}}$.
\newline
(ii) this corepresentation implements an action $(1\underset{N}{_b\otimes_\alpha}\hat{\alpha}, \tilde{\mathfrak a})$ of $(\widehat{\gG})^c$ on $A\rtimes_{\mathfrak a}\gG$, which verifies, for any $x\in A$, $y\in \widehat{M}'$ :
\[\tilde{\mathfrak a}(\mathfrak a(x))=\mathfrak a(x)\underset{N^o}{_{\hat{\alpha}}\otimes_\beta}1\]
\[\tilde{\mathfrak a}(1\underset{N}{_b\otimes_\alpha}y)=1\underset{N}{_b\otimes_\alpha}\widehat{\Gamma}^c(y)\]
This action $(1\underset{N}{_b\otimes_\alpha}\hat{\alpha}, \tilde{\mathfrak a})$ (or $\tilde{\mathfrak a}$ for simplification) of $(\widehat{\gG})^c$ will be called the dual action of $\mathfrak a$. 
\newline
(iii) We have :}
 \[\mathfrak a(A)\subset(A\rtimes_{\mathfrak a}\gG)^{ \tilde{\mathfrak a}}= A\rtimes_{\mathfrak a}\gG\cap A\underset{N}{_b*_\alpha}M\]

\begin{proof}
(i) is given by \ref{Example}. Moreover, as $\mathfrak a(x)$ belongs to $A\underset{N}{_b*_\alpha}M$, we have :
\[(1_\gH\underset{N}{_b\otimes_\alpha}W^{co})(\mathfrak a(x)\underset{N}{_{\hat{\beta}}\otimes_{\hat{\alpha}}}1)(1_\gH\underset{N}{_b\otimes_\alpha}W^{co})^*=\mathfrak a(x)\underset{N^o}{_{\hat{\alpha}}\otimes_\beta}1\]
and, using \ref{thL5}(vi), we get :
\[(1_\gH\underset{N}{_b\otimes_\alpha}W^{co})(1\underset{N}{_b\otimes_\alpha}y\underset{N^o}{_{\hat{\beta}}\otimes_\alpha}1)(1_\gH\underset{N}{_b\otimes_\alpha}W^{co})^*=1\underset{N}{_b\otimes_\alpha}\widehat{\Gamma}^c(y)\]
and, therefore, we have :
\begin{multline*}
(1_\gH\underset{N}{_b\otimes_\alpha}W^{co})A\rtimes_{\mathfrak a}\gG(1_\gH\underset{N}{_b\otimes_\alpha}W^{co})^*\\
\subset (\mathfrak a(A)\underset{N^o}{_{\hat{\alpha}}\otimes_\beta}1\cup 1\underset{N}{_b\otimes_\alpha}\widehat{M}'\underset{N^o}{_{\hat{\alpha}}*_\beta}\widehat{M}')''
\subset A\rtimes_{\mathfrak a}\gG\underset{N^o}{_{\hat{\alpha}}*_\beta}\widehat{M}'
\end{multline*}
which proves (ii). 
\newline
The inclusion $\mathfrak a(A)\subset(A\rtimes_{\mathfrak a}\gG)^{ \tilde{\mathfrak a}}$ is given by (ii); let now $X$ be in $A\rtimes_{\mathfrak a}\gG$, then $X$ belongs to $(A\rtimes_{\mathfrak a}\gG)^{ \tilde{\mathfrak a}}$ if and only if we have :
\[(1\underset{N}{_b\otimes_\alpha}W^{co})(X\underset{N}{_{\hat{\beta}}\otimes_{\hat{\beta}}}1)=
(X\underset{N^o}{_{\hat{\alpha}}\otimes_\beta}1)(1\underset{N}{_b\otimes_\alpha}W^{co})\]
which means that $X$ belongs to $\mathcal L(\gH)\underset{N}{_b*_\alpha}M$. As $A\rtimes_{\mathfrak a}\gG\subset A\underset{N}{_b*_\alpha}\mathcal L(H_\Phi)$, this finishes the proof of (iii). \end{proof}

%%%%dualtri
\subsection{Example}
\label{dualtri}
Let $\gG=(N, M, \alpha, \beta, \Gamma, T, T', \nu)$ be a measured quantum groupoid, and $(id, \beta)$ its trivial action on $N^o$ (\ref{triaction}); then the crossed-product $N^o\rtimes_{\beta}\gG$ is equal to $\widehat{M}'$; moreover, the dual action $\tilde{\beta}$ is equal to $\widehat{\Gamma}^c$, considered as an action of $(\widehat{\gG})^c$ on $\widehat{M}'$ (\ref{Gamma}). 
\newline
Therefore, $(\beta, \Gamma)$, considered as an action of $\gG$ on $M$, is the dual action of the trivial action $(id, \alpha)$ of $\widehat{\gG}^o$ on $N$. 

%%%dualg
\subsection{Example}
\label{dualg}
Let $\mathcal G$ be a measured groupoid, and $\mathfrak a$ an action of $\mathcal G$ on a von neumann algebra $A$, in the sense of \cite{Y1}, \cite{Y2}, \cite{Y3}. We have seen in \ref{gcross} that the crossed-product $A\rtimes_\mathfrak a\mathcal G$ in the sense of \ref{defcross} in equal to the crossed-product in the sense of Yamanouchi; moreover, the dual action $\tilde{\mathfrak a}$ in the sense of \ref{dualaction} is an action of the measured quantum groupoid $\widehat{\mathcal G}^c$ (whose underlying von Neumann algebra is the von Neumann algebra generated by the right convolution algebra of $\mathcal G$). Using (\cite{Y3}, 4.4), we see it is equal to the "dual coaction of $\mathcal G$" introduced in (\cite{Y3}, 4.9). 

%%%%Propdual
\subsection{Proposition}
\label{Propdual}
{\it Let $\gG=(N, M, \alpha, \beta, \Gamma, T, T', \nu)$ be a measured quantum groupoid; let $A$ be a von Neumann algebra acting on a Hilbert space $\gH$, $(b, \mathfrak a)$ an action of $\gG$ on $A$; the dual action $\tilde{\mathfrak a}$ of $(\widehat{\gG})^c$ on the crossed product $A\rtimes_{\mathfrak a}\gG$ satisfies the property (A) defined in \ref{defA}. }

\begin{proof}
Let us consider the von Neumann algebra $B$ generated on the Hilbert space $\gH\underset{\nu}{_b\otimes_\alpha}H_\Phi\underset{\nu^o}{_{\hat{\alpha}}\otimes_\beta}H_\Phi$ by $\tilde{\mathfrak a}(A\rtimes_{\mathfrak a}\gG)$ and $1\underset{N}{_b\otimes_\alpha}1\underset{N^o}{_{\hat{\alpha}}\otimes_\beta}\beta(N)'$. By \ref{dualaction}(ii), it is generated by $\mathfrak a(A)\underset {N^o}{_{\hat{\alpha}}\otimes_\beta}1$, $1\underset{N}{_b\otimes_\alpha}\widehat{\Gamma}'(\widehat{M}')$ and $1\underset{N}{_b\otimes_\alpha}1\underset{N^o}{_{\hat{\alpha}}\otimes_\beta}\beta(N)'$; but, using \ref{GammaA} applied to $(\widehat{\gG})^c$, we get that :
\[(\widehat{\Gamma}'(\widehat{M}')\cup 1\underset{N^o}{_{\hat{\alpha}}\otimes_\beta}\beta(N)')''=\widehat{M}'\underset{N^o}{_{\hat{\alpha}}*_\beta}\mathcal L(H_\Phi)\]
and, therefore :
\[B=(\mathfrak a(A)\underset {N^o}{_{\hat{\alpha}}\otimes_\beta}1\cup 1\underset{N}{_b\otimes_\alpha}\widehat{M}'\underset{N^o}{_{\hat{\alpha}}*_\beta}\mathcal L(H_\Phi))''\]
or :
\begin{eqnarray*}
B'&=&\mathfrak a(A)'\underset {N^o}{_{\hat{\alpha}}*_\beta}\mathcal L(H_\Phi)\cap\mathcal L(\gH)\underset{N}{_b*_\alpha}\widehat{M}\underset{N^o}{_{\hat{\alpha}}\otimes_\beta}1\\
&=&(\mathcal L(\gH)\underset{N}{_b*_\alpha}\widehat{M}\cap\mathfrak a(A)')\underset{N^o}{_{\hat{\alpha}}\otimes_\beta}1
\end{eqnarray*}
which gives that $B=A\rtimes_{\mathfrak a}\gG\underset{N^o}{_{\hat{\alpha}}*_\beta}\mathcal L(H_\Phi)$, and finishes the proof. \end{proof}

%%%%%%Propdual2
\subsection{Theorem}
\label{Propdual2}
{\it Let $\gG=(N, M, \alpha, \beta, \Gamma, T, T', \nu)$ be a measured quantum groupoid; let $A$ be a von Neumann algebra acting on a Hilbert space $\gH$, $(b, \mathfrak a)$ an action of $\gG$ on $A$; the dual action $\tilde{\mathfrak a}$ of $(\widehat{\gG})^c$ on the crossed product $A\rtimes_{\mathfrak a}\gG$ is integrable.}

\begin{proof}
Let $y$ be positive in $\widehat{M}'$; we have :
\[T_{\tilde{\mathfrak a}}(1\underset{N}{_b\otimes_\alpha}y)=(id\underset{\nu^o}{_{\hat{\alpha}}*_\beta}\widehat{\Phi}^c)(1\underset{N}{_b\otimes_\alpha}\widehat{\Gamma}^c(y))=1\underset{N}{_b\otimes_\alpha}\widehat{T}^c(y)\]
As $\widehat{T}^c$ is a semi-finite weight from $\widehat{M}'$ to $\beta(N)$, we get that $T_{\tilde{\mathfrak a}}$ is semi-finite, which finishes the proof. \end{proof}

%%%%defdualw
\subsection{Definition}
\label{defdualw}
Let us take the notations of \ref{Propdual2}. Let $\psi$ be a normal semi-finite faithful weight on $(A\rtimes_\mathfrak a\gG)^{\tilde{\mathfrak a}}$; then, we shall denote $\tilde{\psi}$ the lifted weight $\tilde{\psi}=\psi\circ T_{\tilde{\mathfrak a}}$, which is a normal semi-finite faithful weight on $A\rtimes_\mathfrak a\gG$. 

%%%%lemdualw
\subsection{Lemma}
\label{lemdualw}
{\it With the notations of \ref{defdualw}, we have, for all $x\in (A\rtimes_\mathfrak a\gG)^{\tilde{\mathfrak a}}$ and $a\in\gN_{\widehat{T}^c}$:}
\[\tilde{\psi}(x^*(1\underset{N}{_b\otimes_\alpha}a^*a)x)=\psi(x^*(1\underset{N}{_b\otimes_\alpha}\widehat{T}^c(a^*a))x)\]

\begin{proof}
We have, using \ref{dualaction}(ii) :
\[\tilde{\mathfrak a}(1\underset{N}{_b\otimes_\alpha}a^*a)=
1\underset{N}{_b\otimes_\alpha}\widehat{\Gamma}^c(a^*a)\]
from which we get that :
\[T_{\tilde{\mathfrak a}}(1\underset{N}{_b\otimes_\alpha}a^*a)
=1\underset{N}{_b\otimes_\alpha}\widehat{T}^c(a^*a)\]
and, therefore :
\[T_{\tilde{\mathfrak a}}[x^*(1\underset{N}{_b\otimes_\alpha}a^*a)x]
=x^*[1\underset{N}{_b\otimes_\alpha}\widehat{T}^c(a^*a)]x\]
from which we get the result. \end{proof}

%%%%%propdualw
\subsection{Proposition}
\label{propdualw}
{\it Let's take the hypothesis and notations of \ref{defdualw}; there exists an isometry $V$ from $H_\psi\underset{\nu}{_b\otimes_{\alpha}}H_\Phi$ into $H_{\tilde{\psi}}$, such that, we have, for all $x$ in $\gN_\psi$ and $a$ in $\gN_{\widehat{T}^c}\cap\gN_{{\widehat{\Phi}}^c}$:
\[V(\Lambda_\psi(x)\underset{\nu}{_b\otimes_\alpha}\Lambda_{{\widehat{\Phi}^c}}(a))=
\Lambda_{\tilde{\psi}}((1\underset{N}{_b\otimes_\alpha}a)x)\]
Moreover, we have, for all $b\in \widehat{M}'$, $y\in (A\rtimes_\mathfrak a\gG)^{\tilde{\mathfrak a}}$ :}
\[V(1\underset{N^o}{_b\otimes_\alpha}b)=\pi_{\tilde{\psi}}(1\underset{N}{_b\otimes_\alpha}b)V\]
\[V(J_\psi yJ_\psi\underset{N^o}{_b\otimes_\alpha}1)=J_{\tilde{\psi}}\pi_{\tilde{\psi}}(y)J_{\tilde{\psi}}V\]

\begin{proof}
For any $a$ in $\gN_{\widehat{T}^c}\cap\gN_{{\widehat{\Phi}}^c}$, we have $\Lambda_{{\widehat{\Phi}^c}}(a)=J_{\widehat{\Phi}}\Lambda_{\widehat{\Phi}}(J_{\widehat{\Phi}}aJ_{\widehat{\Phi}})$, which belongs to $D(_{\alpha}H_\Phi, \nu)$ (\ref{basic}), and $R^{\alpha, \nu}(\Lambda_{{\widehat{\Phi}^c}}(a))=J_{\widehat{\Phi}}\Lambda_{\widehat{T}}(J_{\widehat{\Phi}}aJ_{\widehat{\Phi}})J_\nu$, and, therefore, we have:
\[<\Lambda_{{\widehat{\Phi}^c}}(a)), \Lambda_{{\widehat{\Phi}^c}}(a))>_{\alpha, \nu}=
\alpha^{-1}[\widehat{T}(J_{\widehat{\Phi}}a^*aJ_{\widehat{\Phi}})]^o=\beta^{-1}(\widehat{T}^c(a^*a))\]
and we get :
\[\|\Lambda_\psi(x)\underset{\nu}{_b\otimes_\alpha}\Lambda_{{\widehat{\Phi}^c}}(a))\|^2=
(b\circ\beta^{-1}(\widehat{T}^c(a^*a))\Lambda_\psi(x)|\Lambda_\psi(x))\]
which, using \ref{lemdualw}, is equal to $\|\Lambda_{\tilde{\psi}}((1\underset{N}{_b\otimes_\alpha}a)x)\|^2$. Which, by polarisation, proves the existence of the isometry $V$. 
\newline
The formula $V(1\underset{N^o}{_b\otimes_\alpha}b)=\pi_{\tilde{\psi}}(1\underset{N}{_b\otimes_\alpha}b)V$ is then trivial. 
\newline
As, for all $t\in\mathbb{R}$, and $y\in (A\rtimes_\mathfrak a\gG)^{\tilde{\mathfrak a}}$, we have $\sigma_t^{\tilde{\psi}}(y)=\sigma_t^\psi(y)$, we get that if $y$ is analytic with respect to $\psi$, $y$ is also analytic with respect to $\tilde{\psi}$, and, for such an element :
\begin{eqnarray*}
V(J_\psi \sigma_{-i/2}^\psi(y^*)yJ_\psi\underset{N^o}{_b\otimes_\alpha}1)(\Lambda_\psi(x)\underset{\nu}{_b\otimes_\alpha}\Lambda_{{\widehat{\Phi}^c}}(a))
&=&
V(\Lambda_\psi(xy)\underset{\nu}{_b\otimes_\alpha}\Lambda_{{\widehat{\Phi}^c}}(a))\\
&=&
\Lambda_{\tilde{\psi}}((1\underset{N}{_b\otimes_\alpha}a)xy)\\
\end{eqnarray*}
which is equal to :
\[J_{\tilde{\psi}}\pi_{\tilde{\psi}}(\sigma_{-i/2}^{\tilde{\psi}}(y^*))J_{\tilde{\psi}}\Lambda_{\tilde{\psi}}((1\underset{N}{_b\otimes_\alpha}a)x)=
J_{\tilde{\psi}}\pi_{\tilde{\psi}}(\sigma_{-i/2}^\psi(y^*))J_{\tilde{\psi}}\Lambda_{\tilde{\psi}}((1\underset{N}{_b\otimes_\alpha}a)\mathfrak a(x))\]
which, by continuity, gives the result which finishes the proof. \end{proof}

%%%%lemdualw1
\subsection{Lemma}
\label{lemdualw1}
{\it Let's take the hypothesis and notations of \ref{defdualw} and \ref{propdualw}; then, for all $X$ in $A\rtimes_\mathfrak a\gG$, we have $\pi_{\tilde{\psi}}(X)V=VX$. }

\begin{proof}
Let $x\in A$; we have :
\begin{eqnarray*}
(\mathfrak a\underset{N}{_b*_\alpha}id)\mathfrak a(x)
&=&
(id\underset{N}{_b*_\alpha}\Gamma)\mathfrak a(x)\\
&=&(1\underset{N}{_b\otimes_\alpha}\widehat{W^o})(\mathfrak a(x)\underset{N^o}{_{\hat{\alpha}}\otimes_\beta}1)((1\underset{N}{_b\otimes_\alpha}\widehat{W^o})^*
\end{eqnarray*}
and, therefore :
\[(\mathfrak a(x)\underset{N^o}{_{\hat{\alpha}}\otimes_\beta}1)((1\underset{N}{_b\otimes_\alpha}\widehat{W^o})^*=
(1\underset{N}{_b\otimes_\alpha}\widehat{W^o})^*(\mathfrak a\underset{N}{_b*_\alpha}id)\mathfrak a(x)\]
Let $\xi\in D(_\alpha H_\Phi)$, $\eta\in D(H_\beta)$; we have :
\[\mathfrak a(x)(1\underset{N}{_b\otimes_\alpha}(id*\omega_{\xi, \eta}(\widehat{W^o}^*))
=
(\rho_\eta^{\beta, \alpha})^*(1\underset{N}{_b\otimes_\alpha}\widehat{W^o})^*(\mathfrak a\underset{N}{_\beta*_\alpha}id)\mathfrak a(x)\rho_\xi^{b, \alpha}\]
Let $(e_i)_{i\in I}$ be an orthogonal $(\alpha, \nu)$-basis of $H_\Phi$; as $\sum_i\rho_{e_i}^{\beta, \alpha}(\rho_{e_i}^{\beta, \alpha})^*=1$, it is equal to :
\begin{multline*}
\sum_i(\rho_\eta^{\beta, \alpha})^*(1\underset{N}{_b\otimes_\alpha}\widehat{W^o})^*\rho_{e_i}^{\beta, \alpha}(\rho_{e_i}^{\beta, \alpha})^*(\mathfrak a\underset{N}{_\beta*_\alpha}id)\mathfrak a(x)\rho_\xi^{b, \alpha}=\\
\sum_i(1\underset{N}{_b\otimes_\alpha}(id*\omega_{e_i, \eta})(\widehat{W^o}^*))\mathfrak a((id\underset{N}{_b*_\alpha}\omega_{\xi, e_i})\mathfrak a(x))
\end{multline*}
This sum being convergent $\sigma^*$-strongly. 
\newline
Let $y$ be in $\gN_T\cap\gN_\Phi\cap\gN_{RTR}\cap\gN_{\Phi\circ R}$; and let us put $\eta=J_{\widehat{\Phi}}J_\Phi\Lambda_\Phi(y)$; then $J_{\widehat{\Phi}}\eta$ belongs to $J_\Phi\Lambda_\Phi(\gN_\Phi)\cap D(_\alpha H_\Phi, \nu)\cap D((H_\Phi)_\beta, \nu^o)$, thanks to \ref{thL3}(i); therefore, $\eta$ belongs to $D(_\alpha H_\Phi, \nu)\cap D((H_\Phi)_\beta, \nu^o)$; moreover, as $\pi'(J_{\widehat{\Phi}}\eta)$ is a bounded operator, we get (\ref{defI}) that $\omega_{J_{\widehat{\Phi}}\xi, J_{\widehat{\Phi}}\eta}$ belongs to $I_\Phi$, and, therefore, by \ref{Plancherel}, that $(\omega_{J_{\widehat{\Phi}}\xi, J_{\widehat{\Phi}}\eta}*id)(W)$ belongs to $\gN_{\widehat{\Phi}}$; we then get easily that $(id*\omega_{\xi, \eta})(\widehat{W^o}^*)$ belongs to $\gN_{\widehat{\Phi}^c}$, and that :
\begin{eqnarray*}
\Lambda_{\widehat{\Phi}^c}((id*\omega_{\xi, \eta})(\widehat{W^o}^*))
&=&
J_{\widehat{\Phi}}\Lambda_{\widehat{\Phi}}((\omega_{J_{\widehat{\Phi}}\xi, J_{\widehat{\Phi}}\eta}*id)(W))\\
&=&J_{\widehat{\Phi}}a_\Phi(\omega_{J_{\widehat{\Phi}}\xi, J_{\widehat{\Phi}}\eta})\\
&=&J_{\widehat{\Phi}}\pi'(J_{\widehat{\Phi}}\eta)^*J_{\widehat{\Phi}}\xi\\
&=&J_{\widehat{\Phi}}J_\Phi y^*J_\Phi J_{\widehat{\Phi}}\xi
\end{eqnarray*}
\newline
Let $x'$ be in $\gN_\psi$, then, $\mathfrak a(x)(1\underset{N}{_b\otimes_\alpha}(id*\omega_{\xi, \eta})(\widehat{W^o}^*))x'$ belongs to $\gN_{\tilde{\psi}}$, by \ref{propdualw}, and is the $\sigma^*$-strong limit of :
\[\sum_i(1\underset{N}{_b\otimes_\alpha}(id*\omega_{e_i, \eta})(\widehat{W^o}^*))\mathfrak a((id\underset{N}{_b*_\alpha}\omega_{\xi, e_i})\mathfrak a(x))x'\]
Let $J$ finite, $J\subset I$; then :
\[\Lambda_{\tilde{\psi}}[\sum_{i\in J}(1\underset{N}{_b\otimes_\alpha}(id*\omega_{e_i, \eta})(\widehat{W^o}^*))\mathfrak a((id\underset{N}{_b*_\alpha}\omega_{\xi, e_i})\mathfrak a(x))x']\]
is, using \ref{propdualw}, equal to :
\[V(\sum_{i\in J}[(id\underset{N}{_b*_\alpha}\omega_{\xi, e_i})\mathfrak a(x)\Lambda_\psi(x')\underset{N}{_b\otimes_\alpha}\Lambda_{\widehat{\Phi}^c}((id*\omega_{e_i, \eta})(\widehat{W^o}^*))]\]
which is :
\[V(\sum_{i\in J}[(id\underset{N}{_b*_\alpha}\omega_{\xi, e_i})\mathfrak a(x)\Lambda_\psi(x')\underset{N}{_b\otimes_\alpha}J_{\widehat{\Phi}}J_\Phi y^*J_\Phi J_{\widehat{\Phi}}e_i]\]
which can be written :
\[V(1\underset{N}{_b\otimes_\alpha}J_{\widehat{\Phi}}J_\Phi y^*J_\Phi J_{\widehat{\Phi}})\sum_{i\in J}\rho_{e_i}^{\beta, \alpha}(\rho_{e_i}^{\beta, \alpha})^*\mathfrak a(x)(\Lambda_\psi(x')\underset{N}{_b\otimes_\alpha}\xi)\]
which converges in norm, when $J$ is growing, to :
\begin{eqnarray*}
V(1\underset{N}{_b\otimes_\alpha}J_{\widehat{\Phi}}J_\Phi y^*J_\Phi J_{\widehat{\Phi}})\mathfrak a(x)(\Lambda_\psi(x')\underset{N}{_b\otimes_\alpha}\xi)
&=&
V\mathfrak a(x)(\Lambda_\psi(x')\underset{N}{_b\otimes_\alpha}J_{\widehat{\Phi}}J_\Phi y^*J_\Phi J_{\widehat{\Phi}}\xi)\\
&=&
V\mathfrak a(x)(\Lambda_\psi(x')\underset{N}{_b\otimes_\alpha}\Lambda_{\widehat{\Phi}^c}((id*\omega_{\xi, \eta})(\widehat{W^o}^*))
\end{eqnarray*}
As $\Lambda_{\tilde{\psi}}$ is $\sigma^*$-norm closed, we therefore get that :
\begin{eqnarray*}
V\mathfrak a(x)(\Lambda_\psi(x')\underset{N}{_b\otimes_\alpha}\Lambda_{\widehat{\Phi}^c}((id*\omega_{\xi, \eta})(\widehat{W^o}^*))&=&
\Lambda_{\tilde{\psi}}[\mathfrak a(x)(1\underset{N}{_b\otimes_\alpha}(id*\omega_{\xi, \eta})(\widehat{W^o}^*))x']\\
&=&
\pi_{\tilde{\psi}}(\mathfrak a(x))\Lambda_{\tilde{\psi}}[(1\underset{N}{_b\otimes_\alpha}(id*\omega_{\xi, \eta})(\widehat{W^o}^*))x']\\
&=&
\pi_{\tilde{\psi}}(\mathfrak a(x))V(\Lambda_\psi(x')\underset{N}{_b\otimes_\alpha}\Lambda_{\widehat{\Phi}^c}((id*\omega_{\xi, \eta})(\widehat{W^o}^*))
\end{eqnarray*}
and, by linearity and continuity, we get that $\pi_{\tilde{\psi}}(\mathfrak a(x))V=V\mathfrak a(x)$, for any $x\in A$. As we have obtained in \ref{propdualw} that $V(1\underset{N^o}{_b\otimes_\alpha}b)=\pi_{\tilde{\psi}}(1\underset{N}{_b\otimes_\alpha}b)V$ for any $b\in\widehat{M}'$, we obtain that $\pi_{\tilde{\psi}}(X)V=VX$ for any $X$ in the involutive algebra generated by $\mathfrak a(A)$ and $1\underset{N}{_b\otimes_\alpha}\widehat{M}'$; then, by strong limit, we obtain the result. \end{proof}

%%%%aux
\section{An auxilliary weight on the crossed-product}
\label{aux}
Be given an action $\mathfrak a$ of a measured quantum groupoid $\gG$ on a von Neumann algebra $A$, and a normal semi-finite faithful weight $\psi$ on $A$, we construct, in this chapter, (\ref{propo1}) an auxilliary weight $\tilde{\psi}_0$ on $A\rtimes_\mathfrak a\gG$ which will allow us, thanks to chapter \ref{deltainv}, to construct a standard implementation for the dual action (\ref{thdual}). 
We obtain (\ref{thdual}) that the dual action satisfies the saturation property introduced in chapter \ref{tech}; we introduce a third technical property (property (B) in \ref{defB}) which is satisfied by the dual action (\ref{corB}).

%%%%propo1
\subsection{Proposition}
\label{propo1}
{\it Let's take the hypothesis and notations of \ref{defdualw} and \ref{propdualw}; then :
\newline 
(i) there exists a unique normal semi-finite faithful weight $\tilde{\psi}_0$ on $A\rtimes_\mathfrak a\gG$ such that :
\[\gN_{\tilde{\psi}_0}=\{x\in\gN_{\tilde{\psi}}, \Lambda_{\tilde{\psi}}(x)\in Im V\}\]
and such that we may identify $H_{\tilde{\psi}_0}$ with $ImV$, $\Lambda_{\tilde{\psi}_0}$ with the restriction of $\Lambda_{\tilde{\psi}}$ to $\gN_{\tilde{\psi}_0}$, and, for any $y\in A\rtimes_\mathfrak a\gG$, $\pi_{\tilde{\psi}_0}(y)$ to the restriction of $\pi_{\tilde{\psi}}(y)$ to $ImV$. 
\newline
(ii) let us write $\tilde{V}$ the unitary from $H_\psi\underset{\nu}{_b\otimes_\alpha}H_\Phi$ onto $H_{\tilde{\psi}_0}$ defined, for all $\Xi\in H_\psi\underset{\nu}{_b\otimes_\alpha}H_\Phi$, by $\tilde{V}\Xi=V\Xi$; then we have, for all $y$ in $A\rtimes_\mathfrak a\gG$ :}
\[\pi_{\tilde{\psi}_0}(y)=\tilde{V}y\tilde{V}^*\]
{\it (iii) the linear set generated by the elements of the form $(1\underset{N}{_b\otimes_\alpha}a)x$, with $a\in\gN_{\hat{T}^c}\cap\gN_{\hat{\Phi}^c}$ and $x\in\gN_\psi$, is dense in $A\rtimes_\mathfrak a\gG$.}
\begin{proof}
Using \ref{lemdualw}, we get that $\{x\in\gN_{\tilde{\psi}}, \Lambda_{\tilde{\psi}}(x)\in Im V\}$ is an ideal of 
$A\rtimes_\mathfrak a\gG$, which is dense by \ref{propdualw}; then, we can apply (\cite{V2}, 7.4)
to obtain (i); result (ii) is then straightforward. Result (iii) is also a corollary of \ref{lemdualw} and \ref{propdualw}. \end{proof}

%%%%%propo2
\subsection{Proposition}
\label{propo2}
{\it Let's take the hypothesis and notations of \ref{defdualw}, \ref{propdualw} and \ref{propo1}; let us define the unitary $U_\psi$ from $H_\psi\underset{\nu^o}{_a\otimes_\beta}H_\Phi$ onto $H_\psi\underset{\nu}{_b\otimes_\alpha}H_\Phi$ by the formula :
\[U_\psi=\tilde{V}^*J_{\tilde{\psi}_0}\tilde{V}(J_\psi\underset{\nu^o}{_a\otimes_\beta}J_{\widehat{\Phi}})\]
Then :
\newline
(i) for all $y$ in $(A\rtimes_\mathfrak a\gG)^{\tilde{\mathfrak a}}$, we have $y=U_\psi(\pi_\psi(y)\underset{N^o}{_a\otimes_\beta}1)U_\psi^*$; 
\newline
(ii) for all $t\in\mathbb{R}$, we have $\sigma_t^{\tilde{\psi}_0}(y)=\sigma_t^\psi(y)$; 
\newline
(iii) we have $U_\psi(J_\psi\underset{\nu}{_b\otimes_\alpha}J_{\widehat{\Phi}})=(J_\psi\underset{\nu^o}{_a\otimes_\beta}J_{\widehat{\Phi}})U_\psi^*$. }

\begin{proof}
Let us first remark that it is possible to define on elementary tensors an antilinear surjective isometry 
$J_\psi\underset{\nu^o}{_a\otimes_\beta}J_{\widehat{\Phi}}$ from $H_\psi\underset{\nu^o}{_a\otimes_\beta}H_\Phi$ onto $H_\psi\underset{\nu}{_b\otimes_\alpha}H_\Phi$, whose inverse is $J_\psi\underset{\nu}{_b\otimes_\alpha}J_{\widehat{\Phi}}$, which is defined the same way. 
\newline
Let us take now $y\in D(\sigma^\psi_{i/2})$; we have, for all $x\in\gN_\psi$, $a\in\gN_{\widehat{\Phi}^c}\cap\gN_{\hat{T}^c}$ :
\begin{eqnarray*}
\Lambda_{\tilde{\psi}_0}[(1\underset{N}{_b\otimes_\alpha}a)xy]
&=&
\tilde{V}(\Lambda_\psi(xy)\underset{\nu}{_b\otimes_\alpha}\Lambda_{\widehat{\Phi}^c}(a))\\
&=&\tilde{V}(J_\psi\pi_\psi(\sigma_{i/2}^\psi(y))^*J_\psi\underset{N}{_b\otimes_\alpha}1)(\Lambda_\psi(x)\underset{\nu}{_b\otimes_\alpha}\Lambda_{\widehat{\Phi}^c}(a))
\end{eqnarray*}
which is equal to :
\[\tilde{V}(J_\psi\pi_\psi(\sigma_{i/2}^\psi(y))^*J_\psi\underset{N}{_b\otimes_\alpha}1)\tilde{V}^*\Lambda_{\tilde{\psi}_0}[(1\underset{N}{_b\otimes_\alpha}a)x]\]
As the linear set generated by all the elements of the form $(1\underset{N}{_b\otimes_\alpha}a)x$ is, by construction of $\tilde{\psi}_0$, a core for $\Lambda_{\tilde{\psi}_0}$, we get that, for any $X\in\gN_{\tilde{\psi}_0}$, we have :
\[\Lambda_{\tilde{\psi}_0}(Xy)=\tilde{V}(J_\psi\pi_\psi(\sigma_{i/2}^\psi(y))^*J_\psi\underset{N}{_b\otimes_\alpha}1)\tilde{V}^*\Lambda_{\tilde{\psi}_0}(X)\]
from which we deduce that $y$ belongs to $D(\sigma_{i/2}^{\tilde{\psi}_0})$, and that :
\[J_{\tilde{\psi}_0}\pi_{\tilde{\psi}_0}(\sigma_{i/2}^{\tilde{\psi}_0}(y)^*)J_{\tilde{\psi}_0}=\tilde{V}(J_\psi\pi_\psi(\sigma_{i/2}^\psi(y))^*J_\psi\underset{N}{_b\otimes_\alpha}1)\tilde{V}^*\]
which can be written also, thanks to \ref{propo1}(ii) :
\[\sigma_{i/2}^{\tilde{\psi}_0}(y)=U_\psi(\pi_\psi(\sigma_{i/2}^\psi(y))\underset{N^o}{_a\otimes_\beta}1)U_\psi^*\]
Taking the adjoints, we get that, if $y'$ belongs to $D(\sigma^\psi_{-i/2})$, then $y'$ belongs to $D(\sigma^{\tilde{\psi}_0}_{-i/2})$, and :
\[\sigma_{-i/2}^{\tilde{\psi}_0}(y')=U_\psi(\pi_\psi(\sigma_{-i/2}^\psi(y'))\underset{N^o}{_a\otimes_\beta}1)U_\psi^*\]
If now $z$ belongs to $D(\sigma^\psi_{-i})$, then $z$ belongs to $D(\sigma^{\tilde{\psi}_0}_{-i/2})$, and :
\[\sigma_{-i/2}^{\tilde{\psi}_0}(z)=U_\psi(\pi_\psi(\sigma_{-i/2}^\psi(z))\underset{N^o}{_a\otimes_\beta}1)U_\psi^*\]
But, as $\sigma^\psi_{-i}(z)$ belongs to $D(\sigma^\psi_{i/2})$, $\sigma^\psi_{-i}(z)$ belongs to $D(\sigma^{\tilde{\psi}_0}_{i/2})$ and :
\[\sigma_{i/2}^{\tilde{\psi}_0}(\sigma^\psi_{-i}(z))=U_\psi(\pi_\psi(\sigma^\psi_{i/2}(\sigma^\psi_{-i}(z)))\underset{N^o}{_a\otimes_\beta}1)U_\psi^*=U_\psi(\pi_\psi(\sigma^\psi_{-i/2}(z))\underset{N^o}{_a\otimes_\beta}1)U_\psi^*\]
and we get that $\sigma_{-i/2}^{\tilde{\psi}_0}(z)=\sigma_{i/2}^{\tilde{\psi}_0}(\sigma^\psi_{-i}(z))$; therefore, $z$ belongs to $D(\sigma^{\tilde{\psi}_0}_{-i})$, and $\sigma^{\tilde{\psi}_0}_{-i}(z)=\sigma_{-i}^\psi(z)$. Then, using (\cite{H}, 4.3 and 4.4), we obtain (ii). 
Let us return now to the formula obtained for $y\in D(\sigma^\psi_{i/2})$; we had obtained :
\[\sigma_{i/2}^{\tilde{\psi}_0}(y)=U_\psi(\pi_\psi(\sigma_{i/2}^\psi(y))\underset{N^o}{_a\otimes_\beta}1)U_\psi^*\]
Using (ii), we get :
\[\sigma_{i/2}^\psi(y)=U_\psi(\pi_\psi(\sigma_{i/2}^\psi(y))\underset{N^o}{_a\otimes_\beta}1)U_\psi^*\]
from which, by density, one obtains (i). 
\newline
From the definition of $U_\psi$, one gets :
\[U_\psi^*=(J_\psi\underset{\nu}{_b\otimes_\alpha}J_{\widehat{\Phi}})\tilde{V}^*J_{\tilde{\psi}_0}\tilde{V}\]
and the formula (iii) is straightforward. \end{proof}

%%%propo3
\subsection{Proposition}
\label{propo3}
{\it Let's take the hypothesis and notations of \ref{defdualw}, \ref{propdualw} and \ref{propo1}; then,
the weight $\tilde{\psi}_0$ is $\hat{\delta}^{-1}$-invariant with respect to the dual action $\tilde{\mathfrak a}$.}

\begin{proof}
By \ref{Phideltainv} applied to $\widehat{\gG}^c$, we know that $\widehat{\Phi}^c$ is $\hat{\delta}^{-1}$-invariant with respect to the action $\widehat{\Gamma}^c$ of $\widehat{\gG}^c$ on $\widehat{M}'$. Moreover, by \ref{propW} applied to $\widehat{\Gamma}^c$, we get that :
\[\Lambda_{\widehat{\Phi}^c}[((id\underset{N^o}{_{\hat{\alpha}}*_\beta}\omega_{\eta, \xi})\widehat{\Gamma}^c(a)]=(id*\omega_{\hat{\delta}^{-1/2}\eta, \xi})(W^{oc})\Lambda_{\widehat{\Phi}^c}(a)\]
for all $a\in\gN_{\widehat{\Phi}^c}$, $\xi\in D((H_\Phi)_\beta, \nu^o)$ and $\eta\in D((H_\Phi)_\beta, \nu^o)$, such that $\hat{\delta}^{-1/2}\eta$ belongs to $D(_{\hat{\alpha}}H_\Phi, \nu)$. 
\newline
Therefore, for any $x\in\gN_\psi$ and $a\in \gN_{\widehat{T}^c}\cap \gN_{\widehat{\Phi}^c}$, we have :
\[\Lambda_{\tilde{\psi}_0}[(id\underset{N^o}{_{\hat{\alpha}}*_\beta}\omega_{\eta, \xi})\tilde{\mathfrak a}((1\underset{N}{_b\otimes_\alpha}a)x)]=\Lambda_{\tilde{\psi}_0}[(1\underset{N}{_b\otimes_\alpha}(id\underset{N^o}{_{\hat{\alpha}}*_\beta}\omega_{\eta, \xi})\widehat{\Gamma}^c(a)x]\]
which is equal to :
\begin{multline*}
\tilde{V}[\Lambda_\psi(x)\underset{\nu}{_b\otimes_\alpha}\Lambda_{\widehat{\Phi}^c}[(id\underset{N^o}{_{\hat{\alpha}}*_\beta}\omega_{\eta, \xi})\widehat{\Gamma}^c(a)]]\\
=\tilde{V}[1\underset{N}{_b\otimes_\alpha}(id*\omega_{\hat{\delta}^{-1/2}\eta, \xi})(W^{oc})](\Lambda_\psi(x)\underset{\nu}{_b\otimes_\alpha}\Lambda_{\widehat{\Phi}^c}(a))\\
=\tilde{V}[1\underset{N}{_b\otimes_\alpha}(id*\omega_{\hat{\delta}^{-1/2}\eta, \xi})(W^{oc})]\tilde{V}^*\Lambda_{\tilde{\psi}_0}((1\underset{N}{_b\otimes_\alpha}a)\mathfrak a(x))
\end{multline*}
which, by definition of $\tilde{\psi}_0$, and the closedness of $\Lambda_{\tilde{\psi}_0}$, allows us to write, for all $X\in\gN_{\tilde{\psi}_0}$:
\[\Lambda_{\tilde{\psi}_0}((id\underset{N^o}{_{\hat{\alpha}}*_\beta}\omega_{\eta, \xi})\tilde{\mathfrak a}(X))=\tilde{V}[1\underset{N}{_b\otimes_\alpha}(id*\omega_{\hat{\delta}^{-1/2}\eta, \xi})(W^{oc})]\tilde{V}^*\Lambda_{\tilde{\psi}_0}(X)\]
Taking now $(\xi_i)_{i\in I}$ an orthogonal $(\beta, \nu^o)$-basis of $H_\Phi$, we have :
\[\tilde{\psi}_0((id\underset{N^o}{_{\hat{\alpha}}*_\beta}\omega_{\eta})\tilde{\mathfrak a}(X^*X))
=
\sum_i\|\Lambda_{\tilde{\psi}_0}((id\underset{N^o}{_{\hat{\alpha}}*_\beta}\omega_{\eta, \xi_i})\tilde{\mathfrak a}(X))\|^2\]
which is equal to :
\begin{multline*}
\sum_i\|[1\underset{N}{_b\otimes_\alpha}(id*\omega_{\hat{\delta}^{-1/2}\eta, \xi_i})(W^{oc})]\tilde{V}^*\Lambda_{\tilde{\psi}_0}(X)\|^2=\\
=\|(1\underset{N}{_b\otimes_\alpha}W^{oc})\tilde{V}^*(\Lambda_{\tilde{\psi}_0}(X)\underset{\nu}{_{\beta_0}\otimes_{\hat{\alpha}}}\hat{\delta}^{-1/2}\eta)\|^2
=\|\Lambda_{\tilde{\psi}_0}(X)\underset{\nu}{_{\beta_0}\otimes_{\hat{\alpha}}}\hat{\delta}^{-1/2}\eta\|^2
\end{multline*}
where we put, for $n\in N$, $\beta_0(n)=J_{\tilde{\psi}_0}(1\underset{N}{_b\otimes_\alpha}\hat{\alpha}(n^*))J_{\tilde{\psi}_0}$. Which is the result.  \end{proof}

%%%%Stildepsi
\subsection{Lemma}
\label{Stildepsi}
{\it Let $(b, \mathfrak a)$ be an action of a measured quantum groupoid $\gG$ on a von Neumann algebra $A$, and let $\psi$ be a normal semi-finite faithful weight on $A$. Let $A\rtimes_\mathfrak a\gG$ be the crossed product, and $\tilde{\psi}_0$ be the weight defined in \ref{propo1}. Then, the linear span generated by the elements $\pi_{\tilde{\psi}_0}(x^*)\tilde{V}(\Lambda_\psi(y)\underset{\nu}{_b\otimes_\alpha}\Lambda_{\widehat{\Phi}^c}(a))$, where $x$, $y$ belong to $\gN_\psi$ and $a$ belongs to $\gN_{\widehat{T}^c}\cap\gN_{\widehat{\Phi}^c}\cap\gN_{\widehat{T}^c}^*\cap\gN_{\widehat{\Phi}^c}^*$, is a core for $S_{\tilde{\psi}_0}$, and we have :}
\[S_{\tilde{\psi}_0}\pi_{\tilde{\psi}_0}(x^*)\tilde{V}(\Lambda_\psi(y)\underset{\nu}{_b\otimes_\alpha}\Lambda_{\widehat{\Phi}^c}(a))=
\pi_{\tilde{\psi}_0}(y^*)\tilde{V}(\Lambda_\psi(x)\underset{\nu}{_b\otimes_\alpha}\Lambda_{\widehat{\Phi}^c}(a^*))\]

\begin{proof}
By definition of $\tilde{\psi}_0$ (\ref{propo3}(ii)), the linear span generated by the elements $(1\underset{N}{_b\otimes_\alpha}a)x$, where $a$ belongs to $\gN_{\widehat{T}^c}\cap\gN_{\widehat{\Phi}^c}$, and $x$ belongs to $\gN_\psi$ is a core for $\Lambda_{\tilde{\psi}_0}$. We deduce easily the result from this. \end{proof}

%%%Upsi
\subsection{Proposition}
\label{Upsi}
{\it Let $(b, \mathfrak a)$ be an action of a measured quantum groupoid $\gG$ on a von Neumann algebra $A$, and let $\psi$ be a normal semi-finite faithful weight on $A$. Let $A\rtimes_\mathfrak a\gG$ be the crossed product, and $\tilde{\psi}_0$ be the weight defined in \ref{propo1}. Then, the unitary $U_\psi$ defined in \ref{propo2} satifies, for all $b\in M$ :}
\[(1\underset{N}{_b\otimes_\alpha}J_\Phi bJ_\Phi)U_\psi=U_\psi (1\underset{N^o}{_a\otimes_\beta}J_\Phi bJ_\Phi)\]

\begin{proof}
Using the notations of \ref{thleftw} and \ref{Stildepsi}, we have :
\begin{multline*}
(1\underset{N}{_b\otimes_\alpha}J_\Phi bJ_\Phi)\tilde{V}^*S_{\tilde{\psi}_0}\pi_{\tilde{\psi}_0}(x^*)\tilde{V}(\Lambda_\psi(y)\underset{\nu}{_b\otimes_\alpha}\Lambda_{\widehat{\Phi}^c}(a^*))
=\\
(1\underset{N}{_b\otimes_\alpha}J_\Phi bJ_\Phi)\tilde{V}^*\pi_{\tilde{\psi}_0}(y^*)\tilde{V}(\Lambda_\psi(x)\underset{\nu}{_b\otimes_\alpha}\Lambda_{\widehat{\Phi}^c}(a))
\end{multline*}
which is equal to :
\begin{eqnarray*}
(1\underset{N}{_b\otimes_\alpha}J_\Phi bJ_\Phi)y^*(\Lambda_\psi(x)\underset{\nu}{_b\otimes_\alpha}\Lambda_{\widehat{\Phi}^c}(a))
&=&
y^*(1\underset{N}{_b\otimes_\alpha}J_\Phi bJ_\Phi)(\Lambda_\psi(x)\underset{\nu}{_b\otimes_\alpha}\Lambda_{\widehat{\Phi}^c}(a))\\
&=&
y^*(\Lambda_\psi(x)\underset{\nu}{_b\otimes_\alpha}J_\Phi bJ_\Phi\Lambda_{\widehat{\Phi}^c}(a))
\end{eqnarray*}
Using \ref{thL3}(v), we get that, for all $t\in\mathbb{R}$, $\Delta_{\widehat{\Phi}}^{it}J_\Phi bJ_\Phi\Delta_{\widehat{\Phi}}^{-it}$ belongs to $M'$, and we define this way an automorphism group $\mu_t$ of $M'$; if we suppose that $J_\Phi R(b^*)J_\Phi$ belongs to $D(\mu_{-i/2})$, we get that :
\[J_\Phi bJ_\Phi\Lambda_{\widehat{\Phi}^c}(a)=J_{\widehat{\Phi}}\Delta_{\widehat{\Phi}}^{-1/2}\mu_{-i/2}(J_\Phi R(b^*)J_\Phi)\Lambda_{\widehat{\Phi}^c}(a^*)\]
and, therefore, using \ref{Stildepsi}, we get :
\begin{multline*}
(1\underset{N}{_b\otimes_\alpha}J_\Phi bJ_\Phi)\tilde{V}^*S_{\tilde{\psi}_0}\pi_{\tilde{\psi}_0}(x^*)\tilde{V}((\Lambda_\psi(y)\underset{\nu}{_b\otimes_\alpha}\Lambda_{\widehat{\Phi}^c}(a^*))\\
=y^*(\Lambda_\psi(x)\underset{\nu}{_b\otimes_\alpha}J_{\widehat{\Phi}}\Delta_{\widehat{\Phi}}^{-1/2}\mu_{-i/2}(J_\Phi R(b^*)J_\Phi)\Lambda_{\widehat{\Phi}^c}(a^*)\\
=\tilde{V}^*\pi_{\tilde{\psi}_0}(y^*)\tilde{V}(\Lambda_\psi(x)\underset{\nu}{_b\otimes_\alpha}J_{\widehat{\Phi}}\Delta_{\widehat{\Phi}}^{-1/2}\mu_{-i/2}(J_\Phi R(b^*)J_\Phi)\Lambda_{\widehat{\Phi}^c}(a^*)\\
=\tilde{V}^*S_{\tilde{\psi}_0}\pi_{\tilde{\psi}_0}(x^*)\tilde{V}(\Lambda_\psi(y)\underset{\nu}{_b\otimes_\alpha}\mu_{-i/2}(J_\Phi R(b^*)J_\Phi)\Lambda_{\widehat{\Phi}^c}(a^*))\\
=\tilde{V}^*S_{\tilde{\psi}_0}\tilde{V}x^*(\Lambda_\psi(y)\underset{\nu}{_b\otimes_\alpha}\mu_{-i/2}(J_\Phi R(b^*)J_\Phi)\Lambda_{\widehat{\Phi}^c}(a^*))
\end{multline*}
and, by the closedness of $S_{\tilde{\psi}_0}$, we get :
\[(1\underset{N}{_b\otimes_\alpha}J_\Phi bJ_\Phi)\tilde{V}^*S_{\tilde{\psi}_0}\tilde{V}\subset \tilde{V}^*S_{\tilde{\psi}_0}\tilde{V}(1\underset{N}{_b\otimes_\alpha}\mu_{-i/2}(J_\Phi R(b^*)J_\Phi))\]
Taking the adjoints, we have $J_\Phi R(b)J_\Phi\in D(\mu_{i/2})$ and :
\[(1\underset{N}{_b\otimes_\alpha}\mu_{i/2}(J_\Phi R(b)J_\Phi))\tilde{V}^*F_{\tilde{\psi}_0}\tilde{V}\subset
\tilde{V}^*F_{\tilde{\psi}_0}\tilde{V}(1\underset{N}{_b\otimes_\alpha}J_\Phi b^*J_\Phi)\]
Therefore, if we suppose that $J_\Phi bJ_\Phi$ belongs to $\mu_{-i}$, we get :
\[(1\underset{N}{_b\otimes_\alpha}J_\Phi bJ_\Phi)\tilde{V}^*\Delta_{\tilde{\psi}_0}\tilde{V}\subset \tilde{V}^*\Delta_{\tilde{\psi}_0}\tilde{V}((1\underset{N}{_b\otimes_\alpha}\mu_{-i}(J_\Phi b J_\Phi))\]
So, for any $X\in\mathcal L(H_\psi\underset{\nu}{_b\otimes_\alpha}H_\phi)$, let $\epsilon_t(X)=\tilde{V}^*\Delta_{\tilde{\psi}_0}^{it}\tilde{V}X\tilde{V}^*\Delta_{\tilde{\psi}_0}^{-it}\tilde{V}$; we then get that $(1\underset{N}{_b\otimes_\alpha}J_\Phi bJ_\Phi)$ belongs to $D(\epsilon_{-i})$, and, therefore, for any $b$ in $M$, we have :
\[\epsilon_t(1\underset{N}{_b\otimes_\alpha}J_\Phi bJ_\Phi)=1\underset{N}{_b\otimes_\alpha}\mu_t(J_\Phi bJ_\Phi)\]
and :
\[(1\underset{N}{_b\otimes_\alpha}J_\Phi bJ_\Phi)\tilde{V}^*J_{\tilde{\psi}_0}\tilde{V}=\tilde{V}^*J_{\tilde{\psi}_0}\tilde{V}(1\underset{N}{_b\otimes_\alpha}J_\Phi R(b^*)J_\Phi)\]
and, therefore :
\[(1\underset{N}{_b\otimes_\alpha}J_\Phi bJ_\Phi)U_\psi=U_\psi(1\underset{N}{_b\otimes_\alpha}J_\Phi bJ_\Phi)\]
\end{proof}

%%%%corUpsi
\subsection{Corollary}
\label{corUpsi}
{\it Let $(b, \mathfrak a)$ be an action of a measured quantum groupoid $\gG$ on a von Neumann algebra $A$, and let $\psi$ be a normal semi-finite faithful weight on $A$. Let $A\rtimes_\mathfrak a\gG$ be the crossed product, and $\tilde{\psi}_0$ be the weight defined in \ref{propo1}. Then, for all $n\in N$, we have :}
\[J_{\tilde{\psi}_0}\pi_{\tilde{\psi}_0}(1\underset{N}{_b\otimes_\alpha}\hat{\alpha}(n))^*J_{\tilde{\psi}_0}
=\tilde{V}(1\underset{N}{_b\otimes_\alpha}\hat{\beta}(n))\tilde{V}^*\]

\begin{proof}
We have, using \ref{propo1}(ii), 
\[J_{\tilde{\psi}_0}\pi_{\tilde{\psi}_0}(1\underset{N}{_b\otimes_\alpha}\hat{\alpha}(n))^*J_{\tilde{\psi}_0}
=
J_{\tilde{\psi}_0}\tilde{V}(1\underset{N}{_b\otimes_\alpha}\hat{\alpha}(n)^*)\tilde{V}^*J_{\tilde{\psi}_0}\]
which, using \ref{propo2}, is equal to :
\[\tilde{V}U_\psi(J_\psi\underset{N}{_b\otimes_\alpha}J_{\widehat{\Phi}})(1\underset{N}{_b\otimes_\alpha}\hat{\alpha}(n)^*)(J_\psi\underset{N^o}{_a\otimes_\beta}J_{\widehat{\Phi}})U_\psi^*\tilde{V}^*\]
and, using \ref{Upsi}, is equal to :
\[\tilde{V}U_\psi(1\underset{N^o}{_a\otimes_\beta}\hat{\beta}(n))U_\psi^*\tilde{V}^*
=\tilde{V}(1\underset{N}{_b\otimes_\alpha}\hat{\beta}(n))\tilde{V}^*\]
\end{proof}

%%%%thdual
\subsection{Theorem}
\label{thdual}
{\it Let $(b, \mathfrak a)$ be an action of a measured quantum groupoid $\gG$ on a von Neumann algebra $A$; let $A\rtimes_\mathfrak a\gG$ be the crossed product, $\tilde{\mathfrak a}$ the dual action of $\hat{\gG}^c$ on $A\rtimes_\mathfrak a\gG$, and let $\psi$ be a normal semi-finite faithful weight on $(A\rtimes_\mathfrak a\gG)^{\tilde{\mathfrak a}}$. Let $\tilde{\psi}_0$ be the weight defined in \ref{propo1}. Then, 
\newline
(i) the weight $\tilde{\psi}_0$ is $\hat{\delta}^{-1}$-invariant with respect to the dual action $\tilde{\mathfrak a}$, and bears the density property defined in \ref{defdeltainv}. 
\newline
(ii) there exists a standard implementation of the dual action on the Hilbert space $H_{\tilde{\psi}_0}$.
\newline
(iii) the dual action is saturated. }

\begin{proof}
As $V[D((H_\psi)_b, \nu^o)\underset{\nu}{_b\otimes_\alpha}D(_{\hat{\alpha}}H_\Phi, \nu)\cap D((H_\Phi)_{\hat{\beta}}, \nu^o)]$ is dense in $H_{\tilde{\psi}_0}$, we get that $\tilde{\psi}_0$ bears the density property defined in \ref{defdeltainv}; so (i) is then given by \ref{propo3}. Then (ii) is given by \ref{corstandardV}(i), and (iii) is given by \ref{Propdual} and \ref{corstandardV}(ii). \end{proof}

%%%thdual2
\subsection{Theorem}
\label{thdual2}
{\it Let $(b, \mathfrak a)$ be an action of a measured quantum groupoid $\gG$ on a von Neumann algebra $A$; let $A\rtimes_\mathfrak a\gG$ be the crossed product, $\tilde{\mathfrak a}$ the dual action of $\hat{\gG}^c$ on $A\rtimes_\mathfrak a\gG$, and let $\psi$ be a normal semi-finite faithful weight on $(A\rtimes_\mathfrak a\gG)^{\tilde{\mathfrak a}}$. Let $\tilde{\psi}$ be the lifted weight defined in \ref{defdualw},and $\tilde{\psi}_0$ be the weight defined in \ref{propo1}. Then :
\newline
(i) we have $\tilde{\psi}=\tilde{\psi}_0$; 
\newline
(ii) the linear set generated by the elements of the form $(1\underset{N}{_b\otimes_\alpha}a)x$, with $a\in \gN_{\widehat{\Phi}^c}\cap\gN_{\hat{T}^c}$ and $x\in\gN_\psi$, is a core for $\Lambda_{\tilde{\psi}}$. 
\newline
(iii) the weight $\tilde{\psi}$ is $\hat{\delta}^{-1}$-invariant with respect to the dual action $\tilde{\mathfrak a}$, and bears the density property defined in \ref{defdeltainv}. }

\begin{proof}
Thanks to \ref{thdual}(i), we can apply \ref{thV2}(ii) to the action $\tilde{\mathfrak a}$ and to the weight $\tilde{\psi}_0$, and we have, for any $x\in A\rtimes_\mathfrak a\gG$ :
\[\tilde{\mathfrak a}(\sigma_t^{\tilde{\psi}_0}(x))=(\sigma_t^{\tilde{\psi}_0}\underset{N^o}{_{\hat{\alpha}}*_\beta}\hat{\tau}_{-t}^c\sigma_{-t}^{\widehat{\Phi}^c\circ R^c}\sigma_t^{\widehat{\Phi}^c})\tilde{\mathfrak a}(x)\]
and, using \ref{thL2}(vi) and (vii) applied to $\widehat{\Phi}^c$, we get that, for any positive $x$ in $A\rtimes_\mathfrak a\gG$, we have $T_{\tilde{\mathfrak a}}(\sigma_t^{\tilde{\psi}_0}(x))=\sigma_t^{\tilde{\psi}_0}(T_{\tilde{\mathfrak a}}(x))$, which implies, using \ref{propo2}(ii), that :
\[T_{\tilde{\mathfrak a}}(\sigma_t^{\tilde{\psi}_0}(x))=\sigma_t^\psi (T_{\tilde{\mathfrak a}}(x))\]
and, therefore, $\tilde{\psi}(\sigma_t^{\tilde{\psi}_0}(x))=\tilde{\psi}(x)$. 
\newline
By construction (\ref{propo1}(i)), we have $\tilde{\psi}(x)=\tilde{\psi}_0(x)$ for all $x\in\gM_{\tilde{\psi}_0}$; from which we can deduce now (i). Then (ii) is just given by (i) and the definition of $\tilde{\psi}_0$, and (iii) by (i) and \ref{thdual}(i). \end{proof}

%%%%%%propcrossed
\subsection{Proposition}
\label{propcrossed}
{\it Let $\gG=(N, M, \alpha, \beta, \Gamma, T, T', \nu)$ be a measured quantum groupo-id, $_a\gH_b$ a $N-N$ bimodule, $A$ a von Neumann algebra such that $b(N)\subset A\subset a(N)'$, $V$ a corepresentation of $\gG$ on $_a\gH_b$ which implements an action $(b, \mathfrak a)$ of $\gG$ on $A$; let $A\rtimes_{\mathfrak a}\gG$ be the crossed-product, and $\tilde{\mathfrak a}$ the dual action of of $(\widehat{\gG})^c$ on the crossed product; then :
\newline
(i) $A\rtimes_{\mathfrak a}\gG\subset V(\mathcal L(\gH)\underset{N}{_a*_\beta}\widehat{M}')V^*$
\newline
(ii) $(A\rtimes_{\mathfrak a}\gG)^{\tilde{\mathfrak a}}\subset Sat \mathfrak a$. }

\begin{proof}
We have obtained in \ref{lemV}(ii) that $1\underset{N}{_b\otimes_\alpha}\widehat{M}'\subset V(\mathcal L(\gH)\underset{N}{_a*_\beta}\widehat{M}')V^*$; as $\mathfrak a (A)=V(A\underset{N^o}{_a\otimes_\beta}1)V^*$, we obtain (i) by definition of $A\rtimes_{\mathfrak a}\gG$. Then (ii) is a corollary of (i) and \ref{dualaction}(iii). \end{proof}

%%%%%defB
\subsection{Definition}
\label{defB}
Let $\gG=(N, M, \alpha, \beta, \Gamma, T, T', \nu)$ be a measured quantum groupo-id, $A$ a von neumann algebra, $(b, \mathfrak a)$ an action of $\gG$ on $A$; following (\cite{ES1}, II5), we shall say that $\mathfrak a$ satisfies property (B) if we have $\mathfrak a(A)=(A\rtimes_{\mathfrak a}\gG)^{\tilde{\mathfrak a}}$. 

%%%%%satB
\subsection{Proposition}
\label{satB}
{\it  Let $\gG=(N, M, \alpha, \beta, \Gamma, T, T', \nu)$ be a measured quantum grou-poid, $_a\gH_b$ a $N-N$ bimodule, $A$ a von Neumann algebra such that $b(N)\subset A\subset a(N)'$, $V$ a corepresentation of $\gG$ on $_a\gH_b$ which implements an action $(b, \mathfrak a)$ of $\gG$ on $A$; let $A\rtimes_{\mathfrak a}\gG$ be the crossed-product, and $\tilde{\mathfrak a}$ the dual action of of $(\widehat{\gG})^c$ on the crossed product; let us suppose that the action $\mathfrak a$ is saturated, in the sense of \ref{defsat}. Then, the action $\mathfrak a$ satisfies property (B) (in the sense of \ref{defB}). }

\begin{proof}
Using \ref{dualaction}(iii) and \ref{propcrossed}(ii), we have :
\[\mathfrak a(A)\subset (A\rtimes_{\mathfrak a}\gG)^{\tilde{\mathfrak a}}\subset Sat \mathfrak a\]
Therefore, the result is trivial. \end{proof}

%%%%corB
\subsection{Corollary}
\label{corB}
{\it Let $\gG=(N, M, \alpha, \beta, \Gamma, T, T', \nu)$ be a measured quantum groupoid; let $A$ be a von Neumann algebra acting on a Hilbert space $\gH$, $(b, \mathfrak a)$ an action of $\gG$ on $A$; the dual action $\tilde{\mathfrak a}$ of $(\widehat{\gG})^c$ on the crossed product $A\rtimes_{\mathfrak a}\gG$ satisfies property (B) of \ref{defB}.}

\begin{proof}
This is a straightforward corollary of \ref{thdual}(ii) and \ref{satB}. \end{proof}

%%%%%biduality
\section{Biduality}
\label{biduality}
In that chapter, we prove a first version of biduality theorems (\ref{propduality2}, \ref{propduality3}), which implies that any action satisfies properties (A) (\ref{thduality1}(ii)) and (B) (\ref{thduality1}(i)). We then get 
the biduality theorem (\ref{thduality2}) and commutation theorem (\ref{thcom2}).  

%%%lemcrossed
\subsection{Lemma}
\label{lemcrossed}
{\it Let $\gG=(N, M, \alpha, \beta, \Gamma, T, T', \nu)$ be a measured quantum groupoid; let $A$ be a von Neumann algebra acting on a Hilbert space $\gH$, $(b, \mathfrak a)$ an action of $\gG$ on $A$, $A\rtimes_{\mathfrak a}\gG$ the crossed product of $A$ by $\gG$ via this action. Then, we have :}
\[[(\mathfrak a\underset{N}{_b*_\alpha}id)\mathfrak a (A)\cup 1\underset{N}{_b\otimes_\alpha}\widehat{M}'\underset{N}{_\beta*_\alpha}\mathcal L(H_\Phi)]''=A\rtimes_{\mathfrak a}\gG\underset{N}{_\beta*_\alpha}\mathcal L(H_\Phi)\]

\begin{proof}
We have $\mathfrak a(A)\subset A\underset{N}{_b*_\alpha}\mathcal L(H_\Phi)$, and, therefore:  
\[(\mathfrak a\underset{N}{_b*_\alpha}id)\mathfrak a (A)\subset\mathfrak a(A)\underset{N}{_b*_\alpha}\mathcal L(H_\Phi)\]
from which we get :
\[[(\mathfrak a\underset{N}{_b*_\alpha}id)\mathfrak a (A)\cup 1\underset{N}{_b\otimes_\alpha}\widehat{M}'\underset{N}{_\beta*_\alpha}\mathcal L(H_\Phi)]''\subset A\rtimes_{\mathfrak a}\gG\underset{N}{_\beta*_\alpha}\mathcal L(H_\Phi)\]
Conversely, let $X\in\mathcal L(\gH)\underset{N}{_b*_\alpha}\widehat{M}$, such that $X\underset{N}{_\beta\otimes_\alpha}1$ commutes with $(\mathfrak a\underset{N}{_b*_\alpha}id)\mathfrak a (A)$. As, for $x\in A$, we have, using \ref{thL5}(v) :
\[(\mathfrak a\underset{N}{_b*_\alpha}id)\mathfrak a (x)=(id\underset{N}{_b*_\alpha}\Gamma)\mathfrak a(x)
=(1\underset{N}{_b\otimes_\alpha}\sigma_\nu W^o\sigma_\nu)^*(\mathfrak a(x)\underset{N^o}{_{\hat{\alpha}}\otimes_\beta}1)(1\underset{N}{_b\otimes_\alpha}\sigma_\nu W^o\sigma_\nu)\]
and, therefore :
\[\mathfrak a(x)\underset{N^o}{_{\hat{\alpha}}\otimes_\beta}1=
(1\underset{N}{_b\otimes_\alpha}\sigma_\nu W^o\sigma_\nu)(\mathfrak a\underset{N}{_b*_\alpha}id)\mathfrak a (x)(1\underset{N}{_b\otimes_\alpha}\sigma_\nu W^o\sigma_\nu)^*\]
As $X$ belongs to $\mathcal L(\gH)\underset{N}{_b*_\alpha}\widehat{M}$, we get that :
\[(1\underset{N}{_b\otimes_\alpha}\sigma_\nu W^o\sigma_\nu)(X\underset{N}{_\beta\otimes_\alpha}1)=(X\underset{N^o}{_{\hat{\alpha}}\otimes_{\hat{\beta}}}1)(1\underset{N}{_b\otimes_\alpha}\sigma_\nu W^o\sigma_\nu)\]
and, therefore, we get that $X\underset{N^o}{_{\hat{\alpha}}\otimes_{\hat{\beta}}}1$ commutes with 
$\mathfrak a(x)\underset{N^o}{_{\hat{\alpha}}\otimes_\beta}1$, which means that $X$ belongs to $\mathfrak a(A)'$; 
as $X$ belongs to $\mathcal L(\gH)\underset{N}{_b*_\alpha}\widehat{M}$, we finally get that $X$ belongs to $(A\rtimes_{\mathfrak a}\gG)'$. So, we have proved that :
\[[(\mathfrak a\underset{N}{_b*_\alpha}id)\mathfrak a (A)\cup 1\underset{N}{_b\otimes_\alpha}\widehat{M}'\underset{N}{_\beta*_\alpha}\mathcal L(H_\Phi)]'\subset [A\rtimes_{\mathfrak a}\gG\underset{N}{_\beta*_\alpha}\mathcal L(H_\Phi)]'\]
from which we get the result. \end{proof}

%%%%propduality1
\subsection{Proposition}
\label{propduality1}
{\it Let $\gG=(N, M, \alpha, \beta, \Gamma, T, T', \nu)$ be a measured quantum groupo-id; let $A$ be a von Neumann algebra acting on a Hilbert space $\gH$, $(b, \mathfrak a)$ an action of $\gG$ on $A$, $\tilde{\mathfrak a}$ the dual action of $(\widehat{\gG})^c$ on the crossed product $A\rtimes_{\mathfrak a}\gG$, $\tilde{\tilde{\mathfrak a}}$ the bidual action of $({\gG})^{oc}$ on the double crossed-product $(A\rtimes_{\mathfrak a}\gG)\rtimes_{\tilde{\mathfrak a}}\widehat{\gG}^o$; we shall denote :
\[\underline{A}=(\mathfrak a (A)\cup 1\underset{N}{_b\otimes_\alpha}\alpha(N)')''\]
and, for any $X\in\underline{A}$ :
\[\underline{\mathfrak a}(X)=(1\underset{N}{_b\otimes_\alpha}\sigma_{\nu^o} W\sigma_{\nu^o})(id\underset{N}{_b*_\alpha}\varsigma_N)(\mathfrak a\underset{N}{_b*_\alpha}id)(X)(1\underset{N}{_b\otimes_\alpha}\sigma_{\nu^o} W\sigma_{\nu^o})^*\]
Moreover, for any $Y\in\mathcal L(\gH\underset{\nu}{_b\otimes_\alpha}H_\Phi\underset{\nu}{_\beta\otimes_\alpha}H_\Phi)$, we define :
\[\Theta(Y)=(1\underset{\nu}{_b\otimes_\alpha}\sigma_\nu W^o\sigma_\nu)Y(1\underset{N}{_b\otimes_\alpha}\sigma_\nu W^o\sigma_\nu)^*\in\mathcal L(\gH\underset{\nu}{_b\otimes_\alpha}H_\Phi\underset{\nu^o}{_{\hat{\alpha}}\otimes_\beta}H_\Phi)\]
Then, for any $X\in\underline{A}$, we have :
\[\Theta (\mathfrak a\underset{N}{_b*_\alpha}id)(\underline{A})=(A\rtimes_{\mathfrak a}\gG)\rtimes_{\tilde{\mathfrak a}}\widehat{\gG}^o\]
\[(\Theta (\mathfrak a\underset{N}{_b*_\alpha}id)\underset{N}{_b*_\alpha}\mathcal I_{\gG})\underline{\mathfrak a}(X)=\tilde{\tilde{\mathfrak a}}\circ\Theta (\mathfrak a\underset{N}{_b*_\alpha}id)(X)\]
and, for any $Z\in A\rtimes_\mathfrak a\gG$, we get that $
\Theta(\mathfrak a\underset{N}{_b*_\alpha}id)(Z)=\tilde{\mathfrak a}(Z)$.}
\begin{proof}
Let $x\in A$; then $\mathfrak a (x)\in \underline{A}$, $(\mathfrak a\underset{N}{_b*_\alpha}id)\mathfrak a (x)=(id\underset{N}{_b*_\alpha}\Gamma)\mathfrak a (x)$, and, using \ref{thL5}, we have :
\[\Theta(\mathfrak a\underset{N}{_b*_\alpha}id)\mathfrak a (x)=\mathfrak a(x)\underset{N^o}{_{\hat{\alpha}}\otimes_\beta}1\]
Let $y\in \widehat{M}'$; then $1\underset{N}{_b\otimes_\alpha}y\in\underline{A}$, $(\mathfrak a\underset{N}{_b*_\alpha}id)(1\underset{N}{_b\otimes_\alpha}y)=1\underset{N}{_b\otimes_\alpha}1\underset{N}{_\beta\otimes_\alpha}y$, and, using \ref{thL3}(iii) and \ref{thL5}(iv), we get :
\[\Theta(\mathfrak a\underset{N}{_b*_\alpha}id)(1\underset{N}{_b\otimes_\alpha}y)=1\underset{N}{_b\otimes_\alpha}\widehat{\Gamma}^c(y)\]
So, for any $Z\in A\rtimes_\mathfrak a\gG$, we get that $\Theta(\mathfrak a\underset{N}{_b*_\alpha}id)(Z)=\tilde{\mathfrak a}(Z)$. 
Let now $z\in M'$; then $1\underset{N}{_b\otimes_\alpha}z\in\underline{A}$, $(\mathfrak a\underset{N}{_b*_\alpha}id)(1\underset{N}{_b\otimes_\alpha}z)=1\underset{N}{_b\otimes_\alpha}1\underset{N}{_\beta\otimes_\alpha}z$, and :
\[\Theta(\mathfrak a\underset{N}{_b*_\alpha}id)(1\underset{N}{_b\otimes_\alpha}z)=1\underset{N}{_b\otimes_\alpha}1\underset{N^o}{_{\hat{\alpha}}\otimes_\beta}z\]
from which we get that $\Theta (\mathfrak a\underset{N}{_b*_\alpha}id)(\underline{A})=(A\rtimes_{\mathfrak a}\gG)\rtimes_{\tilde{\mathfrak a}}\widehat{\gG}^o$. 
Moreover, using \ref{thL1}, we deduce from $(\mathfrak a\underset{N}{_b*_\alpha} id)\mathfrak a (x)=(id\underset{N^o}{_b*_\alpha}\Gamma)\mathfrak a (x)$, that, for all $x\in A$, we have  $\underline{\mathfrak a}(\mathfrak a (x))=\mathfrak a (x)\underset{N}{_{\hat{\beta}}\otimes_\alpha}1$.
From which we deduce that :
\begin{eqnarray*}
(\Theta (\mathfrak a\underset{N}{_b*_\alpha}id)\underset{N}{_b*_\alpha}\mathcal I_{\gG})\underline{\mathfrak a}(\mathfrak a(x))
&=&\mathfrak a (x)\underset{N^o}{_{\hat{\alpha}}\otimes_\beta}1\underset{N}{_{\hat{\beta}}\otimes_{\hat{\alpha}}}1\\
&=&\tilde{\mathfrak a}(\mathfrak a (x))\underset{N}{_{\hat{\beta}}\otimes_{\hat{\alpha}}}1\\
&=&\tilde{\tilde{\mathfrak a}}(\tilde{\mathfrak a} (\mathfrak a(x)))\\
&=&\tilde{\tilde{\mathfrak a}}\circ\Theta (\mathfrak a\underset{N}{_b*_\alpha}id)\mathfrak a(x)
\end{eqnarray*}
For $y\in\widehat{M}'$, we have $\underline{\mathfrak a}(1\underset{N}{_b\otimes_\alpha}y)=1\underset{N}{_b\otimes_\alpha}y\underset{N}{_{\hat{\beta}}\otimes_\alpha}1$
and, therefore :
\begin{eqnarray*}
(\Theta (\mathfrak a\underset{N}{_b*_\alpha}id)\underset{N}{_b*_\alpha}\mathcal I_{\gG})\underline{\mathfrak a}(1\underset{N}{_b\otimes_\alpha}y)
&=&1\underset{N}{_b\otimes_\alpha}\widehat{\Gamma}^c(y)\underset{N}{_{\hat{\beta}}\otimes_{\hat{\alpha}}}1\\
&=&\tilde{\mathfrak a}(1\underset{N}{_b\otimes_\alpha}y)\underset{N}{_{\hat{\beta}}\otimes_{\hat{\alpha}}}1\\
&=&\tilde{\tilde{\mathfrak a}}(\tilde{\mathfrak a} (1\underset{N}{_b\otimes_\alpha}y))\\
&=&\tilde{\tilde{\mathfrak a}}\circ\Theta (\mathfrak a\underset{N}{_b*_\alpha}id)(1\underset{N}{_b\otimes_\alpha}y)
\end{eqnarray*}
from which we deduce that, for all $Z\in A\rtimes_\mathfrak a\gG$, we have  :
\[(\Theta (\mathfrak a\underset{N}{_b*_\alpha}id)\underset{N}{_b*_\alpha}\mathcal I_{\gG})\underline{\mathfrak a}(Z)=\tilde{\tilde{\mathfrak a}}(\tilde{\mathfrak a} (Z))=\tilde{\tilde{\mathfrak a}}\circ\Theta (\mathfrak a\underset{N}{_b*_\alpha}id)(Z)\]
For $z\in M'$, we have $(\mathfrak a\underset{N}{_b*_\alpha}id)(1\underset{N}{_b\otimes_\alpha}z)=1\underset{N}{_b\otimes_\alpha}1\underset{N}{_\beta\otimes_\alpha}z$ and :
\[\underline{\mathfrak a}(1\underset{N}{_b\otimes_\alpha}z)
=(1\underset{N}{_b\otimes_\alpha}\sigma_{\nu^o})(1\underset{N}{_b\otimes_\alpha}(W(1\underset{N}{_\beta\otimes_\alpha}z)W^*))(1\underset{N}{_b\otimes_\alpha}\sigma_{\nu^o})\]
As $W=(J_{\widehat{\Phi}}\underset{N}{_\beta\otimes_
\alpha} J_\Phi)W^*(J_{\widehat{\Phi}}\underset{N}{_\beta\otimes_
\alpha} J_\Phi)$, 
we get that :
\[\underline{\mathfrak a}(1\underset{N}{_b\otimes_\alpha}z)=
(id\underset{N}{_b*_\alpha}\varsigma_{N^o})(1\underset{N^o}{_b\otimes_\alpha}(J_{\widehat{\Phi}}\underset{N}{_\beta\otimes_
\alpha} J_\Phi)\Gamma(J_\Phi zJ_\Phi)(J_{\widehat{\Phi}}\underset{N}{_\beta\otimes_
\alpha} J_\Phi))\]
and :
\begin{eqnarray*}
(\Theta (\mathfrak a\underset{N}{_b*_\alpha}id)\underset{N}{_b*_\alpha}\mathcal I_{\gG})\underline{\mathfrak a}(1\underset{N}{_b\otimes_\alpha}z)
&=&1\underset{N}{_b\otimes_\alpha}1\underset{N^o}{_{\hat{\alpha}}\otimes_\beta}\Gamma^{oc}(z)\\
&=&\tilde{\tilde{\mathfrak a}}(1\underset{N}{_b\otimes_\alpha}1\underset{N^o}{_{\hat{\alpha}}\otimes_\beta}z)\\
&=&\tilde{\tilde{\mathfrak a}}\circ\Theta (\mathfrak a\underset{N}{_b*_\alpha}id)(1\underset{N}{_b\otimes_\alpha}z)
\end{eqnarray*}
from which we get the result. \end{proof}

%%%%propduality2
\subsection{Proposition}
\label{propduality2}
{\it Let $\gG=(N, M, \alpha, \beta, \Gamma, T, T', \nu)$ be a measured quantum groupo-id; let $A$ be a von Neumann algebra acting on a Hilbert space $\gH$, $(b, \mathfrak a)$ an action of $\gG$ on $A$, $\tilde{\mathfrak a}$ the dual action of $(\widehat{\gG})^c$ on the crossed product $A\rtimes_{\mathfrak a}\gG$, $\tilde{\tilde{\mathfrak a}}$ the bidual action of $({\gG})^{oc}$ on the double crossed-product $(A\rtimes_{\mathfrak a}\gG)\rtimes_{\tilde{\mathfrak a}}\widehat{\gG}^o$; then, with the notations of \ref{propduality1} : 
\newline
(i) $(1\underset{N}{_b\otimes_\alpha}\hat{\beta}, \underline{\mathfrak a})$ is an action of $\gG$ on $\underline{A}$; 
\newline
(ii) this action satisfies property (B);
\newline
(iii) $\underline{A}^{\underline{\mathfrak a}}=A\rtimes_\mathfrak a\gG$}

\begin{proof}
Thanks to \ref{propduality1}, we get that, for any $n\in N$ :
\begin{eqnarray*}
\underline{\mathfrak a}(1\underset{N}{_b\otimes_\alpha}\hat{\beta}(n))
&=&
(id\underset{N}{_b*_\alpha}\varsigma_{N^o})(1\underset{N^o}{_b\otimes_\alpha}(J_{\widehat{\Phi}}\underset{N}{_\beta\otimes_
\alpha} J_\Phi)\Gamma(\alpha(n))(J_{\widehat{\Phi}}\underset{N}{_\beta\otimes_
\alpha} J_\Phi))\\
&=&1\underset{N}{_b\otimes_\alpha}1\underset{N}{_{\hat{\beta}}\otimes_\alpha}\beta(n)
\end{eqnarray*}
On the other hand, we get, using \ref{propduality1}, for any $X\in\underline{A}$, that  :
\[(\Theta (\mathfrak a\underset{N}{_b*_\alpha}id)\underset{N}{_b*_\alpha}\mathcal I_{\gG}\underset{N}{_b*_\alpha}\mathcal I_{\gG})(\underline{\mathfrak a}(X)\underset{N}{_\beta*_\alpha}id)\underline{\mathfrak a}(X)\]
is equal to :
\begin{eqnarray*}
(\tilde{\tilde{\mathfrak a}}\circ\Theta (\mathfrak a\underset{N}{_b*_\alpha}id)\underset{N}{_b*_\alpha}
\mathcal I_{\gG})\underline{\mathfrak a}(X)
&=&
(\tilde{\tilde{\mathfrak a}}\underset{N}{_{\hat{\beta}}\otimes_{\hat{\alpha}}}id)\tilde{\tilde{\mathfrak a}}
\circ\Theta (\mathfrak a\underset{N}{_b*_\alpha}id)(X)\\
&=&
(id\underset{N}{_b\otimes_{\hat{\alpha}}}\Gamma^{oc})\tilde{\tilde{\mathfrak a}}
\circ\Theta (\mathfrak a\underset{N}{_b*_\alpha}id)(X)\\
&=&
(\Theta (\mathfrak a\underset{N}{_b*_\alpha}id)\underset{N}{_b*_\alpha}\Gamma^{oc}\circ\mathcal I_{\gG})\underline{\mathfrak a}(X)
\end{eqnarray*}
which is equal to $(\Theta (\mathfrak a\underset{N}{_b*_\alpha}id)\underset{N}{_b*_\alpha}\mathcal I_{\gG}\underset{N}{_b*_\alpha}\mathcal I_{\gG})(id\underset{N}{_\beta*_\alpha}\Gamma)\underline{\mathfrak a}(X)$, 
from which we get (i). On the other hand, we have :
\begin{eqnarray*}
(\Theta (\mathfrak a\underset{N}{_b*_\alpha}id)\underset{N}{_b*_\alpha}\mathcal I_{\gG})
(\underline{A}\rtimes_{\underline{\mathfrak a}}\gG)
&=&
(\Theta (\mathfrak a\underset{N}{_b*_\alpha}id)\underset{N}{_b*_\alpha}\mathcal I_{\gG})
(\underline{\mathfrak a}(\underline{A})\cup 1\underset{N}{_b\otimes_\alpha}1\underset{N}{_{\hat{\beta}}\otimes_\alpha}\widehat{M}')''\\
&=&
(\tilde{\tilde{\mathfrak a}}[(A\rtimes_\mathfrak a\gG)\rtimes_{\tilde{\mathfrak a}}\widehat{\gG}^o]\cup 1\underset{N}{_b\otimes_\alpha}1\underset{N}{_{\hat{\beta}}\otimes_{\hat{\alpha}}}\widehat{M})''\\
&=&
[(A\rtimes_\mathfrak a\gG)\rtimes_{\tilde{\mathfrak a}}\widehat{\gG}^o]\rtimes_{\tilde{\tilde{\mathfrak a}}}\gG^{oc}
\end{eqnarray*}
Using then \ref{dualaction}(iii), we get that :
\[(\Theta (\mathfrak a\underset{N}{_b*_\alpha}id)\underset{N}{_b*_\alpha}\mathcal I_{\gG})
((\underline{A}\rtimes_{\underline{\mathfrak a}}\gG)^{\tilde{\underline{\mathfrak a}}})
=
(\Theta (\mathfrak a\underset{N}{_b*_\alpha}id)\underset{N}{_b*_\alpha}\mathcal I_{\gG})
((\underline{A}\rtimes_{\underline{\mathfrak a}}\gG\cap\underline{A}\underset{N}{_b*_\alpha}M)\]
which, using \ref{propduality1} and \ref{dualaction}(iii) again, is equal to :
\[[(A\rtimes_\mathfrak a\gG)\rtimes_{\tilde{\mathfrak a}}\widehat{\gG}^o]\rtimes_{\tilde{\tilde{\mathfrak a}}}\gG^{oc}\cap[(A\rtimes_\mathfrak a\gG)\rtimes_{\tilde{\mathfrak a}}\widehat{\gG}^o]\underset{N^o}{_{\hat{\alpha}}*_{\hat{\beta}}}M'
=
([(A\rtimes_\mathfrak a\gG)\rtimes_{\tilde{\mathfrak a}}\widehat{\gG}^o]\rtimes_{\tilde{\tilde{\mathfrak a}}}\gG^{oc})^{\tilde{\tilde{\mathfrak a}}}\]
which, using \ref{corB}, is equal to $\tilde{\tilde{\mathfrak a}}[(A\rtimes_\mathfrak a\gG)\rtimes_{\tilde{\mathfrak a}}\widehat{\gG}^o]$, and, therefore, by \ref{propduality1} again, is equal to 
$(\Theta (\mathfrak a\underset{N}{_b*_\alpha}id)\underset{N}{_b*_\alpha}\mathcal I_{\gG})\underline{\mathfrak a}(\underline{A})$. From which we get that $(\underline{A}\rtimes_{\underline{\mathfrak a}}\gG)^{\tilde{\underline{\mathfrak a}}}=\underline{\mathfrak a}(\underline{A})$
which is (ii). 
\newline
Using \ref{propduality2}, we get that $X$ in $\underline{A}$ belongs to $\underline{A}^{\underline{\mathfrak a}}$ if and only if $\Theta (\mathfrak a\underset{N}{_b*_\alpha}id)(X)$ belongs to $[(A\rtimes_\mathfrak a\gG)\rtimes_{\tilde{\mathfrak a}}\widehat{\gG}^o]^{\tilde{\tilde{\mathfrak a}}}$, which, using  \ref{corB}, is equal to $\tilde{\mathfrak a}(A\rtimes_\mathfrak a\gG)$, or, using \ref{propduality2} again, is equal to $\Theta (\mathfrak a\underset{N}{_b*_\alpha}id)(A\rtimes_\mathfrak a\gG)$; from which we get (iii). \end{proof}

%%%%%propduality3
\subsection{Proposition}
\label{propduality3}
{\it Let $\gG=(N, M, \alpha, \beta, \Gamma, T, T', \nu)$ be a measured quantum groupo-id; let $A$ be a von Neumann algebra acting on a Hilbert space $\gH$, $(b, \mathfrak a)$ an action of $\gG$ on $A$, $\tilde{\mathfrak a}$ the dual action of $(\widehat{\gG})^c$ on the crossed product $A\rtimes_{\mathfrak a}\gG$, $\tilde{\tilde{\mathfrak a}}$ the bidual action of $({\gG})^{oc}$ on the double crossed-product $(A\rtimes_{\mathfrak a}\gG)\rtimes_{\tilde{\mathfrak a}}\widehat{\gG}^o$; let's use the notations of \ref{propduality1} and \ref{propduality2}. Then, we have :}
\[\underline{A}\rtimes_{\underline{\mathfrak a}} \gG=(id\underset{N}{_b*_\alpha}\varsigma_{N^o})[(1\underset{N}{_b\otimes_\alpha}W)(A\rtimes_\mathfrak a\gG\underset{N}{_\beta*_\alpha}\mathcal L(H_\Phi))(1\underset{N}{_b\otimes_\alpha}W)^*]\]

\begin{proof}
By definition of $\underline{A}$, we have :
\[\underline{\mathfrak a}(\underline{A})=(id\underset{N}{_b*_\alpha}\varsigma_{N^o})[(1\underset{N}{_b\otimes_\alpha}W)((\mathfrak a\underset{N}{_b*_\alpha}id)\mathfrak a (A)\cup 1\underset{N}{_b\otimes_\alpha}1\underset{N}{_\beta\otimes_\alpha}\alpha(N)')''(1\underset{N}{_b\otimes_\alpha}W)^*]\]
and $\underline{A}\rtimes_{\underline{\mathfrak a}}\gG=(\underline{\mathfrak a}(\underline{A})\cup 1\underset{N}{_b\otimes_\alpha}1\underset{N}{_{\hat{\beta}}\otimes_\alpha}\widehat{M}')''$. But, on the other hand, using \ref{thL4}(iii), we get that $W^*(\widehat{M}'\underset{N^o}{_\alpha\otimes_{\hat{\beta}}}1)W$ is equal to :
\begin{multline*}
(J_{\widehat{\Phi}}\underset{N^o}{_\alpha\otimes_{\hat{\beta}}} J_\Phi)W(J_{\widehat{\Phi}}\underset{N^o}{_\alpha\otimes_{\hat{\beta}}} J_\Phi)(\widehat{M}'\underset{N^o}{_\alpha\otimes_{\hat{\beta}}}1)(J_{\widehat{\Phi}}\underset{N}{_\beta\otimes_
\alpha} J_\Phi)W^*(J_{\widehat{\Phi}}\underset{N}{_\beta\otimes_
\alpha} J_\Phi)\\
=(J_{\widehat{\Phi}}\underset{N^o}{_\alpha\otimes_{\hat{\beta}}} J_\Phi)\widehat{\Gamma}^o(\widehat{M})(J_{\widehat{\Phi}}\underset{N}{_\beta\otimes_
\alpha} J_\Phi)
\end{multline*}
Therefore, we get that :
\[\widehat{M}'\underset{N^o}{_\alpha\otimes_{\hat{\beta}}}1=W(J_{\widehat{\Phi}}\underset{N^o}{_\alpha\otimes_{\hat{\beta}}} J_\Phi)\widehat{\Gamma}^o(\widehat{M})(J_{\widehat{\Phi}}\underset{N}{_\beta\otimes_
\alpha} J_\Phi)W^*\]
and :
\[1\underset{N}{_{\hat{\beta}}\otimes_\alpha}\widehat{M}'=\varsigma_{N^o}[W(J_{\widehat{\Phi}}\underset{N^o}{_\alpha\otimes_{\hat{\beta}}} J_\Phi)\widehat{\Gamma}^o(\widehat{M})(J_{\widehat{\Phi}}\underset{N}{_\beta\otimes_
\alpha} J_\Phi)W^*]\]
Moreover, we have :
\begin{multline*}
(J_{\widehat{\Phi}}\underset{N^o}{_\alpha\otimes_{\hat{\beta}}} J_\Phi)\widehat{\Gamma}^o(\widehat{M})(J_{\widehat{\Phi}}\underset{N}{_\beta\otimes_
\alpha} J_\Phi)\cup (1\underset{N}{_\beta\otimes_\alpha}\alpha(N)')''=\\
=(J_{\widehat{\Phi}}\underset{N^o}{_\alpha\otimes_{\hat{\beta}}} J_\Phi)[\widehat{\Gamma}^o(\widehat{M})\cup 1\underset{N^o}{_\alpha\otimes_{\hat{\beta}}}\hat{\beta}(N)']''((J_{\widehat{\Phi}}\underset{N}{_\beta\otimes_\alpha} J_\Phi)
\end{multline*}
which, by \ref{GammaA} applied to $\widehat{\gG}^o$, is equal to :
\[(J_{\widehat{\Phi}}\underset{N^o}{_\alpha\otimes_{\hat{\beta}}} J_\Phi)(\widehat{M}\underset{N^o}{_\alpha*_{\hat{\beta}}}\mathcal L(H_\Phi))(J_{\widehat{\Phi}}\underset{N}{_\beta\otimes_\alpha} J_\Phi)=\widehat{M}'\underset{N}{_\beta*_\alpha}\mathcal L(H_\Phi)\]
Therefore, we get that $\underline{A}\rtimes_{\underline{\mathfrak a}}\gG$ is equal to :
\[(id\underset{N}{_b*_\alpha}\varsigma_{N^o})[(1\underset{N}{_b\otimes_\alpha}W)((\mathfrak a\underset{N}{_b*_\alpha}id)\mathfrak a (A)\cup 1\underset{N}{_b\otimes_\alpha}\widehat{M}'\underset{N}{_\beta*_\alpha}\mathcal L(H_\Phi))''(1\underset{N}{_b\otimes_\alpha}W)^*]\]
which, using \ref{lemcrossed}, gives the result. \end{proof}

%%%%thduality1
\subsection{Theorem}
\label{thduality1}
{\it Let $\gG=(N, M, \alpha, \beta, \Gamma, T, T', \nu)$ be a measured quantum groupoid; let $A$ be a von Neumann algebra acting on a Hilbert space $\gH$, $(b, \mathfrak a)$ an action of $\gG$ on $A$, $\tilde{\mathfrak a}$ the dual action of $(\widehat{\gG})^c$ on the crossed product $A\rtimes_{\mathfrak a}\gG$. Then :
\newline
(i) $\mathfrak a$ satisfies property (B), i.e. we have :
\[(A\rtimes_{\mathfrak a}\gG)^{\tilde{\mathfrak a}}=\mathfrak a (A)\]
(ii) $\mathfrak a$ satisfies property (A), i.e. we have :}
\[(\mathfrak a(A)\cup 1\underset{N}{_b\otimes_\alpha}\alpha (N)')''=A\underset{N}{_b*_\alpha}\mathcal L(H_\Phi)\]

\begin{proof} Let's use the constructions and notations of \ref{propduality1}, \ref{propduality2} and \ref{propduality3}. Using \ref{dualaction}(iii), we have $(\underline{A}\rtimes_{\underline{\mathfrak a}}\gG)^{\tilde{\underline{\mathfrak a}}}=\underline{A}\rtimes_{\underline{\mathfrak a}}\gG\cap \mathcal L(\gH\underset{\nu}{_b\otimes_\alpha} H_\Phi)\underset{M}{_\beta*_\alpha}M$, and, therefore, using \ref{propduality3} and \ref{dualaction}(iii) again, we have :
\[(\underline{A}\rtimes_{\underline{\mathfrak a}}\gG)^{\tilde{\underline{\mathfrak a}}}=
(id\underset{N}{_b*_\alpha}\varsigma_{N^o})[(1\underset{N}{_b\otimes_\alpha}W)((A\rtimes_\mathfrak a\gG)^{\tilde{\mathfrak a}}\underset{N}{_\beta*_\alpha}\mathcal L(H_\Phi))(1\underset{N}{_b\otimes_\alpha}W)^*]\]
But, by \ref{propduality2}, we have :
\[(\underline{A}\rtimes_{\underline{\mathfrak a}}\gG)^{\tilde{\underline{\mathfrak a}}}=\underline{\mathfrak a}(\underline{A})
=(id\underset{N}{_b*_\alpha}\varsigma_{N^o})[(1\underset{N}{_b\otimes_\alpha}W)(\mathfrak a\underset{N}{_b*_\alpha}id)(\underline{A})\underset{N}{_\beta*_\alpha}\mathcal L(H_\Phi))(1\underset{N}{_b\otimes_\alpha}W)^*]\]
from which we deduce that 
$(A\rtimes_\mathfrak a\gG)^{\tilde{\mathfrak a}}\underset{N}{_\beta*_\alpha}\mathcal L(H_\Phi)=(\mathfrak a\underset{N}{_b*_\alpha}id)(\underline{A})$, 
from which we get that 
$(A\rtimes_\mathfrak a\gG)^{\tilde{\mathfrak a}}\underset{N}{_\beta*_\alpha}\mathcal L(H_\Phi)\subset \mathfrak a (A)\underset{N}{_\beta*_\alpha}\mathcal L(H_\Phi)$
and we deduce that $(A\rtimes_\mathfrak a\gG)^{\tilde{\mathfrak a}}\subset\mathfrak a (A)$, which, thanks to \ref{dualaction}(iii), gives (i). 
\newline
But we have now :
\begin{eqnarray*}
(\mathfrak a\underset{N}{_b*_\alpha}id)(A\underset{N}{_b*_\alpha}\mathcal L(H_\Phi))
&=&\mathfrak a(A)\underset{N}{_\beta*_\alpha}\mathcal L(H_\Phi)\\
&=&(A\rtimes_\mathfrak a\gG)^{\tilde{\mathfrak a}}\underset{N}{_\beta*_\alpha}\mathcal L(H_\Phi)\\
&=&(\mathfrak a\underset{N}{_b*_\alpha}id)(\underline{A})
\end{eqnarray*}
which gives (ii). \end{proof}

%%%thduality2
\subsection{Theorem}
\label{thduality2}
{\it Let $\gG=(N, M, \alpha, \beta, \Gamma, T, T', \nu)$ be a measured quantum groupoid; let $A$ be a von Neumann algebra acting on a Hilbert space $\gH$, $(b, \mathfrak a)$ an action of $\gG$ on $A$, $\tilde{\mathfrak a}$ the dual action of $(\widehat{\gG})^c$ on the crossed product $A\rtimes_{\mathfrak a}\gG$, $\tilde{\tilde{\mathfrak a}}$ the bidual action of $({\gG})^{oc}$ on the double crossed-product $(A\rtimes_{\mathfrak a}\gG)\rtimes_{\tilde{\mathfrak a}}\widehat{\gG}^o$; let us define, for any $X\in A\underset{N}{_b*_\alpha}\mathcal L(H_\Phi)$ :
\[\underline{\mathfrak a}(X)=(1\underset{N}{_b\otimes_\alpha}\sigma_{\nu^o} W\sigma_{\nu^o})(id\underset{N}{_b*_\alpha}\varsigma_N)(\mathfrak a\underset{N}{_b*_\alpha}id)(X)(1\underset{N}{_b\otimes_\alpha}\sigma_{\nu^o} W\sigma_{\nu^o})^*\]
Moreover, for any $Y\in\mathcal L(\gH\underset{\nu}{_b\otimes_\alpha}H_\Phi\underset{\nu}{_\beta\otimes_\alpha}H_\Phi)$, we define :
\[\Theta(Y)=(1\underset{\nu}{_b\otimes_\alpha}\sigma_\nu W^o\sigma_\nu)Y(1\underset{N}{_b\otimes_\alpha}\sigma_\nu W^o\sigma_\nu)^*\in\mathcal L(\gH\underset{\nu}{_b\otimes_\alpha}H_\Phi\underset{\nu^o}{_{\hat{\alpha}}\otimes_\beta}H_\Phi)\]
Then, for any $X\in A\underset{N}{_b*_\alpha}\mathcal L(H_\Phi)$, we have :
\[\Theta (\mathfrak a\underset{N}{_b*_\alpha}id)(A\underset{N}{_b*_\alpha}\mathcal L(H_\Phi))=(A\rtimes_{\mathfrak a}\gG)\rtimes_{\tilde{\mathfrak a}}\widehat{\gG}^o\]
\[(\Theta (\mathfrak a\underset{N}{_b*_\alpha}id)\underset{N}{_b*_\alpha}\mathcal I_{\gG})\underline{\mathfrak a}(X)=\tilde{\tilde{\mathfrak a}}\circ\Theta (\mathfrak a\underset{N}{_b*_\alpha}id)(X)\]
and, for any $Z\in A\rtimes_\mathfrak a\gG$, we get that $\Theta(\mathfrak a\underset{N}{_b*_\alpha}id)(Z)=\tilde{\mathfrak a}(Z)$. 
\newline
Then $(1\underset{N}{_b\otimes_\alpha}\hat{\beta}, \underline{\mathfrak a})$ is an action of $\gG$ on $A\underset{N}{_b*_\alpha}\mathcal L(H_\Phi)$, and :
\[(A\underset{N}{_b*_\alpha}\mathcal L(H_\Phi))^{\underline{\mathfrak a}}=A\rtimes_{\mathfrak a}\gG\]
Moreover, for any $X\in (A\underset{N}{_b*_\alpha}\mathcal L(H_\Phi))^+$, we have :}
\[T_{\tilde{\tilde{\mathfrak a}}}\Theta(\mathfrak a\underset{N}{_b*_\alpha}id)=\Theta(\mathfrak a\underset{N}{_b*_\alpha}id)T_{\underline{\mathfrak a}}(X)=\tilde{\mathfrak a}(T_{\underline{\mathfrak a}}(X))\]
\begin{proof}
This result is given by \ref{propduality1}, \ref{propduality2} and \ref{thduality1}(ii). \end{proof}

%%%%%thcom2
\subsection{Theorem}
\label{thcom2}
{\it Let $\gG=(N, M, \alpha, \beta, \Gamma, T, T', \nu)$ be a measured quantum groupoid, $_a\gH_b$ a $N-N$ bimodule, $A$ a von Neumann algebra such that $b(N)\subset A\subset a(N)'$, $V$ a corepresentation of $\gG$ on $_a\gH_b$ which implements an action $(b, \mathfrak a)$ of $\gG$ on $A$ and an action $(a, \mathfrak a')$ of $\gG^o$ on $A'$. Then :
\newline
(i) $A\rtimes_\mathfrak a\gG=V(\mathcal L(\gH)\underset{N^o}{_a*_\beta}\widehat{M}')V^*\cap A\underset{N}{_b*_\alpha}\mathcal L(H_\Phi)$
\newline
(ii) $A'\rtimes_{\mathfrak a'}\gG^o=V^*(A\rtimes_\mathfrak a\gG)'V$. }

\begin{proof}
Using \ref{thduality2}, we get that $X\in A\underset{N}{_b*_\alpha}\mathcal L(H_\Phi)$ belongs to $A\rtimes_\mathfrak a\gG$ if and only if $X\in A\underset{N}{_b*_\alpha}{\hat{\beta}}(N)'$, and $\underline{\mathfrak a}(X)=X\underset{N}{_{\hat{\beta}}\otimes_\alpha}1$, or, equivalently :
\[(\mathfrak a\underset{N}{_b*_\alpha}id)\mathfrak a(X)
=(1\underset{N}{_b\otimes_\alpha}W^*)(1\underset{N}{_b\otimes_\alpha}\sigma_\nu)(X\underset{N}{_{\hat{\beta}}\otimes_\alpha}1)(1\underset{N}{_b\otimes_\alpha}\sigma_\nu)(1\underset{N}{_b\otimes_\alpha}W)\]
which means (\ref{lemV}) that $X$ belongs to $A\underset{N}{_b*_\alpha}{\hat{\beta}}(N)'\cap SatV$, and, therefore, by \ref{lemV}(i), we get :
\[A\rtimes_\mathfrak a\gG=V(\mathcal L(\gH)\underset{N^o}{_a*_\beta}\widehat{M}')V^*\cap A\underset{N}{_b*_\alpha} \hat{\beta}(N)'\]
We have seen in \ref{lemV} that $V(\mathcal L(\gH)\underset{N^o}{_a*_\beta}\widehat{M}')V^*$ belongs to $a(N)'\underset{N}{_b*_\alpha}\hat{\beta}(N)'$, which gives (i). 
\newline
From (i), we get that :
\[(A\rtimes_\mathfrak a\gG)'=(V(1\underset{N}{_b\otimes_\alpha}\widehat{M})V^*\cup A'\underset{N}{_b\otimes_\alpha}1)''\]
and :
\[V^*(A\rtimes_\mathfrak a\gG)'V=[V(A'\underset{N}{_b\otimes_\alpha}1)V^*\cup 1\underset{N}{_b\otimes_\alpha}\widehat{M}]''\]
which gives (ii). \end{proof}

%%%%%thduality4
\subsection{Theorem}
\label{thduality4}
{\it Let $\gG=(N, M, \alpha, \beta, \Gamma, T, T', \nu)$ be a measured quantum groupoid; let $A$ be a von Neumann algebra acting on a Hilbert space $\gH$, $(b, \mathfrak a)$ an action of $\gG$ on $A$; then, $\mathfrak a$ is saturated, i.e. we have :}
\[\mathfrak a(A)=\{X\in A\underset{N}{_b*_\alpha}M, (\mathfrak a \underset{N}{_b*_\alpha}id)(X)=(id\underset{N}{_b*_\alpha}\Gamma)(X)\}\]

\begin{proof}
Let us suppose that there exists a corepresentation $V$ on $\gH$ which implements $\mathfrak a$; using \ref{propsat} and \ref{thduality1}(i), we get that $\mathfrak a$ is saturated (\ref{defsat}), which means  that we have the result. Applying this fact to \ref{thimple}, we get that $(\mathfrak a\underset{N}{_b*_\alpha}id)\mathfrak a(A)$ is equal to :
\[\{X\in \mathfrak a(A)\underset{N}{_b*_\alpha}M, (\mathfrak a \underset{N}{_b*_\alpha}id\underset{N}{_\beta*_\alpha}id)(X)=(id\underset{N}{_b*_\alpha}id\underset{N}{_\beta*_\alpha}\Gamma)(X)\}\]
which is the image under $(\mathfrak a\underset{N}{_b*_\alpha}id)$ of :
\[\{Y\in A\underset{N}{_b*_\alpha}M; (\mathfrak a \underset{N}{_b*_\alpha}id\underset{N}{_\beta*_\alpha}id)(\mathfrak a\underset{N}{_b*_\alpha}id)(Y)=(id\underset{N}{_b*_\alpha}id\underset{N}{_\beta*_\alpha}\Gamma)(\mathfrak a\underset{N}{_b*_\alpha}id)(Y)\}\]
As, $(id\underset{N}{_b*_\alpha}id\underset{N}{_\beta*_\alpha}\Gamma)(\mathfrak a\underset{N}{_b*_\alpha}id)(Y)=(\mathfrak a \underset{N}{_b*_\alpha}id\underset{N}{_\beta*_\alpha}id)(id\underset{N}{_b*_\alpha}\Gamma)(Y)$, we get the result.  \end{proof}

%%thdense
\subsection{Theorem}
\label{thdense}
{\it Let $\gG=(N, M, \alpha, \beta, \Gamma, T, T', \nu)$ be a measured quantum groupoid; let $A$ be a von Neumann algebra, $(b, \mathfrak a)$ an action of $\gG$ on $A$ and $A\rtimes_\mathfrak a\gG$ its crossed product. Then, the linear subspace generated by all elements of the form $(1\underset{N}{_b\otimes_\alpha}a)\mathfrak a(x)$, with $a\in \widehat{M}'$ and $x\in A$, is dense in  $A\rtimes_\mathfrak a\gG$. }
\begin{proof}
Using \ref{thduality1}(i), this is a straightforward corollary of \ref{propo1}(iii). \end{proof}
%%%remark
\subsection{Remark}
\label{remark}
These theorems (\ref{thduality1}, \ref{thduality2}, \ref{thduality4}) are the generalization, up to the measured quantum groupoid framework, of biduality theorems ([V2] 2.6, 2.7) obtained for actions of locally compact quantum groups, and we followed the same strategy as [V2] (which is, in fact, the strategy of [ES1]). These theorems generalizes as well the duality theorems for groupoid actions ([Y3], 6.5) and integrable groupoid coactions ([Y3], 7.8). 

%%%%characterization
\section{Characterization of crossed-products}
\label{char}
On this chapter, inspired by (\cite{ES1}, V), a measured quantum groupoid $\gG$ be given, we characterize crossed-products among the von Neumann algebras on which there exists an action of $\widehat{\gG}^c$ (\ref{thchar}). A corollary (\ref{corcross}) will be used in chapter \ref{dualw}. 

%%%%notchar
\subsection{Notations}
\label{notchar}
In this chapter, we shall consider a measured quantum groupoid $\gG$, its dual $\widehat{\gG}$, and a von Neumann algebra $B$ on a Hilbert space $\gH$. We suppose that there exist :
\newline
(i) a normal faithful morphism $\mu$ from $\widehat{M}'$ into $B$;
\newline
(ii) a normal faithful morphism $\mathfrak b$ from $B$ into $B\underset{N^o}{_{\mu\circ\hat{\alpha}}*_\beta}\widehat{M}'$, such that $(\mu\circ\hat{\alpha}, \mathfrak b)$ is an action of $\widehat{\gG}^c$ on $B$; 
\newline
(iii) we have, moreover : $\mathfrak b\circ\mu=(\mu\underset{N^o}{_{\hat{\alpha}}*_\beta}id)\widehat{\Gamma}^c$. 
\newline
This last formula implies that, for any $n\in N$, $\mathfrak b\circ\mu (\beta(n))$ belongs to $(\mu\circ\hat{\alpha}(N))'$ and :
\[\mathfrak b\circ\mu(\beta(n))=(\mu\underset{N^o}{_{\hat{\alpha}}*_\beta}id)\widehat{\Gamma}^c(\beta(n))=
(\mu\underset{N^o}{_{\hat{\alpha}}*_\beta}id)(\beta(n)\underset{\nu^o}{_{\hat{\alpha}}\otimes_\beta}1)=
\mu(\beta(n))\underset{\nu^o}{_{\mu\circ\hat{\alpha}}\otimes_\beta}1\]
and, therefore, by definition, $\mu\circ\beta (N)\subset B^{\mathfrak a}$. 

%%%%Y
\subsection{Lemma}
\label{Y}
{\it Let's take the hypothesis and notations of \ref{notchar}. Then 
\newline
(i) there exists a corepresentation $Y$ of $\gG$ on the bimodule $_{\mu\circ\hat{\alpha}}\gH_{\mu\circ \beta}$. such that, for any $\xi\in D(_\alpha H_\Phi, \nu)$, $\eta\in D((H_\Phi)_\beta, \nu^o)$, we have :
\[(id*\omega_{\xi, \eta})(Y)=\mu[(id*\omega_{\xi, \eta})((\sigma_\nu W^o\sigma_\nu)^*)]\]
(ii) for any $b'\in B'$, we have $Y(b'\underset{N^o}{_{\mu\circ\hat{\alpha}}\otimes_\beta}1)=(b'\underset{N^o}{_{\mu\circ\beta}\otimes_\alpha}1)Y$; 
\newline
(iii) for any $b\in B$, $Y\mathfrak b(x)Y^*$ belongs to $B^\mathfrak b\underset{N}{_{\mu\circ\beta}*_\alpha}\mathcal L(H_\Phi)$; moreover, if $x\in B^\mathfrak b$, then $Y(x\underset{N}{_{\mu\circ\beta}\otimes_\alpha}1)Y^*$ belongs to $B^\mathfrak b\underset{N}{_{\mu\circ\beta}*_\alpha}M$; 
\newline
(iv) the coreprentation $Y$ of $\gG$ implements on $B^\mathfrak b$ an action of $\gG$, we shall denote $\mathfrak d$. More precisely, $(\mu\circ\beta, \mathfrak d)$ is an action of $\gG$ on $B^\mathfrak b$. }

\begin{proof}
By \ref{thL3}(v), $H_\Phi$ can be considered as the standard Hilbert space of $\widehat{M}'$, we can identify with $\widehat{M}^o$. Therefore, the Hilbert space $\gH$ is isomorphic (\ref{rel}) to the relative tensor product $\gH\underset{\widehat{\Phi}}{_\mu\otimes_{\hat{\pi}}}(H_\Phi)$, where $\hat{\pi}$ denotes here the canonical representation of $\widehat{M}$ (resp. $\widehat{M}'$) on $H_\Phi$, and this isomorphism is an isomorphism of right $\widehat{M}$-modules, which means that, modulo the identification of these Hilbert spaces, we have, for all $x\in\widehat{M}'$, $\mu(x)=1\underset{\widehat{M}}{_\mu\otimes_{\hat{\pi}}}x$. 
\newline
Let us put $Y=1_\gH\underset{\widehat{M}}{_\mu\otimes _{\hat{\pi}}}(\sigma_\nu W^o\sigma_\nu)^*$.  So, $Y$ is a unitary from $\gH\underset{\nu^o}{_{\mu\circ\hat{\alpha}}\otimes_\beta}H_\Phi$ onto $\gH\underset{\nu}{_{\mu\circ\beta}\otimes_\alpha}H_\Phi$. As $(\sigma_\nu W^o\sigma_\nu)^*$ is a corepresentation of $\gG$ on the $N-N$ bimodule $_{\hat{\alpha}}(H_\Phi)_\beta$ (\ref{Example}), we get easily that $Y$ is a corepresentation of $\gG$ on the $N-N$ bimodule $_{\mu\circ\hat{\alpha}}\gH_{\mu\circ \beta}$. 
\newline
Moreover, we have then :
\begin{eqnarray*}
(id*\omega_{\xi, \eta})(Y)&=&(id\underset{\widehat{M}}{_\mu*_{\hat{\pi}}}id*\omega_{\xi, \eta})[1\underset{\widehat{M}}{_\mu\otimes_{\hat{\pi}}}(\sigma_\nu W^o\sigma_\nu)^*]\\
&=&1\underset{\widehat{M}}{_\mu\otimes_{\hat{\pi}}}(id*\omega_{\xi, \eta})((\sigma_\nu W^o\sigma_\nu)^*)\\
&=&\mu[(id*\omega_{\xi, \eta})((\sigma_\nu W^o\sigma_\nu)^*)]
\end{eqnarray*}
which gives (i). Then (ii) is  a straightforward corollary of (i). From (ii), we get easily that $Y\mathfrak b(b)Y^*$ belongs to $B\underset{N}{_{\mu\circ\beta}*_\alpha}\mathcal L(H_\Phi)$, and, therefore, that 
$(id\underset{N}{_{\mu\circ\beta}*_\alpha}\omega_\xi)(Y\mathfrak b(b)Y^*)$ belongs to $B$. Let $(\eta_i)_{i\in I}$ be an orthogonal $(\beta, \nu^o)$-basis of $H_\Phi$. We get that the element $(id\underset{N}{_{\mu\circ\beta}*_\alpha}\omega_\xi)(Y\mathfrak b(b)Y^*)$ is equal to the sum :
\begin{multline*}
\sum_{i,j}(id\underset{N}{_{\mu\circ\beta}*_\alpha}\omega_\xi)(Y\theta^{\beta, \nu^o}(\eta_i, \eta_i)\mathfrak b(b)\theta^{\beta, \nu^o}(\eta_j, \eta_j)Y^*)=\\
=\sum_{i,j}(id*\omega_{\xi, \eta_i})(Y)(id\underset{N^o}{_{\mu\circ\hat{\alpha}}*_\beta}\omega_{\eta_i, \eta_j})\mathfrak b(b)[(id*\omega_{\xi, \eta_j})(Y)]^*
\end{multline*}
And, as we know that  $(\sigma_\nu W^o\sigma_\nu)^*=\widehat{W}^c$, we get that :
\[\mathfrak b[(id\underset{N}{_{\mu\circ\beta}*_\alpha}\omega_\xi)(Y\mathfrak b(b)Y^*)]\]
 is equal to the sum over $(i,j)$ of all following elements :
\[\mathfrak b\circ\mu [(id*\omega_{\xi, \eta_i})(\widehat{W}^c)(id\underset{N^o}{_{\mu\circ\hat{\alpha}}*_\beta}id\underset{N^o}{_{\hat{\alpha}}*_\beta}\omega_{\eta_i, \eta_j})(id\underset{N^o}{_{\mu\circ\hat{\alpha}}*_\beta}\widehat{\Gamma}^c)\mathfrak b(b)]\mathfrak b\circ\mu (id*\omega_{\xi, \eta_j})(\widehat{W}^c)^*\]
Using the pentagonal relation for  $\widehat{W}^c$, we get that :
\[\mathfrak b\circ\mu[(id*\omega_{\xi, \eta_i})(\widehat{W}^c)]
=(\mu\underset{N^o}{_{\hat{\alpha}}*_\beta}id)\widehat{\Gamma}^c[(id*\omega_{\xi, \eta_i})(\widehat{W}^c)]\]
is equal to $(\mu\underset{N^o}{_{\hat{\alpha}}*_\beta}id)(id*id*\omega_{\xi, \eta_i})[\sigma_{\beta, \hat{\alpha}}^{2,3}(\widehat{W}^c\underset{N^o}{_\alpha\otimes_\beta}1)(1\underset{N^o}{_{\hat{\alpha}}\otimes_\beta}\sigma_\nu)(1\underset{N^o}{_{\hat{\alpha}}\otimes_\beta}\widehat{W}^c)]$. 
And, as $(id\underset{N^o}{_{\mu\circ\hat{\alpha}}*_\beta}\widehat{\Gamma}^c)\mathfrak b(b)$ is equal to $(1_\gH\underset{\widehat{M}}{_\mu\otimes _{\hat{\pi}}}\widehat{W}^c)^*(1\underset{N^o}{_{\mu\circ\hat{\alpha}}\otimes_\beta}\sigma_{\nu^o})(\mathfrak b(b)\underset{N^o}{_{\hat{\alpha}}\otimes_\beta}1)(1\underset{N^o}{_{\mu\circ\hat{\alpha}}\otimes_\beta}\sigma_{\nu^o})^*(1_\gH\underset{\widehat{M}}{_\mu\otimes _{\hat{\pi}}}\widehat{W}^c)$, 
we obtain finally that this sum is equal to :
\[(\mu\underset{N^o}{_{\hat{\alpha}}*_\beta}id)(id*id*\omega_{\xi})[\sigma_{\beta, \hat{\alpha}}^{2,3}(\widehat{W}^c\mathfrak b(b)(\widehat{W}^c)^*\underset{N^o}{_\alpha\otimes_\beta}1)(\sigma_{\beta, \hat{\alpha}}^{2,3})^*]\]
which can be written :
\[(id\underset{N}{_{\mu\circ\beta}*_\alpha}\omega_\xi)(\sigma_{\beta, \hat{\alpha}}^{2,3}(Y\mathfrak b(b)Y^*\underset{N^o}{_{\hat{\alpha}}\otimes_\beta}1)(\sigma_{\beta, \hat{\alpha}}^{2,3})^*)=
(id\underset{N}{_{\mu\circ\beta}*_\alpha}\omega_\xi)(Y\mathfrak b(b)Y^*)\underset{N^o}{_{\hat{\alpha}}\otimes_\beta}1\]
and, therefore, we get that $Y\mathfrak b(b)Y^*$ belongs to $B^\mathfrak b\underset{N}{_{\mu\circ\beta}*_\alpha}\mathcal L(H_\Phi)$. 
\newline
Let $x\in B^\mathfrak b\subset \mu\circ\hat{\alpha}(N)'$.  We have obtained that $Y(x\underset{N^o}{_{\mu\circ\hat{\alpha}}\otimes_\beta}1)Y^*$ belongs to $B^\mathfrak b\underset{N}{_{\mu\circ\beta}*_\alpha}\mathcal L(H_\Phi)$; as it belongs clearly to $\mathcal L(\gH)\underset{N}{_{\mu\circ\beta}*_\alpha}M$, we get that $Y(x\underset{N^o}{_{\mu\circ\hat{\alpha}}\otimes_\beta}1)Y^*$ belongs to $B^\mathfrak b\underset{N}{_{\mu\circ\beta}*_\alpha}M$, which finishes the proof of (iii). 
\newline
If $n\in N$, we have :
\begin{eqnarray*}
\mathfrak d(\mu\circ\beta(n))&=&(1_\gH\underset{\widehat{M}}{_\mu\otimes _{\hat{\pi}}}\widehat{W}^c)(
1_\gH\underset{\widehat{M}}{_\mu\otimes _{\hat{\pi}}}\beta(n)\underset{N^o}{_{\hat{\alpha}}\otimes_\beta}1)(1_\gH\underset{\widehat{M}}{_\mu\otimes _{\hat{\pi}}}\widehat{W}^c)^*\\
&=&1_\gH\underset{\widehat{M}}{_\mu\otimes _{\hat{\pi}}}\Gamma(\beta(n))\\
&=&1_\gH\underset{\widehat{M}}{_\mu\otimes _{\hat{\pi}}}1\underset{N}{_\beta\otimes_\alpha}\beta(n)\\
&=&1\underset{N}{_{\mu\circ\beta}\otimes_\alpha}\beta(n)
\end{eqnarray*}
which gives (iv). \end{proof}

%%%thchar
\subsection{Theorem}
\label{thchar}
{\it Let $\gG$ be a measured quantum groupoid, $\widehat{\gG}$ its dual measured quantum groupoid, $B$ a von Neumann algebra on a Hilbert space $\gH$.  Then, are equivalent : 
\newline
(i) there exists a von Neumann algebra $A$, an action $(b, \mathfrak a)$ of $\gG$ on $A$, and an isomorphism $\mathcal I$ between the von Neumann algebras $B$ and $A\rtimes_\mathfrak a\gG$. 
\newline
(ii) there exists a normal faithful morphism $\mu$ from $\widehat{M}'$ into $B$, a normal faithful morphism $\mathfrak b$ from $B$ into $B\underset{N^o}{_{\mu\circ\hat{\alpha}}*_\beta}\widehat{M}'$ such that $(\mu\circ\hat{\alpha}, \mathfrak b)$ is an action of $\widehat{\gG}^c$ on $B$ which satisfies
 $\mathfrak b\circ\mu=(\mu\underset{N^o}{_{\hat{\alpha}}*_\beta}id)\widehat{\Gamma}^c$. 
\newline
Moreover, we have then :
\[\mathfrak b\circ \mathcal I=(\mathcal I\underset{N^o}{_{\hat{\alpha}}*_\beta}id)\mathfrak b\]
\[\mathcal I\circ T_\mathfrak b=T_{\tilde{\mathfrak a}}\circ \mathcal I\]
\[\mathcal I(B^\mathfrak b)=\mathfrak a(A)\]
\[B=(B^\mathfrak b\cup\mu(\widehat{M}'))''\]
and, for any $x\in\widehat{M}'$, we have :}
\[\mathcal I\circ\mu(x)=1\underset{N}{_b\otimes_\alpha}x\]

\begin{proof}
Let us suppose (i); then, if we define $\mathfrak b$ by $\mathfrak b=(\mathcal I\underset{N^o}{_{\hat{\alpha}}*_\beta}id)^{-1}\circ
\tilde{\mathfrak a}\circ \mathcal I$, and, for $x\in \widehat{M}'$, if we put $\mu(x)=\mathcal I^{-1}(1\underset{N}{_b\otimes_\alpha}x)$, then, we have (ii). 
\newline
Let us suppose (ii); we can use all notations and results of \ref{Y}. Let us consider the action $(\mu\circ\beta, \mathfrak d)$ of $\gG$ on $B^\mathfrak b$; we have :
\begin{eqnarray*}
B^\mathfrak b\rtimes_\mathfrak d \gG
&=&(\mathfrak d (B^\mathfrak b)\cup 1\underset{N}{_{\mu\circ\beta}\otimes_\alpha}\widehat{M}')''\\
&=&(Y(B^\mathfrak b\underset{N^o}{_{\mu\circ\hat{\alpha}}\otimes_\beta}1)Y^*\cup1\underset{N}{_{\mu\circ\beta}\otimes_\alpha}\widehat{M}')''
\end{eqnarray*}
But, for any $x\in \widehat{M}'$, we have :
\begin{eqnarray*}
Y^*(1\underset{N}{_{\mu\circ\beta}\otimes_\alpha}x)Y
&=&
(1_\gH\underset{\widehat{M}}{_\mu\otimes_{\hat{\pi}}}\widehat{W}^{c*})(1_\gH\underset{\widehat{M}}{_\mu\otimes_{\hat{\pi}}}1\underset{N}{_{\mu\circ\beta}\otimes_\alpha}x)(1_\gH\underset{\widehat{M}}{_\mu\otimes_{\hat{\pi}}}\widehat{W}^{c})\\
&=&1_\gH\underset{\widehat{M}}{_\mu\otimes_{\hat{\pi}}}\widehat{\Gamma}^c(x)\\
&=&(\mu\underset{N^o}{_{\hat{\alpha}}*_\beta} id)\widehat{\Gamma}^c(x)\\
&=&\mathfrak b\circ\mu (x)
\end{eqnarray*}
and, therefore, we have :
\begin{eqnarray*}
B^\mathfrak b\rtimes_\mathfrak d \gG
&=&(Y(B^\mathfrak b\underset{N^o}{_{\mu\circ\hat{\alpha}}\otimes_\beta}1)Y^*\cup1\underset{N}{_{\mu\circ\beta}\otimes_\alpha}\widehat{M}')''\\
&=&(Y\mathfrak b(B^\mathfrak b)Y^*\cup Y\mathfrak b(\mu(\widehat{M}'))Y^*)''\\
&=& Y\mathfrak b((B^\mathfrak b\cup \mu(\widehat{M}')'')Y^*\\
&\subset& Y\mathfrak b(B)Y^*
\end{eqnarray*}
For $y\in B$, let us define $\mathcal I (y)=Y\mathfrak b(y)Y^*$; 
so, we get $B^\mathfrak b\rtimes_\mathfrak d \gG\subset \mathcal I(B)$. 
\newline
Using \ref{Y}(iii), we get that $\mathcal I(B)\subset B^\mathfrak b\underset{N}{_{\mu\circ\beta}*_\alpha}\mathcal L (H_\Phi)$, and, therefore :
\[\mathcal I(B)\subset B^\mathfrak b\underset{N}{_{\mu\circ\beta}*_\alpha}\mathcal L (H_\Phi)\cap Y(\mathcal L(\gH)\underset{N^o}{_{\mu\circ\hat{\alpha}}*_\beta}\widehat{M}')Y^*\]
which, thanks to \ref{thcom2}(i), is equal to $B^\mathfrak b\rtimes_\mathfrak d\gG$.
\newline
So, we have proved that $\mathcal I(B)=B^\mathfrak b\rtimes_\mathfrak d\gG$, and we get (i), with $A=B^\mathfrak b$ and $(b, \mathfrak a)=(\mu\circ\beta, \mathfrak d)$. 
\newline
We have already obtained that  $\tilde{\mathfrak a}\circ \mathcal I=(\mathcal I\underset{N^o}{_{\hat{\alpha}}*_\beta}id)\mathfrak b$, and, for any $x\in\widehat{M}'$, we have 
$\mathcal I\circ\mu(x)=1\underset{N}{_b\otimes_\alpha}x$
from which we get, from the definition of the crossed-product, that $B=(B^\mathfrak b\cup\mu(\widehat{M}')''$. \end{proof}

%%%%corchar
\subsection{Corollary}
\label{corchar}
{\it Let's take the hypothesis and notations of \ref{notchar} and \ref{Y};  then, there exists an isomorphism $\mathcal J$ from $B^\mathfrak b\underset{N}{_{\mu\circ\beta}*_\alpha}\mathcal L(H_\Phi)$ onto $B\rtimes_\mathfrak b\widehat{\gG}^c$ such that, for all $x\in B^\mathfrak b\underset{N}{_{\mu\circ\beta}*_\alpha}\mathcal L(H_\Phi)$, we have :
\[\tilde{\mathfrak b}\mathcal J(x)=(\mathcal J\underset{N}{_{\mu\circ\beta}*_\alpha}\mathcal I_\gG)\mathfrak e(x)\]
where 
\[\mathfrak e(x)=
(1\underset{N}{_{\mu\circ\beta}\otimes_\alpha}\sigma_{\nu^o}W\sigma_{\nu^o})(id\underset{N}{_{\mu\circ\beta}*_\alpha}\varsigma_N)(\mathfrak d\underset{N}{_{\mu\circ\beta}*_\alpha}id)(x)(1\underset{N}{_{\mu\circ\beta}\otimes_\alpha}\sigma_{\nu^o}W\sigma_{\nu^o})^*\]
and $\mathcal I_\gG$ had been defined in \ref{thL5}(vi). }

\begin{proof}
Thanks to \ref{thchar}, there exist an action $(\mu\circ\beta, \mathfrak d)$ of $\gG$ on $B^\mathfrak b$, and an isomorphism $\mathcal I$ from $B$ onto $B^\mathfrak b\rtimes_\mathfrak d\gG$ such that $\tilde{\mathfrak d}\circ\mathcal I=(\mathcal I\underset{N^o}{_{\hat{\alpha}}*_\beta}id)\mathfrak b$. Therefore, we have, using \ref{thduality2}, there exists an isomorphism $\Theta$ such that :
\begin{eqnarray*}
(\mathcal I\underset{N^o}{_{\hat{\alpha}}*_\beta}id)(B\rtimes_\mathfrak b\widehat{\gG}^c)
&=&(B^\mathfrak b\rtimes_\mathfrak d\gG)\rtimes_{\tilde{\mathfrak d}}\widehat{\gG}^c\\
&=&\Theta(\mathfrak d\underset{N}{_{\mu\circ\beta}*_\alpha}id)(B^\mathfrak b\underset{N}{_{\mu\circ\beta}*_\alpha}\mathcal L(H_\Phi))
\end{eqnarray*}
Let's use again \ref{thduality2}; for all $x\in B^\mathfrak b\underset{N}{_{\mu\circ\beta}*_\alpha}\mathcal L(H_\Phi)$, let's define :
\[\mathfrak e(x)=
(1\underset{N}{_{\mu\circ\beta}\otimes_\alpha}\sigma_{\nu^o}W\sigma_{\nu^o})(id\underset{N}{_{\mu\circ\beta}*_\alpha}\varsigma_N)(\mathfrak d\underset{N}{_{\mu\circ\beta}*_\alpha}id)(x)(1\underset{N}{_{\mu\circ\beta}\otimes_\alpha}\sigma_{\nu^o}W\sigma_{\nu^o})^*\]
We have then :
\[\tilde{\tilde{\mathfrak d}}\Theta(\mathfrak d\underset{N}{_{\mu\circ\beta}*_\alpha}id)(x)=\\
(\Theta(\mathfrak d\underset{N}{_{\mu\circ\beta}*_\alpha}id)\underset{N}{_{\mu\circ\beta}*_\alpha}\mathcal I_\gG)\mathfrak e(x)\]
Let us define $\mathcal J=(\mathcal I\underset{N^o}{_{\hat{\alpha}}*_\beta}id)^{-1}\circ\Theta(\mathfrak d\underset{N}{_{\mu\circ\beta}*_\alpha}id)$ which is an isomorphism from $B^\mathfrak b\underset{N}{_{\mu\circ\beta}*_\alpha}\mathcal L(H_\Phi)$ onto $B\rtimes_\mathfrak b\widehat{\gG}^c$; we have then :
\[(\mathcal J\underset{N}{_{\mu\circ\beta}*_\alpha}\mathcal I_\gG)\mathfrak e(x)=\tilde{\tilde{\mathfrak d}}\Theta(\mathfrak d\underset{N}{_{\mu\circ\beta}*_\alpha}id)(x)=\tilde{\mathfrak b}\mathcal J(x)\]
which finishes the proof. \end{proof}

%%%corcross
\subsection{Corollary}
\label{corcross}
{\it Let $\gG$ be a measured quantum groupoid, and $(b, \mathfrak a)$ an action of $\gG$ on a von Neumann algebra $A$; let $A_0$ a von Neumann subalgebra of $A^\mathfrak a$, and let us write $D=A\cap A'_0$, and $\mathfrak d=\mathfrak a_{|D}$; then :
\newline
(i) $(b, \mathfrak d)$ is an action of $\gG$ on $D$; 
\newline
(ii) $D\rtimes_\mathfrak d\gG=A\rtimes_\mathfrak a\gG\cap A'_0\underset{N}{_b*_\alpha}\mathcal L(H_\Phi)$. }

\begin{proof}
Let us define $B=A\rtimes_\mathfrak a\gG\cap A'_0\underset{N}{_b*_\alpha}\mathcal L(H_\Phi)$
For $x\in \widehat{M}'$, let us put $\mu(x)=1\underset{N}{_b\otimes_\alpha}x$, which belongs to $B$, and let us write $\mathfrak b=\tilde{\mathfrak a}_{|B}$; for all $b\in B$, and $a_0\in A_0$, we have :
\begin{eqnarray*}
\mathfrak b(b)(a_0\underset{N}{_b\otimes_\alpha}1\underset{N^o}{_{\hat{\alpha}}\otimes_\beta}1)
&=&\tilde{\mathfrak a}(b)\tilde{\mathfrak a}(\mathfrak a (a_0)\underset{N^o}{_{\hat{\alpha}}\otimes_\beta}1)\\
&=&\tilde{\mathfrak a}(\mathfrak a (a_0)\underset{N^o}{_{\hat{\alpha}}\otimes_\beta}1)\tilde{\mathfrak a}(b)\\
&=&(a_0\underset{N}{_b\otimes_\alpha}1\underset{N^o}{_{\hat{\alpha}}\otimes_\beta}1)\mathfrak b(b)
\end{eqnarray*}
from which we get that $\mathfrak b(b)$ belongs to $A'_0\underset{N}{_b*_\alpha}\mathcal L(H_\Phi)\underset{N^o}{_{\hat{\alpha}}*_\beta}\mathcal L(H_\Phi)$, and, therefore, to $B\underset{N^o}{_{\hat{\alpha}}*_\beta}M'$. So, $(\mu\circ\hat{\alpha}, \mathfrak b)$ is an action of $\widehat{\gG}^c$ on $B$, and we can apply \ref{thchar}, from which we get that $B$ is the crossed product of $B^\mathfrak b$ by an action of $\gG$. 
\newline
But we have $B^\mathfrak b=(A\rtimes_\mathfrak a\gG)^{\tilde{\mathfrak a}}\cap A'_0\underset{N}{_b*_\alpha}\mathcal L(H_\Phi)$, which is equal to $\mathfrak a(A)\cap A'_0\underset{N}{_b*_\alpha}\mathcal L(H_\Phi)$ (\ref{thduality1}(i)), and, therefore, to $\mathfrak a(D)$. We then easily finish the proof. \end{proof}

%%%%%%dualw
\section{Dual weight; bidual weight; depth 2 inclusion associated to an action}
\label{dualw}
Be given an action $(b, \mathfrak a)$ of a measured quantum groupoid $\gG$ on a von Neumann algebra $A$ and a normal semi-finite faithful weight $\psi$ on $A$, we define in \ref{defdualw1} a dual weight $\tilde{\psi}$ on the crossed-product $A\rtimes_\mathfrak a\gG$, using property (B) proved in \ref{thduality1}(i) and the weight constructed in chapter \ref{aux}. We obtain a characterization of these weights (\ref{cardualw}), and the GNS construction associated to this weight (\ref{GNSdual}). Moreover, we study (\ref{Upsi1}) then the unitary $U_\psi$ constructed in \ref{propo2}.  Using then the isomorphism obtained in chapter \ref{biduality} between the double crossed product and $A\underset{N}{_b*_\alpha}\mathcal L(H_\Phi)$, we obtain a characterization of the bidual weight (\ref{spatialder}), which will allow us to construct, using \ref{Upsi1}, Jones' tower associated to the inclusion $\mathfrak a(A)\subset A\rtimes_\mathfrak a\gG$ (\ref{thstandard} and \ref{thD2}(i)). We prove then that this inclusion is depth 2 (\ref{thD2}(iv)) and that the operator-valued weight $T_{\tilde{\mathfrak a}}$ of this inclusion is regular in the sense of \ref{basic}(\ref{threg}).

%%%%defdualw
\subsection{Definition}
\label{defdualw1}
Let $(b, \mathfrak a)$ be an action of a measured quantum groupoid $\gG$ on a von Neumann algebra $A$, and let $\psi$ be a normal semi-finite faithful weight on $A$; using \ref{Propdual2}, we get that the dual action $\tilde{\mathfrak a}$ of $\widehat{\gG}^c$ on the crossed product $A\rtimes_\mathfrak a\gG$ is integrable, and, therefore, by definition, the operator-valued weight $T_{\tilde{\mathfrak a}}$ from $A\rtimes_\mathfrak a\gG$ on $(A\rtimes_\mathfrak a\gG)^{\tilde{\mathfrak a}}$ is semi-finite. On the other hand, using \ref{thduality1}(i), we know that $\mathfrak a(A)=(A\rtimes_\mathfrak a\gG)^{\tilde{\mathfrak a}}$. So, the formula $\tilde{\psi}=\psi\circ\mathfrak a^{-1}\circ T_{\tilde{\mathfrak a}}$ defines a normal semi-finite faithful weight on $A\rtimes_\mathfrak a\gG$, we shall call the dual weight of $\psi$. 

%%%dualwg
\subsection{Example}
\label{dualwg}
Let $\mathcal G$ be a measured groupoid, and $\mathfrak a$ an action of $\mathcal G$ on a von Neumann algebra, as defined in \cite{Y1}, \cite{Y2}, \cite{Y3}. We have seen in \ref{gcross} that the crossed-product defined by Yamanouchi is the same one as ours, and in \ref{dualg} that the dual coaction of $\mathcal G$ introduced by Yamanouchi is our dual action $\tilde{\mathfrak a}$ (of $\widehat{\mathcal G}^c$) on the crossed-product $A\rtimes_\mathfrak a\mathcal G$. Starting from a normal semi-finite faithful weight on $A$, in \cite{Y1}, \cite{Y2}, Yamanouchi defines a dual weight on $A\rtimes_\mathfrak a\mathcal G$ (\cite{Y2}, 4.10), using an operator-valued weight defined in (\cite{Y2}, 3.11), which is equal (\cite{Y2}, bottom of the page 665) to our operator-valued weight $T_{\tilde{\mathfrak a}}$; therefore, Yamanouchi's dual weight is the same as ours. 

%%%%cardualw
\subsection{Theorem}
\label{cardualw}
{\it Let $\gG$ be a quantum measured groupoid, and $(b, \mathfrak a)$ an action of $\gG$ on a von Neumann algebra $A$; let $\varphi$ be a normal semi-finite faithful weight on the crossed product $A\rtimes_\mathfrak a\gG$; then, are equivalent :
\newline
(i) there exists a normal semi-finite faithful weight $\psi$ on $A$ such that $\varphi=\tilde{\psi}$; 
\newline
(ii) the weight $\varphi$ is $\hat{\delta}^{-1}$-invariant, with respect to the dual action $\tilde{\mathfrak a}$ of $\widehat{\gG}^c$ on the crossed product $A\rtimes_\mathfrak a\gG$, and bears the density property defined in \ref{defdeltainv}.}

\begin{proof}
Let us suppose (i): then, we have (ii) by a simple application of \ref{thdual2}(ii).
\newline
 Let us suppose (ii), and let now $\theta$ be a normal semi-finite faithful weight on $A$, and $\tilde{\theta}$ its dual weight; by \ref{cocycle}, we obtain that $(D\varphi: D\tilde{\theta})_t$ belongs to $(A\rtimes_\mathfrak a\gG)^{\tilde{\mathfrak a}}$, which is equal to $\mathfrak a(A)$ by \ref{thduality1}(i). As $\sigma_t^{\tilde{\theta}}\circ\alpha=\alpha\circ\sigma_t^\theta$ by definition of the dual weight, there exists a normal semi-finite faithful weight $\varphi$ on $A$ such that $\mathfrak a[(D\psi:D\theta)_t]=(D\varphi: D\tilde{\theta})_t$. 
 \newline
 On the other hand, by definition of the dual weights, we have :
 \[(D\tilde{\psi}:D\tilde{\theta})_t=\mathfrak a[(D\psi:D\theta)_t]\]
 From which we deduce that $(D\varphi: D\tilde{\theta})_t=(D\tilde{\psi}:D\tilde{\theta})_t$, which gives the result. \end{proof}

%%%%%%GNSdual
\subsection{Theorem}
\label{GNSdual}
{\it Let $\gG$ be a quantum measured groupoid, and $(b, \mathfrak a)$ an action of $\gG$ on a von Neumann algebra $A$; let $\psi$ be a normal faithful semi-finite weight on $A$, and let $\tilde{\psi}$ be the dual weight constructed in \ref{defdualw1} on the crossed product $A\rtimes_\mathfrak a\gG$; then, the linear set generated by all the elements $(1\underset{N}{_b\otimes_\alpha}a)\mathfrak a(x)$, for all $x\in\gN_\psi$, $a\in\gN_{\widehat{\Phi}^c}\cap\gN_{\hat{T}^c}$, is a core for $\Lambda_{\tilde{\psi}}$, and it is possible to identify the GNS representation of $A\rtimes_\mathfrak a\gG$ associated to the weight $\tilde{\psi}$ with the natural representation on $H_\psi\underset{\nu}{_b\otimes_\alpha}H_\Phi$ by writing :
\[\Lambda_\psi(x)\underset{\nu}{_b\otimes_\alpha}\Lambda_{\widehat{\Phi}^c}(a)=\Lambda_{\tilde{\psi}}[(1\underset{N}{_b\otimes_\alpha}a)\mathfrak a(x)]\]
Moreover, using that identification, the linear set generated by the elements of the form $\mathfrak a(y^*)(\Lambda_\psi(x)\underset{\nu}{_b\otimes_\alpha}\Lambda_{\widehat{\Phi}^c}(a))$, for $x, y$ in $\gN_\psi$, and $a$ in $\gN_{\widehat{\Phi}^c}\cap\gN_{\hat{T}^c}\cap\gN_{\widehat{\Phi}^c}^*\cap\gN_{\hat{T}^c}^*$ is a core for $S_{\tilde{\psi}}$, and we have :}
\[S_{\tilde{\psi}}\mathfrak a(y^*)(\Lambda_\psi(x)\underset{\nu}{_b\otimes_\alpha}\Lambda_{\widehat{\Phi}^c}(a))=\mathfrak a(x^*)(\Lambda_\psi(y)\underset{\nu}{_b\otimes_\alpha}\Lambda_{\widehat{\Phi}^c}(a^*))\]

\begin{proof} 
Thanks to \ref{thdual2}(i) and (ii), the unitary $\tilde{V}$ constructed in \ref{propo1}(ii) is equal to the isometry $V$ constructed in \ref{propdualw}. So, using this unitary, we can identify $H_{\tilde{\psi}}$ with $H_\psi\underset{\nu}{_b\otimes_\alpha}H_\Phi$, which leads to the first result, using \ref{thduality1}(i), \ref{propdualw} and \ref{propo1}(ii). The second result comes then from \ref{Stildepsi}. \end{proof}

%%%%Upsi1
\subsection{Proposition}
\label{Upsi1}
{\it Let $\gG$ be a quantum measured groupoid, and $(b, \mathfrak a)$ an action of $\gG$ on a von Neumann algebra $A$; let $\psi$ be a normal faithful semi-finite weight on $A$, and let $\tilde{\psi}$ be the dual weight constructed in \ref{defdualw1} on the crossed product $A\rtimes_\mathfrak a\gG$; let us identify  $H_{\tilde{\psi}}$ with $H_\psi\underset{\nu}{_b\otimes_\alpha}H_\Phi$ as in \ref{GNSdual}. Then, the unitary $U^\mathfrak a_\psi$ from $H_\psi\underset{\nu^o}{_a\otimes_\beta}H_\Phi$ onto $H_\psi\underset{\nu}{_b\otimes_\alpha}H_\Phi$ defined by :
\[U^\mathfrak a_\psi=J_{\tilde{\psi}}(J_\psi\underset{N^o}{_a\otimes_\beta}J_{\widehat{\Phi}})\]
satisfies :
\[U^\mathfrak a_\psi(J_\psi\underset{N}{_b\otimes_\alpha}J_{\widehat{\Phi}})=(U^\mathfrak a_\psi)^*(J_\psi\underset{N^o}{_a\otimes_\beta}J_{\widehat{\Phi}})\]
and we have :
\newline
(i) for all $y\in A$ :
\[\mathfrak a (y)=U^\mathfrak a_\psi(y\underset{N^o}{_a\otimes_\beta}1)(U^\mathfrak a_\psi)^*\]
(ii) for all $b\in M$ :
\[(1\underset{N}{_b\otimes_\alpha}J_\Phi bJ_\Phi)U^\mathfrak a_\psi=U^\mathfrak a_\psi(1\underset{N^o}{_a\otimes_\beta}J_\Phi bJ_\Phi)\]
(iii) for all $n\in N$ :}
\[U_\psi^\mathfrak a(b(n)\underset{N^o}{_a\otimes_\beta}1)=(1\underset{N}{_b\otimes_\alpha}\beta(n))U_\psi^\mathfrak a\]
\[U_\psi^\mathfrak a(1\underset{N^o}{_a\otimes_\beta}\alpha(n))=(a(n)\underset{N}{_b\otimes_\alpha}1)U_\psi^\mathfrak a\]

\begin{proof}
Using again \ref{thdual2}(i) and \ref{thduality1}(i), the first result is given by \ref{propo2}(iii);  (i) is given by  \ref{propo2}(i), (ii) by \ref{Upsi}. Applying (i) to $y=b(n)$, we get the first result of (iii); we get then the second result of (iii) by using $U^\mathfrak a_\psi(J_\psi\underset{N}{_b\otimes_\alpha}J_{\widehat{\Phi}})=(U^\mathfrak a_\psi)^*(J_\psi\underset{N^o}{_a\otimes_\beta}J_{\widehat{\Phi}})$. \end{proof}

%%%lemTb
\subsection{Lemma}
\label{lemTb}
{\it Let $\gG$ be a measured quantum groupoid, and $(b,\mathfrak a)$ an action of $\gG$ on a von Neumann algebra $A$; let $(1\underset{N}{_b\otimes_\alpha}\hat{\beta}, \underline{\mathfrak a})$ be the action of $\gG$ on $A\underset{N}{_b*_\alpha}\mathcal L(H)$ constructed in \ref{thduality2}. Then, for all $a$ in $\gN_{\widehat{\Phi}^c}$ and $X$ in $A^+$, we have :}
\[T_{\underline{\mathfrak a}}(\rho^{b, \alpha}_{\Lambda_{\widehat{\Phi}^c}(a)}X(\rho^{b, \alpha}_{\Lambda_{\widehat{\Phi}^c}(a)})^*)=(1\underset{N}{_b\otimes_\alpha}a)\mathfrak a (X)(1\underset{N}{_b\otimes_\alpha}a^*)\]

\begin{proof}
Let us first remark that $\Lambda_{\widehat{\Phi}^c}(a)=J_{\widehat{\Phi}}\Lambda_{\widehat{\Phi}}(J_{\widehat{\Phi}}aJ_{\widehat{\Phi}})$ which belongs to $D(_\alpha H_\Phi, \nu)$, by \ref{basic}; on the other hand, let us suppose that $A$ is acting on an Hilbert space $\gH$; we verify straightforwardly that $\rho^{b, \alpha}_{\Lambda_{\widehat{\Phi}^c}(a)}X(\rho^{b, \alpha}_{\Lambda_{\widehat{\Phi}^c}(a)})^*$ commutes with all $Y\underset{N}{_b\otimes_\alpha}1$, with $Y\in A'$, and, therefore, belongs to $(A\underset{N}{_b*_\alpha}\mathcal L(H))^+$. 
\newline
Let $\eta$ in $D(\gH_b, \nu^o)$; for all $\xi$ in $D((H_\Phi)_{\hat{\beta}}, \nu^o)$, we have, using \ref{thduality2} :
\begin{multline*}
<T_{\underline{\mathfrak a}}(\rho^{b, \alpha}_{\Lambda_{\widehat{\Phi}^c}(a)}X(\rho^{b, \alpha}_{\Lambda_{\widehat{\Phi}^c}(a)})^*), \omega_{\eta\underset{N}{_b\otimes_\alpha}\xi}>
=\\
\Phi[(\omega_\eta\underset{N}{_b*_\alpha}id)
(1\underset{N}{_b\otimes_\alpha}(\rho_\xi^{\alpha, \hat{\beta}})^* W\rho^{\beta, \alpha}_{\Lambda_{\widehat{\Phi}^c}(a)})\mathfrak a (X)(1\underset{N}{_b\otimes_\alpha}(\rho^{b, \alpha}_{\Lambda_{\widehat{\Phi}^c}(a)})^* W^*\rho_\xi^{\alpha, \hat{\beta}})]
\end{multline*}
We know that $(\rho^{b, \alpha}_{\Lambda_{\widehat{\Phi}^c}(a)})^* W^*\rho_\xi^{\alpha, \hat{\beta}}=(id*\omega_{\xi, \Lambda_{\widehat{\Phi}^c}(a)})(W^*)$
belongs to $\gN_\Phi$, thanks to \ref{propdualw} applied to $\widehat{\gG}$, and that :
\[\Lambda_\Phi((id*\omega_{\xi, \Lambda_{\widehat{\Phi}^c}(a)})(W^*))=a_{\widehat{\Phi}}(\omega_{\xi, \Lambda_{\widehat{\Phi}^c}(a)})=a^*\xi\]
and, therefore, we get that $<T_{\underline{\mathfrak a}}(\rho^{b, \alpha}_{\Lambda_{\widehat{\Phi}^c}(a)}X(\rho^{b, \alpha}_{\Lambda_{\widehat{\Phi}^c}(a)})^*), \omega_{\eta\underset{N}{_b\otimes_\alpha}\xi}>$ is equal to :
\begin{multline*}
\Phi((id*\omega_{\xi, \Lambda_{\widehat{\Phi}^c}(a)})(W^*)^*(\omega_\eta\underset{N}{_b*_\alpha}id)(X)(id*\omega_{\xi, \Lambda_{\widehat{\Phi}^c}(a)})(W^*))=\\
=
((\omega_\eta\underset{N}{_b*_\alpha}id)(X)\Lambda_\Phi((id*\omega_{\xi, \Lambda_{\widehat{\Phi}^c}(a)})(W^*)|\Lambda_\Phi(id*\omega_{\xi, \Lambda_{\widehat{\Phi}^c}(a)})(W^*))\\
=((\omega_\eta\underset{N}{_b*_\alpha}id)(X)a^*\xi|a^*\xi)
\end{multline*}
and, finally, we get that :
\[<T_{\underline{\mathfrak a}}(\rho^{b, \alpha}_{\Lambda_{\widehat{\Phi}^c}(a)}X(\rho^{b, \alpha}_{\Lambda_{\widehat{\Phi}^c}(a)})^*), \omega_{\eta\underset{N}{_b\otimes_\alpha}\xi}>
=((1\underset{N}{_b\otimes_\alpha}a)\mathfrak a (X)(1\underset{N}{_b\otimes_\alpha}a^*)(\eta\underset{N}{_b\otimes_\alpha}\xi)|\eta\underset{N}{_b\otimes_\alpha}\xi)\]
from which we get the result. 
\end{proof}

%%%%spatialder
\subsection{Theorem}
\label{spatialder}
{\it Let $\gG$ be a measured quantum groupoid, and $(b,\mathfrak a)$ an action of $\gG$ on a von Neumann algebra $A$; let $\psi$ be a normal semi-finite faithful weight on $A$; let $A\rtimes_{\mathfrak a}\gG$ be the crossed-product, as defined in \ref{defcross}, and $\tilde{\psi}$ the dual weight, as defined in \ref{defdualw1}; let us consider the dual action, as defined in \ref{dualaction}, the bicrossed-product and the bidual weight $\tilde{\tilde{\psi}}$. We have defined in \ref{thduality2} an isomorphism $\Theta$ such that :
\[\Theta(\mathfrak a\underset{N}{_b*_\alpha}id)(A\underset{N}{_b*_\alpha}\mathcal L(H_\Phi))=(A\rtimes_{\mathfrak a}\gG)\rtimes_{\tilde{\mathfrak a}}\widehat{\gG}^o\]
We have then :}
\[\frac{d\tilde{\tilde{\psi}}\circ \Theta(\mathfrak a\underset{N}{_b*_\alpha}id)}{d\psi^o}=\Delta_{\tilde{\psi}}\]

\begin{proof}
Let us remark that $\tilde{\tilde{\psi}}\circ \Theta(\mathfrak a\underset{N}{_b*_\alpha}id)$ is a normal semi-finite weight on $A\underset{N}{_b*_\alpha}\mathcal L(H_\Phi)$, whose commutant is $A'\underset{N}{_b\otimes_\alpha}1$, which is isomorphic to $A'$, and, therefore to $A^o$. Therefore, the spatial derivative has a meaning. Let us write it $h$ for simplification. By definition of this spatial derivative, for all $\Xi$ in $D(H_\psi\underset{\nu}{_b\otimes_\alpha}H_\Phi, \psi^o)$,  $\Xi$ belongs to $D(h^{1/2})$ if and only if $\theta^{\psi^o}(\Xi, \Xi)$ belongs to $\gM_{\tilde{\tilde{\psi}}\circ \Theta(\mathfrak a\underset{N}{_b*_\alpha}id)}^+$, and we have then :
\[\|h^{1/2}\Xi\|^2=\tilde{\tilde{\psi}}\circ \Theta(\mathfrak a\underset{N}{_b*_\alpha}id)(\theta^{\psi^o}(\Xi, \Xi))\]
which, using \ref{thduality2}, gives that :
\[\|h^{1/2}\Xi\|^2=\tilde{\psi}\circ T_b(\theta^{\psi^o}(\Xi, \Xi))\]
Let now $x$, $y$ in $\gN_\psi$, and $a$ in $\gN_{\hat{T}^c}\cap \gN_{\hat{T}^c}^*\cap\gN_{\widehat{\Phi}}\cap\gN_{\widehat{\Phi}}^*$; we have, for all $z\in\gN_\psi$ :
\begin{eqnarray*}
(J_\psi z^*J_\psi\underset{N}{_b\otimes_\alpha}1)(\mathfrak a(x^*)(\Lambda_\psi(y)\underset{N}{_b\otimes_\alpha}\Lambda_{\widehat{\Phi}^c(a)})
&=&
\mathfrak a (x^*)(J_\psi z^*J_\psi\Lambda_\psi(y)\underset{N}{_b\otimes_\alpha}\Lambda_{\widehat{\Phi}^c(a)})\\
&=&\mathfrak a(x^*)(yJ_\psi\Lambda_\psi(z)\underset{N}{_b\otimes_\alpha}\Lambda_{\widehat{\Phi}^c(a)})\\
&=&\mathfrak a(x^*)\rho^{b, \alpha}_{\Lambda_{\widehat{\Phi}^c(a)}}yJ_\psi\Lambda_\psi(z)
\end{eqnarray*}
from which we get that $\mathfrak a(x^*)(\Lambda_\psi(y)\underset{N}{_b\otimes_\alpha}\Lambda_{\widehat{\Phi}^c(a)})$ belongs to $D(H_\psi\underset{\nu}{_b\otimes_\alpha}H_\Phi, \psi^o)$, and that :
\[R^{\psi^o}(\mathfrak a(x^*)(\Lambda_\psi(y)\underset{N}{_b\otimes_\alpha}\Lambda_{\widehat{\Phi}^c(a)}))=\mathfrak a(x^*)\rho^{b, \alpha}_{\Lambda_{\widehat{\Phi}^c(a)}}y\]
Therefore, for $x, y$ in $\gN_\psi$, $a$ in $\gN_{\hat{T}^c}\cap \gN_{\hat{T}^c}^*\cap\gN_{\widehat{\Phi}}\cap\gN_{\widehat{\Phi}}^*$, we get that :
\[\theta^{\psi^o}(\mathfrak a(x^*)(\Lambda_\psi(y)\underset{N}{_b\otimes_\alpha}\Lambda_{\widehat{\Phi}^c(a)}), \mathfrak a(x^*)(\Lambda_\psi(y)\underset{N}{_b\otimes_\alpha}\Lambda_{\widehat{\Phi}^c(a)})=\mathfrak a(x^*)\rho^{b, \alpha}_{\Lambda_{\widehat{\Phi}^c(a)}}yy^*(\rho^{b, \alpha}_{\Lambda_{\widehat{\Phi}^c(a)}})^*\mathfrak a(x)\]
which, thanks to \ref{lemTb}, belongs to $\gM_{T_{\underline{\mathfrak a}}}^+$, and we have :
\begin{multline*}
T_{\underline{\mathfrak a}}(\theta^{\psi^o}(\mathfrak a(x^*)(\Lambda_\psi(y)\underset{N}{_b\otimes_\alpha}\Lambda_{\widehat{\Phi}^c(a)}), \mathfrak a(x^*)(\Lambda_\psi(y)\underset{N}{_b\otimes_\alpha}\Lambda_{\widehat{\Phi}^c(a)}))=\\
\mathfrak a(x^*)(1\underset{N}{_b\otimes_\alpha}a)\mathfrak a(yy^*)(1\underset{N}{_b\otimes_\alpha}a^*)\mathfrak a(x)
\end{multline*}
Using now \ref{GNSdual}, we get that this last operator belongs to $\gM_\psi^+$, and, therefore, that 
$\mathfrak a(x^*)(\Lambda_\psi(y)\underset{N}{_b\otimes_\alpha}\Lambda_{\widehat{\Phi}^c(a)})$ belongs to $D(h^{1/2})$. And, by \ref{GNSdual}, we get that $\mathfrak a(y^*)(1\underset{N}{_b\otimes_\alpha}a^*)\mathfrak a(x)$ belongs to $\gN_{\tilde{\psi}}\cap\gN_{\tilde{\psi}}^*$, and that :
\[\tilde{\psi}\circ T_{\underline{\mathfrak a}}(\theta^{\psi^o}(\mathfrak a(x^*)(\Lambda_\psi(y)\underset{N}{_b\otimes_\alpha}\Lambda_{\widehat{\Phi}^c}(a)), \mathfrak a(x^*)(\Lambda_\psi(y)\underset{N}{_b\otimes_\alpha}\Lambda_{\widehat{\Phi}^c}(a)))\]
is equal to :
\begin{eqnarray*}
\|\Lambda_{\tilde{\psi}}(\mathfrak a(y^*)(1\underset{N}{_b\otimes_\alpha}a^*)\mathfrak a(x))\|^2
&=&
\|S_{\tilde{\psi}}\Lambda_{\tilde{\psi}}[\mathfrak a(x^*)(1\underset{N}{_b\otimes_\alpha}a)\mathfrak a(y)]\|^2\\
&=&\|\Delta_{\tilde{\psi}}^{1/2}\mathfrak a(x^*)(\Lambda_\psi(y)\underset{\nu}{_b\otimes_\alpha}\Lambda_{\widehat{\Phi}^c}(a))\|^2
\end{eqnarray*}
Therefore, we get that :
\[\|h^{1/2}\mathfrak a(x^*)(\Lambda_\psi(y)\underset{\nu}{_b\otimes_\alpha}\Lambda_{\widehat{\Phi}^c}(a))\|=\|\Delta_{\tilde{\psi}}^{1/2}\mathfrak a(x^*)(\Lambda_\psi(y)\underset{\nu}{_b\otimes_\alpha}\Lambda_{\widehat{\Phi}^c}(a))\|\]
and, using again \ref{GNSdual}, as the vectors $\mathfrak a(x^*)(\Lambda_\psi(y)\underset{\nu}{_b\otimes_\alpha}\Lambda_{\widehat{\Phi}^c}(a))$ are a core for $\Delta_{\tilde{\psi}}^{1/2}$, we have $\|h^{1/2}\xi\|=\|\Delta_{\tilde{\psi}}^{1/2}\xi\|$ for all $\xi$ in $D(h^{1/2})$, from which we get that $\Delta_{\tilde{\psi}}\subset h$, and, as they are self-adjoint operators, we get the result. 
\end{proof}

%%%%thstandard
\subsection{Theorem}
\label{thstandard}
{\it Let $\gG$ be a measured quantum groupoid, and $(b,\mathfrak a)$ an action of $\gG$ on a von Neumann algebra $A$, $A\rtimes_\mathfrak a \gG$ its crossed product; let $(1\underset{N}{_b\otimes_\alpha}\hat{\beta}, \underline{\mathfrak a})$ be the action of $\gG$ on $A\underset{N}{_b*_\alpha}\mathcal L(H_\Phi)$ introduced in \ref{thduality2}; then, the inclusion :
\[\mathfrak a (A)\subset A\rtimes_\mathfrak a \gG\subset A\underset{N}{_b*_\alpha}\mathcal L(H_\Phi)\]
is standard, and the operator-valued weight $T_{\underline{\mathfrak a}}$ from $A\underset{N}{_b*_\alpha}\mathcal L(H_\Phi)$ to $A\rtimes_\mathfrak a \gG$ is the basic construction made from the operator-valued weight $T_{\tilde{\mathfrak a}}$ from $A\rtimes_\mathfrak a \gG$ to $\mathfrak a (A)$.}

\begin{proof}
Let $\psi$ be a normal semi-finite faithful weight on $A$, and $\tilde{\psi}$ its dual weight on $A\rtimes_\mathfrak a \gG$; let us represent the inclusion $\mathfrak a (A)\subset A\rtimes_\mathfrak a \gG$ on the Hilbert space $H_{\tilde{\psi}}$, which had been identified (\ref{GNSdual}) with $H_\psi\underset{\nu}{_b\otimes_\alpha}H_\Phi$, equipped with the natural representation of the inclusion. 
We have then, using \ref{Upsi1} :
\begin{eqnarray*}
J_{\tilde{\psi}}\mathfrak a(A)J_{\tilde{\psi}}&=&J_{\tilde{\psi}}U^\mathfrak a_\psi(A\underset{N^o}{_a\otimes_\beta}1)(U^\mathfrak a_\psi)^*J_{\tilde{\psi}}\\
&=&(J_\psi\underset{N^o}{_a\otimes_\beta}J_{\hat{\Phi}})(A\underset{N^o}{_a\otimes_\beta}1)
(J_\psi\underset{N}{_b\otimes_\alpha}J_{\hat{\Phi}})\\
&=&A'\underset{N}{_b\otimes_\alpha}1
\end{eqnarray*}
and, therefore, we get :
\[J_{\tilde{\psi}}\mathfrak a(A)'J_{\tilde{\psi}}=(A'\underset{N}{_b\otimes_\alpha}1)'=A\underset{N}{_b*_\alpha}\mathcal L(H_\Phi)\]
which gives that $A\underset{N}{_b*_\alpha}\mathcal L(H_\Phi)$ is the standard construction made from the inclusion $\mathfrak a (A)\subset A\rtimes_\mathfrak a \gG$. 
\newline
Let $T$ be the standard construction made from the operator-valued weight $T_{\tilde{\mathfrak a}}$. By definition (\ref{basic}), we have :
\[\frac{d\tilde{\psi}\circ T}{d\psi^o}=\Delta_{\tilde{\psi}}\]
which, thanks to \ref{spatialder}, gives that $\tilde{\psi}\circ T=\tilde{\tilde{\psi}}\circ\Theta(\mathfrak a\underset{N}{_b*_\alpha}id)$, and, using \ref{thduality2}, we have :
\[\tilde{\psi}\circ T=\tilde{\psi}\circ T_{\underline{\mathfrak a}}\]
from which we get that $T=T_{\underline{\mathfrak a}}$, which finishes the proof. \end{proof}

%%%%%thD2
\subsection{Theorem}
\label{thD2}
{\it Let $\gG$ be a measured quantum groupoid, and $(b,\mathfrak a)$ an action of $\gG$ on a von Neumann algebra $A$, $A\rtimes_\mathfrak a \gG$ its crossed product; then :
\newline
(i) the Jones' tower associated to the inclusion $\mathfrak a(A)\subset A\rtimes_\mathfrak a\gG$ is :
\[\mathfrak a(A)\underset{N^o}{_{\hat{\alpha}}\otimes_\beta}1\subset \tilde{\mathfrak a}(A\rtimes_\mathfrak a\gG)\subset (A\rtimes_\mathfrak a\gG)\rtimes_{\tilde{\mathfrak a}}\widehat{\gG}^c\subset (A\rtimes_\mathfrak a\gG)\underset{N^o}{_{\hat{\alpha}}*_\beta}\mathcal L(H_\Phi)\]
(ii) moreover, the operator-valued weight  $T_{\underline{(\tilde{\mathfrak a})}}$ from $(A\rtimes_\mathfrak a\gG)\underset{N^o}{_{\hat{\alpha}}*_\beta}\mathcal L(H_\Phi)$ to $(A\rtimes_\mathfrak a\gG)\rtimes_{\tilde{\mathfrak a}}\widehat{\gG}^c$ is the operator-valued weight constructed by successive basic constructions in this tower. 
\newline
(iii) the derived tower is :
\begin{multline*}
\tilde{\mathfrak a}(A\rtimes_\mathfrak a\gG\cap\mathfrak a(A)')\subset
(A\rtimes_\mathfrak a\gG)\rtimes_{\tilde{\mathfrak a}}\widehat{\gG}^c\cap
\mathfrak a(A)'\underset{N^o}{_{\hat{\alpha}}*_\beta}\mathcal L(H_\Phi)\\\subset
[A\rtimes_\mathfrak a\gG\cap\mathfrak a(A)']\underset{N^o}{_{\hat{\alpha}}*_\beta}\mathcal L(H_\Phi)
\end{multline*}
(iv) let us write $B=A\rtimes_{\mathfrak a}\gG\cap \mathfrak a(A)'$, and $\mathfrak b=\tilde{\mathfrak a}_{|B}$. Then $(1\underset{N}{_b\otimes_\alpha}\hat{\alpha}, \mathfrak b)$ is an action of $\widehat{\gG}^c$ on $B$, and we have :
\[B\rtimes_\mathfrak b\widehat{\gG}^c=(A\rtimes_\mathfrak a\gG)\rtimes_{\tilde{\mathfrak a}}\widehat{\gG}^c\cap
\mathfrak a(A)'\underset{N^o}{_{\hat{\alpha}}*_\beta}\mathcal L(H_\Phi)\]
(v) the inclusion $\mathfrak a(A)\subset A\rtimes_\mathfrak a\gG$ is depth 2.}

\begin{proof}
We have got in \ref{thstandard} that the inclusion :
\[\mathfrak a (A)\subset A\rtimes_\mathfrak a \gG\subset A\underset{N}{_b*_\alpha}\mathcal L(H_\Phi)\]
is standard. Let's apply the isomorphism $\Theta(\mathfrak a\underset{N}{_b*_\alpha}id)$ from $A\underset{N}{_b*_\alpha}\mathcal L(H_\Phi)$ onto $(A\rtimes_\mathfrak a\gG)\rtimes_{\tilde{\mathfrak a}}\widehat{\gG}^c$ defined in \ref{thduality2}; we then get the inclusion :
\[\mathfrak a(A)\underset{N^o}{_{\hat{\alpha}}\otimes_\beta}1\subset \tilde{\mathfrak a}(A\rtimes_\mathfrak a\gG)\subset (A\rtimes_\mathfrak a\gG)\rtimes_{\tilde{\mathfrak a}}\widehat{\gG}^c\]
of von Neumann algebras on $H_\psi\underset{\nu}{_b\otimes_\alpha}H_\Phi\underset{\nu^o}{_{\hat{\alpha}}\otimes_\beta}H_\Phi$; let's apply now \ref{thstandard} to the inclusion :
\[\tilde{\mathfrak a}(A\rtimes_\mathfrak a\gG)\subset (A\rtimes_\mathfrak a\gG)\rtimes_{\tilde{\mathfrak a}}\widehat{\gG}^c\]
We get that the basic construction made from this inclusion gives the von Neumann algebra $(A\rtimes_\mathfrak a\gG)\underset{N^o}{_{\hat{\alpha}}*_\beta}\mathcal L(H_\Phi)$, which gives (i). 
\newline
By \ref{thduality2} again, we see that the isomorphism $\Theta(\mathfrak a\underset{N}{_b*_\alpha}id)$ sends the operator-valued weight $T_{\underline{\mathfrak a}}$ on $T_{\tilde{\tilde{\mathfrak a}}}$. So, 
applying again \ref{thstandard} to this operator-valued weight, we get (ii). 
\newline
Now, (iii) is just a corollary from (i), and (iv) is just an application of \ref{corcross}. So, we get that 
the derived tower obtained in (iii) can be written as :
\[\mathfrak b(B)\subset B\rtimes_\mathfrak b\widehat{\gG}^c\subset B\underset{N^o}{_{\hat{\alpha}}*_\beta}\mathcal L(H_\Phi)\]
which is standard, by \ref{thstandard} applied to the action $\mathfrak b$; which finishes the proof. \end{proof}

%%%%%threg
\subsection{Theorem}
\label{threg}
{\it Let $\gG$ be a measured quantum groupoid, and $(b,\mathfrak a)$ an action of $\gG$ on a von Neumann algebra $A$, $A\rtimes_\mathfrak a \gG$ its crossed product, $\underline{\mathfrak a}$ the action of $\gG$ on $A\underset{N}{_b*_\alpha}\mathcal L(H_\Phi)$ introduced in \ref{thduality2}. Then :
\newline
(i) for any $x\in M'^+$, $1\underset{N}{_b\otimes_\alpha}x$ belongs to $A\underset{N}{_b*_\alpha}\mathcal L(H_\Phi)\cap \mathfrak a(A)'$, and we have :
\[T_{\underline{\mathfrak a}}(1\underset{N}{_b\otimes_\alpha}x)=1\underset{N}{_b\otimes_\alpha}T^{co}(x)\]
(ii) the restriction of $T_{\underline{\mathfrak a}}$ to $A\underset{N}{_b*_\alpha}\mathcal L(H_\Phi)\cap \mathfrak a(A)'$ is semi-finite;
\newline
(iii) the operator-valued weight $T_{\tilde{\mathfrak a}}$ from $A\rtimes_\mathfrak a \gG$ to $\mathfrak a(A)$ is regular in the sense of \ref{basic}. }

\begin{proof}
Let $x\in M'$; we have, using \ref{thL4} :
\begin{eqnarray*}
W(1\underset{N}{_\beta\otimes_\alpha}x)W^*
&=&
(J_{\widehat{\Phi}}\underset{N}{_\beta\otimes_\alpha}J_\Phi)W^*(1\underset{N^o}{_\alpha\otimes_{\hat{\beta}}}J_\Phi xJ_\Phi)W(J_{\widehat{\Phi}}\underset{N}{_\beta\otimes_\alpha}J_\Phi)\\
&=&(J_{\widehat{\Phi}}\underset{N}{_\beta\otimes_\alpha}J_\Phi)\Gamma(J_\Phi xJ_\Phi)(J_{\widehat{\Phi}}\underset{N}{_\beta\otimes_\alpha}J_\Phi)
\end{eqnarray*}
and, therefore, using \ref{thduality2} :
\[\underline{\mathfrak a}(1\underset{N}{_b\otimes_\alpha}x)=
(id\underset{N}{_b*_\alpha}\varsigma_N)(1\underset{N}{_b\otimes_\alpha}[(J_{\widehat{\Phi}}\underset{N}{_\beta\otimes_\alpha}J_\Phi)\Gamma(J_\Phi xJ_\Phi)(J_{\widehat{\Phi}}\underset{N}{_\beta\otimes_\alpha}J_\Phi)])\]
from which we get, for $x$ positive :
\begin{eqnarray*}
T_{\underline{\mathfrak a}}(1\underset{N}{_b\otimes_\alpha}x)
&=&
1\underset{N}{_b\otimes_\alpha}J_\Phi(\Phi\circ R\underset{N}{_\beta*_\alpha}id)\Gamma(J_\Phi xJ_\Phi)J_\Phi\\
&=&1\underset{N}{_b\otimes_\alpha}J_\Phi RTR(J_\Phi xJ_\Phi)J_\Phi\\
&=&1\underset{N}{_b\otimes_\alpha}T^{co}(x)
\end{eqnarray*}
which is (i).  Then, (ii) is just a corollary of (i).
\newline
Let us apply now (ii) to the action $\tilde{\mathfrak a}$, we get that the restriction of $T_{\underline{(\tilde{\mathfrak a})}}$ to the von Neumann algebra $(A\rtimes_\mathfrak a\gG)\underset{N^o}{_{\hat{\alpha}}*_\beta}\mathcal L(H_\Phi)\cap \tilde{\mathfrak a}(A\rtimes_\mathfrak a\gG)'$ is semi-finite, which, with (ii), gives (iii). \end{proof}

%%%remark2
\subsection{Remark}
\label{remark2}
The fact that the inclusion $\mathfrak a(A)\subset A\rtimes_\mathfrak a\gG$ is depth 2 (\ref{thD2}(v)) had been obtained in ([V2], 5.10) for actions of locally compact quantum groups in a somehow different way. 

%%%%mqga
\section{The measured quantum groupoid associated to an action}
\label{mqga}

In this section, we apply the results of (\cite{EV}, \cite{E4}), recalled in \ref{d2}, \ref{d2action} and \ref{d2cross}, to the depth 2 inclusion (\ref{thD2}) $\mathfrak a(A)\subset A\rtimes_\mathfrak a\gG$, where $(b, \mathfrak a)$ is an action of $\gG$ on a von Neumann algebra $A$. From such data, we therefore obtain a new quantum measured groupoid $\gG(\mathfrak a)$ (\ref{th2mqga}), and we show that the dual action $\tilde{\mathfrak a}$ of $\widehat{\gG}^c$ can be considered as an action of $\gG(\tilde{\mathfrak a})$ (\ref{cor1mqga}). As the underlying von Neumann algebra of $\gG(\mathfrak a)$ is a crossed-product of the relative commutant $A\rtimes_\mathfrak a\gG\cap\mathfrak a(A)'$ by the restriction of the action $\tilde{\mathfrak a}$ of $\widehat{\gG}^c$, we show that the initial measured quantum groupoid can be naturally sent into $\gG(\mathfrak a)$ (\ref{prop1mqga}), and, moreover, that, in the particular case when we are starting with an outer action of a locally compact quantum group, we recover the initial locally compact quantum group (\ref{outer}). 

%%%%th1mqga
\subsection{Theorem}
\label{th1mqga}
{\it Let $\gG$ be a measured quantum groupoid, and $(b, \mathfrak a)$ an action of $\gG$ on a von Neumann algebra $A$. Let $\psi$ be a normal faithful semi-finite weight on $A$, and $\tilde{\psi}$ be the dual weight on the crossed product $A\rtimes_\mathfrak a\gG$. Let us write $B$ for $A\rtimes_\mathfrak a\gG\cap\mathfrak a(A)'$, $\mathfrak b$ for $\tilde{\mathfrak a}_{|B}$, and $\tilde{M}$ for $A\underset{N}{_b*_\alpha}\mathcal L(H_\Phi)\cap \mathfrak a(A)'$; then :
\newline
(i) let us write $id$ for the inclusion of $B$ into $\tilde{M}$, and $j$ for the anti-$*$-homomorphism $j(x)=J_{\tilde{\psi}}x^*J_{\tilde{\psi}}$, which sends $B$ into $A\underset{N}{_b*_\alpha}\mathcal L(H_\Phi)\cap (A\rtimes_\mathfrak a\gG)'$ (and, therefore, into $\tilde{M}$). Then, there exists $\tilde{\Gamma}$, such that :
\[(B, \tilde{M}, id, j, \tilde{\Gamma})\]
is a Hopf-bimodule. 
\newline
(ii) let us write also $j(x)=J_{\tilde{\psi}}x^*J_{\tilde{\psi}}$ the $*$-anti-automorphism of $\tilde{M}$; then $j$ is a co-inverse for the coproduct constructed in (i). 
\newline
(iii) let us write $\tilde{T}$ for the restriction of $T_{\underline{\mathfrak a}}$ to $\tilde{M}$; then, $\tilde{T}$ is a normal, semi-finite, faithful operator-valued weight from $\tilde{M}$ onto $B$, which is left-invariant with respect to the coproduct $\tilde{\Gamma}$.
\newline
(iv) the isomorphism $\Theta (\mathfrak a\underset{N}{_b*_\alpha}id)$ sends $\tilde{M}$ onto $B\rtimes_\mathfrak b\widehat{\gG}^c$ and $B$ onto $\mathfrak b(B)$. Let us denote $j_1$ the $*$-anti-automorphism of $B\rtimes_\mathfrak b\widehat{\gG}^c$ obtained by transporting $j$ via this isomorphism, and $\tilde{\Gamma}_1$ the coproduct obtained the same way. Then :
\[ (B, B\rtimes_\mathfrak b\widehat{\gG}^c, \mathfrak b, j_1\circ\mathfrak b, \tilde{\Gamma}_1)\]
is a Hopf-bimodule, isomorphic to the one obtained in (i), $j_1$ is a co-inverse, and $T_{\tilde{\mathfrak b}}$ is a left-invariant operator-valued weight. }

\begin{proof}
We had got in \ref{thD2} that the inclusion $\mathfrak a(A)\subset A\rtimes_\mathfrak a\gG$ is depth 2, and, in \ref{threg}, that the operator-valued weight $T_{\underline{\mathfrak a}}$ is regular. So, we can apply \ref{d2}(i) and (ii). As, by \ref{thstandard}, the basic construction made from this inclusion gives the algebra $A\underset{N}{_b*_\alpha}\mathcal L(H_\Phi)$, and the beginning of the derived tower is $B\subset \tilde{M}$, we get (i), (ii) and (iii). Using \ref{thD2}(iv), we get the beginning of (iv); moreover, using \ref{propduality1}, we get that this isomorphism sends the restriction of $T_{\underline{\mathfrak a}}$ on the restriction of $T_{\tilde{\tilde{\mathfrak a}}}$, which is $T_{\tilde{\mathfrak b}}$.
\end{proof}

%%%%th2mqga
\subsection{Theorem}
\label{th2mqga}
{\it Let us take the notations of \ref{th1mqga}. Let us suppose that there exists a normal semi-finite faithful weight $\chi$ on $B$, invariant under the modular automorphism group $\sigma_t^{T_{\tilde{\mathfrak a}}}$; then :
\[(B, \tilde{M}, id, j, \tilde{\Gamma}, \tilde{T}, j\circ\tilde{T}\circ j, \chi)=(B, B\rtimes_\mathfrak b\widehat{\gG}^c, \mathfrak b, j_1\circ\mathfrak b, \tilde{\Gamma}_1, T_{\tilde{\mathfrak b}}, j_1\circ T_{\tilde{\mathfrak b}}\circ j_1, \chi)\]
is a measured quantum groupoid, we shall denote by $\gG(\mathfrak a)$. Moreover, the dual action $\tilde{\mathfrak a}$ satisfies the same hypothesis, and $\gG(\tilde{\mathfrak a})=\widehat{\gG(\mathfrak a)}^o$. }
\begin{proof}
This is just \ref{d2} and \ref{thD2}. \end{proof}

%%%%%%th3mqga
\subsection{Theorem}
\label{th3mqga}
{\it Let us take the notations of \ref{th1mqga} and \ref{th2mqga}; there exists an action $(id, \overline{\mathfrak a})$ of $\gG(\tilde{\mathfrak a})$ on $A\rtimes_\mathfrak a\gG$ (where $id$ means the inclusion of $B$ into $A\rtimes_\mathfrak a\gG$, which is a anti-$*$-homomorphism of the basis $B^o$ of $\gG(\tilde{\mathfrak a})$), such that $(A\rtimes_\mathfrak a\gG)^{\overline{\mathfrak a}}=\mathfrak a (A)$, and such that the crossed-product $(A\rtimes_\mathfrak a\gG)\rtimes_{\overline{\mathfrak a}}\gG(\tilde{\mathfrak a})$ is isomorphic to $A\underset{N}{_b*_\alpha}\mathcal L(H_\Phi)$. More precisely, there exists an isomorphism $I_{\overline{\mathfrak a}}$ from $(A\rtimes_\mathfrak a\gG)\rtimes_{\overline{\mathfrak a}}\gG(\tilde{\mathfrak a})$ onto $A\underset{N}{_b*_\alpha}\mathcal L(H_\Phi)$, such that :
\newline
(i) for any $x\in A\rtimes_\mathfrak a\gG$, we have $I_{\overline{\mathfrak a}}(\overline{\mathfrak a}(x))=x$; 
\newline
(ii) for any $y\in A\underset{N}{_b*_\alpha}\mathcal L(H_\Phi)\cap \mathfrak a (A)'$, we have $I_{\overline{\mathfrak a}}(1\underset{B^o}{_{id}\otimes_{\hat{\jmath}\circ \mathfrak b}}y)=y$, where $\hat{\jmath}(z)=J_{\tilde{\chi}}z^*J_{\tilde{\chi}}$, for all $z\in B\rtimes_\mathfrak b\widehat{\gG}^c$.}

\begin{proof}
This is just \ref{d2cross} applied to \ref{thD2} and \ref{th2mqga}. \end{proof}

%%%%%th4mqga
\subsection{Theorem}
\label{th4mqga}
{\it Let us take the notations of \ref{th1mqga}, \ref{th2mqga} and \ref{th3mqga}; then, 
the application defined, for any $x\in\gN_{\tilde{\psi}}$, $y\in\gN_\chi$, $b\in\gN_{\Phi^{oc}}\cap\gN_{T^{oc}}$ by :
\[\mathcal U[\Lambda_{\tilde{\psi}}(x)\underset{\chi^o}{_{id}\otimes _{\hat{\jmath}\circ b}}(\Lambda_\chi(y)\underset{\nu^o}{_{\hat{\alpha}}\otimes_\beta}\Lambda_{\Phi^{oc}}(b))]=
\Lambda_{\tilde{\psi}}(yx)\underset{\nu^o}{_{\hat{\alpha}}\otimes_\beta}\Lambda_{\Phi^{oc}}(b)\]
is a unitary from $H_{\tilde{\psi}}\underset{\chi^o}{_{id}\otimes _{\hat{\jmath}\circ b}}H_{\tilde{\chi}}$ onto $H_{\tilde{\psi}}\underset{\nu^o}{_{\hat{\alpha}}\otimes_\beta}H_\Phi$, such that, for all $X\in (A\rtimes_\mathfrak a\gG)\rtimes_{\tilde{\mathfrak a}}\widehat{\gG}^c$, we have :
\[\mathcal UX\mathcal U^*=\Theta\circ(\mathfrak a\underset{N}{_b*_\alpha}id)I_{\overline{\mathfrak a}}(X)\]
and, in particular, for all $Y\in A\rtimes_\mathfrak a\gG$, $c\in M'$, we have :}
\[\mathcal U\overline{\mathfrak a}(Y)\mathcal U^*=\tilde{\mathfrak a}(Y)\]
\[\mathcal U(1\underset{B^o}{_{id}\otimes_{\hat{\jmath}\circ\mathfrak b}}(1\underset{N}{_b\otimes_\alpha}1\underset{N^o}{_{\hat{\alpha}}\otimes_\beta}c))\mathcal U^*=1\underset{N}{_b\otimes_\alpha}1\underset{N^o}{_{\hat{\alpha}}\otimes_\beta}c\]
\begin{proof}
Thanks to \ref{th3mqga} and \ref{thduality2}, we get that $\Theta\circ(\mathfrak a\underset{N}{_b*_\alpha}id)I_{\overline{\mathfrak a}}$ is an isomorphism between $(A\rtimes_\mathfrak a\gG)\rtimes_{\overline{\mathfrak a}}\gG(\tilde{\mathfrak a})$ and $(A\rtimes_\mathfrak a\gG)\rtimes_{\tilde{\mathfrak a}}\widehat{\gG}^c$, which verifies, for all $Y$ and $c$ :
\[\Theta\circ(\mathfrak a\underset{N}{_b*_\alpha}id)I_{\overline{\mathfrak a}}(\overline{\mathfrak a}(Y))=\tilde{\mathfrak a}(Y)\]
\[\Theta\circ(\mathfrak a\underset{N}{_b*_\alpha}id)I_{\overline{\mathfrak a}}(1\underset{B^o}{_{id}\otimes_{\hat{\jmath}\circ\mathfrak b}}(1\underset{N}{_b\otimes_\alpha}1\underset{N^o}{_{\hat{\alpha}}\otimes_\beta}c))=1\underset{N}{_b\otimes_\alpha}1\underset{N^o}{_{\hat{\alpha}}\otimes_\beta}c\]
On the other hand, let us put, for all $n\in N$, $\tilde{b}(n)=J_{\tilde{\chi}}\hat{\alpha}(n^*)J_{\tilde{\chi}}$. We have, by definition of $U_\chi^\mathfrak b$, for all $y\in B$ :
\begin{eqnarray*}
\hat{\jmath}(\mathfrak b(y))
&=&
J_{\tilde{\chi}}\mathfrak b(y^*)J_{\tilde{\chi}}\\
&=&
J_{\tilde{\chi}}U_\chi^\mathfrak b (y^*\underset{\nu}{_{\tilde{b}}\otimes_{\hat{\alpha}}}1)(U_\chi^\mathfrak b)^*J_{\tilde{\chi}}\\
&=& J_\chi y^*J_\chi\underset{\nu^o}{_{\hat{\alpha}}\otimes_\beta}1\\
&=& s(y)\underset{\nu^o}{_{\hat{\alpha}}\otimes_\beta}1
\end{eqnarray*}
where $s(y)=J_\chi y^*J_\chi$. 
Therefore, the identification of $\Lambda_{\tilde{\psi}}(x)\underset{\chi^o}{_{id}\otimes_s}\Lambda_\chi (y)$ with $\Lambda_{\tilde{\psi}}(yx)$ (\ref{rel}) gives the definition of $\mathcal U$. 
\newline
Let us write $\widetilde{(\tilde{\psi})}$ for the dual weight of $\tilde{\psi}$ on $(A\rtimes_\mathfrak a\gG)\rtimes_{\overline{\mathfrak a}}\gG(\tilde{\mathfrak a})$. Using \ref{GNSdual} applied to the weight $\tilde{\chi}$, then to the weight $\widetilde{(\tilde{\psi})}$, we get that :
\begin{eqnarray*}
\Lambda_{\tilde{\psi}}(x)\underset{\chi^o}{_{id}\otimes _{\hat{\jmath}\circ b}}(\Lambda_\chi(y)\underset{\nu^o}{_{\hat{\alpha}}\otimes_\beta}\Lambda_{\Phi^{oc}}(b))
&=&
\Lambda_{\tilde{\psi}}(x)\underset{\chi^o}{_{id}\otimes _{\hat{\jmath}\circ b}}\Lambda_{\tilde{\chi}}((1\underset{N}{_b\otimes_\alpha}1\underset{N^o}{_{\hat{\alpha}}\otimes_\beta}b)\mathfrak b(y))\\
&=&
\Lambda_{\widetilde{(\tilde{\psi})}}[(1\underset{B^o}{_{id}\otimes_{\hat{\jmath}\circ \mathfrak b}}(1\underset{N}{_b\otimes_\alpha}1\underset{N^o}{_{\hat{\alpha}}\otimes_\beta}b)\mathfrak b(y))\overline{\mathfrak a}(x)]
\end{eqnarray*}
As we have $\overline{\mathfrak a}(y)=1\underset{B^o}{_{id}\otimes_{\hat{\jmath}\circ \mathfrak b}}\mathfrak b(y)$, we finally get, using again \ref{GNSdual} applied to the bidual weight $\tilde{\tilde{\psi}}$ :
\[\mathcal U\Lambda_{\widetilde{(\tilde{\psi})}}[(1\underset{B^o}{_{id}\otimes_{\hat{\jmath}\circ \mathfrak b}}(1\underset{N}{_b\otimes_\alpha}1\underset{N^o}{_{\hat{\alpha}}\otimes_\beta}b))\overline{\mathfrak a}(yx)]=
\Lambda_{\tilde{\tilde{\psi}}}((1\underset{N}{_b\otimes_\alpha}1\underset{N^o}{_{\hat{\alpha}}\otimes_\beta}b)\tilde{\mathfrak a}(yx))\]
which we can write, using the isomorphism $\mathcal I=\Theta\circ(\mathfrak a\underset{N}{_b*_\alpha}id)I_{\overline{\mathfrak a}}$ :
\[\mathcal U\Lambda_{\widetilde{(\tilde{\psi})}}[(1\underset{B^o}{_{id}\otimes_{\hat{\jmath}\circ \mathfrak b}}(1\underset{N}{_b\otimes_\alpha}1\underset{N^o}{_{\hat{\alpha}}\otimes_\beta}b))\overline{\mathfrak a}(yx)]=
\Lambda_{\tilde{\tilde{\psi}}}(\mathcal I[(1\underset{B^o}{_{id}\otimes_{\hat{\jmath}\circ \mathfrak b}}(1\underset{N}{_b\otimes_\alpha}1\underset{N^o}{_{\hat{\alpha}}\otimes_\beta}b))\overline{\mathfrak a}(yx)])\]
Using \ref{thdual2}(ii) applied to $\widetilde{(\tilde{\psi})}$ and to $\tilde{\tilde{\psi}}$, we see that these elements are a core for, respectively, $\Lambda_{\widetilde{(\tilde{\psi})}}$ and $\Lambda_{\tilde{\tilde{\psi}}}$, and we finaly get that $\mathcal U\Lambda_{\widetilde{(\tilde{\psi})}}(Z)=\Lambda_{\tilde{\tilde{\psi}}}(\mathcal I(Z))$, for all $Z\in\gN_{\widetilde{(\tilde{\psi})}}$. From which one gets that :
\[\mathcal UX\mathcal U^*=\Theta\circ(\mathfrak a\underset{N}{_b*_\alpha}id)I_{\overline{\mathfrak a}}(X)\]
which finishes the proof. \end{proof}

%%%%cor1mqga
\subsection{Corollary}
\label{cor1mqga}
{\it Let us take the notations of \ref{th1mqga}, \ref{th2mqga} and \ref{th3mqga}. If we make the identification of $H_{\tilde{\psi}}\underset{B^o}{_{id}\otimes_{\hat{\jmath}\circ\mathfrak b}}H_{\tilde{\chi}}$ with $H_{\tilde{\tilde{\psi}}}$ by writing, for any $x\in\gN_{\tilde{\psi}}$, $y\in\gN_\chi$, $b\in\gN_{\Phi^{oc}}\cap\gN_{T^{oc}}$:
\[\Lambda_{\tilde{\psi}}(x)\underset{\chi^o}{_{id}\otimes _{\hat{\jmath}\circ b}}(\Lambda_\chi(y)\underset{\nu^o}{_{\hat{\alpha}}\otimes_\beta}\Lambda_{\Phi^{oc}}(b))=
\Lambda_{\tilde{\psi}}(yx)\underset{\nu^o}{_{\hat{\alpha}}\otimes_\beta}\Lambda_{\Phi^{oc}}(b)\] 
we then identify $\tilde{\mathfrak a}$ with $\overline{\mathfrak a}$, and, therefore, $(id, \tilde{\mathfrak a})$ (where $id$ is the inclusion of $B^o$ into $A\rtimes_\mathfrak a\gG$) is an action of $\gG(\tilde{\mathfrak a})$ on $A\rtimes_\mathfrak a\gG$. }

\begin{proof} Clear by \ref{th4mqga}. \end{proof}

%%%%cor2mqga
\subsection{Corollary}
\label{cor2mqga}
{\it Let us use the notations of \ref{th1mqga}, \ref{th2mqga}, \ref{th3mqga}. Then, the dual action $(j_1\circ b, \tilde{\overline{\mathfrak a}})$ of $\widehat{\gG(\tilde{\mathfrak a})}^c=\gG(\mathfrak a)$ on $(A\rtimes_\mathfrak a\gG)\rtimes_{\tilde{\mathfrak a}}\widehat{\gG}^c$ can be, as well, be identified with  the action $(1\underset{N}{_b\otimes_\alpha}1\underset{N^o}{_{\hat{\alpha}}\otimes_{\beta}}\hat{\beta}, \tilde{\tilde{\mathfrak a}})$ of $\gG^{oc}$, and the action $(1\underset{N}{_b\otimes_\alpha}\hat{\beta}, \underline{\mathfrak a})$ of $\gG$ on $A\underset{N}{_b*_\alpha}\mathcal L(H_\Phi)$ can be considered as well as an action $(j, \underline{\mathfrak a})$ of $\gG(\mathfrak a)^{oc}$. }

\begin{proof} Left to the reader. \end{proof}

%%%%prop1mqga
\subsection{Proposition}
\label{prop1mqga}
{\it Let us use the notations introduced in \ref{th1mqga} and \ref{th4mqga}; let us define, for all $x\in M'$, $\mu(x)=1\underset{N}{_b\otimes_\alpha}x$. We define this way an injective $*$-homomorphism from $M'$ into $\tilde{M}$. Moreover, we have :
\newline
(i) for all $x\in M'$, $j_1(\mu(x))=\mu(R^c(x))$; in particular, for all $n\in N$, $\mu(\hat{\alpha}(n))$ belongs to $\tilde{N}$ and $j_1(\mu(\hat{\alpha}(n)))=1\underset{N}{_b\otimes_\alpha}\hat{\beta}(n)$. 
\newline
(ii) for all positive $x\in M'^+$, we have $T_{\underline{\mathfrak a}}(\mu(x))=\mu(T^{oc}(x))$.
\newline
(iii) for all $x\in M'$, we have :}
\[\widetilde{\Gamma}(\mu(x))\simeq (\mu*\mu)\Gamma^{oc}(x)\]

\begin{proof}
We have, using successively \ref{Upsi1}, \ref{thL4}(v) and \ref{Upsi1}(ii) :
\begin{eqnarray*}
j_1(1\underset{N}{_b\otimes_\alpha}x)
&=&J_{\tilde{\psi}}(1\underset{N}{_b\otimes_\alpha}x^*)J_{\tilde{\psi}}\\
&=&U_\psi^\mathfrak a(J_\psi\underset{N}{_b\otimes_\alpha}J_{\widehat{\Phi}})(1\underset{N}{_b\otimes_\alpha}x^*)(J_\psi\underset{N^o}{_a\otimes_\beta}J_{\widehat{\Phi}})(U_\psi^\mathfrak a)^*\\
&=&U_\psi^\mathfrak a(1\underset{N^o}{_a\otimes_\beta}R^c(x))(U_\psi^\mathfrak a)^*\\
&=&1\underset{N}{_b\otimes_\alpha}R^c(x)
\end{eqnarray*}
which gives (i). Result (ii) is just a rewriting of \ref{threg}(i). 
\newline
Let us apply now \ref{th4mqga} and \ref{cor2mqga} to $1\underset{N}{_b\otimes_\alpha} y$, with $y\in \widehat{M}'$; we get :
\begin{eqnarray*}
1\underset{N}{_b\otimes_\alpha}1\underset{N^o}{_{\hat{\alpha}}\otimes_\beta}\Gamma^{oc}(y)
&=&
\tilde{\tilde{\mathfrak a}}(1\underset{N}{_b\otimes_\alpha}1\underset{N^o}{_{\hat{\alpha}}\otimes_\beta}y)\\
&=&\tilde{\overline{\mathfrak a}}(1\underset{N}{_b\otimes_\alpha}1\underset{N^o}{_{\hat{\alpha}}\otimes_\beta}y)\\
&=&1\underset{N}{_b\otimes_\alpha}1\underset{B}{_{id}\otimes_j}\tilde{\Gamma}(1\underset{N}{_b\otimes_\alpha}y)
\end{eqnarray*}
which, up to the identifications, gives the result. \end{proof}

%%%lcqgaction
\subsection{Example}
\label{lcqgact}
Let ${\bf G}$ be a locally compact quantum group,and let $\mathfrak a$ be an action of ${\bf G}$ on a von Neumann algebra $A$; by applying \ref{lcqgaction} and \ref{cor2mqga} to $\mathfrak a$, we can construct a measured quantum groupoid $\gG(\mathfrak a)$, whose underlying von Neumann algebra is the relative commutant $A\otimes\mathcal L(H_\Phi)\cap \mathfrak a (A)'$, the basis being $B=A\rtimes_\mathfrak a{\bf G}\cap \mathfrak a (A)'$;  it can be described also by looking at the restriction $\mathfrak b$ of the dual action $\tilde{\mathfrak a}$ of $\widehat{\bf G}^c$ to $B$ and looking at the inclusion of $\mathfrak b(B)$ into the crossed product $B\rtimes_\mathfrak b\widehat{\bf G}^c$, which is isomorphic to the inclusion of $B$ into $A\otimes\mathcal L(H_\Phi)\cap \mathfrak a (A)'$. This example will be studied later on in another article, and we need to compare this example with Vainerman's construction made in \cite{Va}. 

%%%%%outer
\subsection{Remark}
\label{outer}
Let us recall that an action $\mathfrak a$ of a locally compact quantum group ${\bf G}$ on a von Neumann algebra $A$ is called outer when $A\rtimes_\mathfrak a{\bf G}\cap \mathfrak a (A)'=\mathbb{C}$ (which means that the inclusion $\mathfrak a(A)\subset A\rtimes_\mathfrak a{\bf G}$ is irreducible). It was proven in (\cite{EN}, \cite{E2}) that any depth 2 irreducible inclusion of factors $M_0\subset M_1$ generates a locally compact quantum group ${\bf G}$, and an outer action $\mathfrak a$ of $\widehat{\bf G}^c$ on $M_1$, such that $M_1^{\mathfrak a}=M_0$, and that the crossed product $M_1\rtimes_\mathfrak a\widehat{\bf G}^c$ is isomorphic to the basic construction made from the inclusion $M_0\subset M_1$; conversely, Vaes proved (\cite{V2}) that, for any action $\mathfrak a$ of a locally compact quantum group ${\bf G}$ on a von Neumann algebra $A$, the inclusion $\mathfrak a(A)\subset A\rtimes_\mathfrak a{\bf G}$ is depth 2. 
\newline
Here, we have constructed from this depth 2 inclusion, a measured quantum groupoid $\gG (\mathfrak a)$; when the action $\mathfrak a$ is outer, this measured quantum groupoid is a locally compact quantum group (\cite{EN}, \cite{E2}), and we get (\ref{lcqgact}) that this locally compact quantum group is $\bf G^{oc}$. This result was certainly known by specialists, but we were not able to find a proof anywhere. 
\newline
As Vaes (\cite{V3}) proved that any locally compact quantum group has an outer action (on some type $III_1$ factor), this proves that any locally compact quantum group comes from an irreducible depth 2 inclusion of factors, by the construction made in (\cite{EN}, \cite{E2}). 

\newpage

%%%%appendice

\begin{appendix}
\centerline{\bf Appendix}
%%%%%%antipod
\section{Coinverse and scaling group of a measured quantum groupoid}
\label{coinverse}
In this chapter, we are dealing with a Hopf-bimodule $(N, \alpha, \beta, M, \Gamma)$, equipped with a left-invariant operator-valued weight $T$, and a right-invariant operator-valued weight $T'$. If $\nu$ denotes a normal semi-finite faithful weight on the basis, let $\Phi$ (resp. $\Psi$) be the lifted normal faithful semi-finite weight on $M$ by $T$ (resp. $T'$). Then, with the additional hypothesis that the two modular automorphism groups associated to the two weight $\Phi$ and $\Psi$ commute (we then say that $\nu$ is relatively invariant with respect to $T$ and $T'$ (\ref{defMQG}), we can construct a co-inverse, a scaling group and an antipod, using slight generalizations of the constructions made in (\cite{L2},9) for "adapted measured quantum groupoids".

%%%%%lemE2
\subsection{Lemma}
\label{lemE2}
{\it Let $(N, \alpha, \beta, M, \Gamma)$ be a Hopf-bimodule; let us suppose that there exist a left-invariant operator-valued weight $T$, a right-invariant valued weight $T'$ and $\nu$ a normal semi-finite faithful weight on $N$, relatively invariant with respect to $T$ and $T'$ in the sense of \ref{defMQG}; we shall denote $\Phi=\nu\circ\alpha^{-1}\circ T$ and $\Psi=\nu\circ\beta^{-1}\circ T'$ the two lifted normal semi-finite weights on $M$. Let us denote $\delta$ the modulus of $\Psi$ with respect to $\Phi$ and $\lambda$ the scaling operator (\ref{Vaes}). We shall use the notations of \ref{propbasic}. Then :
\newline
(i) let $x\in \mathcal T_{\Psi, T'}$ and $n\in\mathbb{N}$ and $y=e_nx$, with the notations of \ref{Vaes}; then $y$ belongs to $\gN_\Psi\cap\gN_{T'}$, is analytical with respect to $\Psi$, and the operator $\sigma_{-i/2}^\Psi(y^*)\delta^{1/2}$ is bounded, and its closure $\overline{\sigma_{-i/2}^\Psi(y^*)\delta^{1/2}}$ belongs to $\gN_{\Phi}$; moreover, with the identifications made in \ref{Vaes}, we have :
\[\Lambda_\Phi(\overline{\sigma_{-i/2}^\Psi(y^*)\delta^{1/2}})=J_\Psi\Lambda_\Psi(y)\]
(ii) let $E$ be the linear space generated by all such elements of the form $\overline{\sigma_{-i/2}^\Psi(y^*)\delta^{1/2}}$, for all $x\in \mathcal T_{\Psi, T'}$ and $n\in\mathbb{N}$; then $E$ is a weakly dense subspace of $\gN_\Phi$, and, for all $z\in E$, $\Lambda_\Phi(z)\in D((H_\Phi)_\beta, \nu^o)$;
\newline
(iii) the linear set of all products $<\Lambda_\Phi(z), \Lambda_\Phi(z')>_{\beta, \nu^o}$ (for $z$, $z'$ in $E$) is a dense subspace of $N$. }

\begin{proof}
As $e_n$ is analytical with respect to $\Psi$, $y$ belongs to $\gN_\Psi\cap\gN_{T_R}$, is analytical with respect to $\Psi$, and $\sigma_{-i/2}^\Psi(y^*)\delta^{1/2}$
is bounded (\cite{V1}, 1.2); as $\delta^{-1}$ is the modulus of $\Phi$ with respect to $\Psi$, we get that $\overline{\sigma_{-i/2}^\Psi(y^*)\delta^{1/2}}$ belongs to $\gN_{\Phi}$; we identify $\Lambda_\Phi(\overline{\sigma_{-i/2}^\Psi (y^*)\delta^{1/2}})$ with $\Lambda_\Psi(\sigma_{-i/2}^\Psi(y^*))=J_\Psi\Lambda_\Psi(y)$, which is (i). 
\newline
The subspace $E$ contains all elements of the form $\sigma_{-i/2}^\Psi (x^*)\overline{\delta^{1/2}\sigma_{-i/2}^\Psi (e_n)}$ ($x\in\mathcal T_{\Psi, T'}$), and, by density of $\mathcal T_{\Psi, T'}$ in $M$, we get that the closure of $E$ contains all elements of the form $a\overline{e_n\delta^{-1/2}}\overline{\delta^{1/2}\sigma_{-i/2}^\Psi (e_n)}=ae_n\sigma_{-i/2}^\Psi(e_n)$, for all $a\in M$; now, as $e_n\sigma_{-i/2}^\Psi(e_n)$ is converging to $1$, we finally get that $E$ is dense in $M$; as $\Lambda_\Phi(E)\subset J_\Psi\Lambda_\Psi(\gN_\psi\cap\gN_{T'})$, we get, by \ref{basic}, that, for all $z$ in $E$, $\Lambda_\Phi(z)$ belongs to $D((H_\Phi)_\beta, \nu^o)$; more precisely, we have :
\[R^{\beta, \nu^o}(\Lambda_\Phi(\sigma_{-i/2}^\Psi (x^*)\overline{\delta^{1/2}\sigma_{-i/2}^\Psi (e_n)}))=R^{\beta, \nu^o}(J_\psi\Lambda_\psi(e_nx))=\Lambda_{T'}(e_nx)\]
Therefore, the set of elements of the form $<\Lambda_\Phi(z), \Lambda_\Phi(z')>_{\beta, \nu^o}$ contains all elements of the form $\beta^{-1}\circ T'(x^*e_ne_nx)$, for all $x$ in $\mathcal T_{\Psi, T'}$ and $n\in \mathbb{N}$; as we have :
\[T'(x^*e_ne_nx)=\Lambda_{T'}(e_nx)^*\Lambda_{T'}(e_nx)=\Lambda_{T'}(x)^*e_n^*e_n\Lambda_{T'}(x)\]
we get that its closure contains all elements of the form $\beta^{-1}\circ T'(x^*x)$, and, therefore, it contains $\beta^{-1}\circ T'(\gM_{T'}^+)$, which finishes the proof. \end{proof}

%%%%DefS
\subsection{Definition}
\label{defS}
As in (\cite{L2}, 9.2), we can define, for all $\lambda\in\mathbb{C}$, a closed operator $\Delta_\Phi^\lambda\underset{N^o}{_\alpha\otimes_{\hat{\beta}}}\Delta_\Phi^\lambda$, with natural values on elementary tensor products; it is possible also to define a unitary antilinear operator $J_\Phi\underset{N^o}{_\alpha\otimes_{\hat{\beta}}}J_\Phi$ from $H_\Phi\underset{N^o}{_\alpha\otimes_{\hat{\beta}}}H_\Phi$ onto $H_\Phi\underset{N}{_{\hat{\beta}}\otimes_\alpha}H_\Phi$ (whose inverse will be $J_\Phi \underset{N}{_{\hat{\beta}}\otimes_\alpha}J_\Phi$); by composition, we define then a closed antilinear operator $S_\Phi\underset{N^o}{_\alpha\otimes_{\hat{\beta}}}S_\Phi$, with natural values on elementary tensor products, whose adjoint will be $F_\Phi \underset{N}{_{\hat{\beta}}\otimes_\alpha}F_\Phi$. 

%%%%%%PropS
\subsection{Proposition}
\label{PropS}
{\it For all $a$, $c$ in $(\gN_\Phi\cap\gN_{T})^*(\gN_\Psi\cap\gN_{T'})$, $b$, $d$ in $\mathcal T_{\Psi, T'}$ and $g$, $h$ in $E$, the following vector :
\[U^*_{H_\Phi}\Gamma(g^*)[\Lambda_\Phi(h)\underset{\nu}{_\beta\otimes_\alpha}(\lambda^{\beta, \alpha}_{\Lambda_\Psi(\sigma_{-i}^\Psi(b^*))})^*U_{H_\Psi}(\Lambda_\Psi(a)\underset{\nu^o}{_\alpha\otimes_{\hat{\beta}}}\Lambda_\Phi((cd)^*))]\]
belongs to $D(S_\Phi\underset{\nu^*}{_\alpha\otimes_{\hat{\beta}}}S_\Phi)$, and the value of $\sigma_\nu(S_\Phi\underset{\nu^*}{_\alpha\otimes_{\hat{\beta}}}S_\Phi)$ on this vector is equal to :}
\[U^*_{H_\Phi}\Gamma(h^*)[\Lambda_\Phi(g)\underset{\nu}{_\beta\otimes_\alpha}(\lambda^{\beta, \alpha}_{\Lambda_\Psi(\sigma_{-i}^\Psi(d^*))})^*U_{H_\Psi}(\Lambda_\Psi(c)\underset{\nu^o}{_\alpha\otimes_{\hat{\beta}}}\Lambda_\Phi((ab)^*))]\]

\begin{proof}
The proof is identical to (\cite{L2},9.9), thanks to \ref{lemE2}(ii). \end{proof}

%%%%%%%propG
\subsection{Proposition}
\label{propG}
{\it There exists a closed densely defined anti-linear operator $G$ on $H_\Phi$ such that the linear span of :
\[(\lambda^{\beta, \alpha}_{\Lambda_\Psi(\sigma_{-i}^\Psi(b^*))})^*U_{H_\Psi}(\Lambda_\Psi(a)\underset{\nu^o}{_\alpha\otimes_{\hat{\beta}}}\Lambda_\Phi((cd)^*))\]
with $a$, $c$ in $(\gN_\Phi\cap\gN_{T})^*(\gN_\Psi\cap\gN_{T'})$, $b$, $d$ in $\mathcal T_{\Psi, T'}$, is a core of $G$, and we have :}
\[G[(\lambda^{\beta, \alpha}_{\Lambda_\Psi(\sigma_{-i}^\Psi(b^*))})^*U_{H_\Psi}(\Lambda_\Psi(a)\underset{\nu^o}{_\alpha\otimes_{\hat{\beta}}}\Lambda_\Phi((cd)^*))]=
(\lambda^{\beta, \alpha}_{\Lambda_\Psi(\sigma_{-i}^\Psi(d^*))})^*U_{H_\Psi}(\Lambda_\Psi(c)\underset{\nu^o}{_\alpha\otimes_{\hat{\beta}}}\Lambda_\Phi((ab)^*))\]

\begin{proof}
The proof is identical to (\cite{L2},9.10), thanks to \ref{lemE2}(iii). \end{proof}

%%%%%%%%tau
\subsection{Theorem}
\label{tau}
{\it Let $(N, \alpha, \beta, M, \Gamma)$ be a Hopf-bimodule; let us suppose that there exist a left-invariant operator-valued weight $T$, a right-invariant valued weight $T'$ and $\nu$ a normal semi-finite faithful weight on $N$, relatively invariant with respect to $T$ and $T'$ in the sense of \ref{defMQG}; we shall denote $\Phi=\nu\circ\alpha^{-1}\circ T$ and $\Psi=\nu\circ\beta^{-1}\circ T'$ the two lifted normal semi-finite weights on $M$. Let $G$ be the closed densely defined antilinear operator defined in \ref{propG}, and let $G=ID^{1/2}$ its polar decomposition. Then, the operator $D$ is positive self-adjoint and non singular; there exists a one-parameter automorphism group $\tau_t$ on $M$ defined, for $x\in M$, by :
\[\tau_t(x)=D^{-it}xD^{it}\]
We have, for all $n\in N$ and $t\in\mathbb{R}$ :
\[\tau_t(\alpha(n))=\alpha(\sigma^\nu_t(n))\]
\[\tau_t(\beta(n))=\beta(\sigma^\nu_t(n))\]
which allows us to define $\tau_t\underset{N}{_\beta*_\alpha}\tau_t$, $\tau_t\underset{N}{_\beta*_\alpha}\sigma_t^\Phi$ and $\sigma_t^\Psi\underset{N}{_\beta*_\alpha}\tau_{-t}$ on $M\underset{N}{_\beta*_\alpha}M$; moreover, we have :}
\[\Gamma\circ\tau_t=(\tau_t\underset{N}{_\beta*_\alpha}\tau_t)\Gamma\]
\[\Gamma\circ\sigma^\Phi_t=(\tau_t\underset{N}{_\beta*_\alpha}\sigma_t^\Phi)\Gamma\]
\[\Gamma\circ\sigma^\Psi_t=(\sigma_t^\Psi\underset{N}{_\beta*_\alpha}\tau_{-t})\Gamma\]

\begin{proof}
The proof is identical to \cite{L2}, 9.12 to 9.28. \end{proof}

%%%%%%R
\subsection{Theorem}
\label{R}
{\it Let $(N, \alpha, \beta, M, \Gamma)$ be a Hopf-bimodule; let us suppose that there exists a left-invariant operator-valued weight $T$, a right-invariant valued weight $T'$ and $\nu$ a normal semi-finite faithful weight on $N$, relatively invariant with respect to $T$ and $T'$ in the sense of \ref{defMQG}; we shall denote $\Phi=\nu\circ\alpha^{-1}\circ T$ and $\Psi=\nu\circ\beta^{-1}\circ T'$ the two lifted normal semi-finite weights on $M$. Let $G$ be the closed densely defined antilinear operator defined in \ref{propG}, and let $G=ID^{1/2}$ its polar decomposition. Then, the operator $I$ is antilinear, isometric, surjective, and we have $I=I^*=I^2$; there exists a $*$-antiautomorphism $R$ on $M$ defined, for $x\in M$, by :
\[R(x)=Ix^*I\]
such that, for all $t\in\mathbb{R}$, we get $R\circ\tau_t=\tau_t\circ R$ and $R^2=id$. 
\newline
For any $a$, $b$ in $\gN_\Psi\cap\gN_{T'}$ we have :
\[R((\omega_{J_\Psi\Lambda_\Psi(a)}\underset{N}{_\beta*_\alpha}id)\Gamma(b^*b))=
(\omega_{J_\Psi\Lambda_\Psi(b)}\underset{N}{_\beta*_\alpha}id)\Gamma(a^*a)\]
and for any $c$, $d$ in $\gN_\Phi\cap\gN_{T}$, we have :
\[R((id\underset{N}{_\beta*_\alpha}\omega_{J_\Phi\Lambda_\Phi(c)})\Gamma(d^*d))=
(id\underset{N}{_\beta*_\alpha}\omega_{J_\Phi\Lambda_\Phi(d)})\Gamma(c^*c))\]
\newline
For all $n\in N$, we have $R(\alpha(n))=\beta(n)$, which allows us to define $R\underset{N}{_\beta*_\alpha}R$ from $M\underset{N}{_\beta*_\alpha}M$ onto $M\underset{N^o}{_\alpha*_\beta}M$ (whose inverse will be $R\underset{N^o}{_\alpha*_\beta}R$), and we have :}
\[\Gamma\circ R=\varsigma_{N^o}(R\underset{N}{_\beta*_\alpha}R)\Gamma\]

\begin{proof}
The proof is identical to \cite{L2}, 9.38 to 9.42. \end{proof}

%%%%%dens
\subsection{Theorem}
\label{dens}
{\it Let $(N, \alpha, \beta, M, \Gamma)$ be a Hopf-bimodule; let us suppose that there exists a left-invariant operator-valued weight $T$, a right-invariant valued weight $T'$ and $\nu$ a normal semi-finite faithful weight on $N$, relatively invariant with respect to $T$ and $T'$ in the sense of \ref{defMQG}; we shall denote $\Phi=\nu\circ\alpha^{-1}\circ T_L$; then :
\newline
(i) $M$ is the weak closure of the linear span of all elements of the form $(\omega\underset{N}{_\beta*_\alpha}id)\Gamma(x)$, for all $x\in M$ and $\omega\in M_*$ such that there exists $k>0$ with $\omega\circ\beta\leq k\nu$. 
\newline
(ii) $M$ is the weak closure of the linear span of all elements of the form $(id\underset{N}{_\beta*_\alpha}\omega)\Gamma(x)$, for all $x\in M$ and $\omega\in M_*$ such that there exists $k>0$ with $\omega\circ\alpha\leq k\nu$. 
\newline
(iii) $M$ is the weak closure of the linear span of all elements of the form $(id*\omega_{v, w})(W)$, where $v$ belongs to $D(_\alpha H_\Phi, \nu)$ and $w$ belongs to $D((H_\Phi)_{\hat{\beta}}, \nu^o)$. }

\begin{proof}
The proof is identical to \cite{L2}, 9.25. \end{proof}

%%%%%def
\subsection{Definition}
\label{def}
Let $(N, \alpha, \beta, M, \Gamma)$ be a Hopf-bimodule; let us suppose that there exists a left-invariant operator-valued weight $T$, a right-invariant valued weight $T'$ and $\nu$ a normal semi-finite faithful weight on $N$, relatively invariant with respect to $T$ and $T'$ in the sense of \ref{defMQG}; we shall denote $\Phi=\nu\circ\alpha^{-1}\circ T$ and $\Psi=\nu\circ\beta^{-1}\circ T'$ the two lifted normal semi-finite weights on $M$; let $\tau_t$ the one-parameter automorphism group constructed in \ref{tau} and let $R$ be the involutive $*$-antiautomorphism constructed in \ref{R}. We shall call $\tau_t$ the scaling group of $(N, \alpha, \beta, M, \Gamma, T, T', \nu)$ and $R$ the coinverse of $(N, \alpha, \beta, M, \Gamma, T, T', \nu)$. Thanks to \ref{R} and \ref{dens}, we see that, $T$ and $\nu$ being given, $R$ does not depend on the choice of the right-invariant operator-valued weight $T'$. 
\newline
Similarly, from \ref{tau}, one gets that, for all $x$ in $M$, $\omega\in M_*$ such that there exists $k>0$ with $\omega\circ\alpha\leq k\nu$, $\omega'\in M_*$ such that there exists $k>0$ with $\omega\circ\beta\leq k\nu$, one has :
\[\tau_t((id\underset{N}{_\beta*_\alpha}\omega)\Gamma(x))=(id\underset{N}{_\beta*_\alpha}\omega\circ\sigma_{-t}^\Phi)\Gamma\sigma^\Phi_t(x)\]
\[\tau_t((\omega'\underset{N}{_\beta*_\alpha}id)\Gamma(x))=(\omega'\circ\sigma_{t}^\Psi\underset{N}{_\beta*_\alpha}id)\Gamma\sigma_{-t}^\Psi(x)\]
So, $T$ and $\nu$ being given, $\tau_t$ does not depend on the choice of the right-invariant operator-valued weight $T'$. 

%%%thS
\subsection{Theorem}
\label{thS}
{\it Let $(N, \alpha, \beta, M, \Gamma)$ be a Hopf-bimodule; let us suppose that there exists a left-invariant operator-valued weight $T$, a right-invariant valued weight $T'$ and $\nu$ a normal semi-finite faithful weight on $N$, relatively invariant with respect to $T$ and $T'$ in the sense of \ref{defMQG}; we shall denote $\Phi=\nu\circ\alpha^{-1}\circ T$; then, for any $\xi$, $\eta$ in $D(_\alpha H_\Phi, \nu)\cap D((H_\Phi)_{\hat{\beta}}, \nu^o)$, $(id*\omega_{\xi, \eta})(W)$ belongs to $D(\tau_{i/2})$, and, if we define $S=R\tau_{i/2}$, we have :
\[S((id*\omega_{\xi, \eta})(W))=(id*\omega_{\eta, \xi})(W)^*\]
More generally, for any $x$ in $D(S)=D(\tau_{i/2})$, we get that $S(x)^*$ belongs to $D(S)$ and $S(S(x)^*)^*=x$; 
$S$ will be called the antipod of the measured quantum groupoid, and, therefore, the co-inverse and the scaling group, given by polar decomposition of the antipod, rely only upon the pseudo-multiplicative $W$. }

\begin{proof}
It is proved similarly to \cite{L2} 9.35 and 9.36. \end{proof}

%%%%RTR
\subsection{Proposition}
\label{RTR}
{\it Let $(N, \alpha, \beta, M, \Gamma)$ be a Hopf-bimodule, equipped with a left-invariant operator-valued weight $T$, and a right-invariant valued weight $T'$; let $\nu$ be a normal semi-finite faithful weight on $N$, relatively invariant with respect to $T$ and $T'$ in the sense of \ref{defMQG}; let $\tau_t$ be the scaling group of $(N, \alpha, \beta, M, \Gamma, T, T', \nu)$ and $R$ the coinverse of $(N, \alpha, \beta, M, \Gamma, T_L, T_R, \nu)$; then :
\newline
(i) the operator-valued weight $RT'R$ is left-invariant, the operator valued-weight $RTR$ is right-invariant, and $\nu$ is relatively invariant with respect to $RT'R$ and $RTR$. 
\newline
(ii) $\tau_t$ is the scaling group of $(N, \alpha, \beta, M, \Gamma, RT'R, RTR, \nu)$ }

\begin{proof}
Let $\Phi=\nu\circ\alpha^{-1}\circ T$ and $\Psi=\nu\circ\beta^{-1}\circ T'$ the two lifted normal semi-finite weights on $M$ by $T$ and $T'$;  the lifted weight by $RT'R$ (resp. $RTR$) is then $\Psi\circ R$ (resp. $\Phi\circ R$). As $\sigma_t^{\Psi\circ R}=R\circ\sigma_{-t}^\Psi\circ R$ and $\sigma_s^{\Phi\circ R}=R\circ\sigma_{-s}^\Phi\circ R$, we get that $\sigma^{\Psi\circ R}$ and $\sigma^{\Phi\circ R}$ commute, which is (i). 
\newline
From \ref{tau} and \ref{R}, we get that :
\begin{multline*}
\Gamma\circ\sigma_t^{\Psi\circ R}=\Gamma\circ R\circ\sigma_{-t}^\Psi\circ R=
\varsigma_{N^o}(R\underset{N}{_\beta*_\alpha}R)\Gamma\circ\sigma_{-t}^\Psi\circ R\\
=\varsigma_{N^o}(R\circ\sigma_{-t}^\Psi\circ R\underset{N^o}{_\alpha*_\beta}R\circ\tau_t\circ R)\varsigma_{N}\Gamma=(\tau_t\underset{N}{_\beta*_\alpha}\sigma_t^{\Psi\circ R})\Gamma
\end{multline*}
from which we get that, for all $x\in M$ and $\omega\in M_*$ such that there exists $k>0$ such that $\omega\circ\alpha<k\nu$, we have :
\[\tau_t((id\underset{N}{_\beta*_\alpha}\omega)\Gamma(x))=(id\underset{N}{_\beta*_\alpha}\omega\circ\sigma_{-t}^{\Psi\circ R})\Gamma(\sigma_t^{\Psi\circ R}(x))\]
from which we get, by \ref{dens}, that $\tau_t$ is the scaling group associated to $RT_RR$, $RT_LR$ and $\nu$. \end{proof}

%%%%%QMG
\section{Automorphism groups on the basis of a measured quantum groupoid}
\label{AGB}
In this section, with the same hypothesis as in appendix \ref{coinverse}, we construct two one-parameter automorphism groups on the basis $N$ (\ref{gamma}), and we prove (\ref{thcentral}) that these automorphisms leave invariant the quasi-invariant weight $\nu$.  We prove also in \ref{thcentral} that the weight $\nu$ is also quasi-invariant with respect to $T$ and $RTR$. We finish by proving (\ref{th}) that our axioms and Lesieur's axioms are equivalent.

%%%%%lembeta
\subsection{Lemma}
\label{lembeta}
{\it Let $(N, \alpha, \beta, M, \Gamma)$ be a Hopf-bimodule, equipped with a left-invariant operator-valued weight $T$, and a right-invariant valued weight $T'$; let $\nu$ be a normal semi-finite faithful weight on $N$, relatively invariant with respect to $T$ and $T'$ in the sense of \ref{defMQG}. Let $x\in M\cap\alpha(N)'$ and $y\in M\cap\beta(N)'$. Then :
\newline
(i) $x$ belongs to $\beta(N)$ if and only if we have :
\[\Gamma (x)=1\underset{N}{_\beta\otimes_\alpha}x\]
(ii) $y$ belongs to $\alpha(N)$ if and only if we have :
\[\Gamma(y)=y\underset{N}{_\beta\otimes_\alpha}1\]
More generally, if $x_1$, $x_2$ are in $M\cap\alpha(N)'$ and such that $\Gamma (x_1)=1\underset{N}{_\beta\otimes_\alpha}x_2$, then $x_1=x_2\in\beta(N)$. }
\begin{proof}
The proof is given in \cite{L2}, 4.4. \end{proof}

%%%%%%gamma
\subsection{Proposition}
\label{gamma}
{\it Let $(N, \alpha, \beta, M, \Gamma)$ be a Hopf-bimodule; let us suppose that there exists a left-invariant operator-valued weight $T$, a right-invariant valued weight $T'$ and $\nu$ a normal semi-finite faithful weight on $N$, relatively invariant with respect to $T$ and $T'$ in the sense of \ref{defMQG}.  Then, there exists a unique one-parameter group of automorphisms $\gamma_t^{L}$ of $N$ such that, for all $t\in\mathbb{R}$ and $n\in N$, we have :
\[\sigma_t^{T}(\beta(n))=\beta(\gamma_t^{L}(n))\]
\[\sigma_t^{RTR}(\alpha (n))=\alpha(\gamma_{-t}^L(n))\]
Moreover, the automorphism groups $\gamma^{L}$ and $\sigma^\nu$ commute, and there exists a positive self-adjoint non-singular operator $h_{L}$ $\eta$ $Z(N)\cap N^{\gamma^L}$ such that, for any $x\in N^+$ and $t\in \mathbb{R}$, we have :
\[\nu\circ\gamma^L_t(x)=\nu(h_{L}^tx)\]
Using the operator-valued weights $RT'R$ and $RTR$, we obtain another one-parameter group of automorphisms $\gamma_t^R$ of $N$, such that we have :
\[\sigma_t^{RT'R}(\beta(n))=\beta(\gamma_t^R(n))\]
\[\sigma_t^{T'}(\alpha(n))=\alpha(\gamma_{-t}^R(n))\]
and a positive self-adjoint non-singular operator $h_R$ $\eta$ $Z(N)\cap N^{\gamma^R}$ such that we have :} 
\[\nu\circ\gamma^R_t(x)=\nu(h_R^tx)\]

\begin{proof}
The existence of $\gamma_t^{L}$ is given by \cite{L2}, 4.5; moreover, from the formula $\sigma^\Phi_t\circ\sigma^\Psi_s(\beta(n))=\sigma^\Psi_s\circ\sigma^\Phi_t(\beta(n))$, we obtain :
\[\beta(\gamma_t^{L}\circ\sigma_{-s}^\nu(n))=\beta(\sigma_{-s}^\nu\circ\gamma_t^{L}(n))\]
which gives the commutation of $\gamma_t^{L}$ and $\sigma_{-s}^\nu$. The existence of $h_{L}$ is then straightforward. The construction of $\gamma^R$ and $h_R$ is just the application of the preceeding results to $RT'R$, $RTR$ and $\nu$. \end{proof}

%%%%%proph
\subsection{Proposition}
\label{proph}
{\it Let $(N, \alpha, \beta, M, \Gamma)$ be a Hopf-bimodule; let us suppose that there exists a left-invariant operator-valued weight $T$, a right-invariant valued weight $T'$ and $\nu$ a normal semi-finite faithful weight on $N$, relatively invariant with respect to $T$ and $T'$ in the sense of \ref{defMQG}. Let $T_L$ (resp. $T_R$) be another left (resp. right)-invariant operator-valued weight; we shall denote $\Phi=\nu\circ\alpha^{-1}\circ T$, $\Phi'=\nu\circ\alpha^{-1}\circ T_L$, $\Psi=\nu\circ\beta^{-1}\circ T'$ and $\Psi'=\nu\circ\beta^{-1}\circ T_R$ the lifted normal semi-finite weights on $M$; then, we have :
\[\beta(h_{L}^{ist})=(D\Psi'\circ\sigma_t^\Phi:D\Psi'\circ\tau_t)_s\]
\[\alpha(h_R^{ist})=(D\Phi'\circ\sigma_{-t}^\Psi:D\Phi'\circ\tau_t)_s\]
where $\tau_s$ is the scaling group constructed from $T$, $T'$ and $\nu$ as well from $RT'R$, $RTR$ and $\nu$ (\ref{tau} and \ref{RTR}). }

\begin{proof}
From \ref{tau}, we get, for all $t\in\mathbb{R}$, $\Gamma\circ\sigma_t^\Phi\tau_{-t}=(id\underset{N}{_\beta*_\alpha}\sigma_t^\Phi\tau_{-t})\Gamma$, and, therefore, by the right-invariance of $T_R$, we get, for all $x\in\gM_{T_R}^+$, that $\tau_t\sigma_{-t}^\Phi T_R\sigma_t^\Phi\tau_{-t}(x)= T_R(x)$; let now $x\in\gM_{\Psi'}^+$; $T_R(x)$ is an element of the positive extended part of $\beta(N)$ which can be written :
\[\int_0^\infty \lambda de_\lambda +(1-p)\infty\]
where $p$ is a projection in $\beta(N)$, and $e_\lambda$ is a resolution of $p$. As $x$ belongs to $\gM_{\Psi'}^+$, it is well known that $p=1$, and $T_R(x)=\int_0^\infty \lambda de_\lambda$. There exists also a projection $q$ and a resolution of $q$ such that :
\[\tau_t\sigma_{-t}^\Phi T_R\sigma_t^\Phi\tau_{-t}(x)=\int_0^\infty \lambda df_\lambda +(1-q)\infty\]
and, for all $\mu\in \mathbb{R}^+$, we have, because $e_\mu x e_\mu$ belongs to $\gM_{T_R}^+$ :
\begin{eqnarray*}
e_\mu(\int_0^\infty \lambda df_\lambda)e_\mu +e_\mu(1-q)e_\mu\infty 
&=&e_\mu\tau_t\sigma_{-t}^\Phi T_R\sigma_t^\Phi\tau_{-t}(x)e_\mu\\
&=&\tau_t\sigma_{-t}^\Phi T_R\sigma_t^\Phi\tau_{-t}(e_\mu xe_\mu)\\
&=&T_R(e_\mu xe_\mu)\\
&=&\int_0^\mu \lambda de_\lambda
\end{eqnarray*}
from which we infer that $(1-q)e_\mu=0$, and, therefore, that $q=1$; then, we get that $e_\mu\tau_t\sigma_{-t}^\Phi T_R\sigma_t^\Phi\tau_{-t}(x)e_\mu$ is increasing with $\mu$ towards $T_R(x)$. Therefore, we get 
that :
\[\tau_t\sigma_{-t}^\Phi T_R\sigma_t^\Phi\tau_{-t}(x)\subset T_R(x)\]
and, finally, the equality, for all $x\in\gM_{\Psi'}^+$ :
\[\tau_t\sigma_{-t}^\Phi T_R\sigma_t^\Phi\tau_{-t}(x)=T_R(x)\]
Moreover, as we have, for all $n\in N$ :
\[\tau_t\sigma_{-t}^\Phi(\beta(n))=\beta(\sigma_{t}^\nu\gamma_{-t}^{L}(n))\]
we get, using \ref{gamma}, that, for all $x\in \gM_{\Psi'}^+$ :
\[\Psi'(\beta(h_{L}^{-t/2})\sigma_t^\Phi\tau_{-t}(x)\beta(h_L^{-t/2}))=\Psi'(x)\]
and, therefore, that, for all $x\in M^+$ :
 \[\Psi'(\beta(h_{L}^{-t/2})\sigma_t^\Phi\tau_{-t}(x)\beta(h_L^{-t/2}))\leq\Psi'(x)\]
 A similar calculation (with $\tau_t\sigma_{-t}^\Phi$ instead of $\sigma_t^\Phi\tau_{-t}$) leads to :
 \[\Psi'(\beta(h_L^{t/2})\tau_t\sigma_{-t}^\Phi(x)\beta(h_L^{t/2}))\leq\Psi'(x)\]
 which leads to the equality, from which we get the first result. 
\newline
Applying this result to $RT_RR$, $RT_LR$ and $\nu$, we get, using again \ref{RTR} :
\begin{eqnarray*}
\beta(h_R^{ist})&=&(D\Phi'\circ R\circ\sigma_t^{\Psi\circ R}:D\Phi'\circ R\circ\tau_t)_s\\
&=&(D\Phi'\circ\sigma_{-t}^\Psi\circ R: D\Phi'\circ\tau_t\circ R)_s\\
&=&R[((D\Phi'\circ\sigma_{-t}^\Psi : D\Phi'\circ\tau_t)_{-s})^*]
\end{eqnarray*}
which leads to the result. \end{proof}

%%%%%%corh
\subsection{Corollary}
\label{corh}
{\it Let $(N, \alpha, \beta, M, \Gamma)$ be a Hopf-bimodule; let us suppose that there exists a left-invariant operator-valued weight $T$, a right-invariant valued weight $T'$ and $\nu$ a normal semi-finite faithful weight on $N$, relatively invariant with respect to $T$ and $T'$ in the sense of \ref{defMQG}. We shall denote $\Phi=\nu\circ\alpha^{-1}\circ T$ and $\Psi=\nu\circ\beta^{-1}\circ T'$ the two lifted normal semi-finite weights on $M$, $R$ the coinverse and $\tau_t$ the scaling group constructed in \ref{R} and \ref{tau}; we shall denote $\lambda$ the scaling operator of $\Psi$ with respect to $\Phi$ (\ref{Vaes}), $h_L$ and $h_R$ the operators constructed in \ref{gamma}. Then, for all $s$, $t$ in $\mathbb{R}$ :
\newline
(i) $(D\Psi:D\Psi\circ\tau_t)_s=\lambda^{ist}\beta(h_{L}^{ist})$
\newline
(ii) $(D\Phi:D\Phi\circ\tau_t)_s=\lambda^{ist}\alpha(h_R^{ist})$
\newline
(iii) $(D\Phi:D\Phi\circ\sigma_{-t}^{\Phi\circ R})_s=\lambda^{ist}\alpha(h_R^{ist})\alpha(h_L^{-ist})$
\newline
(iv) $(D\Psi:D\Psi\circ\sigma_t^{\Psi\circ R})_s=\lambda^{ist}\beta(h_L^{ist})\beta(h_R^{-ist})$.}
\begin{proof}
Applying \ref{proph} with $T_R=T'$, as $(D\Psi\circ\sigma^\Phi_t:D\Psi)_s=\lambda^{-ist}$ (\ref{Vaes}), we obtain (i). Applying \ref{proph} with $T_L=T$, as $(D\Phi:D\Phi\circ\sigma_{-t}^\Psi)_s=\lambda^{ist}$, we obtain (ii). Applying \ref{proph} with $T_R=RTR$, we obtain :
\begin{eqnarray*}
\beta(h_L^{ist})&=&(D\Phi\circ R\circ\sigma_t^\Phi : D\Phi\circ R\circ \tau_t)_s\\
&=&(D\Phi\circ\sigma_{-t}^{\Phi\circ R}\circ R:D\Phi\circ\tau_t\circ R)_s\\
&=&R((D\Phi\circ\sigma_{-t}^{\Phi\circ R}:D\Phi\circ\tau_t)_{-s}^*)
\end{eqnarray*}
and, therefore $\alpha(h_L^{ist})=(D\Phi\circ\sigma_{-t}^{\Phi\circ R}:D\Phi\circ\tau_t)_{-s}^*$ from which one gets :
\[\alpha(h_L^{ist})=(D\Phi\circ\sigma_{-t}^{\Phi\circ R}:D\Phi\circ\tau_t)_s\]
Using (ii), we get :
\[(D\Phi:D\Phi\circ\sigma_{-t}^{\Phi\circ R})_s=\lambda^{ist}\alpha(h_R^{ist})\alpha(h_L^{-ist})\]
which is (iii). And applying \ref{proph} with $T_L=RT'R$, we obtain (iv). 
\end{proof}

%%%%%%%lem1tau
\subsection{Lemma}
\label{lem1tau}
{\it Let $M$ be a von Neumann algebra, $\Phi$ a normal semi-finite faithful weight on $M$, $\theta_t$ a one parameter group of automorphisms of $M$. Let us suppose that there exists a positive non singular operator $\mu$ affiliated to $M^\Phi$ such that, for all $s$, $t$ in $\mathbb{R}$, we have 
\[(D\Phi\circ\theta_t:D\Phi)_s=\mu^{ist}\]
We have then, for all $t\in\mathbb{R}$, $\theta_t(\mu)=\mu$. Let us write $\mu=\int_0^\infty \lambda de_\lambda$ the spectral decomposition of $\mu$, and let us define $f_n=\int_{1/n}^n de_\lambda$. We have then, for all $a$ in $\gN_\Phi$, $t$ in $\mathbb{R}$, $n$ in $\mathbb{N}$ :}
\[\omega_{J_\Phi\Lambda_\Phi(af_n)}\circ\theta_t=\omega_{J_\Phi\Lambda_\Phi(\theta_{-t}(a)f_n\mu^{t/2})}\]

\begin{proof}
Let us remark first that $\theta_t(\mu)=\mu$, and, therefore, $\theta_t(f_n)=f_n$. On the other hand, for any $a$ in $M$, we have :
\[\theta_{-t}\sigma_s^\Phi\theta_t(x)=\sigma_s^{\Phi\circ\theta_t}(x)=\mu^{ist}\sigma_s^\Phi(x)\mu^{-ist}\]
and then :
\[\theta_{-t}\sigma_s^\Phi(x)=\mu^{ist}\sigma_s^\Phi\theta_{-t}(x)\mu^{-ist}\]
If now $x$ is analytic with respect to $\Phi$, we get that $\theta_{-t}(f_nxf_m)$ is analytic with respect to $\Phi$ and that :
\[f_n\theta_{-t}\sigma_{i/2}^\Phi(x)f_m=\mu^{-t/2}f_n\sigma_{i/2}^\Phi(\theta_{-t}(x))f_m\mu^{t/2}\]
Let us take now $a$ in $\gN_\Phi$, analytic with respect to $\Phi$; we have, for any $y$ in $M$ :
\begin{eqnarray*}
\omega_{J_\Phi\Lambda_\Phi(f_naf_m)}\circ\theta_t(y)
&=&(\theta_t(y)J_\Phi\Lambda_\Phi(f_naf_m)|J_\Phi\Lambda_\Phi(f_naf_m))\\
&=&(\theta_t(y)\Lambda_\Phi(f_m\sigma_{-i/2}^\Phi(a^*)f_n)|\Lambda_\Phi(f_m\sigma_{-i/2}^\Phi(a^*)f_n))\\
&=&\Phi(f_n\sigma_{i/2}^\Phi(a)f_m\theta_t(y)f_m\sigma_{-i/2}^\Phi(a^*)f_n)
\end{eqnarray*}
which, using the preceeding remarks, is equal to :
\[\Phi\circ\theta_t(\mu^{-t/2}f_n\sigma_{i/2}^\Phi(\theta_{-t}(a))f_m\mu^{t/2}y\mu^{t/2}f_m\sigma_{-i/2}^\Phi(\theta_{-t}(a^*))f_n\mu^{-t/2})\]
and, making now $f_n$ increasing to $1$, we get that $\omega_{J_\Phi\Lambda_\Phi(af_m)}\circ\theta_t(y)$ is equal to :
\begin{multline*}
\Phi(\sigma_{i/2}^\Phi(\theta_{-t}(a))f_m\mu^{t/2}y\mu^{t/2}f_m\sigma_{-i/2}^\Phi(\theta_{-t}(a^*)))\\
=(y\Lambda_\Phi(f_m\mu^{t/2}\sigma_{-i/2}^\Phi(\theta_{-t}(a^*)))|\Lambda_\Phi(f_m\mu^{t/2}\sigma_{-i/2}^\Phi(\theta_{-t}(a^*)))\\
=(yJ_\Phi\Lambda_\Phi(\theta_{-t}(a)f_m\mu^{t/2})|J_\Phi\Lambda_\Phi(\theta_{-t}(a)f_m\mu^{t/2}))
\end{multline*}
from which we get the result. \end{proof}

%%%%%lem2tau
\subsection{Lemma}
\label{lem2tau}
{\it Let $(N, \alpha, \beta, M, \Gamma)$ be a Hopf-bimodule, equipped with a left-invariant operator-valued weight $T$, and a right-invariant valued weight $T'$; let $\nu$ be a normal semi-finite faithful weight on $N$, relatively invariant with respect to $T$ and $T'$ in the sense of \ref{defMQG}. We shall denote $\Phi=\nu\circ\alpha^{-1}\circ T$ and $\Psi=\nu\circ\beta^{-1}\circ T'$ the two lifted normal semi-finite weights on $M$, $R$ the coinverse and $\tau_t$ the scaling group constructed in \ref{R} and \ref{tau}. Then, we have :
\newline
(i) there exists a positive non singular operator $\mu_1$ affiliated to $M^\Phi$ and invariant under $\tau_t$, such that $(D\Phi\circ\tau_t:D\Phi)_s=\mu_1^{ist}$; let us write $\mu_1=\int_0^\infty \lambda de_\lambda$ and $f_n=\int_{1/n}^nde_\lambda$; we have then, for all $a$ in $\gN_\Phi$, $t$ in $\mathbb{R}$, $n$ in $\mathbb{N}$ and $x$ in $M^+$ :
\[\omega_{J_\Phi\Lambda_\Phi(\tau_t(a)f_n)}=\omega_{J_\Phi\Lambda_\Phi(af_n\mu_1^{t/2})}\circ\tau_{-t}\]
\[T\circ\tau_t(x)=\alpha\circ\sigma_t^\nu\circ\alpha^{-1}(T(\mu_1^{t/2}x\mu_1^{-t/2}))\]
(ii) there exists a positive non singular operator $\mu_2$ affiliated to $M^\Phi$ and invariant under $\sigma_t^{\Phi\circ R}$, such that $(D\Phi\circ\sigma_{-t}^{\Phi\circ R}:D\Phi)_s=\mu_2^{ist}$; let us write $\mu_2=\int_0^\infty \lambda de'_\lambda$ and $f'_n=\int_{1/n}^nde'_\lambda$; we have then, for all $b$ in $\gN_\Phi$, $t$ in $\mathbb{R}$ and $n$ in $\mathbb{N}$ :
\[\omega_{J_\Phi\Lambda_\Phi(bf'_n)}\circ\sigma_t^{\Phi\circ R}=\omega_{J_\Phi\Lambda_\Phi(\sigma_{-t}^{\Phi\circ R}(b)f'_n\mu_2^{-t/2})}\]
\[T(\sigma_{-t}^{\Phi\circ R}(\mu_1^{-t/2}x\mu_1^{t/2}))=\alpha\circ\gamma_{t}^L\circ\alpha^{-1}(T(x))\]
Moreover, we have $\mu_1^{is}=\lambda^{-is}\alpha(h_R^{-is})$, $\mu_2^{is}=\mu_1^{is}\alpha(h_L^{is})$, and $\mu_1^{is}$, $\mu_2^{is}$, $\alpha(h_L^{is})$ belong to $\alpha(N)'\cap M^\Phi$. The non-singular operators $\mu_1$, $\mu_2$ and $\alpha(h_L)$ commute two by two. }
\begin{proof} 
By \ref{corh}(ii), we get that $(D\Phi\circ\tau_t:D\Phi)_s=\lambda^{-ist}\alpha(h_R^{-ist})$, as $\lambda$ is positive non singular, affiliated to the center $Z(M)$, and $h_R$ is positive non singular affiliated to the center of $N$, we get there exists $\mu_1$ positive non singular, affiliated to $M^\Phi$ such that :
\[\mu_1^{ist}=\lambda^{-ist}\alpha(h_R^{-ist})=(D\Phi\circ\tau_t:D\Phi)_s\]
We can then apply \ref{lem1tau} to $\tau_t$ and $\tau_t(a)f_n$ (which belongs to $\gN_\Phi$) to get the first formula of (i). On the other hand, we get that $\alpha\circ\sigma_{-t}^\nu\circ\alpha^{-1}\circ T\circ \tau_t$ is a normal semi-finite operator-valued weight which verify, for all $x\in M^+$ 
\[\alpha\circ\sigma_{-t}^\nu\circ \alpha^{-1}\circ T\circ \tau_t(x)=T(\mu_1^{t/2}x\mu_1^{t/2})\]
from which we get the second formula of (i). 
\newline
By \ref{corh}(iii), we get that $(D\Phi\circ\sigma_{-t}^{\Phi\circ R}:D\Phi)_s=\lambda^{-ist}\alpha(h_R^{-ist})\alpha(h_L^{ist})$; with the same arguments, we get that there exists $\mu_2$ positive non singular, affiliated to $M^\Phi$ such that:
\[\mu_2^{ist}=\lambda^{-ist}\alpha(h_R^{-ist})\alpha(h_L^{ist})=(D\Phi\circ\sigma_{-t}^{\Phi\circ R}:D\Phi)_s\]
and we get the first formula of (ii) by applying again \ref{lem1tau} with $\sigma_{-t}^{\Phi\circ R}$. 
\newline
On the other hand, using \ref{gamma}, we get that $\alpha\circ\gamma_{-t}^L\circ\alpha^{-1}\circ T\circ\sigma_{-t}^{\Phi\circ R}$ is an operator-valued weight which verify, for all $x\in M^+$ :
\begin{eqnarray*}
\nu\circ\gamma_{-t}^L\circ\alpha^{-1}\circ T\circ\sigma_{-t}^{\Phi\circ R}(x)
&=&
\nu(h_L^{-t/2}\alpha^{-1} (T\sigma_{-t}^{\Phi\circ R}(x))h_l^{-t/2})\\
&=&\Phi(\alpha(h_L^{-t/2}\sigma_{-t}^{\Phi\circ R}(x)\alpha(h_L^{-t/2}))\\
&=&\Phi\circ\sigma_{-t}^{\Phi\circ R}[\alpha(h_L^{-t/2})x\alpha(h_L^{-t/2})]\\
&=&\Phi(\mu_2^{t/2}\alpha(h_L^{-t/2})x\alpha(h_L^{-t/2})\mu_2^{t/2})
\end{eqnarray*}
from which we get, because $\mu_2^{t/2}\alpha(h_L^{-t/2})$ commutes with $\alpha(N)$ :
\[\alpha\circ\gamma_{-t}^L\circ\alpha^{-1}\circ T\circ\sigma_{-t}^{\Phi\circ R}(x)=
T(\mu_2^{t/2}\alpha(h_L^{-t/2})x\alpha(h_L^{-t/2})\mu_2^{t/2})\]
or :
\[ T(\sigma_{-t}^{\Phi\circ R}(x))=\alpha\circ\gamma_t^L\circ\alpha^{-1} (T(\mu_1^{t/2}x\mu_1^{t/2}))\]
from which we finish the proof. \end{proof}

%%%%%thcentral
\subsection{Proposition}
\label{thcentral}
{\it Let $(N, \alpha, \beta, M, \Gamma)$ be a Hopf-bimodule, equipped with a left-invariant operator-valued weight $T$, and a right-invariant valued weight $T'$; let $\nu$ be a normal semi-finite faithful weight on $N$, relatively invariant with respect to $T$ and $T'$ in the sense of \ref{defMQG}. We shall denote $\Phi=\nu\circ\alpha^{-1}\circ T$ and $\Psi=\nu\circ\beta^{-1}\circ T'$ the two lifted normal semi-finite weights on $M$, $R$ the coinverse and $\tau_t$ the scaling group constructed in \ref{R} and \ref{tau}; let $\lambda$ be the scaling operator of $\Psi$ with respect to $\Phi$ (\ref{Vaes}), $\gamma^L$ and $\gamma^R$ the two one-parameter automorphism groups of $N$ introduced in \ref{gamma}  ; then :
\newline
(i) for all $t\in\mathbb{R}$, we have $\Gamma\circ\tau_t=(\sigma_t^{\Phi}\underset{N}{_\beta*_\alpha}\sigma_{-t}^{\Phi\circ R})\Gamma=
(\sigma_t^{\Psi\circ R}\underset{N}{_\beta*_\alpha}\sigma_{-t}^{\Psi})\Gamma$.
\newline
(ii) we have $h_L=h_R=1$, and $\nu\circ\gamma^L=\nu\circ\gamma^R=\nu$.
\newline
(iii) for all $s$, $t$ in $\mathbb{R}$, we have $(D\Phi:D\Phi\circ\tau_t)_s=(D\Psi:D\Psi\circ\tau_t)_s=\lambda^{ist}$. 
\newline
(iv) for all $s$, $t$ in $\mathbb{R}$, we have $(D\Phi\circ\sigma_t^{\Phi\circ R}: D\Phi)_s=\lambda^{ist}$. 
\newline
Therefore, the modular automorphism groups $\sigma^\Phi$ and $\sigma^{\Phi\circ R}$ commute, the weight $\nu$ is relatively invariant with respect to $\Phi$ and $\Phi\circ R$ and $\lambda$ is the scaling operator of $\Phi\circ R$ with respect to $\Phi$; and we have $\tau_t(\lambda)=\lambda$, $R(\lambda)=\lambda$; 
\newline
(v) there exists a non singular positive operator $q$ affiliated to $Z(N)$ such that $\lambda=\alpha(q)=\beta(q)$. }

\begin{proof}
As, for all $n\in N$, we have :
\[\sigma_{-t}^{\Phi\circ R}(\alpha(n))=R\sigma_t^\Phi R(\alpha(n))=\alpha(\gamma_t^L(n))\]
and, by definition, $\sigma_t^\Phi (\beta(n))=\beta(\gamma_t^L(n))$, using a remark made in \ref{fiber}, we may consider the automorphism $\sigma_{-t}^\Phi\underset{N}{_\beta*_\alpha}\sigma_{t}^{\Phi\circ R}$ on $M\underset{N}{_\beta*_\alpha}M$; let's take $a$ and $b$ in $\gN_\Phi\cap\gN_{T}$; let's write $h_L=\int_O^\infty \lambda de^L_\lambda$ and let us write $h_p=\int_{1/p}^p de^L_\lambda$; moreover, let's use the notations of \ref{lem2tau}; we get that :
\[(id\underset{N}{_\beta*_\alpha}\omega_{J_\Phi\Lambda_\Phi(b\alpha(h_p)f'_m)})(\sigma_{-t}^{\Phi}\underset{N}{_\beta*_\alpha}\sigma_{t}^{\Phi\circ R})\Gamma\circ\tau_t(f_na^*af_n)\]
is equal to :
\[\sigma_{-t}^\Phi(id\underset{N}{_\beta*_\alpha}\omega_{J_\Phi\Lambda_\Phi(b\alpha(h_p)f'_m)}\circ\sigma_t^{\Phi\circ R})\Gamma\circ\tau_t(f_na^*af_n)\]
which, thanks to \ref{lem2tau}(ii), can be written, because $\alpha(h_p)$ belongs to $\alpha(N)'\cap M^\Phi$, and therefore $b\alpha(h_p)$ belongs to $\gN_\Phi$ : 
\[\sigma_{-t}^\Phi(id\underset{N}{_\beta*_\alpha}\omega_{J_\Phi\Lambda_\Phi(\sigma_{-t}^{\Phi\circ R}(b\alpha(h_p))f'_m\mu_2^{-t/2})})\Gamma\circ\tau_t(f_na^*af_n)\]
or :
\[R\sigma_{t}^{\Phi\circ R}R(id\underset{N}{_\beta*_\alpha}\omega_{J_\Phi\Lambda_\Phi(\sigma_{-t}^{\Phi\circ R}(b\alpha(h_p))f'_m\mu_2^{-t/2})})\Gamma\circ\tau_t(f_na^*af_n)\]
By \ref{lem2tau} and \ref{lemT}, we know that $af_n\mu_1^{t/2}$ belongs to $\gN_\Phi\cap\gN_{T}$; using now \ref{lem2tau}(i), we get that $\tau_t(af_n)=\tau_t(a)f_n$ belongs to $\gN_\Phi\cap\gN_{T}$.
\newline
On the other hand, by \ref{lem2tau} and \ref{lemT}, we know that $b\alpha(h_p)f'_m$ belongs to $\gN_\Phi\cap\gN_{T}$; using now \ref{lem2tau}(ii), we get that :
\[\sigma_{-t}^{\Phi\circ R}(b\alpha(h_p)f'_m\mu_1^{-t/2})=\sigma_{-t}^{\Phi\circ R}(b)f'_m\mu_2^{-t/2}\alpha(h_p)\alpha(h_L^{t/2})\]
belongs to $\gN_\Phi\cap\gN_{T}$, and so, using again \ref{lemT}, 
\[\sigma_{-t}^{\Phi\circ R}(b)f'_m\mu_2^{-t/2}\alpha(h_p)=\sigma_{-t}^{\Phi\circ R}(b)f'_m\mu_2^{-t/2}\alpha(h_p)\alpha(h_L^{t/2})\alpha(h_p)\alpha(h_L^{-t/2})\] 
 belongs also to $\gN_\Phi\cap\gN_{T}$; therefore, we can use \ref{R}, and we get it is equal to : 
\[R\sigma_{t}^{\Phi\circ R}(id\underset{N}{_\beta*_\alpha}\omega_{J_\Phi\Lambda_\Phi(\tau_t(a)f_n)})\Gamma(\mu_2^{-t/2}f'_m\alpha(h_p)\sigma_{-t}^{\Phi\circ R}(b^*b)\alpha(h_p)f'_m\mu_2^{-t/2})\]
which can be written, thanks to \ref{lem2tau}(i) : 
\[R\sigma_{t}^{\Phi\circ R}(id\underset{N}{_\beta*_\alpha}\omega_{J_\Phi\Lambda_\Phi(af_n\mu_1^{t/2})}\circ\tau_{-t})\Gamma(\mu_2^{-t/2}f'_m\alpha(h_p)\sigma_{-t}^{\Phi\circ R}(b^*b)\alpha(h_p)f'_m\mu_2^{-t/2})\]
or, $\alpha(h_p)$, as well as $\mu_2^{-t/2}f'_m$, being invariant under $\sigma_t^{\Phi\circ R}$ : 
\[R(id\underset{N}{_\beta*_\alpha}\omega_{J_\Phi\Lambda_\Phi(af_n\mu_1^{t/2})})(\sigma_{t}^{\Phi\circ R}\underset{N}{_\beta*_\alpha}\tau_{-t})\Gamma\circ\sigma_{-t}^{\Phi\circ R}
(\mu_2^{-t/2}f'_m\alpha(h_p)b^*b\alpha(h_p)f'_m\mu_2^{-t/2})\]
and using \ref{tau}, and again \ref{R}, we get it is equal to : 
\begin{multline*}
R[(id\underset{N}{_\beta*_\alpha}\omega_{J_\Phi\Lambda_\Phi(af_n\mu_1^{t/2})})\Gamma(\mu_2^{-t/2}f'_m\alpha(h_p)b^*b\alpha(h_p)f'_m\mu_2^{-t/2})]\\
=
(id\underset{N}{_\beta*_\alpha}\omega_{J_\Phi\Lambda_\Phi(b\alpha(h_p)f'_m\mu_2^{-t/2})})\Gamma(\mu_1^{t/2}f_na^*af_n\mu_1^{t/2})
\end{multline*}
Finally, we have proved that, for all $a$, $b$ in $\gN_\Phi\cap \gN_{T}$, $m,n,p$ in $\mathbb{N}$, we have :
\begin{multline*}
(id\underset{N}{_\beta*_\alpha}\omega_{J_\Phi\Lambda_\Phi(b\alpha(h_p)f'_m)})(\sigma_{-t}^{\Phi}\underset{N}{_\beta*_\alpha}\sigma_{t}^{\Phi\circ R})\Gamma\circ\tau_t(f_na^*af_n)=\\
(id\underset{N}{_\beta*_\alpha}\omega_{J_\Phi\Lambda_\Phi(b\alpha(h_p)f'_m\mu_2^{-t/2})})\Gamma(\mu_1^{t/2}f_na^*af_n\mu_1^{t/2})
\end{multline*}
But, for all $x,y\in M$, we have :
\[\omega_{J_\Phi\Lambda_\Phi(b\alpha(h_p)f'_m)}(x)=\omega_{J_\Phi\Lambda_\Phi(b)}(\alpha(h_p)f'_mxf'_m\alpha(h_p))\]
\[\omega_{J_\Phi\Lambda_\Phi(b\alpha(h_p)f'_m\mu_2^{-t/2})}(y)=\omega_{J_\Phi\Lambda_\Phi(b)}(\alpha(h_p)f'_m\mu_2^{-t/2}x\mu_2^{-t/2}f'_m\alpha(h_p))\]
and, therefore, we get that :
\[(id\underset{N}{_\beta*_\alpha}\omega_{J_\Phi\Lambda_\Phi(b)})[(1\underset{N}{_\beta\otimes_\alpha}\alpha(h_p) f'_m)(\sigma_{-t}^{\Phi}\underset{N}{_\beta*_\alpha}\sigma_{t}^{\Phi\circ R})\Gamma\circ\tau_t(f_na^*af_n)(1\underset{N}{_\beta\otimes_\alpha}f'_m\alpha(h_p))]\]
is equal to :
\[(id\underset{N}{_\beta*_\alpha}\omega_{J_\Phi\Lambda_\Phi(b)})[(1\underset{N}{_\beta\otimes_\alpha}\alpha(h_p) f'_m\mu_2^{-t/2})\Gamma(\mu_1^{t/2}f_na^*af_n\mu_1^{t/2})
(1\underset{N}{_\beta\otimes_\alpha}\mu_2^{-t/2}f'_m\alpha(h_p))]\]
and, by density, we get that :
\[(1\underset{N}{_\beta\otimes_\alpha}\alpha(h_p) f'_m)(\sigma_{-t}^{\Phi}\underset{N}{_\beta*_\alpha}\sigma_{t}^{\Phi\circ R})\Gamma\circ\tau_t(f_na^*af_n)(1\underset{N}{_\beta\otimes_\alpha}f'_m\alpha(h_p))\]
is equal to :
\[(1\underset{N}{_\beta\otimes_\alpha}\alpha(h_p) f'_m\mu_2^{-t/2})\Gamma(\mu_1^{t/2}f_na^*af_n\mu_1^{t/2})
(1\underset{N}{_\beta\otimes_\alpha}\mu_2^{-t/2}f'_m\alpha(h_p))\]
and, after making $p$ going to $\infty$, we obtain that :
\[(1\underset{N}{_\beta\otimes_\alpha} f'_m)(\sigma_{-t}^{\Phi}\underset{N}{_\beta*_\alpha}\sigma_{t}^{\Phi\circ R})\Gamma\circ\tau_t(f_na^*af_n)(1\underset{N}{_\beta\otimes_\alpha}f'_m)\]
is equal to $(*)$:
\[(1\underset{N}{_\beta\otimes_\alpha}\ f'_m\mu_2^{-t/2})\Gamma(\mu_1^{t/2}f_na^*af_n\mu_1^{t/2}))(1\underset{N}{_\beta\otimes_\alpha}\mu_2^{-t/2}f'_m)
\]
Let's now take a file $a_i$ in $\gN_\Phi\cap\gN_{T}$ weakly converging to $1$; going to the limit, we get that :
\[(1\underset{N}{_\beta\otimes_\alpha} f'_m)(\sigma_{-t}^{\Phi}\underset{N}{_\beta*_\alpha}\sigma_{t}^{\Phi\circ R})\Gamma\circ\tau_t(f_n)(1\underset{N}{_\beta\otimes_\alpha} f'_m)=
(1\underset{N}{_\beta\otimes_\alpha} f'_m\mu_2^{-t/2})\Gamma(\mu_1^{t/2}f_n\mu_1^{t/2})(1\underset{N}{_\beta\otimes_\alpha}\mu_2^{-t/2}f'_m)\]
When $n$ goes to $\infty$, then $f_n$ is increasing to $1$, the first is increasing to $1\underset{N}{_\beta\otimes_\alpha}f'_m$, and the second is increasing to $(1\underset{N}{_\beta\otimes_\alpha} f'_m\mu_2^{-t/2})\Gamma(\mu_1^{t})(1\underset{N}{_\beta\otimes_\alpha}\mu_2^{-t/2}f'_m)$
which is therefore bounded.
\newline
Taking now $m$ going to $\infty$, we get that the two non-singular operators $\Gamma(\mu_1^t)$ and $1\underset{N}{_\beta\otimes_\alpha}\mu_2^t$ are equal. Using \ref{lembeta}, we get then that $\mu_1$ is equal to $\mu_2$ (and is affiliated to $\beta(N)$), from which we get, using \ref{lem2tau}, that $h_L=1$. Applying all these calculations to $(N, \alpha, \beta, M, \Gamma, RT'R, T', \nu)$, we get that $h_R=1$, which is (ii). 
\newline
Let's come back to the equality $(*)$ above; we obtain that : 
\[(1\underset{N}{_\beta\otimes_\alpha} f'_m)(\sigma_{-t}^{\Phi}\underset{N}{_\beta*_\alpha}\sigma_{t}^{\Phi\circ R})\Gamma\circ\tau_t(f_na^*af_n)(1\underset{N}{_\beta\otimes_\alpha} f'_m)\]
is equal to :
\[(1\underset{N}{_\beta\otimes_\alpha}f'_m)\Gamma(f_na^*af_n)(1\underset{N}{_\beta\otimes_\alpha} f'_m)\]
So, when $n$ and $m$ go to $\infty$, we obtain :
\[(\sigma_{-t}^{\Phi}\underset{N}{_\beta*_\alpha}\sigma_{t}^{\Phi\circ R})\Gamma\circ\tau_t(a^*a)
=\Gamma(a^*a)\]
which, by density, gives the first formula of (i), the secong being given then by \ref{RTR}. 
\newline
From (ii) and \ref{corh} (i) and (ii), we get (iii). 
\newline
From (ii) and \ref{corh}(iii), we get that $(D\Phi\circ\sigma_t^{\Phi\circ R}: D\Phi)_s=\lambda^{ist}$; therefore, as $\lambda$ is affiliated to $Z(M)$, we get the commutation of the modular groups $\sigma^\Phi$ and $\sigma^{\Phi\circ R}$. Using \ref{Vaes}, we get that there exists $\lambda_R$ positive non singular affiliated to $Z(M)$ and $\delta_R$ positive non singular affiliated to $M$ such that $(D\Phi\circ R:D\Phi)_t=\lambda_R^{it^2/2}\delta_R^{it}$, and the properties of $R$ allows us to write that $R(\lambda_R)=\lambda_R$. But, on the other hand, the formula $(D\Phi\circ\sigma_t^{\Phi\circ R}: D\Phi)_s=\lambda_R^{ist}$ (\ref{Vaes}), gives that $\lambda_R=\lambda$ and, therefore, we get that $R(\lambda)=\lambda$. The formula $\tau_t(\lambda)=\lambda$ comes from (iii), which finishes the proof of (iv). 
\newline
By (i), we have $\lambda=\mu_1=\mu_2$, and, as we had proved that $\mu_1$ is affiliated to $\beta(N)$, we get that $\lambda$ is affilated to $\beta(N)$; as $R(\lambda)=\lambda$ by (iv), we get (v). \end{proof}

%%%th
\subsection{Theorem}
\label{th}
{\it Let $(N, M, \alpha, \beta, \Gamma, T, T', \nu)$ be a measured quantum groupoid in the sense of \ref{defMQG}, and let us denote $R$ (resp. $\tau_t$) the co-inverse (resp. the scaling group) constructed in \ref{R} (resp. \ref{tau}). Then :
\newline
(i) $(N, M, \alpha, \beta, \Gamma, T, RTR, \nu)$ (resp. $(N, M, \alpha, \beta, \Gamma, RT'R, T', \nu)$) is a measured quantum groupoid . Moreover, $R$ (resp. $\tau_t$) remains the co-inverse (resp. the scaling group) of this measured quantum groupoid. 
\newline
(ii) $(N, M, \alpha, \beta, \Gamma, T, R, \tau, \nu)$ is a measured quantum groupoid in the sense of \cite{L2}, 4.1. 
\newline
(iii) conversely if $(N, M, \alpha, \beta, \Gamma, T, R, \tau, \nu)$ is a measured quantum groupoid in the sense of \cite{L2}, 4.1, then $(N, M, \alpha, \beta, \Gamma, T, RTR, \nu)$ is a measured quantum groupoid in the sense of \ref{defMQG}. }

\begin{proof}
Result (i) is given by \ref{thcentral}(iv); (ii) is given by \ref{R}, \ref{def} and \ref{thcentral}(ii); and (iii) is given by \cite{L2}, 5.3. \end{proof}

\end{appendix}

\newpage

%%%%%%bibli

\end{document}